\documentclass[11pt]{amsart}
\usepackage{amsfonts,amsthm,amsmath,amssymb,etoolbox,comment,graphicx,array,tabu,commath,latexsym,environ,tikz,tikz-cd,stackengine,diagbox,caption,enumitem}
\usepackage{soul}

\usepackage[colorlinks, linkcolor = black, anchorcolor = black, citecolor = black, urlcolor=black]{hyperref}

\newcommand{\cC}{\mathcal{C}}
\newcommand{\cD}{\mathcal{D}}
\newcommand{\cE}{\mathcal{E}}

\newcommand{\cI}{\mathcal{I}}

\newcommand{\cS}{\mathcal{S}}

\newcommand{\cV}{\mathcal{V}}

\newcommand{\cY}{\mathcal{Y}}

\newcommand{\E}{E}

\newcommand{\R}{\mathbb R}

\newcommand{\Z}{\mathbb Z}

\newcommand{\Hess}{\operatorname{Hess}}
\newcommand{\dist}{\operatorname{dist}}

\renewcommand{\div}{\mathrm{div}}
\renewcommand{\tilde}{\widetilde}

\newcommand{\vol}{\operatorname{vol}}

\newcommand{\x}{\times}

\newcommand{\Ric}{\operatorname{Ric}}

\newcommand{\Diff}{\mathrm {Diff}}

\usepackage[style=alphabetic, backend=biber, sorting=nyt, url=false ]{biblatex}
\newtheorem{thm}{Theorem}[section]

\newtheorem{prop}[thm]{Proposition}
\newtheorem{lem}[thm]{Lemma}
\newtheorem{cor}[thm]{Corollary}

\theoremstyle{definition}
\newtheorem{definition}[thm]{Definition}
\newtheorem{example}[thm]{Example}

\theoremstyle{claim}
\newtheorem{claim}[thm]{Claim}

\theoremstyle{remark}
\newtheorem{remark}[thm]{Remark}

\title[]{
A variational construction for  
Einstein metrics on connected sums of $S^2\x S^2$}

\author{Adrian Chun-Pong Chu}
\author{Daniel Stern}
\address{Cornell University, Department of Mathematics, 310 Malott Hall, Ithaca, NY 14853}
\email{cc2938@cornell.edu}
\email{daniel.stern@cornell.edu}

\begin{document}
\maketitle
\begin{abstract}
For each positive integer $n$, we give a variational construction of a positive Einstein metric on the connected sum of $n$ copies of $S^2\times S^2$, with coupled toric and dihedral symmetries. More precisely, we obtain these metrics via a min-max construction for the total scalar curvature, which is in some sense dual to the usual max-min formulation of the smooth Yamabe invariant within this symmetry class, with more favorable analytic properties. 
\end{abstract}

\setcounter{tocdepth}{1}

\tableofcontents

\section{Introduction}
Understanding the existence and structure of canonical metrics on a given manifold is among the oldest problems in differential geometry. In dimension 2, the classical uniformization theorem supplies every conformal class with a metric of constant Gaussian curvature. Likewise, for manifolds with dimension 3 or above, the solution to the Yamabe problem shows that every conformal class contains at least one metric of constant scalar curvature \cite{Yamabe1960Deformation,Trudinger1968Conformal,aubin1976yamabe,schoen1984conformal}, obtained by minimizing the total scalar curvature functional. 

A more restrictive notion of canonical metric is the {\it Einstein metric}: a metric $g$ whose Ricci curvature $\Ric_g$ satisfies 
$$\Ric_g=\Lambda g$$
for some constant $\Lambda$. Depending on the sign of $\Lambda$, Einstein metrics fall into three types: {\it positive Einstein}, {\it Ricci-flat}, and {\it negative Einstein}. In dimensions two and three, where the Ricci curvature determines the full Riemannian curvature tensor, every  Einstein manifold is a space form. In dimension four, Einstein metrics diverge from the constant-curvature models, and many nontrivial examples have been discovered since their introduction over a century ago. For compact negative Einstein metrics, in addition to the homogeneous examples, many K\"ahler examples were first produced by the Aubin-Yau theorem \cite{Aub78,Yau1978}. Subsequent examples were constructed  in \cite{Anderson2006,Bamler2012,FinePremoselli2020,HamenstadtJackel2024,GuenanciaHamenstadt2025}. 

In the Ricci-flat setting, there are famously few compact examples, with the only ones currently known in dimension four being quotients of the flat $T^4$ and K3 surfaces with Calabi-Yau metrics \cite{Yau1978}. However, the library of complete \emph{non-compact} four-dimensional Ricci-flat manifolds, also known as {\it gravitational instantons} \cite{Hawking1977GravitationalInstantons}, is much larger, and growing in recent years. In addition to several classical models with origins in general relativity, like the (Riemannian) Schwarzschild, Kerr, Taub--NUT, and Eguchi--Hanson metrics, in recent decades, several other families have been discovered; we refer readers to Li-Sun's excellent survey \cite{liSun2026geometryInstantons} for further discussion.

As for positive Einstein 4-manifolds, the classical simply-connected examples are $S^4$, $S^2\x S^2$, and  the Fubini-Study metric on $\mathbb{CP}^2$. In \cite{Page1978}, physicist D. Page constructed a positive Einstein metric on  $\mathbb{CP}^2\#\overline{\mathbb{CP}^2}$, Tian and Yau  constructed K\"ahler-Einstein examples on  $\mathbb{CP}^2\# k\overline{\mathbb{CP}^2}$  for $3\leq k\leq 8$ \cite{TianYau1987,Tian1990}, and Chen-LeBrun-Weber also found a conformally K\"ahler one on $\mathbb{CP}^2\# 2\overline{\mathbb{CP}^2}$ \cite{ChenLeBrunWeber2008}.  

In this paper, we use variational methods for the total scalar curvature to show that the connected sum of arbitrarily many copies of $S^2\times S^2$ admits a positive Einstein metric.

\begin{thm}\label{thm:main}
There exists a positive Einstein metric on $\#n(S^2\x S^2)$ for every  $n$.
\end{thm}

    While we were completing this manuscript, three separate preprints by Ozuch, Li-Sun, and Broder-Lafuente appeared on the arXiv, which also construct positive Einstein metrics on $\#n(S^2\times S^2)$, among other spaces \cite{Ozuch2026PositiveEinstein,liSun2026ToricEinstein,broderLafuente2026t2Einstein}. Their methods for constructing these metrics appear to differ substantially from our variational approach, which we hope will be of independent interest.

 \medskip

We arrive at these metrics via a min-max scheme for the volume-normalized counterpart
$$\tilde\cE(g):=\frac 1{\vol(M,g)^{1/2}}\int_M R_gdV_g$$
of the Einstein-Hilbert functional $\cE(g):=\int_M R_g dV_g$, within a class $\cD$ of metrics on $\#n(S^2\times S^2)$ invariant under coupled actions of $T^2$ and the dihedral group $D_{4n+4}$, which we describe in more detail in Section \ref{sect:geom.ein}. Modulo gauge fixing, we show that metrics in $\cD$ can be parametrized by a trio of scalar functions $\phi,v,q$, with $v,q$ parametrizing all  \emph{conformal classes} in $\cD$, and $\phi$ a positive function which determines the area of the $T^2$-fibers, and serves as a universal conformal factor. For each $\phi$, we then consider the (possibly infinite) supremum of $\tilde{\cE}(g_{\phi,v,q})$ over all conformal types $(v,q)$, then minimize over all $\phi$.

In principle, one could also consider a dual max-min scheme, restricting the classic smooth Yamabe invariant of Kobayashi and Schoen  \cite{Kobayashi1987,schoen1987variational} to metrics in $\cD$: that is, for each conformal class $(v,q)$, one could minimize over all $\phi$ as in \cite{HebeyVaugon1993Equivariant} to obtain a constant scalar curvature representative, then take the supremum over all conformal classes $(v,q)$. 

However, the min-max path taken here has two key advantages over the classic Yamabe approach. First, for any fixed $\phi$, the functional $(v,q)\mapsto \tilde{\cE}(g_{\phi,v,q})$ enjoys a kind of convexity property, which forces \emph{uniqueness} of a maximizing pair $(v,q)$ when one exists, in contrast to the Yamabe minimization problem, where Yamabe metrics may be non-unique in conformal classes of positive Yamabe type. One important caveat is that, even when the supremum $\sup_{(v,q)} \tilde{\cE}(g_{\phi,v,q})$ is finite, the maximizing pair $(v_{\phi},q_{\phi})$ typically gives rise to a non-smooth metric $g_{\phi,v_{\phi},q_{\phi}}$, with edge-cone type singularities along $2$-spheres preserved by the $T^2$ action, and to interpret $(v_{\phi},q_{\phi})$ as  a maximum, the Einstein-Hilbert functional must be extended suitably to metrics of this type.

The second advantage is a straightforward lower semi-continuity property for the functional $\phi\mapsto \sup_{(v,q)}\tilde{\cE}(g_{\phi,v,q})$ with respect to bounded weak convergence of $\phi$ in $W^{1,2}(M)$ and strong convergence in $L^4(M)$ for some fixed reference metric. As a result, showing existence of some $\phi\in W^{1,2}(M)$ realizing the infimum is analytically identical to the classic conformally-constrained Yamabe problem, and can be reduced to a single estimate as in \cite{HebeyVaugon1993Equivariant} ruling out bubbling at the fixed points of the $T^2$ action. In particular, the weak existence question is reduced to a single numerical estimate, leaving only the substantial technical challenge of establishing regularity of this weak minimizer. 

In the presence of additional geometric structure, we note that there have been many other successful applications of variational methods to the construction of Einstein metrics, e.g. on homogeneous spaces \cite{Jensen1971,WangZiller1986,Bohm2004,BohmWangZiller2004}, and in the K\"ahler setting \cite{BergerCherrier1988,Lauret2002,BermanBerndtsson2013,BermanBoucksomGuedjZeriahi2013,BermanGuenancia2014,DarvasRubinstein2017,BermanBoucksomEyssidieuxGuedjZeriahi2019,BermanBoucksomJonsson2021}.

\subsection{Geometry of the Einstein manifolds}\label{sect:geom.ein}
Let us  describe  the geometry of the  Einstein metrics in Theorem \ref{thm:main}. Roughly speaking, each of them is a {\it diagonal toric metric with a certain dihedral symmetry}.

Let $D_k$ denote the dihedral group of order $k.$ For each $n\geq 0$, let $(\Sigma,h)$ denote any smooth Riemannian surface with $2n+2$ geodesic boundary edges and $2n+2$ right-angled corners. Write $\partial\Sigma$ as a union  $\Gamma_1\cup\Gamma_2$, where each $\Gamma_i$ consists of $n+1$ closed, pairwise disjoint edges.  We impose on $\Sigma$ a dihedral symmetry $D_{4n+4}$. To be explicit, we can write 
$$ D_{4n+4}=\langle r,s|r^{2n+2}=s^2=(sr)^2=1\rangle,$$ where $r$ is a  rotation by angle $ \frac{\pi}{n+1}$ and $s$ a reflection about the axis perpendicular to two arbitrarily chosen opposite edges. Now, let $S^1=\R/2\pi\Z$. We consider on $\Sigma\x S^1\x S^1$ warped metrics of the form 
\begin{equation}\label{eq:h.f1.f2}
g=h+f_1^2d\theta_1^2+f_2^2d\theta_2^2,
\end{equation}
where $f_1,f_2: \Sigma\to\R$ are smooth functions satisfying the following.
\begin{itemize}
    
\item {\it Boundary conditions}:  Let $\nu$ denote the outward unit normal along the boundary of $(\Sigma,h)$. Then
\begin{align*}
    f_1&=0, \quad\partial_\nu f_1=-1,\quad \partial_\nu f_2=0\quad \textrm{ on }\Gamma_1;\\
    f_2&=0, \quad\partial_\nu f_2=-1,\quad \partial_\nu f_1=0\quad \textrm{ on }\Gamma_2.
\end{align*}
\item {\it Positivity}: $f_1>0$ on $\Sigma\backslash\Gamma_1$, and $f_2>0$ on $\Sigma\backslash\Gamma_2$.
    \item {\it Symmetry}:  
\begin{align*}
     f_1\circ r=  f_2,\quad f_2\circ r=f_1, \quad f_1\circ s=f_1,\quad f_2\circ s=f_2.
\end{align*}
\end{itemize}
It is elementary to check that the above boundary conditions close up the boundary of  $\Sigma\x S^1\x S^1$  {\it smoothly}  and produce a 4-manifold $\#n(S^2\x S^2)$. The Einstein metrics in Theorem \ref{thm:main} take the form \eqref{eq:h.f1.f2} in which $h,$ $f_1$, and $f_2$ satisfy all  conditions above.

\subsection{A min-max scheme for the Einstein--Hilbert functional}\label{sect:intro.min.max}
We now describe in more detail the min-max problem for the Einstein--Hilbert functional, restricted to the collection of diagonal toric metrics with dihedral symmetry on $M:=\# n(S^2\x S^2)$, which satisfy the description in \S \ref{sect:geom.ein}. We denote this collection of metrics by $\cD$.

\subsubsection{Gauge fixing}
On the $(2n+2)$-gon $\Sigma$, consider the collection of all metrics $h$ such that $(\Sigma,h)$ has geodesic edges and right-angled corners, with dihedral symmetry $D_{4n+4}$ (recall \S \ref{sect:geom.ein}). Then, this collection of metrics all lies in a single conformal class, modulo $D_{4n+4}$-equivariant diffeomorphisms. In particular, $(\Sigma,h)$ can be realized as one quarter of a closed surface $S$ of genus $n$ with an action of $\Z^2_2\rtimes D_{4n+4}$, which determines a unique symmmetric conformal structure on $S$.

As a corollary of this elementary fact, we can fix on $\Sigma$ a background metric $h_0$, such that every metric in $\cD$, modulo diffeomorphisms, can be expressed in the form
\begin{equation}\label{eq:intro.g.u.f1.f2}
    g=e^{2u}h_0 + f_1^2 d\theta_1^2+f_2^2 d\theta_2^2
\end{equation}
for some $u:\Sigma\to\R$ and $f_1,f_2:\Sigma\to [0,\infty)$ satisfying the appropriate symmetry conditions, with boundary conditions that ensure smoothness on  $M$. This gauge-fixing will be crucial for the variational problem for $\tilde \cE$, and it is the main benefit of imposing the dihedral symmetry $D_{4n+4}$.

Before   describing the variational setup, it will be convenient to switch to an alternative parametrization for metrics in $\cD$. Namely, we first fix on $M$ an arbitrary smooth background metric
\begin{equation}\label{eq:intro.g0}
    g_0=h_0+\varphi^2_1 d\theta_1^2+\varphi_2^2d\theta_2^2\in\cD.
\end{equation}
It is easy to see that any metric of the form \eqref{eq:intro.g.u.f1.f2} can also be expressed as
$$g=\phi^2(e^{2v}h_0+e^{2q}\varphi_1^2d\theta_1^2+e^{-2q}\varphi_2^2d\theta_2^2),$$
for some $\phi\in C^\infty(M)$ and $v,q\in C^\infty(\Sigma)$ that satisfy  the natural symmetry conditions and boundary conditions to ensure smoothness. Note that the pair $(v,q)$ parametrizes the conformal classes in $\cD$, modulo diffeomorphism, while $\phi$ serves as a universal conformal factor. 

\subsubsection{Min-max for $\tilde\cE$} Fixing on $M$ an arbitrary background metric $g_0$ of the form \eqref{eq:intro.g0}, we consider the following min-max problem. First, for any $(v,q)$ satisfying the requisite boundary conditions, let
$g_{v,q}:=e^{2v}h_0+e^{2q}\varphi_1^2d\theta_1^2+e^{-2q}\varphi_2^2d\theta_2^2$. For any fixed $\phi\not\equiv 0$,  we take supremum over the conformal classes:
\begin{equation}\label{eq:intro.S.g0}
    \cS_{g_0}(\phi):=\sup_{v,q} \tilde \cE(\phi^2g_{v,q}),
\end{equation}
then take infimum over all $\phi$:
\begin{equation}\label{sect:intro.I.g0}
    \cI_{g_0}:=\inf_{\phi} \cS_{g_0}(\phi).
\end{equation}

There are several nice features of this min-max problem. First, while the supremum \eqref{eq:intro.S.g0} for any fixed $\phi$ depends on the background metric $g_0$, the min-max value  \eqref{sect:intro.I.g0} is {\it independent of $g_0$}: This is elementary to prove, based on the fact that total scalar curvature and volume are invariant under diffeomorphisms.

Moreover, it is not difficult to show that the functional $F(\phi,v,q):= \cE(\phi^2g_{v,q})$ can be extended continuously to arbitrary smooth triples $(\phi,v,q)$ satisfying the requisite positivity and symmetry conditions, allowing for edge-cone type singularities in the metrics $\phi^2g_{v,q}$. Within this class, it is straightforward to find many functions $\phi$ for which the maximum is achieved.

\begin{remark}
    The extension of $\cE$ to metrics with edge-cone singularities should be compared with the modified mass of Alaee, Khuri, and Kunduri in the setting of graviational instantons \cite{AlaeeKhuriKunduri2026MassToricManifolds}. 
\end{remark}

For any fixed $\phi$, the pair $(v,q)$ is a maximizer for \eqref{eq:intro.S.g0} if and only if $\Ric_g=\Lambda g$ holds {\it in the $T^2$-fiber direction} and $\Lambda=\frac 1{4\vol(g)}\int_MR_g$ for the metric $g=\phi^2g_{v,q}$, which typically will have edge-cone singularities; moreover, if the maximizer exists, it is not difficult to show that it must be unique. 

Finally, for the full min-max problem, a smooth triple $(\phi,v,q)$ attains the min-max value \eqref{sect:intro.I.g0}, meaning $$\tilde \cE(\phi^2g_{v,q})=\cS_{g_0}(\phi)=\cI_{g_0},$$
if and only if $\phi^2 g_{v,q}$ is  Einstein, and $\phi^2g_{v,q}$ has no edge-cone singularities. This crucially uses that the triples $(\phi,v,q)$  parametrize all metrics in $\cD$, modulo diffeomorphism.

\subsubsection{Existence and regularity}

The functional 
$$\cS_{g_0}: C^{\infty}(M,(0,\infty))\to (-\infty,+\infty]$$
defined by \eqref{eq:intro.S.g0} admits a simple, though typically infinite, extension to $W^{1,2}_{g_0}(M)$, which is easily seen to be lower semi-continuous with respect to $L^4$-convergence along minimizing sequences. As a result, the existence of a minimizer $\phi\in W^{1,2}_{g_0}(M)$ can be reduced, much as in \cite{HebeyVaugon1993Equivariant}, to verifying the strict inequality
$$\cI_{g_0}<\sqrt{2n+2}Y_0$$
for the infimum, where $Y_0$ is the Yamabe invariant of the standard $4$-sphere, to rule out bubbling at the $2n+2$ fixed points of the $T^2$ action. For each $n$, an explicit $\phi_n\in W^{1,2}(\#n(S^2\times S^2))$ satisfying $\cS_{g_0}(\phi)<\sqrt{2n+2}Y_0$ is constructed in Proposition \ref{prop:bound.I} below from a kind of branched cover of $S^4$, confirming this estimate.

With the existence of an $\cS_{g_0}$-minimizing $\phi\in W^{1,2}_{g_0}(M)$ in place, we can then establish a weak Euler-Lagrange equation, showing that the (possibly very singular) metric $g=\phi^2 g_{v_{\phi},q_{\phi}}$ satisfies a weak version of $Ric_g=\Lambda g$ along the $T^2$ fibers, a weak version of the constant scalar curvature condition $R_g=4\Lambda$, and a weak version of the boundary conditions for $v,q$ needed to rule out edge-cone singularities.

From the weak Euler-Lagrange equations, regularity of the triple $(\phi,v,q)$ in the interior of the polygon $\Sigma$ can be proved in a relatively straightforward way, and using the gauge-invariance of the infimum $\cI_{g_0}$ and a simple variational argument, we can check that the resulting metric $g$ is critical for $\tilde{\cE}$ among all perturbations in the interior of $\Sigma$, hence Einstein.

From this point, as discussed in Section \ref{sect:boot.v.leq.C}, regularity across the remaining $2$-spheres and fixed points in $M$ corresponding to edges and vertices of $\partial\Sigma$ can be reduced to proving an upper bound for the scalar function $v$. The remainder of the argument exploits the Einstein equations, the weak Euler-Lagrange equations for $\cS_{g_0}$, and the a priori global bounds following from the finiteness $\cS_{g_0}(\phi)<\infty$ to provide such a bound for $v$. 

Away from the fixed points, the boundedness of $v$ near invariant spheres is obtained from a combination of lower bounds for $\max\{f_1,f_2\}^{-2}e^{-2u}|d\rho|_{h_0}^2$, $\max\{f_1,f_2\}$, and $\phi$, and upper bounds for $|d\rho|_{h_0}^2$, where $\rho=f_1f_2$ is the fiber area function. The Einstein equation is used here to supply several useful Bochner identities and a crucial ODE for $f_2^{-2}e^{-2u}|d\rho|_{h_0}^2$ along the level set flow of $\rho$ in the region $\{f_2>f_1\}$, while the weak Euler-Lagrange equations and global energy bounds are needed to justify some global a priori bounds and a weak version of the boundary conditions for $v$ and $q$ required to close the arguments. 

Finally, regularity up to the fixed points is proved by establishing uniform bounds for $v$ on small (quarter-)annuli around the vertices of $\Sigma$, with many of the key estimates resembling careful rescaled versions of those used to establish regularity away from the fixed point set. These bounds are used to justify a maximum principle argument supplying a global upper bound for the Gauss curvature $K_h$ of the base metric $h=e^{2u}h_0$, from which the desired bound for $v$ follows relatively easily.

\subsection{Organization} In \S \ref{sect:dia.toric.metric}, we describe the class of diagonal toric metrics with dihedral symmetry on $\# n(S^2\x S^2)$, and in \S \ref{sect:variation.Ein.Hil} describe the relevant min-max problem for the Einstein--Hilbert functional on this class of metrics. In \S \ref{sect:existenceMin}, we prove the existence of a minimizing $\phi$ that attains the min-max value $\cI$. In \S \ref{sect:weak.E-L}, we derive various a priori estimates from the finiteness $\cS_{g_0}(\phi)<\infty$ and derive the weak form of the Euler-Lagrange equations for the minimizer. In \S \ref{sect:int.reg}, we confirm that the associated metric that is smooth and Einstein in the interior of $\Sigma\x S^1\x S^1$. To show that the metric is fully smooth on $M$, after collecting some useful identities and a priori estimates in \S \ref{sect:comp.sec}, we show in section \S \ref{sect:boot.v.leq.C} that boundedness of $v$ implies the full regularity result. Finally, we obtain the desired bounds for $v$ across the invariant two-spheres, away from the fixed points, in \S \ref{sect.reg.rods}, and finally across the $2n+2$ fixed points in \S \ref{sect:full.reg}, finishing the proof of Theorem \ref{thm:main}.

\subsection*{Acknowledgments} 
The authors would like to thank Michael Anderson, Shih-Kai Chiu, Misha Karpukhin, Alan Tsz-Lok Lam, Vincent S.H. Lee, and Mingyang Li for helpful and insightful discussions. During the completion of this project, AC was partially supported by the AMS-Simons travel grant, and DS was partially supported by NSF grant DMS 2404992, the Simons Foundation, and the Alfred P. Sloan foundation.

\subsection{AI use}\label{ai:log}

In the first phase of this project, from October 2025-June 2026, AI played only a limited role, aiding with literature searches and standard curvature calculations.

\begin{itemize}
    \item In October 2025, after completing \cite{chuStern2025minimalDoublings}, we began to explore Anderson's proposed gluing construction for negative Einstein metrics described in \cite[\S 7]{anderson2008survey}, and the possibility of adapting it to construct positive Einstein metrics on $\# n(S^2\x S^2)$. 
    \item In January 2026, we started to focus on the diagonal toric metrics with dihedral symmetry described in \S \ref{sect:dia.toric.metric}, aiming to construct a positive Einstein analog of the periodic Schwarzschild metrics of Myers and Korotkin-Nicolai. At this time, we also started investigating the possibility of a variational approach to the construction. 

    \item In March-April 2026, we formulated the min-max construction used here, and began to focus exclusively on this approach.
    \item From April-June 2026, we proved conditional existence of the weak minimizer (cf. Proposition \ref{prop:exist.minimizer.conditional}), derived a version of the weak Euler-Lagrange equation (cf. Proposition \ref{weak.el}), and began approaching regularity in the interior of $\Sigma$ via a strategy similar to that of Section \ref{sect:int.reg}, but with an (ultimately unsuccessful) Moser-Trudinger argument in place of the Brezis-Merle estimates. 
\end{itemize}

In July 2026, we began to use AI in essential ways.

\begin{itemize}
\item In July 2026, ChatGPT 5.5 found the candidate in Proposition \ref{prop:bound.I} for establishing the upper bound on $\cI_{g_0}$, and explained how to use Brezis-Merle estimates to complete the proof of regularity in the interior of $\Sigma$.
\item In early September 2026, ChatGPT 6 supplied an initial proof of regularity across the invariant spheres and fixed points, which was refined in subsequent exchanges between the authors and AI. At this point, we focused on digesting these arguments, and began to write the manuscript. We have carefully checked and expanded on those arguments, and take responsibility for their correctness, but cannot claim any real credit for this portion of the proof. 
\end{itemize}

In the final write-up, ChatGPT once again aided with some standard computations, streamlining arguments, and further literature review, as well as identifying many typos and minor gaps in earlier drafts. 

{\bf In short:} the content of Sections 1-5, except for Proposition \ref{prop:bound.I} and Section \ref{l0.sec}, is due primarily to the authors. The main ideas of Section 6 contain a mix of AI and human input. For the content of Sections 7-10, our role is largely expository, expanding on outlines suggested by AI and polishing them into a more coherent narrative. 

This article contains no AI-generated text.

\section{Diagonal toric metrics with dihedral symmetry}\label{sect:dia.toric.metric}
Fix any non-negative integer $n$, and denote by $\# n(S^2\x S^2)$ the connected sum of $n$ copies of $S^2\x S^2$ with the standard smooth structure; the $n=0$ case corresponds to the standard  $S^4$. We are going to consider metrics on it with a certain diagonal toric action and dihedral symmetry. Results in this section are all elementary. Readers may refer to  \cite[Chapter 4]{dammerman2006metricsWithSymmetry} and \cite{liu2025toricEinstein}  for some background on toric metrics, and also Li-Sun's survey \cite{liSun2026geometryInstantons} for previous works on  toric Ricci-flat metrics.
\subsection{The smooth manifold $\#n(S^2\x S^2)$}\label{sect:topo.mfd}
Let $\Sigma$ be a smooth surface with corners, homeomorphic to a disc, with  $2n+2$ boundary edges and $2n+2$ corners.   Write $\partial\Sigma$ as a union  $\Gamma_1\cup\Gamma_2$, where each $\Gamma_i$ consists of $n+1$ closed, pairwise disjoint edges.  Denote by $\mathring\Gamma_i$ the interior of $\Gamma_i$.

Let $S^1:= \R/2\pi\Z$ and $T^2:= S^1\x S^1$. We consider the  smooth manifold with corners $\Sigma\x T^2$, and  denote its points by $(p,\theta_1,\theta_2)$. We close up its boundary as follows.
\begin{itemize}
    \item For any $p\in \mathring\Gamma_1$ and $\theta_2\in S^1$, collapse $\{p\}\x S^1\x\{\theta_2\}$ into a point.
    \item For any $p\in \mathring\Gamma_2$ and $\theta_1\in S^1$, collapse $\{p\}\x\{\theta_1\}\x S^1$ into a point.
    \item For any vertex $p\in  \Gamma_1\cap \Gamma_2$,  collapse  $\{p\}\x S^1\x S^1$ into a point.
\end{itemize}
Denote by $M$ the quotient space.

To equip $M$ with a smooth structure,  we first define some coordinate charts\footnote{We abuse notation here, as the  charts will not map into $\R^4$ and will be  non-injective near $\partial\Sigma\x T^2$.} on $\Sigma\x T^2$. For any $p\in \Sigma$, let $U_p\subset \Sigma$ be a neighborhood around $p$.

\begin{itemize}
    \item If $p\in \mathring\Gamma_1\cup\mathring\Gamma_2$, we  identify $U_p\cong \{(x,y):x\geq 0,y\in\R\}$. Define a  chart $U_p\x T^2\to \R^3\x S^1$  by    
    \begin{align*}
        &(x\cos\theta_1,x\sin\theta_1,y,\theta_2)\quad \textrm{ if }p\in \mathring\Gamma_1,\\
        &(x\cos\theta_2,x\sin\theta_2,y,\theta_1)\quad \textrm{ if }p\in \mathring\Gamma_2.
    \end{align*}
 
    \item If $p$ is a vertex, we identify $U_p\cong \{(x,y):x,y\geq 0\}$.
Define a  chart 
$U_p\x T^2\to \R^2\x\R^2$ by
\begin{equation}\label{eq:vertex.chart}
(x\cos\theta_1,x\sin\theta_1,y\cos\theta_2,y\sin\theta_2)
\end{equation} 
\item If $p\in\mathring\Sigma$, we just identify $U_p\x T^2\cong \R^2\x T^2.$
\end{itemize}
It is a straightforward exercise to show that the  above charts on $\Sigma\x T^2$ descend naturally onto $M$, turning $M$ into a closed smooth 4-manifold.

\begin{example}\label{example:S4}
The case $n=0$ can be written down explicitly. Denote by $S^k\subset\R^{k+1}$ the unit sphere. Consider the bi-gon
\begin{equation}\label{eq:bi-gon}
    \Sigma=\{(x,y,z)\in S^2:x,y\geq 0\}.
\end{equation}
Define a map $F:\Sigma\x T^2\to S^4$ by
$$F((x,y,z),\theta_1,\theta_2)=(x\cos\theta_1,x\sin\theta_1,y\cos\theta_2,y\sin\theta_2,z).$$
It is elementary to check that $F$ induces a smooth diffeomorphism from $M$ to $S^4$, under the above smooth charts on $M$.
\end{example}

\begin{prop}\label{prop:M.homeo.nS2xS2} For each  $n\geq 0$, $M$ is diffeomorphic to $\#n(S^2\x S^2)$.
\end{prop}
\begin{proof} The case $n=0$ was already explained in Example \ref{example:S4}. For $n\geq 1$, we prove by induction. For clarity we write $\Sigma_n$ and $M_n$ to specify the dependence on $n$.     When $n=1$, we can take $\Sigma_1=[0,1]\x [0,1]$ without loss of generality. Then  writing $\Sigma_1\x T^2=([0,1]\x S^1)^2$, the quotienting process  above closes up each $[0,1]\x S^1$ factor into an $S^2$. Hence, $M_1$ is diffeomorphic to $S^2\x S^2.$

Assume $M_{n-1}$  is diffeomorphic to $\#(n-1)(S^2\x S^2)$. Fix a vertex $p$ of $\Sigma_{n-1}$ and  a neighborhood $U\subset\Sigma_{n-1}$ around it. We can assume  $U\cong\{(x,y):x,y\geq 0,x^2+y^2\leq 1\}$, so that  the chart \eqref{eq:vertex.chart} sends   $U\x T^2$ to the unit 4-ball. Thus, after quotienting, $(\Sigma_{n-1}\backslash U)\x T^2$ becomes diffeomorphic to $\#(n-1)(S^2\x S^2)$ with a 4-ball removed. Similarly, we fix a vertex $q$ of $\Sigma_1$ and neighborhood $V\subset \Sigma_1$ of $q$ so that after quotienting, $(\Sigma_{1}\backslash V)\x T^2$ becomes diffeomorphic to $S^2\x S^2$ with a ball removed. Finally, noting  $\Sigma_n$ can be obtained by gluing $\Sigma_{1}\backslash V$ to $\Sigma_{n-1}\backslash U$, it follows  that $M_n$ is the connected sum of $\#(n-1)(S^2\x S^2)$ and $S^2\x S^2$, as desired.
\end{proof}

The $T^2$-action on $\Sigma\x T^2$ naturally gives a smooth $T^2$-action on $M$, with $\Sigma=M/T^2$. Let $S_1\subset M$ (resp. $S_2$) be the fixed point set of the action by the first (resp. second) $S^1$-factor. Note, each $S_i$ is a collection of $n+1$ disjoint smooth spheres, with $S_i/T^2=\Gamma_i$.

\subsection{Diagonal toric metrics}\label{sect:diagonal.toric}
We can glue  four copies of $\Sigma$ together along their boundaries appropriately to obtain a smooth closed  surface $S$ of genus $n$. In fact, we can fix an inclusion $\Sigma \hookrightarrow S$ such that there exist ``reflection" maps $R_1,R_2:S\to S$ for which $S=\Sigma \cup R_1(\Sigma)\cup R_2(\Sigma)\cup R_1R_2(\Sigma)$. We may assume $R_i$ is the reflection that fixes in particular $\Gamma_i$, which sits in $S$ under $\Sigma \hookrightarrow S$. 

Let $h$ be any smooth Riemannian metric on $S$ that is invariant under $R_1$ and $R_2$, so that $\Z^2_2\cong\langle R_1, R_2\rangle<\textrm{Isom}(S,h)$. This makes each edge of $\partial\Sigma\subset S$ geodesic and each vertex  right-angled. We then consider smooth functions $f_1,f_2:S\to \R$ with the following properties:
\begin{itemize}
    \item {\it $\Z_2^2$-symmetry}:  
    \begin{align}
    \label{eq:f1.sym}    f_1\circ R_1=-f_1, \quad f_1\circ R_{2}=f_1,\\
    \label{eq:f2.sym}      f_2\circ R_1=f_2, \quad f_2\circ R_{2}=-f_2.
    \end{align}
    \item {\it Positivity}: for $i=1,2$, 
    \begin{equation}\label{eq:positivity}
        f_i>0\quad\textrm{ on }\Sigma\backslash \Gamma_i
    \end{equation} 
\end{itemize}
In particular, letting  $\nu_h$ be the outward unit normal for $(\Sigma,h)$ along $\mathring\Gamma_1\cup\mathring\Gamma_2$, we can deduce  that  
$$f_1=0 \textrm{ and }\partial_{\nu_h}f_2=0\;\;\textrm{ on } \Gamma_1,\quad f_2=0 \textrm{ and }\partial_{\nu_h}f_1=0\;\;\textrm{ on } \Gamma_2,$$
using their  smoothness across $\partial\Sigma\subset S$.

Under the above setting, we consider metrics on  $\Sigma\x T^2$ of the form
\begin{equation}\label{eq:g}
    g := h + f_1^2d\theta_1^2+ f_2^2d\theta_2^2.
\end{equation}
Since $f_i=0$ on $\Gamma_i$, the metric $g$ is  singular on $\partial\Sigma\x T^2$. Moreover, it can be straightforwardly checked that $g$ induces on $M$ a metric with conical singularities along $S_1\cup S_2$. In fact, due to the toric action, at each point on $(S_1\cup S_2)\backslash (S_1\cap S_2)$ the tangent cone takes the form $C_\beta\x\R^2$, where $C_\beta$ is the cone with angle $\beta$ defined by
$$dr^2+(\beta/2\pi)^2 r^2 d\theta^2,\quad (r,\theta)\in [0,\infty)\x(\R/2\pi\Z);$$
and at each point on $S_1\cap S_2$, the tangent cone has the form $C_{\beta_1}\x C_{\beta_2}$. 

But there is a criterion on $(f_1,f_2)$ that ensures  $g$ induces a fully smooth metric on $M$. Namely, we say that $(f_1,f_2)$ is {\it admissible} with respect to $h$ if \begin{equation}\label{eq:admiss.f1f2}
    \partial_{\nu_h} f_i=-1\quad\textrm{ on } \mathring \Gamma_i.
\end{equation}

\begin{prop}\label{prop:induce.smooth.metric.M}
The metric $g$ on $\Sigma\x T^2$ given by \eqref{eq:g} induces a smooth metric on $M$ if and only if   $(f_1,f_2)$ is admissible with respect to $h$.
\end{prop}

\begin{proof}The calculation is elementary; see for example \cite[\S 2]{liu2025toricEinstein} for details.
\end{proof}

\subsection{Dihedral symmetry}\label{sect:dihedral.sym}
Let $D_{k}$ denote the dihedral group of order $k$. On top of the $\Z_2^2$-symmetry, we now impose in addition a $D_{4n+4}$-symmetry on $S$, by requiring $\Sigma$ to be a {\it regular} $(2n+2)$-gon. We can express $$ D_{4n+4}=\langle r,s|r^{2n+2}=s^2=(sr)^2=1\rangle<\mathrm{Isom}(S),$$ in which, when restricted onto $\Sigma$, $r$ is a  rotation by angle $ \frac{\pi}{n+1}$ and $s$ a reflection about the axis perpendicular to two arbitrarily chosen opposite edges. Note,  $r$  sends $\Gamma_1$ to $\Gamma_2$ and vice versa, while $s$ preserves each of  $\Gamma_1$ and $\Gamma_2$. Let $L$ be the group  generated by $r,s,R_1,R_2\in\mathrm{Isom}(S)$. It can be readily checked that $L$ is a semi-direct product $\Z^2_2\rtimes D_{4n+4}$. Note, if $n\geq 2$, $L=\textrm{Isom}(S)$.

For any closed smooth  $(S,h)$ with $L$-symmetry, we consider $f_1,f_2\in C^\infty(S)$ that satisfy the {\it $L$-symmetry} condition:
\begin{align}
    \label{eq:f1.L.sym}    f_1\circ R_1=-f_1, \quad f_1\circ R_{2}=f_1, \quad f_1\circ r=f_2,\quad f_1\circ s=f_1,\\
    \label{eq:f2.L.sym}      f_2\circ R_1=f_2, \quad f_2\circ R_{2}=-f_2, \quad f_2\circ r=f_1,\quad f_2\circ s=f_2.
\end{align}
Now, for any $f_1,f_2\in C^\infty(S)$ that satisfy this $L$-symmetry, positivity \eqref{eq:positivity}, and admissibility \eqref{eq:admiss.f1f2}, we consider metrics on $\Sigma\x T^2$  of the form
$$g=h+f_1^2d\theta_1^2+f^2_2d\theta_2^2.$$
By Proposition \ref{prop:induce.smooth.metric.M}, these metrics induce smooth metrics on $M=\# n(S^2\x S^2)$. We denote by $\cD$ the collection of all  metrics on $M$ that arise this way.

The maps $r$ and $s$ on $\Sigma$ induce the following  diffeomorphisms on $M$, which will still be denoted respectively $r$ and $s$:
$$r(p,\theta_1,\theta_2):=(r(p),\theta_2,\theta_1),\quad s(p,\theta_1,\theta_2):=(s(p),\theta_1,\theta_2).$$
Note, every  metric $g\in\cD$ is invariant under  $r,s\in \textrm{Diff}(M)$ and the  $T^2$-action. Denote by $G<\mathrm{Diff}(M)$ the group generated by $r,s$ and the $T^2$-action, which in fact is a semi-direct product $T^2\rtimes D_{4n+4}$. 

\begin{example} 
When $n=0$, $M$ is the standard 4-sphere.  As follows from  Example \ref{example:S4}, denoting by $h_0$ the round metric on the bi-gon $\Sigma$ given by \eqref{eq:bi-gon}, the  round metric on $S^4$ can be written in the form
$$h_0+x^2d\theta_1^2+y^2d\theta_2^2\in\cD.$$
When $n=1$, letting $\Sigma=\{(x,y)\in [0,\pi]^2\}$ and $h_0$ be the flat metric, the standard metric on $S^2\x S^2$ can be written as
$$h_0+\sin^2 x d\theta_1^2+\sin^2 y d\theta_2^2\in \cD.$$
\end{example}

\subsection{Gauge-fixing} With the $L$-action on $S$ defined  above, we let $\Diff^L_0(S)$ denote the set of all diffeomorphisms on $S$ that are $L$-equivariantly isotopic to the identity map. The following proposition illustrates a key advantage of imposing dihedral symmetry: 
\begin{prop}\label{prop:unique.conf.class} For each $n\geq 2$,  
the $L$-invariant smooth metrics on $S$ all lie in one conformal class, modulo   diffeomorphisms in $\Diff^L_0(S)$.
\end{prop}
\begin{proof}
By the uniformization theorem, it suffices to show there is a unique $L$-invariant hyperbolic metric on $S$, up to isometries. It therefore suffices to show, for the $(2n+2)$-gon $\Sigma$, the uniqueness of a hyperbolic metric with geodesic edges, right angles, and $D_{4n+4}$-symmetry. 

On the hyperbolic plane, we choose an origin $O$ and a $D_{4n+4}$-action by isometries that fixes $O$. It is easy to see that there is a 1-parameter family of $D_{4n+4}$-invariant $(2n+2)$-gons with geodesic edges, parametrized by the distance of the vertexes from $O$. In fact, by the Gauss-Bonnet formula, there is a unique one with right angles. 
\end{proof}

This proposition allows gauge-fixing on $M$. Namely, we fix an $L$-invariant background metric $h_0$ on $S$. We consider $u\in C^\infty(S)$ with the symmetry
\begin{equation}\label{eq:u.sym}
   u\circ R_1=u\circ R_2= u\circ r=u\circ s=u.
\end{equation} 
Note, for any  $f_1,f_2,u\in C^\infty(S)$ satisfying the $L$-symmetry   \eqref{eq:f1.L.sym},   \eqref{eq:f2.L.sym}, and \eqref{eq:u.sym},  and  positivity \eqref{eq:positivity}, $(f_1,f_2)$ is admissible with respect to $e^{2u}h_0$ if and only if $$\partial_{\nu_{h_0}} f_i=-e^u\quad\textrm{ on }\mathring\Gamma_i.$$ In this case, we say the triple  $(u,f_1,f_2)$ is {\it admissible} with respect to the background $h_0$. Then by Proposition \ref{prop:unique.conf.class}, we can fix a gauge on $M$:
\begin{cor}\label{cor:guage.fix.S} On $S$ we fix a background metric $h_0$ with $L$-symmetry. 
For any   metric  $g_1\in \cD$ on $M$, there exist an   admissible triple $(u,f_1,f_2)$ with respect to $h_0$ and a  $\Phi\in\Diff^L_0(S)$ such that 
$$g_1=\Phi^*(e^{2u}h_0) + (f_1\circ\Phi)^2d\theta_1^2+ (f_2\circ\Phi)^2d\theta_2^2.$$
\end{cor}

\subsection{Curvature formulas}\label{sect:curv.formula} 
We record some curvature formulas for metrics $g=h+f_1^2d\theta_1^2+f^2_2d\theta^2_2$ on $M$. From now on, we will denote $\rho:=f_1f_2$. Consider on $M$ the metric 
$$g=h+f_1^2d\theta_1^2+f^2_2d\theta^2_2.$$
The Ricci curvature of $g$ can be readily computed \cite[\S 9.106]{besse2007einstein}: 
\begin{align}\label{eq:Ric.g}
    \nonumber\Ric_g=&\;\left(K_hh-\frac{\Hess_h f_1}{f_1}-\frac{\Hess_h f_2}{f_2}\right)\\
    &-\left(f_1\Delta_h f_1+\frac{f_1}{f_2}\langle df_1,df_2\rangle_h\right)d\theta^2_1-\left(f_2\Delta_h f_2+\frac{f_2}{f_1}\langle df_1,df_2\rangle_h\right)d\theta^2_2,
\end{align}
where $K_h$ is the Gaussian curvature of $h$. Hence, the Einstein equation $\Ric_g-\lambda g=0$ is equivalent to the following system of equations:
\begin{itemize}
    \item From the fiberwise equations $\Ric_g(\partial_{\theta_i},\partial_{\theta_i})=\lambda f_i^2$, we obtain
    \begin{align}
    \label{eq:fiberEin1} -\Delta_h\rho&=2\lambda \rho\\
    \label{eq:fiberEin2}  \div_h(f_1df_2-f_2df_1)&=0.
    \end{align}
\item Tracing over the base directions $T\Sigma$ gives 
\begin{equation}\label{eq:pureTracePart}
2K_h\rho-f_1\Delta_h f_2-f_2\Delta_h f_1=2\lambda \rho.
    \end{equation}
    \item And the trace-free component $(\Ric_g|_{T\Sigma})^{\circ}$ of the Ricci curvature restricted to $T\Sigma$ vanishes.
\end{itemize}
 
\begin{prop}\label{prop:adm.pure.trace}Let $(S,h)$ be any smooth surface with $L$-symmetry, and $f_1,f_2\in C^\infty(S)$ satisfy the $L$-symmetry \eqref{eq:f1.L.sym} and \eqref{eq:f2.L.sym}, and  positivity \eqref{eq:positivity}.  Let $$g=h+f_1^2d\theta_1^2+f_2^2d\theta_2^2.$$ 
If  $(\Ric_g|_{T\Sigma})^\circ =0$ on   $\mathring\Sigma\x T^2$, then $\partial_{\nu_h} f_i=c$ on $\Gamma_i$ for some constant $c\leq 0.$ 
\end{prop}
\begin{proof}If $(\Ric_g|_{T\Sigma})^\circ=0$, then  by \eqref{eq:Ric.g}
$$f_1(\Hess_h f_2)^\circ+f_2(\Hess_h f_1)^\circ=0,$$
with $^\circ$ denoting the trace-free part. On $\mathring\Gamma_1$, $f_1=0$ and $f_2>0$, therefore $(\Hess_h f_1)^\circ=0$. In fact, since  $\Delta_h f_1=0$ along $\Gamma_1$ due to $f_1\circ R_1=-f_1$,   we actually have $\Hess_h f_1=0$. Hence, denoting by $\tau_h$ the unit tangent along $\Gamma_1$, and using the fact that $\Gamma_1$ is geodesic in $S$, we have
$$\partial_{\tau_h}(\partial_{\nu_h} f_1)=\Hess_h f_1(\tau_h,\nu_h)=0.$$
This implies $\partial_{\nu_h} f_1=c$ along $\Gamma_1$ for some constant $c$, with $c\leq0$ by the assumption $f_1>0$ on $\mathring\Sigma$. Arguing similarly and using the dihedral symmetry, we know $\partial_{\nu_h} f_2=c$ along $\Gamma_2$ with the same constant $c$. 
\end{proof}

Now, let $h_0$ be any background metric on $S$.  If $h=e^{2u}h_0$, then we have the following relation for Gaussian curvatures:
$$K_h=e^{-2u}(K_{h_0}-\Delta_{h_0}u).$$
For the metric $g=e^{2u}h_0+f_1^2d\theta_1^2+f_2^2d\theta_2^2$, we can express $\Ric_g$ with respect to the reference metric $h_0$ as 
\begin{align*}
   \Ric_g =&\;(K_{h_0} -\Delta_{h_0}  u)h_0-\sum^2_{i=1}\frac{1}{f_i}(\Hess_{h_0} f_i-du\otimes df_i-df_i\otimes du+\langle du,df_i\rangle_{h_0} h_0)\\
    &-e^{-2u}\left( f_1\Delta_{h_0} f_1+\frac {f_1}{f_2} \langle df_1,df_2\rangle_{h_0}\right) d\theta_1^2-e^{-2u}\left( f_2\Delta_{h_0} f_2+\frac {f_2}{f_1} \langle df_1,df_2\rangle_{h_0}\right) d\theta_2^2,
\end{align*}
Taking trace, the scalar curvature $R_g$ can be written as 
\begin{equation}\label{eq:R_g}
    R_g=2 e^{-2u}\left(K_{h_0}-\Delta_{h_0} u-\frac{\Delta_{h_0} f_1}{f_1}-\frac{\Delta_{h_0} f_2}{f_2}-\frac{\langle df_1,df_2\rangle_{h_0}}{f_1f_2}\right).
\end{equation} 
To simplify notation,  we sometimes omit the subscripts $h_0$.

\section{A variational setup for the Einstein--Hilbert functional}\label{sect:variation.Ein.Hil}
For the remainder of the paper, we   assume $n\geq 2$.
\subsection{Notations}  Throughout this section,  we choose   on $M$  some   diagonal toric background metric with dihedral symmetry
\begin{equation}\label{eq:g0}
g_0=h_0+\varphi_1^2d\theta_1^2+\varphi_2^2d\theta_2^2\in\cD
\end{equation}
(recall the definition of $\cD$  in  \S \ref{sect:dihedral.sym}).  
For the variational theory of the Einstein--Hilbert functional, instead of expressing elements  $g\in\cD$ in the form $e^{2u}h_0+f_1^2d\theta_1^2+f_2^2d\theta_2^2$, it is more natural to consider the  alternative parametrization \begin{equation}\label{eq:g.phi.v.q}
    g=\phi^2(e^{2v}h_0+e^{2q}\varphi_1^2d\theta_1^2+e^{-2q}\varphi_2^2d\theta_2^2),
\end{equation}
where we view $\phi\in C^\infty(M)$ and $v,q\in C^\infty(S)$. 

To be more precise, let $\mathcal C$ be the set of all pairs $(v,q)\in C^\infty(S)\x C^\infty(S)$ satisfying the following assumptions.
\begin{itemize}
    \item {\it $L$-symmetry}:
\begin{align}
    &v\circ R_1=v\circ R_2=v\circ r=v\circ s=v,\\
   \label{eq:L.sym.q} & q\circ R_1=q\circ R_2=-q\circ r=q\circ s=q.
\end{align}
\item {\it Admissibility}:
$$v=q\quad\textrm{ on }\Gamma_1,\quad v=-q\quad\textrm{ on }\Gamma_2.$$
\end{itemize}
Moreover, denote by $C^\infty(M,\R_{>0})^G$ the set of all $G$-invariant positive functions $\phi\in C^\infty(M)$.

\begin{lem}\label{lem:change.param} 
The following relation gives a one-to-one correspondence between admissible triples $(u,f_1,f_2)$ with respect to $h_0$ and  triples $(\phi,v,q)\in C^\infty(M,\R_{>0})^G\x \cC$:
$$e^{2u}h_0+f_1^2d\theta_1^2+f_2^2d\theta_2^2=\phi^2(e^{2v}h_0+e^{2q}\varphi_1^2d\theta_1^2+e^{-2q}\varphi_2^2d\theta_2^2).$$
\end{lem}
\begin{proof}
    The proof is straightforward.
\end{proof}

From Corollary \ref{cor:guage.fix.S} and Lemma \ref{lem:change.param}, it follows:
\begin{cor}\label{cor:guage.fix.M}
For any   metric  $g_1\in \cD$ on $M$, there exist some   triple  $(\phi,v,q)\in C^\infty(M,\R_{>0})^G\x \cC$   and      $\Phi\in\Diff^L_0(S)$ such that 
$$g_1=\Phi^*(\phi^2(e^{2v}h_0+e^{2q}\varphi_1^2d\theta_1^2+e^{-2q}\varphi_2^2d\theta_2^2)).$$
\end{cor}
Note, we abused notation by viewing $\Phi$ as a diffeomorphism on $M$ too.

In other words, fixing a background $g_0$, the set $\cC$ can be  identified with the set of all conformal classes of metrics in $\cD$ modulo diffeomorphisms. For convenience, we will denote $\rho_0:=\varphi_1\varphi_2$, and 
$$g_{v,q}:=e^{2v}h_0+e^{2q}\varphi_1^2d\theta_1^2+e^{-2q}\varphi_2^2d\theta_2^2.$$
Moreover, to simplify notation, we often omit the subscripts $g_0$ and $h_0$ in various geometric quantities or operators. We may also omit measures $dV_{g_0}$ and $dA_{h_0}$ in integrals.

\begin{remark}
For each  $(v,q)\in\cC$, 
the area of each $T^2$-fiber under the metric $g_{v,q}$ is   the same as that under the background metric $g_0$; in fact, the areas are given by the function $4\pi^2\rho_0$. As for the metric $\phi^2 g_{v,q}$, the areas of the $T^2$-fibers are given by $4\pi^2\phi^2\rho_0.$
\end{remark}

\subsection{The min-max scheme for the restricted Einstein--Hilbert action}
On the 4-manifold $M$, the {\it Einstein--Hilbert functional} $\cE$ on the set of all smooth metrics is given by the total  scalar curvature:
$$\cE(g):=\int_M R_g dV_g.$$
The {\it normalized} Einstein--Hilbert functional is given by
$$\tilde \cE(g):=\frac 1{\vol(g)^{1/2}}\int_M R_g dV_g.$$
In this section, we will restrict to the class of metrics   $\cD$, and introduce a min-max scheme.

First, fixing a background $g_0\in \cD$, we know from Corollary \ref{cor:guage.fix.M} that every  metric in $ \cD$ can be written in the form $\phi^2 g_{v,q}$ modulo diffeomorphisms.
Recall that, the total scalar curvature of $g:=\phi^2 g_{v,q}$ can be computed  by
$$\cE(\phi^2 g_{v,q})=\int_M(6|d\phi|^2_{g_{v,q}}+R_{g_{v,q}}\phi^2)dV_{g_{v,q}},$$
and the volume is $\int_M \phi^4dV_{g_{v,q}}$ (see, e.g., the classic surveys by Lee-Parker \cite{leeParker1987yamabe} and Schoen \cite{schoen1987variational} on the Yamabe problem).  
\begin{definition}
For  any
$(\phi,v,q)\in W^{1,2}(M)^G\x\cC$,  define
$$E_{g_0}(\phi,v,q):=\int_M (6|d\phi|^2_{g_{v,q}}+R_{g_{v,q}}\phi^2)dV_{g_{v,q}}.$$
If in addition $\phi$ is not a.e. zero, we may also normalize by volume and define
$$\tilde E_{g_0}(\phi,v,q):=\frac{E_{g_0}(\phi,v,q)}{(\int_M \phi^4dV_{g_{v,q}})^{1/2}}.$$ 
\end{definition}
In other words, $E_{g_0}$ corresponds to the Einstein--Hilbert functional, while $\tilde E_{g_0}$ corresponds to the normalized  functional. Next, fixing any  $\phi$, we take supremum over all conformal classes modulo diffeomorphisms.

\begin{definition}
Define $\cS_{g_0}:\{\phi\in W^{1,2}(M)^G|\phi \textrm{ not a.e. zero}\}\to (-\infty,+\infty]$ by taking a supremum:
$$\cS_{g_0}(\phi):=\sup\left\{\tilde E_{g_0}(\phi,v,q)\middle| (v,q)\in \cC\right\}.$$  
Then we take infimum over all $\phi$:
$$\cI_{g_0}:=\inf\{\cS_{g_0}(\phi)|\phi\in W^{1,2}(M)^G \textrm{ not a.e. zero}\}.$$
\end{definition}

In the next \S \ref{sect:basic.prop}, we will show that this infimum is finite. But first, let us  also introduce a variant of the usual symmetric smooth Yamabe invariant over $\cD$. Namely, for any  metric $g\in\cD$, define 
$$\cY_\cD(g):=\inf\left\{\frac{\int_M (6|d\phi|^2_{g }+R_{g }\phi^2)dV_{g }}{(\int_M \phi^4dV_{g })^{1/2}}\middle|\phi\in W^{1,2}(M)^G\textrm{ not a.e. zero}\right\},$$ 
and define the {\it smooth Yamabe invariant} for the class of metrics $\cD$ by 
$$\sigma_\cD:=\sup\{\cY_\cD(g):g\in \cD\},$$  
which is easily seen to satisfy $\sigma_\cD>-\infty$.

\subsection{Basic properties}\label{sect:basic.prop}

Next, we note that a lower bound $\cI_{g_0}>-\infty$ follows easily from the following comparison.
 \begin{prop}\label{prop:sigmaD.leq.Emin}For any background metric $g_0\in\cD$, we have
    $\sigma_\cD\leq \cI_{g_0}.$
\end{prop}
\begin{proof}
For any $\phi\in W^{1,2}(M)^G$ not a.e. zero and   $(v,q)\in\cC$, the definitions of $\mathcal{Y}_\cD$ and $\cS_{g_0}$ immediately give
    $$\cY_\cD(g_{v,q})\leq \tilde E_{g_0}(\phi,v,q) \leq \cS_{g_0}(\phi).$$
Taking infimum over $\phi$ and supremum over $(v,q)$, and using Corollary  \ref{cor:guage.fix.M}, the desired claim follows.
\end{proof}

Next, we show that $\cI_{g_0}<\infty$. Let us first express $E_{g_0}$ more transparently with respect to the background $g_0$.  By a direct computation from \eqref{eq:R_g}, 
$$e^{2v} R_{g_{v,q}}=R_{g_0}-2\Delta_{h_0} v-2\langle d\log(\varphi_1/\varphi_2),dq\rangle_{h_0}-2|dq|_{h_0}^2.$$
Hence, noting the relation $|d\phi|^2_{g_{v,q}} dV_{g_{v,q}}=|d\phi|^2_{g_0} dV_{g_0}$,
\begin{align}\label{eq:E.in.g0}
E_{g_0}(\phi,v,q)=&\;\int_M(6|d\phi|^2_{g_0}+R_{g_0}\phi^2)dV_{g_0}-2\int_M \phi^2\Delta_{h_0}v dV_{g_0}\\
\nonumber &-2\int_M\phi^2(\langle d\log(\varphi_1/\varphi_2),dq\rangle_{h_0}+|dq|_{h_0}^2)dV_{g_0} 
\end{align}
In particular, the first integral depends only on $\phi$, the second on $\phi$ and $v$, and the third on $\phi$ and $q$. To simplify notation, we sometimes omit the subscripts $g_0$ and $h_0$. 

If $\phi$ is in addition smooth, then  integrating by parts we have
\begin{align*}
 E_{g_0}(\phi,v,q)=&\;\int_M(6|d\phi|^2+R_{g_0}\phi^2)dV-8\pi^2\int_\Sigma v\Delta\rho dA\\
&\;+8\pi^2\int_\Sigma[-q\div(\phi^2(\varphi_1d\varphi_2-\varphi_2d\varphi_1))-\rho|dq|^2 ] dA\\
&\;+8\pi^2\int_{\partial\Sigma}[v\partial_\nu\rho+q\phi^2(\varphi_1\partial_\nu\varphi_2-\varphi_2\partial_\nu\varphi_1)]ds
\end{align*}
Using $\varphi_i=0$ on $\Gamma_i$ and the admissibility of $(v,q)$, the boundary integral vanishes, giving
\begin{align} 
E_{g_0}(\phi,v,q)=&\;\int_M(6|d\phi|^2+R_{g_0}\phi^2)dV-8\pi^2\int_\Sigma v\Delta\rho dA \label{eq:E.wrt.g0}  \\ 
 &\;+8\pi^2\int_\Sigma[-q\mathrm{div}(\phi^2(\varphi_1d\varphi_2-\varphi_2d\varphi_1))-\rho|dq|^2 ] dA. \nonumber
\end{align}

\begin{prop} For any background metric $g_0\in\cD$, we have $\cI_{g_0}<\infty.$
\end{prop} 
\begin{proof}

To prove $\cI_{g_0}<\infty$, let $\rho$ be a first eigenfunction for $-\Delta_{h_0}$ on $\Sigma$, with Dirichlet boundary condition $\rho|_{\partial\Sigma}=0$, and denote the first eigenvalue $\lambda>0$. We then take the candidate $\phi:=\sqrt{\rho/\rho_0}$, which can be easily checked to be smooth on $M$. It suffices to establish an upper bound for   $\tilde E_{g_0}(\phi,v,q)$ uniformly in $(v,q)$.  

Below, we let $C$ denote any expressions that depend only on $\phi$ and the background $g_0$. For the above choice of $\phi$, 
 \eqref{eq:E.wrt.g0} becomes
\begin{align}\label{eq:E.first.efn}
 &E_{g_0}(\phi,v,q)\\
\nonumber =&\;C+8\pi^2\lambda\int_\Sigma v\rho +8\pi^2\int_\Sigma[-q\div(\phi^2(\varphi_1d\varphi_2-\varphi_2d\varphi_1))-\rho|dq|^2 ]\\
\nonumber=&\; C+8\pi^2\lambda\int_\Sigma v\rho +2\int_M[-q\frac{\div(\phi^2(\varphi_1d\varphi_2-\varphi_2d\varphi_1))}{\rho_0}-\phi^2|dq|^2 ]
\end{align}

In order to bound the integral term with divergence, observe that 
$$\frac{\div(\phi^2(\varphi_1d\varphi_2-\varphi_2d\varphi_1))}{\rho_0}=\frac{\phi^2\Delta\varphi_2}{\varphi_2}-\frac{\phi^2\Delta\varphi_1}{\varphi_1}+ 2{\phi}\left(\frac{\langle d\phi,d\varphi_2\rangle}{\varphi_2}-\frac{\langle d\phi,d\varphi_1\rangle}{\varphi_1}\right),$$
and each term on the right-hand side is smooth on $S$, using the symmetry of $\varphi_i$ across $\partial\Sigma$. 
Therefore, using the Cauchy-Schwarz inequality, we see that
\begin{align}\label{eq.q.div/rho0.bdd}
    &\left|\int_M   q  \frac{\div(\phi^2(\varphi_1d\varphi_2-\varphi_2d\varphi_1))}{\rho_0}\right|\\
\nonumber    &\leq \left(\int_M  q^2\right)^{1/2}\left(\int_M \left|\frac{\div(\phi^2(\varphi_1d\varphi_2-\varphi_2d\varphi_1))}{\rho_0}\right|^2\right)^{1/2}\\
\nonumber     &\leq K_\phi\left(\int_M q^2\right)^{1/2}
\end{align} 
for some constant $K_\phi>0$. Now, since $\int_Mq=0$ by the symmetry $q\circ r=-q$, we can  apply the Poincar\'e inequality to $\|q\|_{L^2}$, with the $L^2$ norm taken in $(M,g_0)$. Together with the fact that $\min\phi>0$,  we have
\begin{equation}\label{eq:q.L2.bdd}
    K_\phi\|q\|_{L^2} \leq C\|dq\|_{L^2}\leq C\|\phi dq\|_{L^2}\leq C+\frac 12\int_M\phi^2|dq|^2.
\end{equation}
Hence, putting \eqref{eq.q.div/rho0.bdd} and \eqref{eq:q.L2.bdd} back into \eqref{eq:E.first.efn}, we obtain 
\begin{equation}\label{eq:bdd.E.int.vrho}
    E_{g_0}(\phi,v,q)\leq   C+8\pi^2\lambda\int_\Sigma v\rho -\int_M \phi^2|dq|^2\leq   C+2 \lambda\int_M v\phi^2. 
\end{equation}

Now, we proceed to bound $\tilde E_{g_0}(\phi,v,q)$. We first rewrite the right-hand side as
$${C} +2 \lambda\int_M v\phi^2=2\lambda\int_M\Bigr(1+\bigl(v+ \frac {C}{2\lambda\int_M \phi^2}-1 \bigr) \Bigr)\phi^2.$$
This allows us to apply the elementary inequality $1+t\leq e^t$, followed by Cauchy--Schwarz inequality, to obtain
\begin{align*}
E_{g_0}(\phi,v,q)
&\leq 2\lambda \exp\Bigr(\frac {C}{2\lambda\int_M \phi^2}-1\Bigr) \int_Me^v\phi^2\\
&\leq 2\lambda \exp\Bigr(\frac {C}{2\lambda\int_M \phi^2}-1\Bigr) \vol(g_0)^{1/2}\left(\int_M e^{2v}\phi^4 \right)^{1/2}.
\end{align*}
Thus, $\tilde E_{g_0}(\phi,v,q)$ is bounded from above independently of $(v,q)$. \end{proof}

Another simple but important feature of $\cI_{g_0}$ is its independence from the choice of background $g_0\in \cD$. 
\begin{prop}\label{prop:Emin.indep.background}
For any  metrics   $g_0,g_1\in\cD$, we have $ \cI_{g_0}= \cI_{g_1}.$
\end{prop}
\begin{proof} Write $g_0=h_0+\varphi_1^2d\theta_1^2+\varphi_2^2d\theta_2^2$.
Using Corollary \ref{cor:guage.fix.M}, there exist some positive $\phi_1\in C^\infty(M)^G$,    $(v_1,q_1)\in \cC$, and    $\Phi\in \mathrm{Diff}^L_0(S)$ such that
$$g_1=\Phi^*(\phi_1^2(e^{2v_1}h_0+e^{2q_1}\varphi_1^2 d\theta_1^2+e^{-2q_1}\varphi_2^2 d\theta_2^2))$$

Take any $\phi\in W^{1,2}(M)^G$ and any   $(v,q)\in \cC$.
Since  total scalar curvature and volume are invariant under diffeomorphisms, it is easy to check that  
\begin{align*}
    \E_{{g}_1}((\phi\phi_1^{-1})\circ \Phi,v\circ \Phi,q\circ \Phi)&=\E_{{g}_0}(\phi ,v+v_1,q+q_1),\\ \tilde E_{{g}_1}((\phi\phi_1^{-1})\circ \Phi,v\circ \Phi,q\circ \Phi)&=\tilde E_{{g}_0}(\phi ,v+v_1,q+q_1).
\end{align*}
Taking supremum over $(v,q)$, we obtain 
$\cS_{g_1}((\phi\phi_1^{-1})\circ\Phi)=\cS_{g_0}(\phi)$. Taking infimum over $\phi$,  the proposition follows.
\end{proof}
 

For the usual Yamabe problem, for any metric $g$, the normalized total scalar curvature functional $\phi\mapsto \tilde \cE(\phi^2 g)$ may have non-unique critical points, and in particular non-unique minimizers (cf. \cite{viaclovsky2010monopole}). 
In contrast, the maximization problem for $\tilde E_{g_0}(\phi,v,q)$ over $(v,q)$ has a simpler structure, as the proposition below shows. In what follows, denote by $\partial_{\theta_i}$ the coordinate vector field of the $i$-th $S^1$ factor, for $i=1,2.$  
 
\begin{prop}\label{prop:fiber.ein.max} Take an arbitrary  triple  $(\phi,v,q)\in C^\infty(M,\R_{>0})^G\x\cC$, and write $g=\phi^2g_{v,q}$. The following conditions are equivalent.
\begin{enumerate}
    \item\label{item:fiber.ein}  $\Ric_g(\partial_{\theta_i},\partial_{\theta_i})=\lambda g(\partial_{\theta_i},\partial_{\theta_i})$ for $i=1,2$ and $\lambda=\frac 1{4\vol(g)}  \int_MR_g dV_g$.
    \item \label{item:crit} $(v,q)$ is a critical point for $\tilde E_{g_0}(\phi,\cdot,\cdot)$.
    \item \label{item:max} $(v,q)$ is the unique maximizer for    $\tilde E_{g_0}(\phi,\cdot,\cdot)$.
\end{enumerate}
In addition, if any item holds, then $\lambda>0$.
\end{prop} 
\begin{proof}
Fix an arbitrary $(\phi,v,q)$, and let $g=\phi^2 g_{v,q}$. Let $r_i:=\frac{\Ric_g(\partial_{\theta_i},\partial_{\theta_i})}{g(\partial_{\theta_i},\partial_{\theta_i})}$, $i=1,2$. Let us first show that \eqref{item:fiber.ein} implies $\lambda>0$. Namely, by \eqref{eq:fiberEin1} in \S \ref{sect:curv.formula}, the condition $r_1=r_2=\lambda$ implies    $\rho:=\phi^2\varphi_1\varphi_2$ is a first eigenfunction of $-\Delta_{\phi^2 e^{2v}h_0}$ on $\Sigma$ with Dirichlet boundary condition, and $\rho$ has eigenvalue $2\lambda$. Hence, $\lambda>0$.

To finish the proof, we just need to show   \eqref{item:fiber.ein}, \eqref{item:crit}, and \eqref{item:max} are equivalent.

Since $\cC$ is a vector space, any element of $\cC$ can be expressed as $(v+v_1,q+q_1)$.
Using \eqref{eq:E.in.g0}, and also  the curvature formulas in \S  \ref{sect:curv.formula},  we can compute
\begin{align*}
&E_{g_0}(\phi,v+v_1,q+q_1)\\
&=E_{g_0}(\phi,v,q)-2\int_M\phi^2\Delta v_1-2\int_M\phi^2(|dq_1|^2+\langle d\log(\varphi_1/\varphi_2)+2dq,dq_1\rangle)\\
&= E_{g_0}(\phi,v,q)-2\int_M\phi^2|dq_1|^2+2\int_M (v_1(r_1+r_2)-q_1(r_1-r_2))\phi^4 e^{2v}.
\end{align*}
Moreover, 
\begin{align}\label{eq:vol.new}
    \vol(\phi^2g_{v+v_1,q+q_1})=\int_M \phi^4 e^{2(v+v_1)}=\int_M e^{2v_1} dV_g.
\end{align}
As a consequence,  we can compute that 
$$\frac d{dt} \tilde E_{g_0}(\phi,v+tv_1,q+tq_1)|_{t=0}=\frac 2{\vol(g)^{1/2}}\int_M (v_1(r_1+r_2-2\lambda)-q_1(r_1-r_2))dV_g.$$
Now, the equivalence of  items \eqref{item:fiber.ein} and   \eqref{item:crit} is clear from the above, as $(v_1,q_1)\in\cC$ is arbitrary. Also,  \eqref{item:max} implies \eqref{item:crit} trivially, so it remains to show \eqref{item:fiber.ein} implies \eqref{item:max}.

Suppose $r_1=r_2=\lambda$. The above computation then gives
\begin{align*}
    E_{g_0}(\phi,v+v_1,q+q_1)&= 4\lambda\vol(g)+4\lambda\int_M  v_1\phi^4 e^{2v}-2\int_M\phi^2|dq_1|^2\\
    &\leq 4\lambda\int_M(1+v_1)dV_g\leq 4\lambda\int_M e^{v_1} dV_g\\
    &\leq 4\lambda \vol(g)^{1/2}\left(\int_M e^{2v_1} dV_g\right)^{1/2},
\end{align*}
in which the last two inequalities used the elementary fact $1+v_1\leq e^{v_1}$ and the Cauchy-Schwarz inequality respectively. Applying \eqref{eq:vol.new}, we obtain  
$$\tilde E_{g_0}(\phi,v+v_1,q+q_1)\leq \tilde E_{g_0}(\phi,v,q),$$
meaning $(v,q)$ is a maximizer. In addition, if equality holds, then necessarily $v_1=0$ and $dq_1=0$, which implies $q_1=0$ by the symmetry $q\circ r=-q.$ This means the maximizer $(v,q)$ is indeed unique, so \eqref{item:max} holds as desired. 
\end{proof}

The following proposition illustrates the relation between the min-max structure of $\tilde E_{g_0}$ and the Einstein condition, and shows that \emph{any} Einstein metric in $\cD$ must realize both $\cI_{g_0}$ and $\sigma_{\cD}$. 

\begin{prop} Take an arbitrary  triple  $(\phi,v,q)\in C^\infty(M,\R_{>0})^G\x\cC$, and write $g=\phi^2g_{v,q}$. If 
$\Ric_g=\lambda g$, then $\lambda>0$ and 
    $$\tilde E_{g_0}(\phi,v,q)=\cS_{g_0}(\phi)=\cI=\sigma_\cD.$$
\end{prop} 
\begin{proof}   Suppose $\Ric_g=\lambda g$. Then immediately from Proposition \ref{prop:fiber.ein.max}, we know   $\lambda>0$. In fact, using  the implication \eqref{item:fiber.ein} $\Rightarrow$ \eqref{item:max} therein, and then  Obata's theorem \cite{obata1971conformal,schoen1987variational}, we have  
    $$\cS_{g_0}(\phi)=\tilde E_{g_0}(\phi,v,q)=\cY_\cD(g_{v,q}).$$
So for any  $\phi_1\in W^{1,2}(M)^G$,
$$\cS_{g_0}(\phi)=\tilde E_{g_0}(\phi,v,q) \leq \tilde E_{g_0}(\phi_1,v,q)\leq \cS_{g_0}(\phi_1).$$
Taking infimum over $\phi_1$, we obtain  $\cS_{g_0}(\phi)=\tilde E_{g_0}(\phi,v,q)=\cI$. Finally, as we have shown $\cY_\cD(g_{v,q})=\cI$,  we   immediately obtain  $\sigma_\cD=\cI$, since $\sigma_\cD\leq \cI$ by Proposition \ref{prop:sigmaD.leq.Emin}.  This completes the proof.
\end{proof}

\section{Existence of the minimizer}\label{sect:existenceMin}
The goal of this section is to prove the existence of $\phi\in W^{1,2}(M)^G$ minimizing $\cS_{g_0}$, via the direct method. The analysis turns out to be quite similar to that involved in the standard Yamabe or Hebey-Vaugon problem: As we explain below, it is easy to establish lower-semicontinuity of $\cS_{g_0}$ along sequences which are bounded in $W^{1,2}(M)$ and convergent in $L^4(M)$, so the main challenge is to rule out bubbling for a minimizing sequence. By essentially the same concentration-compactness analysis used in the Hebey-Vaugon problem \cite{HebeyVaugon1993Equivariant}, this allows us to reduce the existence question to a single inequality, as follows.

\subsection{Conditional existence of minimizer}\label{sect: cond.exis}
Let $Y_0$ denote the Yamabe invariant of the standard  metric on $S^4$. In this section we prove:
\begin{prop}\label{prop:exist.minimizer.conditional} Assume that $\cI<\sqrt{2n+2}Y_0$. Then for any  $g_0\in\cD$ there exists some minimizer $\phi\in W^{1,2}(M)^G$ such that
$\cS_{g_0}(\phi)=\cI.$
\end{prop}
Let $\phi_j$ be a minimizing sequence for $\cS_{g_0}$ in $W^{1,2}(M)^G$ so that $\cS_{g_0}(\phi_j)$ approaches $ \inf\cS_{g_0}=\cI$. We normalize them such that $\|\phi_j\|_{L^4(M,g_0)}=1.$ Hence, putting $v=q=0$ into \eqref{eq:E.in.g0}, we have
$$\int_M(6|d\phi_j| ^2+R_{g_0}\phi_j^2)dV =\tilde E_{g_0}(\phi_j,0,0)\left(\int_M\phi^4\right)^{1/2}\leq\cS_{g_0}(\phi_j).$$
Note the integrals taken are with respect to the measure $dV_{g_0}$.
Therefore, $\{\phi_j\}$ is a bounded sequence in $W^{1,2}$, and by passing to a subsequence we can assume $\phi_j$ converges weakly in $W^{1,2}$ to some $\phi.$ Our goal is to prove $\cS_{g_0}(\phi)=\cI$.

By  the concentration-compactness in Lions' \cite[Lemma~I.1]{Lions1985LimitI} (see also \cite[\S 1 Theorem 9]{evans1990weakConv}), after passing to a subsequence, we have a weak$^*$-convergence of measures
$$\phi_j^4dV_{g_0}\rightharpoonup^*  \phi^4 dV_{g_0}+\sum_{k=1}^\infty\alpha_k\delta_{p_k},$$
for some $\alpha_k\geq 0$ and a sequence $p_1,p_2,\ldots\in M$. In fact, since $\phi_j\in W^{1,2}(M)^G$, the only option is 
\begin{equation}\label{eq:weak*}
    \phi_j^4dV_{g_0}\rightharpoonup^*  \phi^4 dV_{g_0}+\alpha\sum_{p\in S_1\cap S_2} \delta_{p};
\end{equation}
recall from \S \ref{sect:topo.mfd} that $S_1\cap S_2\subset M$ denotes the fixed point set of the $T^2$-action, which consists of $2n+2$ points. Moreover, by Rellich compactness, boundedness of $\{\phi_j\}$ in $W^{1,2}$ implies a strong $L^2$-convergence $\phi_j\to\phi.$

Now, we take any arbitrary admissible $(v,q)$, and study the limit of $ \tilde E_{g_0}(\phi_j,v,q)$. By the boundedness of $\{\phi_j\}$ in $W^{1,2}$, we can assume the Dirichlet energy $\int_M|d\phi_j|^2$ converges by passing to a subsequence. Consider the energy loss
$$D:=\lim_{j\to\infty}6\int_M(|d\phi_j|^2-|d\phi|^2).$$
\begin{lem}\label{lem:lim.E.vol}  As $j\to\infty$, we have the following.
\begin{enumerate}
    \item\label{item:lim.E} $\lim_j E_{g_0}(\phi_j,v,q)= E(\phi,v,q)+D$.
    \item\label{item:lim.phi_j^4} $\lim_j\int_M\phi_j^4 e^{2v} = \int_M\phi^4 e^{2v}+({2n+2})\alpha.$
\end{enumerate}
\end{lem}
\begin{proof}
First, noting the relations
$dV_{g_{v,q}}=e^{2v}dV_{g_0}$ and $|d\phi|^2_{g_{v,q}}=e^{-2v}|d\phi|^2_{g_0}$, 
we  express 
$$E(\phi,v,q)=\int_{M}(6|d\phi|_{g_0}^2+ R_{g_{v,q}}e^{2v}\phi^2)dV_{g_0}.$$
We can analogously rewrite $E(\phi_j,v,q)$, and see 
$$E_{g_0}(\phi_j,v,q)-E_{g_0}(\phi,v,q)=6\int_M(|d\phi_j|^2-|d\phi|^2 ) +\int_M R_{g_{v,q}} e^{2v}(\phi_j^2-\phi^2).$$
Taking $j\to\infty$, the first integral goes to $D$, while the second goes to zero by the $L^2$-convergence $\phi_j\to\phi$. So item \eqref{item:lim.E} follows.

Integrating \eqref{eq:weak*} against $e^{2v}$, we have
$$\lim_j\int_M\phi_j^4e^{2v}= \int_M\phi^4 e^{2v}+\alpha \sum_{p\in S_1\cap S_2} e^{2v(p)}.$$
By the admissibility of $(v,q)$, at every $p\in S_1\cap S_2$ we have $v(p)=q(p)=-q(p)$, which means $v(p)=0$. So we obtain item \eqref{item:lim.phi_j^4} too.
\end{proof}

Next, following the estimate for  concentration-compactness  in  \cite[Lemma~I.1]{Lions1985LimitI},  we can bound the energy loss $D$ from below in terms of $\alpha$.
\begin{lem}\label{lem:energy.loss}
$D\geq (2n+2)\sqrt{\alpha} Y_0.$
\end{lem}
\begin{proof}
For simplicity, denote by $F$ the fixed point set $S_1\cap S_2\subset M$ of the $T^2$-action. Fix an arbitrary $\epsilon>0$. Since $|d\phi|^2, \phi^4\in L^1(M,g_0)$, we can find $\delta=\delta({\epsilon})\in (0,\epsilon)$ such that 
$$\int_{B_{\delta}(F)}(|d\phi|^2+\phi^4) <\epsilon.$$
Let $\sigma(\delta)$ be the   Sobolev constant
$$\sigma(\delta):=\inf_{0\ne f\in W_0^{1,2}(B_{\delta}(p))}\frac{\|df\|_{L^2}^2}{\|f\|_{L^4}^2}$$
in the Sobolev embedding $W_0^{1,2}(B_{\delta}(p))\to L^4(B_{\delta}(p))$. It is elementary to check that as $\delta\to 0$, $\sigma(\delta)$ tends to the   Sobolev constant on $\R^4$
$$\inf_{0\ne f\in W_0^{1,2}(\R^4)}\frac{\|df\|_{L^2}^2}{\|f\|_{L^4}^2},$$
which is known to be  $Y_0/{6}$ using stereographic projection and the conformal invariance of the Yamabe quotient \cite{Aubin1976Sobolev, Talenti1976Sobolev}. Therefore,  $\sigma(\delta)\to {Y_0}/{6}$ as $\delta\to 0$. 

For each $p\in F$, choose a cut-off $\chi_p\in C^\infty_c(B_\delta(p))$ with
$$0\leq \chi_p\leq 1,\quad \chi_p=1 \textrm{ on } B_{\delta/2}(p), \quad |d\chi_p|\leq C(g_0)/\delta.$$
Let $\chi=\sum_{p\in F}\chi_p.$
Applying the Sobolev inequality to each $\chi_p\phi_j$, we have
\begin{align}\label{eq:sigma.chi.phi}
 \sigma(\delta)\sum_{p\in F}\|\chi_p \phi_j\|^2_{L^4}&\leq \sum_{p\in F}\int_M|d(\chi_p\phi_j)|^2\\
\nonumber&=\int_M \chi^2|d\phi_j|^2+2\int_M\chi \phi_j\langle d\phi_j,d\chi\rangle+\int_M\phi_j^2|d\chi|^2\\
\nonumber &\leq \int_M\chi^2|d\phi_j|^2+2 \|d\phi_j\|_{L^2}\|\phi_jd\chi\|_{L^2} +\| \phi_jd\chi\|_{L^2}^2,
\end{align}
in which the last inequality used the Cauchy-Schwarz inequality.

Now, keeping $\epsilon,\delta$ fixed, we let $j\to\infty$ and examine each term. Below, $C$ denotes any constant that depends only on $n$ and $g_0$.
\begin{itemize}
\item For the left-hand side, as in the proof of Lemma \ref{lem:lim.E.vol} \eqref{item:lim.phi_j^4}, we have $$\lim_j \|\chi_p \phi_j\|^2_{L^4}= \left(\int_M(\chi_p \phi)^4 +  \alpha\right)^{1/2}\geq \sqrt\alpha.$$
\item For the first term on the right-hand side, by   the lower semi-continuity
$$\liminf_j\int_M(1-\chi^2)|d\phi_j|^2\geq \int_M(1-\chi^2)|d\phi|^2,$$
we have
$$\limsup_j\int_M\chi^2|d\phi_j|^2\leq \frac D6+\int_M\chi^2|d\phi|^2\leq \frac D6+\epsilon.$$
\item For the second term,  we observe
$$\lim_j \|\phi_j d\chi\|_{L^2} =\|\phi d\chi\|_{L^2}  \leq \frac{C}\delta\vol(B_\delta(F))^{1/4}\|\phi \|_{L^4(B_\delta(F))}\leq C\epsilon^{1/4}.$$
Moreover,  $\|d\phi_j\|_{L^2} $  is bounded independently of $j$, as $\{\phi_j\}$ is bounded in $W^{1,2}$.
\item For the third term, $ \lim_j\|\phi_j d\chi\|_{L^2}^2\leq    C\epsilon^{1/2}.$
\end{itemize}
Applying these to \eqref{eq:sigma.chi.phi}, we obtain 
$$\sigma(\delta)(2n+2)\sqrt\alpha\leq D/6+\epsilon+C\epsilon^{1/4}+C\epsilon^{1/2}.$$
Finally, we let $\epsilon\to 0$, which means $\delta\to 0$. Together with the previous fact that $\sigma(\delta)\to Y_0/6$ as $\delta\to 0$,
we conclude   $Y_0(2n+2)\sqrt\alpha\leq D$.\end{proof}

Now, from definition,
$$E_{g_0}(\phi_j,v,q)\leq \cS_{g_0}(\phi_j)\left(\int_M\phi_j^4e^{2v}\right)^{1/2}.$$
Take $j\to\infty$.
Using Lemma \ref{lem:lim.E.vol} and the fact that $\{\phi_j\}$ is a minimizing sequence,
\begin{equation}\label{eq:E+D}
    E_{g_0}(\phi,v,q)+D\leq \cI\left(\int_M\phi^4 e^{2v}+({2n+2})\alpha\right)^{1/2}.
\end{equation}

{\it At this point, we start using the assumption  $\cI<\sqrt{2n+2}Y_0$.} We first argue that $\phi$ is not a.e. zero. Suppose for contradiction $\phi\equiv 0$. Then by putting $v=0$ into Lemma \ref{lem:lim.E.vol} \eqref{item:lim.phi_j^4} and using the normalization $\|\phi_j\|_{L^4}=1$, we obtain  $1=(2n+2)\alpha$. Thus, Lemma \ref{lem:energy.loss} becomes $D\geq \sqrt{2n+2} Y_0$. On the other hand, applying $\phi\equiv 0$ and $\cI<\sqrt{2n+2}Y_0$ to \eqref{eq:E+D} gives $$D\leq \cI<\sqrt{2n+2}Y_0.$$ 
A contradiction arises. Therefore $\phi$ is not a.e. zero. 

Next, we treat the two cases $\mathcal I\leq 0$ and $\mathcal I>0$ separately. In the first case, because $D\geq 0$,  \eqref{eq:E+D}  immediately implies
$$E_{g_0}(\phi,v,q)\leq \cI\left(\int_M\phi^4 e^{2v}\right)^{1/2},$$
so that $\tilde E_{g_0}(\phi,v,q)\leq\cI$. And since $(v,q)$ is arbitrary, this gives $\cS_{g_0}(\phi)\leq\cI$, which implies the desired claim $\cS_{g_0}(\phi)= \cI$.

If instead $\mathcal I>0$, we apply Lemma \ref{lem:energy.loss} to \eqref{eq:E+D} and then use the fact that the square root function is concave:  
\begin{align}
\nonumber    E_{g_0}(\phi,v,q)&\leq\cI\left(\int_M\phi^4 e^{2v}+({2n+2})\alpha\right)^{1/2}-(2n+2)\sqrt{\alpha} Y_0\\
\label{eq:E.upper.bdd}    &\leq \cI\left(\int_M\phi^4 e^{2v}\right)^{1/2}+(\cI-\sqrt{2n+2}Y_0)\sqrt{(2n+2)\alpha}.
\end{align}
Hence, from the condition $\cI<\sqrt{2n+2}Y_0$,  
we  immediately deduce $\tilde E_{g_0}(\phi,v,q)\leq \cI$, as desired.  This proves Proposition \ref{prop:exist.minimizer.conditional}.

\subsection{Upper bound on $\cI$} 
\begin{prop}\label{prop:bound.I}
For each $n\geq 1$, $\cI<\sqrt{n+1}Y_0$.
\end{prop}
We remark that $Y_0=8\sqrt 6\pi.$
This proposition, together with Proposition \ref{prop:exist.minimizer.conditional}, immediately gives:
\begin{thm}\label{wk.ex}
For any  $g_0\in\cD$ there exists some minimizer $\phi\in W^{1,2}(M)^G$ such that
$\cS_{g_0}(\phi)=\cI.$
\end{thm}

\begin{proof}[Proof of Proposition \ref{prop:bound.I}]
Consider the bi-gon $$\Sigma_0:=\{(x,y,z)\in S^2:x,y\geq 0\}$$ in the unit 2-sphere $S^2$, and let $h_{\mathrm{std}}$ be the induced standard metric on $\Sigma_0$; recall from Example \ref{example:S4} that 
\begin{equation}\label{eq:round.S4.metric}
g_{\mathrm{std}}=h_{\mathrm{std}}+x^2d\theta_1^2+y^2d\theta_2^2    
\end{equation}
gives the standard metric on $S^4$.

Let $F:\Sigma\to\Sigma_0$ be the $(n+1)$-fold covering map with a single branch point $p\in\Sigma$ such that $F(p)=(\frac 1{\sqrt 2},\frac 1{\sqrt 2},0)\in\Sigma_0$. Consider the pull-back metric $F^*h_{\mathrm{std}}$ on $\Sigma\backslash\{p\}$. Under this metric, $\Sigma$ admits a dihedral symmetry $D_{4n+4}$, and has  geodesic edges and right-angled corners. We choose a smooth $L$-invariant metric $h_1$ on $\Sigma$ that lies in the same conformal class as $F^*h_{\mathrm{std}}$ on $\Sigma\backslash\{p\}$, and  such that $h_1=F^*h_{\mathrm{std}}$ on some neighborhood around $\partial\Sigma$. This way, we can write $F^*h_{\mathrm{std}}=e^{2w}h_1$ for some function $w\in C^\infty(\Sigma\backslash \{p\})$ with $w=0$ near $\partial\Sigma$. Now, define on $M$ the smooth metric
$$g_1:=h_1+ (x\circ F)^2d\theta_1^2+(y\circ F)^2d\theta_2^2\in\cD.$$
As before, we denote $\varphi_1:=x\circ F,\varphi_2:=y\circ F,$ and $\rho_1=\varphi_1\varphi_2$ for simplicity.
To prove Proposition \ref{prop:bound.I}, we are going to show that $\cS_{g_1}(1)<\sqrt nY_0$.

Let $(v,q)\in\cC$ be arbitrary. Using \eqref{eq:E.wrt.g0},   \begin{align} 
E_{g_1}(1,v,q)=&\;4\pi^2\int_\Sigma R_{g_1} \rho_1 dA_{h_1}-8\pi^2\int_\Sigma v\Delta_{h_1}\rho_1 dA_{h_1} \label{eq:E.1.wrt.g0}  \\ 
 &\;+8\pi^2\int_\Sigma[-q\mathrm{div}_{h_1}(\varphi_1d\varphi_2-\varphi_2d\varphi_1)-\rho_1|dq|^2_{h_1} ] dA_{h_1}. \nonumber
\end{align}
As before, to ease notation, we will often omit subscripts of $h_1$.
Now, we use the fact that  the pullback metric $F^*g_{\mathrm{std}}$ is smooth and Einstein over $\Sigma\backslash\{p\}$.  
\begin{itemize}
\item By \eqref{eq:fiberEin1}, $-\Delta_{F^*h_\mathrm{std}}\rho_1=6\rho_1$, so
\begin{equation}\label{eq:Delta.h1.rho1}
    -\Delta\rho_1=6e^{2w}\rho_1.
\end{equation}
\item By \eqref{eq:fiberEin2}, $\div_{F^*h_\mathrm{std}}(\varphi_1d\varphi_2-\varphi_2d\varphi_1)=0$, so \begin{equation}\label{eq:div.h1.rho1}
    \div(\varphi_1d\varphi_2-\varphi_2d\varphi_1)=0.
\end{equation}
\item  Applying \eqref{eq:R_g} to calculate the scalar curvature of $F^*g_{\textrm{std}}$ and $g_1$ respectively, we derive \begin{align*}
    R_{F^*g_{\textrm{std}}}&=2 e^{-2w}\left(K_{h_1}-\Delta w-\frac{\Delta \varphi_1}{\varphi_1}-\frac{\Delta \varphi_2}{\varphi_2}-\frac{\langle d\varphi_1,d\varphi_2\rangle }{\varphi_1\varphi_2}\right)\\
    &=e^{-2w}(R_{g_1}-2\Delta w).
\end{align*} 
Hence, since  $R_{F^*g_{\textrm{std}}}=R_{g_{\textrm{std}}}=12$, we have the pointwise identity $$R_{g_1}=12 e^{2w}+2\Delta w.$$
But due to the logarithmic singularity of $w$ at $p$, we should also consider the  distributional identity
\begin{equation}\label{eq.R1.dAh1}
    R_{g_1}dA_{h_1}=(12e^{2w}+2\Delta  w) dA_{h_1}-4\pi n\delta_p.
\end{equation} 
\end{itemize}
Putting \eqref{eq:Delta.h1.rho1},  \eqref{eq:div.h1.rho1}, and \eqref{eq.R1.dAh1} back into \eqref{eq:E.1.wrt.g0}, we have
\begin{align}
\nonumber E_{g_1}(1,v,q)
 &=4\pi^2\int_\Sigma  (12e^{2w}+2\Delta  w)\rho_1  -4\pi^2\cdot 4\pi n\rho_1(p)+48\pi^2\int_\Sigma ve^{2w}\rho_1\\
\nonumber&\quad -8\pi^2\int_\Sigma\rho_1|dq|^2\\
&\leq 8\pi^2 \int_\Sigma \rho_1\Delta w+48\pi^2\int_\Sigma(1+v)e^{2w}\rho_1-16\pi^3n\rho_1(p).\label{eq:E.g1.1.v.q.ineq}
\end{align}

Next, we deal with the $\int_\Sigma \rho_1\Delta w$ term. By Green's formula,
$$\int_\Sigma(\rho_1\Delta w-w\Delta \rho_1) =\int_{\partial\Sigma}(\rho_1\partial_\nu w-w\partial_\nu\rho_1)  .$$
Together with \eqref{eq:Delta.h1.rho1} and  the fact that $w=0$ near $\partial\Sigma$, we have
\begin{equation*}
    \int_\Sigma\rho_1\Delta  w=-6\int_\Sigma w e^{2w}\rho_1.
\end{equation*}
Putting this back into \eqref{eq:E.g1.1.v.q.ineq}, and expressing the integrals over $M$, 
\begin{align} 
\nonumber E_{g_1}(1,v,q)&\leq 12\int_M (1+v-w)e^{2w}-16\pi^3n\rho_1(p)\\
&= 12\int_M (1+v-w)e^{2w}-8\pi^3n.\label{eq:E.g1.1.v.q.ineq2}
\end{align}
The last equality used that $\rho_1(p)=(xy)\circ F(p)=(\frac 1{\sqrt 2})^2=\frac 12.$

Now, recalling that $F$ induces an $(n+1)$-fold branched covering map $M\to S^4$, we see that $\int_M e^{2w}$ is just $n+1$ times the volume of $(S^4,g_{\mathrm{std}})$. Namely,  $\int_M e^{2w}=\frac 83\pi^2(n+1).$ Thus, we can rewrite \eqref{eq:E.g1.1.v.q.ineq2} as
\begin{align*}
E_{g_1}(1,v,q)&\leq12 \int_M  \left (1+v-w-\frac{\pi n}{4(n+1)}\right)e^{2w}.
\end{align*}
Applying the elementary inequality $1+t\leq e^t$ with $t=v-w-\frac{\pi n}{4(n+1)}$ and then the Cauchy-Schwarz inequality, 
\begin{align*}
E_{g_1}(1,v,q)&\leq12 \int_M   e^{v+w-\frac{\pi n}{4(n+1)}}\leq 12 e^{-\frac{\pi n}{4(n+1)}}\left(\int_Me^{2v}\right)^{1/2}\left(\int_Me^{2w}\right)^{1/2}.
\end{align*}
Together with $\int_M e^{2w}=\frac 83\pi^2(n+1)$ and the fact that $Y_0=8\sqrt 6\pi$, 
$$\tilde E_{g_1}(1,v,q)\leq e^{-\frac{\pi n}{4(n+1)}}\sqrt{n+1} Y_0.$$
Taking supremum over all $(v,q)$, we obtain the desired bound. \end{proof}

\section{A priori estimates and weak Euler--Lagrange equations}\label{sect:weak.E-L} 
\subsection{Four reference spaces and some preliminary estimates}

As above, we equip $M=\# n(S^2\times S^2)$ with a fixed, but arbitrary reference metric
$$g_0=h_0+\varphi_1^2d\theta_1^2+\varphi_2^2d\theta_2^2\in \mathcal{D}$$
in the desired symmetry class. Recall that the quotient by the given $T^2$ action yields a metric $h_0$ on the $2(n+1)$-gon $\Sigma$, which extends to a smooth $L$-invariant metric on the closed surface $S$ of genus $n$. Denote the function
$$\rho_0:=\varphi_1\varphi_2$$
on $\Sigma$, so that the $g_0$-area of the $T^2$ fiber above each point is given by $4\pi^2\rho_0$, and note that it extends smoothly to the closed surface $S$ as an odd function with respect to the reflections $R_1$ and $R_2$.

For our analysis, it will be useful to carry out some estimates with respect to additional reference geometries of dimension three and, perhaps less obviously, dimension six, defined as follows. 

For $i\in \{1,2\}$, we denote by $Y_i$ the three-dimensional manifold obtained as the quotient of $M$ by the $S^1$ action whose fibers are parametrized by $\theta_i$, equipped with the corresponding quotient metrics
$$g_0^1=h_0+\varphi_2^2d\theta_2^2\text{ on }Y_1\text{ and }g_0^2=h_0+\varphi_1^2d\theta_1^2\text{ on }Y_2.$$
Then $Y_i$ is a $3$-manifold with smooth boundary $\{\varphi_i=0\}$ consisting of $n+1$ disjoint two-spheres. By the symmetry assumption $g_0\in \mathcal{D}$, it is easy to see that $Y_1$ and $Y_2$ are isometric.

\begin{example}
If $n=1$ and $(M,g_0)$ is the standard product of unit spheres $S^2_1\times S^2_1$, then we have $Y_1=[-1,1]\times S^2_1$ and $Y_2=S^2_1\times [-1,1]$. 
\end{example}

At a few key steps in the proof of regularity for the Einstein metric associated to the $\cS$-minimizer $\phi$, it is also useful to work with the $6$-dimensional $O(3)\times O(3)$-symmetric reference metric
$$G_0:=h_0+\varphi_1^2g_{S^2_1}+\varphi_2^2g_{S^2_1}$$
on the $6$-dimensional manifold $P=(\Sigma\times S^2\times S^2)/\sim$, where $g_{S^2_1}$ is the standard unit-sphere metric, and $(x,p,q)\sim (x,p',q')$ when $\varphi_1(x)=0$ and $q=q'$, $\varphi_2(x)=0$ and $p=p'$, or $\varphi_1(x)=\varphi_2(x)=0$ with $p,p',q,q'$ arbitrary. It is straightforward to check that the conditions which ensure smoothness of $g_0$ on $M$ also ensure smoothness of $G_0$ on $P$.

\begin{example}
If $n=1$ and $(M,g_0)\cong S^2_1\times S^2_1$, then $(P,G_0)\cong S^3_1\times S^3_1$ is the product of standard unit $3$-spheres.
\end{example}

In many of the arguments that follow, we will need to pass estimates for $T^2$-invariant (or $O(3)\times O(3)$-invariant) functions between these different reference models. For now, we begin with the following simple lemma.

\begin{lem}\label{mtos.w11}
    For every nonnegative, $T^2$-invariant function $\chi\in W^{1,1}(M)$, we have $\chi\rho_0\in W_0^{1,1}(\Sigma)$, with
    \begin{equation}
        \|\chi\rho_0\|_{W^{1,1}(\Sigma)}\leq C\|\chi\|_{W^{1,1}(M)}.
    \end{equation}
    In particular, $\chi\rho_0$ extends to a function $\chi\rho_0\in W^{1,1}(S)$ which is odd with respect to $R_1$ and $R_2$.
\end{lem}
\begin{proof}
    It suffices to prove the estimate in the case $\chi$ is smooth, and conclude by density of smooth functions in $W^{1,1}(M)$. At the level of $L^1$ estimates, we have the simple equality
    $$\|\chi\|_{L^1(M,g_0)}=4\pi^2\int_{\Sigma}\chi \rho_0 dA_{h_0}=4\pi^2\|\chi\rho_0\|_{L^1(\Sigma,h_0)},$$
    so we only need to obtain the desired control for $\|d(\chi\rho_0)\|_{L^1({\Sigma})}$. Moreover, writing $d(\chi\rho_0)=\rho_0d\chi+\chi d\rho_0$, we clearly have
    $$4\pi^2\|\rho_0|d\chi|\|_{L^1(\Sigma)}=\|d\chi\|_{L^1(M)}\leq \|\chi\|_{W^{1,1}(M)},$$
    so all that remains is to bound $\|\chi d\rho_0\|_{L^1(\Sigma)}$ by $\|\chi\|_{W^{1,1}(M)}$. 

    To this end, write
    \begin{eqnarray*}
        \int_{\Sigma}\chi \varphi_1|d\varphi_2|^2&=&\int_{\Sigma}\chi [div_{h_0}(\rho_0d\varphi_2)-\varphi_2div_{h_0}(\varphi_1d\varphi_2)]dA_{h_0}\\
        &=&-\int_{\Sigma}\langle d\chi,\rho_0d\varphi_2\rangle-\varphi_2\chi div_{h_0}(\varphi_1d\varphi_2).
    \end{eqnarray*}
    Since $div_{h_0}(\varphi_1d\varphi_2)=\varphi_1\Delta_{h_0}\varphi_2+\langle d\varphi_1,d\varphi_2\rangle_{h_0}$ is odd with respect to both $R_1$ and $R_2$ on $S$, it must be bounded by $C\rho_0$ for some constant $C$, and it follows that
    $$\int_{\Sigma}\chi\varphi_1|d\varphi_2|^2\leq C\int_{\Sigma}\rho_0(|d\chi|+\chi).$$
    In particular, by construction, $\varphi_2$ satisfies $|d\varphi_2|=1$ on $\{\varphi_2=0\}$ and $|d\varphi_2|\geq \frac{1}{2}$ on $\{\varphi_2<\delta_0\}$ for some fixed $\delta_0>0$. Thus, we can estimate
    \begin{eqnarray*}
        \int_{\Sigma}\chi \varphi_1&\leq &C_0\int_{\varphi_2<\delta_0}\chi\varphi_1|d\varphi_2|^2+\int_{\varphi_2>\delta_0}\chi\varphi_1\\
        &\leq & C'\int_{\Sigma}\rho_0(|d\chi|+\chi)+\frac{1}{\delta_0}\int_{\Sigma}\chi\rho_0.
    \end{eqnarray*}
    This gives $\int_{\Sigma}\chi\varphi_1\leq C\|\chi\|_{W^{1,1}(M)}$, and symmetrically $\int_{\Sigma}\chi \varphi_2\leq C\|\chi\|_{W^{1,1}(M)}$ as well. Since $|d\rho_0|\leq C(\varphi_1+\varphi_2)$, this gives the desired estimate.
\end{proof}

\subsection{Consequences of $\cS(\phi)<\infty$}\label{s.fin.sec}

To begin analyzing the regularity of the minimizing $\phi\in W^{1,2}(M)^G$ supplied by Theorem \ref{wk.ex}, we first record some general properties that hold for any $\phi\in W^{1,2}(M)^G$ with $\cS_{g_0}(\phi)<\infty$. In what follows, let $0\leq \phi\in W^{1,2}(M)^G$ satisfy $\cS(\phi)<\infty$, and let 
\begin{equation}\label{rho.def}
\rho:=\phi^2\rho_0
\end{equation}
be the associated fiber area function.

Recall from Section \ref{sect:variation.Ein.Hil} that
$$\mathcal{S}_{g_0}(\phi):=\sup\{\tilde{E}_{g_0}(\phi,v,q)\mid (v,q)\in \mathcal{C}\},$$
where
$$\tilde{E}_{g_0}(\phi,v,q):=\frac{E_{g_0}(\phi,v,q)}{\left(\int_M\phi^4 dV_{g_{v,q}}\right)^{1/2}},$$
and we interpret the numerator in terms of the expansion \eqref{eq:E.in.g0}, so that
\begin{equation}
    E_{g_0}(\phi,v,q):=E_{g_0}(\phi,0,0)+8\pi^2[-\int_{\Sigma}\rho\Delta_{h_0}v dA_{h_0}+T(\phi,q)],
\end{equation}
where
\begin{equation}\label{tq.def}
T(\phi,q):=-\int_\Sigma\rho[|dq|_{h_0}^2+\langle d\log(\varphi_1/\varphi_2),dq\rangle_{h_0}]dA_{h_0}.
\end{equation}
In particular, noting that $dV_{g_{v,q}}=e^{2v}dV_{g_0}$, it follows that $\cS(\phi)<\infty$ if and only if
\begin{equation}\label{fin.char}
-\int_{\Sigma}\rho \Delta_{h_0}v+T(\phi,q)\leq C_{\phi}\left(1+\int_{\Sigma}\phi^4 e^{2v}\rho_0dA_{h_0}\right)
\end{equation}
for every $(v,q)\in \mathcal{C}$, for some $C_{\phi}<\infty$. 

We derive some initial consequences of \eqref{fin.char} by testing the estimate in the $q$ and $v$ variables separately, starting with the $v$ variable. Namely, fixing $q=0$ in \eqref{fin.char} and applying the estimate to all pairs in $\mathcal{C}$ of the form $(v,0)$ gives the following.

\begin{lem}\label{rho.lap.lem}
    Let $0\leq \phi\in W^{1,2}(M)^G$ be nonzero, satisfying $\cS_{g_0}(\phi)<\infty$. Then $\rho\in W^{1,p}(S)$ for every $p\in [1,2)$, and there exists a nonnegative function $f_{\phi}\in L^1(M,\phi^4dV_{g_0})$ such that the weak Laplacian $-\Delta_{h_0}\rho$ is given by
    $$-\Delta_{h_0}\rho=f_{\phi}\phi^4\rho_0\text{ on }S,$$
    and 
    $$\int_M f_{\phi}|\log f_{\phi}|\phi^4 dV_{g_0}<\infty.$$
\end{lem}
\begin{proof}
    To begin, note that since $\phi\in W^{1,2}(M)$, we have
    $$\int_M |d(\phi^2)|_{g_0}dV_{g_0}\leq 2\|\phi\|_{L^2(M)}\|d\phi\|_{L^2(M)}<\infty;$$
    hence, $\phi^2\in W^{1,1}(M)$, and it follows from Lemma \ref{mtos.w11} that $\rho=\phi^2\rho_0\in W_0^{1,1}(\Sigma)$, and extends naturally to a function $\rho\in W^{1,1}(S)$ which is odd with respect to $R_1$ and $R_2$.
    
    Next, it follows from \eqref{fin.char} that 
    $$-\int_{\Sigma}\rho\Delta_{h_0}v\leq C_{\phi}\left(1+\int_{\Sigma}\phi^4 e^{2v}\rho_0 dA_{h_0}\right)$$
    whenever $(v,0)\in \mathcal{C}$, i.e. for any $L$-invariant $v\in C^{\infty}(S)$ such that $v|_{\partial\Sigma}=0$. Moreover, since $\rho\in W_0^{1,1}(\Sigma)$ and $\rho$ is dihedrally symmetric in $\Sigma$, it follows that
    \begin{equation}\label{sig.pair}
    \int_{\Sigma}\langle d\rho,dv\rangle_{h_0}dA_{h_0}\leq C_{\phi}\left(1+\int_{\Sigma}\phi^4 e^{2v}\rho_0 dA_{h_0}\right)
    \end{equation}
    for any Lipschitz function $v\in Lip(\Sigma)$ such that $v|_{\partial\Sigma}=0$. As a consequence, we deduce that
    \begin{equation}\label{c0.star}
    |\langle d\rho,dv\rangle_{L^2(S)}|\leq C\|v\|_{C^0}
    \end{equation}
    for any $v\in Lip(S)$ such that $v|_{\partial\Sigma}=0$, and since the symmetries of $\rho$ force
    $$\langle d\rho,dv\rangle_{L^2(S)}=\langle d\rho,d\bar{v}\rangle_{L^2(S)},$$
    where $\bar{v}=\frac{1}{4}(v-v\circ R_1-v\circ R_2+v\circ R_1R_2)$ satisfies $\bar{v}|_{\partial\Sigma}=0$, it follows that \eqref{c0.star} holds for all $v\in Lip(S)$. In particular, 
    $$\mu:=-\Delta_{h_0}\rho\in C^0(S)^*$$
    defines a signed Radon measure on $S$, odd with respect to the actions of $R_1$ and $R_2$, satisfying
    \begin{equation}\label{mu.bd}
    \int_{\Sigma} v d\mu\leq C_{\phi}\left(1+\int_{\Sigma}\phi^4 e^{2v}\rho_0 dA_{h_0}\right)
    \end{equation}
    for every $v\in C^0(\Sigma)$, which already suffices to show that $\rho\in W^{1,p}(S)$ for every $p\in [1,2)$.

    For any fixed, nonnegative $0\leq \chi\in C^0(S)$, testing \eqref{mu.bd} against $v=-t\chi\leq 0$ with $t>0$ yields
    $$-\int_{\Sigma}\chi d\mu\leq \frac{C_{\phi}}{t}\left(1+\int_{\Sigma}\phi^4 \rho_0dA_{h_0}\right),$$
    and taking $t\to\infty$ shows that $-\int_{\Sigma}\chi d\mu\leq 0$ for any nonnegative $\chi$; thus, $\mu\geq 0$ on $\Sigma$. Moreover, it follows readily from \eqref{mu.bd} that $\mu$ is absolutely continuous with respect to the measure $\phi^4\rho_0dA_{h_0}$, so that $\mu$ has the form
    $$\mu=f_{\phi}\phi^4\rho_0 dA_{h_0},$$
    where $f_{\phi}$ is necessarily nonnegative and $L$-symmetric on $S$. In particular, since $\phi^4\rho_0\in L^1(\Sigma)$, it follows that $\mu$ is absolutely continuous with respect to $dA_{h_0}$ as well, and \eqref{mu.bd} holds for any $v\in L^{\infty}(\Sigma)$.

    To complete the proof, it remains to confirm that
    $$\int_Mf_{\phi}|\log f_{\phi}|\phi^4 dV_{g_0}=4\pi^2\int_{\Sigma}f_{\phi}|\log f_{\phi}|\phi^4\rho_0 dA_{h_0}<\infty.$$
    Testing \eqref{mu.bd} against the $L^{\infty}$ truncations $v_k:=\frac{1}{2}\max\{-k,\min\{ \log f_{\phi},k\}\}$ and letting $k\to\infty$, we see that
    $$\int_{\Sigma}\log(f_{\phi})f_{\phi}\phi^4\rho_0dA_{h_0}\leq 2C_{\phi}\left(1+\int_{\Sigma}\phi^4f_{\phi}\rho_0 dA_{h_0}\right)<\infty,$$
    giving the desired bound.
\end{proof}

As a consequence of the superharmonicity $-\Delta_{h_0}\rho\geq 0$ in $\Sigma$ supplied by Lemma \ref{rho.lap.lem}, we also have a positive lower bound for $\phi$.

\begin{lem}\label{phi.below}
    For $\phi$ as in Lemma \ref{rho.lap.lem}, the essential infimum satisfies 
    $$\inf_M\phi=c_{\phi}>0.$$
\end{lem}
\begin{proof}
    Since $|\Delta_{h_0}\rho|\in L^1(S)$, there is some $\epsilon>0$ such that
$$\int_{K_{\epsilon}}|\Delta_{h_0}\rho|>\frac{1}{2}\int_{\Sigma}|\Delta_{h_0}\rho|,$$
where $K_{\epsilon}:=\{y\in \Sigma\mid \rho_0(y)\geq\epsilon\}.$ On the other hand, for any $y\in K_{\epsilon}$ observe that the Dirichlet Green's function
$$-\Delta_{h_0}G_y=\delta_y,\text{ }G_y|_{\partial\Sigma}=0$$
is harmonic in a neighborhood of $\partial\Sigma$ and extends smoothly to a harmonic function across the boundary which is odd with respect to $R_1$ and $R_2$. As a consequence, it's easy to see that $\frac{G_y(x)}{\rho_0(x)}$ has a positive lower bound, which must be uniform as $y$ varies over the compact set $K_{\epsilon}\subset\subset \Sigma\setminus \partial\Sigma$. In other words, we have
$$\inf_{x\in \Sigma}\inf_{y\in K_{\epsilon}}\frac{G_y(x)}{\rho_0(x)}\geq c_{\epsilon}>0,$$
while Lemma \ref{rho.lap.lem} gives
$$\rho(x)=\int_{\Sigma}G_y(x)(-\Delta_{h_0}\rho)(y)dy\geq \int_{K_{\epsilon}}G_y(x)|\Delta_{h_0}\rho|(y)dy$$
for any $x\in \Sigma$, so that combining our estimates gives
\begin{eqnarray*}
    \rho(x)&\geq & \int_{K_{\epsilon}}G_y(x)|\Delta_{h_0}\rho|(y)dy\\
    &\geq & c_{\epsilon}\rho_0(x)\cdot \frac{1}{2}\int_{\Sigma}|\Delta_{h_0}\rho|dA_{h_0}
\end{eqnarray*}
Since $\rho$ is $R_i$-odd and nonconstant, clearly $\int_{\Sigma}|\Delta_{h_0}\rho|dA_{h_0}>0$, so setting
$$c_{\phi}^2:=c_{\epsilon}\cdot \frac{1}{2}\int_{\Sigma}|\Delta_{h_0}\rho|>0,$$
we arrive at the desired estimate
$$\phi^2(x)\geq \frac{\rho(x)}{\rho_0(x)}\geq c_{\phi}^2.$$
\end{proof}

Next, we set $v=0$ in \eqref{fin.char} and consider the consequences of the bound with respect to the variable $q$. To this end, first consider the space
$$Z_0:=\{q\in C^{\infty}(M)\mid (0,q)\in \mathcal{C}\},$$
and note that, since every $q\in Z_0$ satisfies $q\circ r=-q$, we must have $\int_M q dV_{g_0}=0$. By the usual Poincar\'e inequality on $(M,g_0)$ and the lower bound in Lemma \ref{phi.below} for $\phi$, it follows that
\begin{equation}\label{zphi.def}
\langle q,q'\rangle_{Z_{\phi}}:=\int_{\Sigma}\rho\langle dq,dq'\rangle dA_{h_0}=\frac 1{4\pi^2}\int_M\phi^2\langle dq,dq'\rangle_{g_0}dV_{g_0}
\end{equation}
defines an inner product on $Z_0$, such that
$$\|q\|_{W^{1,2}(M,g_0)}\leq C_{\phi}\|q\|_{Z_{\phi}}$$
for every $q\in Z_0$. Then, define
$$Z_{\phi}:=\overline{Z_0}^{\|\cdot\|_{Z_{\phi}}}\subset W^{1,2}(M,g_0)$$
to be the closure of $Z_0$ with respect to $\|\cdot\|_{Z_{\phi}}$.

Then we have the following.

\begin{lem}\label{q.char}
    For $\phi$ as in Lemma \ref{rho.lap.lem}, the linear functional $\ell:Z_0\to \mathbb{R}$ given by
    $$\ell(q):=\int_{\Sigma}\rho\langle d\log(\varphi_1/\varphi_2),dq\rangle_{h_0}$$
    extends to a bounded linear functional on $Z_{\phi}$. Defining
    $$\hat{T}(\phi,q):=-\|q\|_{Z_{\phi}}^2-\ell(q),$$
    for $q\in Z_{\phi}$, the supremum
    \begin{equation}\label{t.def}
        T(\phi):=\sup\{\hat{T}(\phi,q)\mid q\in Z_{\phi}\}
    \end{equation}
    is finite, and realized by a unique $q_{\phi}\in Z_{\phi}$ satisfying
    $$\hat{T}(\phi,q_{\phi})=\int_{\Sigma} \rho|dq_{\phi}|_{h_0}^2dA_{h_0}$$
    and 
    $$\int_{\Sigma}\rho\langle 2 dq_\phi,d\eta\rangle_{h_0}dA_{h_0}+\ell(\eta)=0$$
    for all $\eta\in Z_{\phi}$.
\end{lem}
\begin{proof}
    Denoting by $\ell:Z_0\to \mathbb{R}$ the linear function
    $$\ell(q):=\int_{\Sigma}\rho\langle d\log(\varphi_1/\varphi_2),dq\rangle_{h_0}dA_{h_0},$$
    we can write
    $$T(\phi,q)=-\|q\|_{Z_{\phi}}^2-\ell(q).$$
    Restricting to $v=0$ in \eqref{fin.char} clearly implies that
    $$\sup\{T(\phi,q)\mid q\in Z_0\}\leq C_{\phi}'<\infty,$$
    which is evidently equivalent to the boundedness
    $$|\ell(q)|\leq C\|q\|_{Z_{\phi}}$$
    of $\ell$ with respect to the $Z_{\phi}$ norm. 

    It follows that $\ell$ extends continuously to a bounded linear functional $\ell: Z_{\phi}\to \mathbb{R}$ on the Hilbert space $Z_\phi$, on which $T(\phi,q)$ is given by
    $$T(\phi,q)=-\|q\|_{Z_{\phi}}^2-\ell(q),$$
    so that 
    $$T(\phi):=\sup\{T(\phi,q)\mid q\in Z_{\phi}\}=\sup\{T(\phi,q)\mid q\in Z_0\}\leq C_{\phi}'<\infty.$$
    Elementary functional analysis shows that the minimization problem for $q\mapsto -T(\phi,q)=\|q\|_{Z_{\phi}}^2+\ell(q)$ on the Hilbert space $Z_{\phi}$ is well-posed, and the desired conclusions follow easily.    
\end{proof}

To show that the maximization problems for the $v$ and $q$ variables in $\tilde{E}_{g_0}(\phi,v,q)$ can be truly decoupled, we will also use the following lemma, extending the analysis above to arbitrary $(v,q)\in \cC$.

\begin{lem}\label{zphi.density}
    For any $(v,q)\in \cC$, we have $q\in Z_{\phi}$, and 
    $$-\int_{\Sigma}\rho\Delta_{h_0}v+T(\phi,q)=-\int_{\Sigma}v f_{\phi}\phi^4\rho_0 dA_{h_0}+\hat{T}(\phi,q).$$
\end{lem}
\begin{proof}
    Fix $(v,q)\in \cC$. For any $L$-invariant $\chi\in C_c^{\infty}(\Sigma)$, it is easy to see that $\chi q\in Z_0$, while
    $$\|q-\chi q\|_{Z_{\phi}}^2\leq 2\int_{\Sigma}\rho((1-\chi)^2|dq|_{g_0}^2+q^2|d\chi|_{g_0}^2)dA_{h_0}.$$
    We will apply this for a family of logarithmic cutoff-type functions $0\leq \chi_{\delta}\in C_c^{\infty}(\Sigma)$ satisfying
    $$0\leq \chi_{\delta}\leq 1\text{ on }\Sigma,$$
    $$\chi_{\delta}=0\text{ on }\{\rho_0<\delta^2\},\text{ }\chi_{\delta}=1\text{ on }\{\rho_0>\delta\},$$
    and
    $$|d\chi_{\delta}|\leq 2\frac{|d\rho_0|}{\rho_0\log(1/\delta)}\text{ on }\{\delta^2<\rho_0<\delta\}.$$
    For this family, we then have
    \begin{eqnarray*}
        \|q-\chi_{\delta}q\|_{Z_{\phi}}^2&\leq &2\|q\|_{C^1}^2\int_{\Sigma}\rho (1-\chi_{\delta})^2dA_{h_0}+2\|q\|_{C^0}^2\int_{\Sigma}\rho|d\chi_\delta|^2_{h_0}dA_{h_0}\\
        &\leq &2\|q\|_{C^1}^2\int_{\{\rho_0<\delta\}}\rho dA_{h_0}+8\|q\|_{C^0}^2\int_{\{\delta^2<\rho_0<\delta\}}\frac{\phi^2|d\rho_0|^2}{\rho_0\log(1/\delta)^2}dA_{h_0}.
    \end{eqnarray*}
    Since $\rho\in L^1(\Sigma)$, it's clear that the first term on the last line vanishes as $\delta\to 0$. To prove that $\lim_{\delta\to 0}\|q-\chi_{\delta}q\|_{Z_{\phi}}=0$, we therefore just need to show that
    \begin{equation}\label{log.err.van}
        \lim_{\delta\to 0}\int_{\{\delta^2<\rho_0<\delta\}}\frac{\phi^2|d\rho_0|^2}{\rho_0\log(1/\delta)^2}dA_{h_0}=0.
    \end{equation}

    To confirm \eqref{log.err.van}, note that for any $0<\epsilon<1/10$, the function $w_{\epsilon}$ given by
    $$w_{\epsilon}:=\frac{1}{\log \epsilon}-\frac{1}{\log(\rho_0+\epsilon)}\text{ where }\rho_0<\frac{1}{2}$$
    and
    $$w_{\epsilon}:=\frac{1}{\log\epsilon}-\frac{1}{\log(1/2+\epsilon)}\text{ where }\rho_0\geq \frac{1}{2}$$
    is Lipschitz and vanishes on $\partial\Sigma$, with $w_{\epsilon}\leq C$, where $C$ is independent of $\epsilon$. Thus, we can apply \eqref{sig.pair} with $w_{\epsilon}$ to obtain an estimate of the form
    $$\int_{\Sigma}\langle d\rho,dw_{\epsilon}\rangle_{h_0}dA_{h_0}\leq C_{\phi}'.$$
    On the other hand, since $dw_{\epsilon}=\frac{d\rho_0}{(\rho_0+\epsilon)\log(\rho_0+\epsilon)^2}$, we see that
    \begin{eqnarray*}
        \int_{\Sigma}\langle d\rho,dw_{\epsilon}\rangle&=&\int_{\{\rho_0<1/2\}}\frac{\phi^2|d\rho_0|^2+2\rho_0\phi\langle d\phi,d\rho_0\rangle_{h_0}}{(\rho_0+\epsilon)\log(\rho_0+\epsilon)^2}dA_{h_0}\\
        &\geq &\int_{\{\rho_0<1/2\}}\frac{1}{2}\frac{\phi^2|d\rho_0|^2}{(\rho_0+\epsilon)\log(\rho_0+\epsilon)^2}-2\frac{\rho_0^2|d\phi|^2}{(\rho_0+\epsilon)\log(\rho_0+\epsilon)^2}dA_{h_0}\\
        &\geq &\frac{1}{2}\int_{\{\rho_0<1/2\}}\frac{\phi^2|d\rho_0|^2}{(\rho_0+\epsilon)\log(\rho_0+\epsilon)^2}-C\|\phi\|_{W^{1,2}(M)}^2,
    \end{eqnarray*}
    so taking $\epsilon\to 0$ in the preceding estimate yields
    \begin{equation}\label{lcut.l1}
    \int_{\{\rho_0<1/2\}}\frac{\phi^2|d\rho_0|^2}{\rho_0\log(\rho_0)^2}<\infty.
    \end{equation}
    Finally, to confirm \eqref{log.err.van}, note that
    $$\int_{\{\delta^2<\rho_0<\delta\}}\frac{\phi^2|d\rho_0|^2}{\rho_0\log(1/\delta)^2}dA_{h_0}\leq 4\int_{\{\delta^2<\rho_0<\delta\}}\frac{\phi^2|d\rho_0|^2}{\rho_0\log(\rho_0)^2},$$
    and since the annuli $\{\delta^2<\rho_0<\delta\}$ are contained in $\{\rho_0<1/2\}$ with vanishing measure as $\delta\to 0$, the desired vanishing \eqref{log.err.van} follows from the integrability statement \eqref{lcut.l1}. In particular, any sequence $\delta_j\to 0$ gives a Cauchy sequence $\chi_{\delta_j}q$ in $Z_{\phi}$ converging to the initial $q$ in $Z_{\phi}$ norm, giving the desired inclusion.

    Finally, for the given smooth pair $(v,q)\in \mathcal{C}$, we can define a function $\bar{q}\in Lip(\Sigma)$ by setting $\bar{q}=q$ on each of the $n+1$ triangles $U_1,\ldots,U_{n+1}$ in $\Sigma$ with an edge in $\Gamma_1$, and $\bar{q}=-q$ on the $(n+1)$ triangles $V_1,\ldots,V_{n+1}$ with an edges in $\Gamma_2$. Since $(v,q)\in \mathcal{C}$, it follows that $v-\bar{q}\in W_0^{1,\infty}(\Sigma)$, and by Lemma \ref{rho.lap.lem}, we must have
    \begin{eqnarray*}
    -\int_{\Sigma}\rho\Delta_{h_0}v&=&\int_{\Sigma}\langle d\rho,dv\rangle_{h_0}\\
    &=&\int_{\Sigma}\langle d\rho,d(v-\bar{q})\rangle_{h_0}+\int_{\Sigma}\langle d\rho,d\bar{q}\rangle_{h_0}dA_{h_0}\\
    &=&\int_{\Sigma}f_{\phi}\phi^4(v-\bar{q})\rho_0 dA_{h_0}+\int_{\Sigma}\langle d\rho,d\bar{q}\rangle_{h_0}dA_{h_0}.
    \end{eqnarray*}
    Expanding the latter term and exploiting symmetries, we see that
    \begin{eqnarray*}
        \int_{\Sigma}\langle d\rho,d\bar{q}\rangle_{h_0}dA_{h_0}&=&\sum_{j=1}^{n+1}\int_{U_j}\langle d\rho,dq\rangle_{h_0}dA_{h_0}-\sum_{j=1}^{n+1}\int_{V_j}\langle d\rho,dq\rangle_{h_0}\\
        &=&\int_{\Sigma}\phi^2\langle \varphi_1d\varphi_2-\varphi_2d\varphi_2,dq\rangle_{h_0}dA_{h_0},
    \end{eqnarray*}
    so that
    \begin{eqnarray*}    \int_{\Sigma}\rho\Delta_{h_0}v+T(\phi,q)&=&\int_{\Sigma}f_{\phi}\phi^4(v-\bar{q})\rho_0 dA_{h_0}-\|q\|_{Z_{\phi}}^2.
    \end{eqnarray*}
    Finally, since $f_{\phi}\phi^4\rho_0\in L^1(\Sigma)$, it's easy to see that the quantities on the right-hand side vary continuously under the approximation
    $$(\chi_{\delta}v,\chi_{\delta}q)\to (v,q)$$
    constructed above, so we must have
    $$\hat{T}(\phi,q)=-\|q\|_{Z_{\phi}}^2-\int_{\Sigma}f_{\phi}\phi^4 \bar{q}\rho_0 dA_{h_0}.$$
    
\end{proof}

\begin{remark}\label{ell.rk}
    In particular, it follows from the proof above that, for any smooth $q\in Z_{\phi}$, or more generally any $q\in Z_{\phi}\cap L^{\infty}(M)$, the linear functional $\ell$ in Lemma \ref{q.char} must be given by    $$\ell(q)=\int_{\Sigma}f_{\phi}\phi^4\bar{q}\rho_0dA_{h_0}=\int_{\Sigma}q div_{h_0}(\phi^2[\varphi_1d\varphi_2-\varphi_2d\varphi_1]),$$
    where the divergence in the latter term should be interpreted in the pointwise sense in the interior of $\Sigma$, rather than distributionally.
\end{remark}

With these basic consequences of $\cS_{g_0}(\phi)<\infty$ established, we can give the following useful closed-form expression for $\cS_{g_0}(\phi)$.

\begin{prop}\label{alt.s.char}
    Let $0\leq \phi\in W^{1,2}(M)^G\setminus\{0\}$ satisfy $\cS_{g_0}(\phi)<\infty$.
    Letting $T(\phi)$ be defined as in \eqref{t.def}, and writing $-\Delta_{h_0}\rho=f_{\phi}\phi^4\rho_0$ as in Lemma \ref{rho.lap.lem}, it follows that 
    \begin{eqnarray*}
        \log(\cS_{g_0}(\phi))&=&\frac{E_{g_0}(\phi,0,0)+8\pi^2T(\phi)+\int_M\log(f_{\phi})f_{\phi}\phi^4 dV_{g_0}}{2\int_M f_{\phi}\phi^4 dV_{g_0}}\\
        &&-1+{\log2}+\frac{1}{2}\log\left(\int_Mf_{\phi}\phi^4 dV_{g_0}\right).
    \end{eqnarray*}
\end{prop}
\begin{proof}
    It follows from Lemmas \ref{rho.lap.lem}, \ref{q.char} and \ref{zphi.density} that, whenever $\cS_{g_0}(\phi)<\infty$, the function
    $$\mathcal{C}\ni (v,q)\mapsto E_{g_0}(\phi,v,q)$$
    has a continuous extension to $L^{\infty}(M)^G\times Z_{\phi}$, given by
    $$E_{g_0}(\phi,v,q)=E_{g_0}(\phi,0,0)+2\int_M v f_{\phi}\phi^4 dV_{g_0}+8\pi^2 \hat{T}(\phi,q),$$
    so that
    $$\cS_{g_0}(\phi)=\sup\{\frac{E_{g_0}(\phi,v,q)}{\left(\int_M\phi^4 e^{2v}dV_{g_0}\right)^{1/2}}\mid v\in L^{\infty}(M)^G,\text{ }q\in Z_{\phi}\},$$
    or equivalently,
    $$\cS_{g_0}(\phi)=\sup\{\frac{E_{g_0}(\phi,0,0)+8\pi^2T(\phi)+2\int_M v f_{\phi}\phi^4 dV_{g_0}}{\left(\int_M e^{2v} \phi^4 dV_{g_0}\right)^{1/2}}\mid v\in L^{\infty}(M)^G\}.$$
    Temporarily setting $A_{\phi}:=E_{g_0}(\phi,0,0)+8\pi^2 T(\phi)$, we observe that for any fixed $v\in L^{\infty}(M)^G$, the maximum of
    $$\mathbb{R}\ni t\mapsto \frac{A_{\phi}+2\int_M (v+t) f_{\phi}\phi^4 dV_{g_0}}{\left(\int_M e^{2(v+t)}\phi^4 dV_{g_0}\right)^{1/2}}$$
    occurs  at the unique value $t_v$ for which $A_{\phi}+2\int_M (v+t_v)f_{\phi}\phi^4 dV_{g_0}=2\int_M f_{\phi}\phi^4 dV_{g_0}$, so restricting to those $v$ for which $t_v=0$, we can recast $\cS_{g_0}(\phi)$ as the supremum of 
    \begin{equation}\label{norm.quot}
    {2}\int_M f_{\phi}\phi^4 dV_{g_0}\cdot \left(\int_M e^{2v}\phi^4 dV_{g_0}\right)^{-1/2}
    \end{equation}
    among all $v\in L^{\infty}(M)^G$ satisfying
    \begin{equation}\label{scal.norm}
        \frac{A_{\phi}}{2}+\int_M v f_{\phi}\phi^4 dV_{g_0}=\int_M f_{\phi}\phi^4 dV_{g_0}.
    \end{equation}
    
    By Lemma \ref{rho.lap.lem}, it follows that 
    $$v_{\phi}:=\frac{1}{2}(\log(f_{\phi})-t_{\phi}),$$
    where
    $$t_{\phi}:=\frac{A_{\phi}+\int_M \log(f_{\phi})f_{\phi}\phi^4 dV_{g_0}}{\int_M f_{\phi}\phi^4 dV_{g_0}}-2,$$
    belongs to $L^1(M,f_{\phi}\phi^4 dV_{g_0})^G$, and it is straightforward to find an approximating sequence $v_k\in L^{\infty}(M)^G$ such that $|v_k|\leq |v|$, $v_k\to v_{\phi}$ almost everywhere in $(M,g_0)$, and $v_k$ satisfy \eqref{scal.norm}. Testing $v_k$ in \eqref{norm.quot} and letting $k\to\infty$, we deduce that
    $$\cS_{g_0}(\phi)\geq {2}\left(\int_M f_{\phi}\phi^4 dV_{g_0}\right)^{1/2} e^{t_{\phi}/2}.$$
    On the other hand, any other test function $v$ for the supremum can be written $v=v_k+w_k$, where $\int_M w f_{\phi}\phi^4=0$, and using the elementary pointwise inequality $e^{2w}\geq 1+2w$, it follows that
    \begin{eqnarray*}
        \int_M e^{2v}\phi^4 dV_{g_0}&=&\int_M e^{2v_k}e^{2w_k}\phi^4 dV_{g_0}\\
        &\geq &\int_M e^{2v_k}(1+2w_k)\phi^4 dV_{g_0},
    \end{eqnarray*}
    and letting $k\to\infty$ gives
    \begin{eqnarray*}
        \int_M e^{2v}\phi^4 dV_{g_0}&\geq &\int_M e^{2v_{\phi}}\phi^4 dV_{g_0}+2 e^{-t_{\phi}}\int_M w f_{\phi}\phi^4 dV_{g_0}\\
        &=&\int_M e^{2v_{\phi}}\phi^4 dV_{g_0}.
    \end{eqnarray*}
    It follows that the supremum of \eqref{norm.quot} is achieved precisely at $v=v_{\phi}$, so that
    $$\cS_{g_0}(\phi)={2}\left(\int_M f_{\phi}\phi^4 dV_{g_0}\right)^{1/2}e^{t_{\phi}/2}.$$
    Taking logarithms and expanding $t_{\phi}$ then gives the desired identity.
\end{proof}

\subsection{Regularity theory for a weighted problem via the six-dimensional lift}\label{l0.sec}

Both for deriving the weak Euler-Lagrange condition for our minimizing $\phi$ and establishing regularity, it is convenient to have some regularity theory in place for the operator
$$\chi\mapsto \frac{-\Delta_{h_0}(\chi\rho_0)}{\rho_0}$$
on $C^{\infty}(M)^G$. 

To this end, recall the six-dimensional reference manifold $(P,G_0)$, and observe that any $\chi\in C^{\infty}(M)^G$ has a smooth $O(3)\times O(3)$-invariant lift to $P$, for which
$$\int_P|d\chi|_{G_0}^2dV_{G_0}=16\pi^2\int_{\Sigma}|d\chi|_{h_0}^2\rho_0^2dA_{h_0}.$$
In particular, it follows that the Laplacian $\Delta_{G_0}\chi$ in this case has the form
$$\Delta_{G_0}\chi=\rho_0^{-2}div_{h_0}(\rho_0^2d\chi)=\Delta_{h_0}\chi+2\rho_0^{-1}\langle d\rho_0,d\chi\rangle_{h_0}.$$
On the other hand, since
$$\Delta_{h_0}(\chi \rho_0)=\rho_0\Delta_{h_0}\chi+2\langle d\rho_0,d\chi\rangle_{h_0}+\chi \Delta_{h_0}\rho_0,$$
defining an operator on $C^{\infty}(P)$ by
\begin{equation}\label{l0.def}
L_0\chi:=\Delta_{G_0}\chi+\frac{\Delta_{h_0}\rho_0}{\rho_0}\chi,
\end{equation}
it follows that
\begin{equation}\label{l0.char}
L_0\chi=\frac{\Delta_{h_0}(\chi\rho_0)}{\rho_0}
\end{equation}
whenever $\chi\in C^{\infty}(M)^G$.

Since the quotient $\frac{\Delta_{h_0}\rho_0}{\rho_0}$ is smooth, the operator $L_0$ is simply a smooth $0$th-order perturbation of the Laplacian on $(P,G_0)$. Moreover, since $\rho_0^{-1}\Delta_{h_0}\rho_0$ is $O(3)\times O(3)$-invariant, we see that the first eigenfunction for $-L_0$ must belong to $C^{\infty}(M)^G$, which together with \eqref{l0.char} implies that $-L_0$ is positive-definite. As a consequence, standard elliptic regularity theory applied to $L_0$ gives the following.

\begin{lem}\label{l0.reg}
    For any $q\in (1,\infty)$ and $V\in L^q(P)$, there exists a unique $\chi\in W^{2,q}(P)$ such that
    $$-L_0\chi=V,$$
    with $\|\chi\|_{W^{2,q}(P)}\leq C_q\|V\|_{L^q(P)}$. In particular, if $V\in L^{\infty}(M)^G$, there exists a unique $\chi\in C^{1,\alpha}(M)^G$ for every $\alpha\in (0,1)$ such that
    $$-\Delta_{h_0}(\chi\rho_0)=V\rho_0,$$
    satisfying
    $$\|\chi\|_{C^{1,\alpha}(\Sigma)}\leq C_{\alpha}\|V\|_{L^{\infty}}.$$
\end{lem}
\begin{proof}
    The first statement is a standard consequence of $L^q$ regularity theory for elliptic operators, as in e.g. \cite[\S 9]{GilbargTrudinger2001Elliptic}. For the second statement, we simply observe that the symmetries of the operator $L_0$ imply $L_0^{-1}V\in C^{\infty}(M)^G$ whenever $V\in C^{\infty}(M)^G$, and combine the $L^q$ regularity result on $P$ with the embedding $W^{2,q}(P)\to C^{1,1-6/q}$ for $q>6$, noting that $\|\chi\|_{C^{1,\alpha}(\Sigma)}=\|\chi\|_{C^{1,\alpha}(P)}$ whenever $\chi\in C^{\infty}(M)^G$.    
\end{proof}

\subsection{A weak Euler-Lagrange equation for the minimizer}

Now, let $\phi\in W^{1,2}(M)^G$ be the minimizer of $\cS_{g_0}$ supplied by Theorem \ref{wk.ex}, and recall from Lemma \ref{phi.below} that $\phi\geq c_{\phi}>0$ on $M$. By Proposition \ref{alt.s.char}, it then follows that, along any path $t\mapsto \phi_t\in W^{1,2}(M)^G$ such that $\cS(\phi_t)<\infty$, $\phi_0=\phi$, the function
\begin{equation}\label{t.fam}
    F(t):=\frac{E_{g_0}(\phi_t,0,0)+8\pi^2T(\phi_t)+\int_M\log(f_{\phi_t})f_{\phi_t}\phi_t^4 dV_{g_0}}{\int_M f_{\phi_t}\phi_t^4 dV_{g_0}}+\log\left(\int_Mf_{\phi_t}\phi_t^4 dV_{g_0}\right)
\end{equation}
is minimized at $t=0$. Using this, we first derive the following conditional variational inequality for the minimizer $\phi$, which will serve as the basis for the full Euler-Lagrange equations.

In the following, from the proof of Proposition \ref{alt.s.char}, recall the notation
$$A_{\phi}:=E_{g_0}(\phi,0,0)+8\pi^2T(\phi)$$
and
\begin{equation}\label{v.def}
    v_{\phi}:=\frac{1}{2}\log(f_{\phi})-\frac{A_{\phi}+\int_M\log(f_{\phi})f_\phi\phi^4 dV_{g_0}}{2\int_Mf_{\phi}\phi^4 dV_{g_0}}+1.
\end{equation}

\begin{lem}\label{el.ineq}
    Let $0<\phi\in W^{1,2}(M)^G$ be a minimizer for $\cS_{g_0}$, and let $\zeta\in L^{\infty}(M)^G$ be a bounded function such that 
    $$\inf \frac{\zeta}{f_{\phi}\phi^4}>-\infty.$$
    Letting $\chi\in C^{1,\alpha}(\Sigma)$ be the solution of $-\Delta_{h_0}(\chi\rho_0)=\zeta \rho_0,$ we then have
    \begin{eqnarray*}
        2\int_M\zeta v_{\phi} dV_{g_0}&\geq & 2\int_M\chi \phi^2 f_{\phi} dV_{g_0}-\int_M(6\langle d\phi,d(\chi/\phi)\rangle_{g_0}+R_{g_0}\chi) dV_{g_0}\\
        &&+8\pi^2\int_{\Sigma}(\chi\rho_0|dq_{\phi}|_{h_0}^2+q_{\phi} \cdot div_{h_0}(\chi(\varphi_1d\varphi_2-\varphi_2d\varphi_1)))dA_{h_0}.
\end{eqnarray*}
\end{lem}

\begin{proof}
To begin, fix $\zeta\in L^{\infty}(M)^G$ and, appealing to Lemma \ref{l0.reg}, let $\chi\in C^{1,\alpha}(\Sigma)$ be the unique solution of 
$$-\Delta_{h_0}(\chi\rho_0)=\zeta\rho_0.$$
Since $\inf \phi>0$, defining $\phi_t:=\sqrt{\phi^2+t\chi}$ then gives a differentiable path in $W^{1,2}(M)^G$ such that the functions $$\rho^t:=\phi_t^2\rho_0=\rho+t\chi\rho_0$$
satisfy
$$-\Delta_{h_0}\rho^t=-\Delta_{h_0}\rho+t\zeta \rho_0.$$

For any such family, note that the ratio $\frac{\phi_t}{\phi}$ satisfies $\frac{\phi_t^2}{\phi^2}\in [1-C|t|,1+C|t|]$ for a suitable constant $C=C(\zeta,\phi).$ As a consequence, for $t$ sufficiently small, the norms $Z_{\phi_t}$ must be equivalent, and the functionals $T(\phi_t,q)$ satisfy
\begin{eqnarray*}
    T(\phi_t,q)-T(\phi,q)&=&\int_{\Sigma}(\rho-\rho^t)[|dq|_{h_0}^2+\langle d\log(\varphi_1/\varphi_2),dq\rangle_{h_0}]dA_{h_0}\\
    &=&-t\int_{\Sigma}\chi\rho_0[|dq|_{h_0}^2+\langle d\log(\varphi_1/\varphi_2),dq\rangle_{h_0}]dA_{h_0}.
\end{eqnarray*}
Note that for any $q\in C_c^{\infty}(\Sigma)\cap Z_{\phi}$,  we can integrate by parts on $(M,g_0)$ to write
\begin{eqnarray*}
    \int_{\Sigma}\chi\rho_0\langle d\log(\varphi_1/\varphi_2),dq\rangle_{h_0}&=&-\int_{\Sigma} q div_{h_0}(\rho_0 \chi d\log(\varphi_1/\varphi_2)) dA_{h_0}\\
    &=&\int_{\Sigma}q div_{h_0}(\chi (\varphi_1d\varphi_2-\varphi_2d\varphi_1)) dA_{h_0}.
\end{eqnarray*}
Since $\chi\in C^{1,\alpha}(\Sigma)$, it is straightforward to see that the last line extends to a continuous function on $Z_{\phi}$,  and we write
$$T(\phi_t,q)-T(\phi,q)=-t\int_{\Sigma}\left(\chi \rho_0|dq|_{h_0}^2+q div_{h_0}(\chi(\varphi_1d\varphi_2-\varphi_2d\varphi_1))\right)dA_{h_0}.$$
Using this, we easily deduce that $\lim_{t\to 0}T(\phi_t)=T(\phi)$, and the associated maximizers $q_{\phi_t}\in Z_{\phi}$ define a differentiable path in $Z_{\phi}$. In particular, noting that $\frac{d}{dt}|_{t=0}T(\phi,q_{\phi_t})=0$ by maximality of $q_{\phi}$ for $T(\phi,q)$, it follows that
\begin{eqnarray*}
    \frac{d}{dt}|_{t=0}T(\phi_t)&=&\lim_{t\to 0}\frac{T(\phi_t,q_{\phi_t})-T(\phi,q_{\phi})}{t}\\
    &=&\lim_{t\to 0}\frac{T(\phi_t,q_{\phi_t})-T(\phi,q_{\phi_t})}{t}\\
    &=&-\int_{\Sigma}(\chi\rho_0|dq_{\phi}|_{h_0}^2+q_{\phi} div_{h_0}(\chi(\varphi_1d\varphi_2-\varphi_2d\varphi_1)))dA_{h_0}.
\end{eqnarray*}

Next, assume that 
\begin{equation}\label{supsol}
    -\Delta_{h_0}\rho+t\zeta\rho_0\geq 0\text{ for }t>0\text{ sufficiently small}.
\end{equation}
Then $-\Delta_{h_0}\rho^t\geq 0$ for small $t>0$, and writing
$$-\Delta_{h_0}\rho^t=f_{\phi_t}\phi_t^4\rho_0,$$
observe that
\begin{equation}\label{ft.exp}
f_{\phi_t}=\frac{\phi^4f_{\phi}}{\phi_t^4}+t\zeta\phi_t^{-4}.
\end{equation}
By continuity of $x\mapsto x\log x$, it's then easy to see that $\int_M \log(f_{\phi_t})f_{\phi_t}\phi_t^4dV_{g_0}<\infty$ for $t>0$ sufficiently small, and an application of the elementary inequality 
$$(x+y)\log(x+y)\leq x\log(x)+y+y\log(x+y)$$
gives
\begin{equation}\label{xlogx.cp}   
\log(f_{\phi_t})f_{\phi_t}\phi_t^4\leq\log(f_{\phi}\phi^4/\phi_t^{4})f_{\phi}\phi^4+t\zeta [1+\log(f_{\phi_t})].
\end{equation}

In particular, this together with the nonnegativity \eqref{supsol}, the $W^{1,2}(M)$-convergence $\phi_t\to\phi$, and $T(\phi_t)\to T(\phi)$ allows us to conclude that $\cS_{g_0}(\phi_t)<\infty$ for $t>0$ sufficiently small. This is enough to ensure that the function $F(t)$ given by \eqref{t.fam} is well-defined for small $t>0$ and minimized at $t=0$. Since
$$\frac{d}{dt}|_{t=0}\phi_t=\frac{\chi}{2\phi},$$
standard computations give
$$\frac{d}{dt}|_{t=0}E_{g_0}(\phi_t,0,0)=\int_M6\langle d\phi,d(\chi/\phi)\rangle_{g_0}+R_{g_0}\chi,$$
and we have already computed
$$\frac{d}{dt}|_{t=0}T(\phi_t)=-\int_{\Sigma}(\chi\rho_0|dq_{\phi}|_{h_0}^2+q_{\phi} \cdot div_{h_0}(\chi(\varphi_1d\varphi_2-\varphi_2d\varphi_1)))dA_{h_0}.$$

For the entropy term $\int_M\log(f_{\phi_t})f_{\phi_t}\phi_t^4 dV_{g_0}$, we can apply \eqref{xlogx.cp} to see that
\begin{eqnarray*}
    \int_M[\log(f_{\phi_t})f_{\phi_t}\phi_t^4-\log(f_{\phi})f_{\phi}\phi^4] dV_{g_0}&\leq &\int_M 4\log(\phi/\phi_t)f_{\phi}\phi^4 dV_{g_0}\\
    &&+t\int_M\zeta[1+\log(f_{\phi_t})]dV_{g_0}.
\end{eqnarray*}
Note that $\frac{d}{dt}|_{t=0}4\log(\phi/\phi_t)=-2\frac{\chi}{\phi^2}$ holds with limits interpreted in the $L^{\infty}$ sense; moreover, using \eqref{ft.exp} and the fact that $f_{\phi}\log(f_{\phi})\in L^1(M)$, we can use the monotone convergence theorem to justify the limit
$$\lim_{t\to 0}\int_M\zeta[1+\log(f_{\phi_t})]dV_{g_0}=\int_M\zeta[1+\log(f_{\phi})]dV_{g_0},$$
where the integral $\int_M\zeta \log(f_{\phi})dV_{g_0}$ takes values in $[-\infty,\infty)$ a priori. Putting this together, it follows that
\begin{eqnarray*}
    &&\limsup_{t\downarrow 0}\frac{1}{t}\int_M[\log(f_{\phi_t})f_{\phi_t}\phi_t^4-\log(f_{\phi})f_{\phi}\phi^4] dV_{g_0}\\&\leq &-2\int_M\chi \phi^2 f_{\phi} dV_{g_0}+\int_M\zeta[1+\log(f_{\phi})]dV_{g_0}.
\end{eqnarray*}

Finally, since $\int_Mf_{\phi_t}\phi_t^4dV_{g_0}=\int_Mf_{\phi}\phi^4dV_{g_0}+t\int_M\zeta dV_{g_0}$, it's easy to see that
$$\frac{d}{dt}|_{t=0}\log\left(\int_Mf_{\phi_t}\phi_t^4dV_{g_0}\right)=\frac{\int_M\zeta dV_{g_0}}{\int_Mf_{\phi}\phi^4 dV_{g_0}}.$$
Combining the preceding computations, and using that $F(0)\leq F(t)$ for all $t>0$, it follows that
\begin{eqnarray*}
    0&\leq &\int_Mf_{\phi}\phi^4 dV_{g_0}\cdot\limsup_{t\downarrow 0}\frac{F(t)-F(0)}{t}\\
    &\leq &\int_M6\langle d\phi,d(\chi/\phi)\rangle_{g_0}+R_{g_0}\chi\\
    &&-8\pi^2\int_{\Sigma}(\chi\rho_0|dq_{\phi}|_{h_0}^2+q_{\phi} \cdot div_{h_0}(\chi(\varphi_1d\varphi_2-\varphi_2d\varphi_1)))dA_{h_0}\\
    &&-2\int_M\chi\phi^2 f_{\phi} dV_{g_0}+\int_M\zeta[1+\log(f_{\phi})]dV_{g_0}+\int_M\zeta dV_{g_0}\\
    &&-\frac{A_{\phi}+\int_M\log(f_{\phi})f_\phi\phi^4 dV_{g_0}}{\int_M f_{\phi}\phi^4dV_{g_0}}\int_M\zeta dV_{g_0},
\end{eqnarray*}
for any $\zeta\in L^{\infty}(M)^G$ satisfying \eqref{supsol}, or equivalently, $\inf\frac{\zeta}{f_{\phi}\phi^4}>-\infty$. With a little rearrangement, and recalling \eqref{v.def}, we obtain the stated inequality. 
\end{proof}

By the results of Section \ref{s.fin.sec}, we see that the right-hand side of the final inequality in Lemma \ref{el.ineq} has a lower bound of the form $-C_{\phi}\|\chi\|_{C^{1,\alpha}}\geq -C'\|\zeta\|_{L^{\infty}}$, so that
$$\int_M \zeta \log(f_{\phi})dV_{g_0}\geq -C'\|\zeta\|_{L^{\infty}}$$
whenever $\inf \frac{\zeta}{f_{\phi}\phi^4}>-\infty$. In particular, this holds whenever $\zeta\geq 0$, and together with the upper bound $\int_M f_{\phi}\phi^4 dV_{g_0}<\infty$, this leads to an $L^1$ bound
\begin{equation}\label{v.l1}
\int_M|\log(f_{\phi})|dV_{g_0}<\infty,
\end{equation}
so that $v_{\phi}\in L^1(M,dV_{g_0})$.

With this, we can now upgrade Lemma \ref{el.ineq} to the full weak Euler-Lagrange equation, as follows.

\begin{prop}\label{weak.el}
    Let $0<\phi\in W^{1,2}(M)^G$ be a minimizer for $\cS_{g_0}$. For any $\zeta\in L^{\infty}(M)^G$ , letting $\chi\in C^{1,\alpha}(\Sigma)$ be the solution of $-\Delta_{h_0}(\chi\rho_0)=\zeta \rho_0,$ we have
    \begin{eqnarray*}
        2\int_M\zeta v_{\phi}dV_{g_0}&= & 2\int_M\chi \phi^2 f_{\phi} dV_{g_0}-\int_M(6\langle d\phi,d(\chi/\phi)\rangle_{g_0}+R_{g_0}\chi) dV_{g_0}\\
        &&+8\pi^2\int_{\Sigma}(\chi\rho_0|dq_{\phi}|_{h_0}^2+q_{\phi} \cdot div_{h_0}(\chi(\varphi_1d\varphi_2-\varphi_2d\varphi_1)))dA_{h_0}.
\end{eqnarray*}
\end{prop}
\begin{proof}
Given $\zeta\in L^{\infty}(M)^G$ and $k\in \mathbb{N}$, define $\zeta_k\in L^{\infty}(M)^G$ by
$$\zeta_k:=\zeta\text{ on }\{f_{\phi}>1/k\}$$
and
$$\zeta_k=0\text{ on }\{f_{\phi}\leq 1/k\}.$$
We then see that $\|\zeta_k\|_{L^{\infty}}$ is uniformly bounded and $\zeta_k\to \zeta$ almost everywhere as $k\to\infty$, thanks to \eqref{v.l1}. Moreover, each $\zeta_k$ satisfies the hypotheses of Lemma \ref{el.ineq}, so letting $\chi_k$ solve
$$-\Delta_{h_0}(\chi_k\rho_0)=\zeta_k\rho_0,$$
it follows that

\begin{eqnarray*}
        2\int_M\zeta_k v_{\phi}dV_{g_0}&\geq & 2\int_M\chi_k \phi^2 f_{\phi} dV_{g_0}-\int_M(6\langle d\phi,d(\chi_k/\phi)\rangle_{g_0}+R_{g_0}\chi_k) dV_{g_0}\\
        &&+8\pi^2\int_{\Sigma}(\chi_k\rho_0|dq_{\phi}|_{h_0}^2+q_{\phi} \cdot div_{h_0}(\chi_k(\varphi_1d\varphi_2-\varphi_2d\varphi_1)))dA_{h_0}.
\end{eqnarray*}

On the other hand, since $\|\zeta_k\|_{L^{\infty}}\leq C$ and $\zeta_k\to\zeta$ almost everywhere as $k\to\infty$, we must have $\|\zeta_k-\zeta\|_{L^p(P)}\to 0$ as $k\to\infty$ for each $p\in [1,\infty)$. By Lemma \ref{l0.reg}, it follows that
$$\lim_{k\to\infty}\|\chi_k-\chi\|_{C^{1,\alpha}(\Sigma)}=0$$
for every $\alpha\in [0,1)$. As a consequence of this $C^{1,\alpha}$ convergence, it's not hard to see that the right-hand side of the inequality above must converge as $k\to\infty$ to the corresponding expression with $\chi$ in place of $\chi_k$, so passing the inequality to the limit gives
\begin{eqnarray*}
    2\int_M\zeta v_{\phi}dV_{g_0}&\geq & 2\int_M\chi\phi^2 f_{\phi} dV_{g_0}-\int_M(6\langle d\phi,d(\chi/\phi)\rangle_{g_0}+R_{g_0}\chi) dV_{g_0}\\
        &&+8\pi^2\int_{\Sigma}(\chi\rho_0|dq_{\phi}|_{h_0}^2+q_{\phi} \cdot div_{h_0}(\chi(\varphi_1d\varphi_2-\varphi_2d\varphi_1)))dA_{h_0}.
\end{eqnarray*}
Finally, since $\zeta\in L^{\infty}(M)^G$ was arbitrary, replacing $\zeta$ with $-\zeta$ gives the opposite inequality, completing the proof.
\end{proof}

\section{Interior regularity and recovery of the Einstein equations}\label{sect:int.reg}

As in the previous sections, fix a minimizer $0<\phi\in W^{1,2}(M)^G$ for the energy $\cS_{g_0}$, and let $q_{\phi}\in Z_{\phi}\subset W^{1,2}(M)$ and $v_{\phi}\in L^1(M)^G$ be the associated functions given by Lemma \ref{q.char} and \eqref{v.def}, so that for $\rho=\phi^2\rho_0$, the results of the previous section give:
\begin{equation}\label{rho.wk}
    -\Delta_{h_0}\rho= f_{\phi}\phi^4\rho_0=C_{\phi}e^{2v_{\phi}}\phi^4\rho_0
\end{equation}
weakly on $S$ for a suitable constant $C_{\phi}$,
\begin{equation}\label{q.wk}
    \int_{\Sigma}\rho\langle 2dq+d\log(\varphi_1/\varphi_2),d\eta\rangle_{h_0}dA_{h_0}=0
\end{equation}
for $\eta\in C_c^{\infty}(M)$, and
\begin{eqnarray*}
    2\int_M v_{\phi}\frac{-\Delta_{h_0}(\chi \rho_0)}{\rho_0}&=& 2\int_M\chi \phi^2 f_{\phi} dV_{g_0}-\int_M(6\langle d\phi,d(\chi/\phi)\rangle_{g_0}+R_{g_0}\chi) dV_{g_0}\\
        &&+8\pi^2\int_{\Sigma}(\chi\rho_0|dq_{\phi}|_{h_0}^2+q_{\phi} \cdot div_{h_0}(\chi(\varphi_1d\varphi_2-\varphi_2d\varphi_1)))dA_{h_0}.
\end{eqnarray*}
whenever $\frac{\Delta_{h_0}(\chi\rho_0)}{\rho_0}\in L^{\infty}(M)^G$.

In this section, we will show that $\phi$, $v_{\phi},$ and $q_{\phi}$ are smooth in the interior of $\Sigma$, i.e. in the region of $M$ on which $T^2$ acts freely, then use this together with the gauge-invariance of the minimization problem for $\cS_{g_0}$ to deduce that the associated metric is Einstein in the interior of $\Sigma$.

\subsection{Smoothness of the triple $(\phi,v_{\phi},q_{\phi})$ away from $\partial\Sigma$}

To begin, let's apply Proposition \ref{weak.el} with $\chi\in C_c^{\infty}(\Sigma)$ supported away from $\partial\Sigma$. Using that $\phi,q_{\phi}\in W^{1,2}(M)$, we  integrate by parts in some of the terms involving $d\chi$ on the right-hand side, to see that
\begin{eqnarray*}
    8\pi^2\int_{\Sigma}v_{\phi}(-\Delta_{h_0}(\chi\rho_0))dA_{h_0}&=&2\int_M\chi \phi^2 f_{\phi}dV_{g_0}\\
    &&+\int_M(6|d\log(\phi)|_{g_0}^2\chi-6\langle d\log(\phi),d\chi\rangle_{g_0}-R_{g_0}\chi)dV_{g_0}\\
    &&+8\pi^2\int_{\Sigma}\chi(\rho_0|dq_{\phi}|_{h_0}^2-\langle dq,[\varphi_1d\varphi_2-\varphi_2d\varphi_1]\rangle_{h_0})dA_{h_0}.
\end{eqnarray*}
As a consequence, setting 
\begin{equation}\label{uphi.def}
u_{\phi}:=v_{\phi}+\log(\phi)\in L^1(M)^G,
\end{equation}
by adding $8\pi^2\int_{\Sigma}\langle d\log(\phi),d(\chi\rho_0)\rangle_{h_0}dA_{h_0}$ to both sides above, we find that
\begin{eqnarray*}
8\pi^2\int_{\Sigma}u_{\phi}(-\Delta_{h_0}(\chi\rho_0))&=&2\int_M\chi \phi^2 f_{\phi}dV_{g_0}+8\pi^2\int_{\Sigma}\chi\langle d\log(\phi),d\rho_0\rangle_{h_0}dA_{h_0}\\
    &&+\int_M(6|d\log(\phi)|_{g_0}^2\chi-4\langle d\log(\phi),d\chi\rangle_{g_0}-R_{g_0}\chi)dV_{g_0}\\
    &&+8\pi^2\int_{\Sigma}\chi(\rho_0|dq_{\phi}|_{h_0}^2-\langle dq,[\varphi_1d\varphi_2-\varphi_2d\varphi_1]\rangle_{h_0})dA_{h_0}.
\end{eqnarray*}
To express the right-hand side solely via $0$th order terms in $\chi$, note that since $\chi\in C_c^{\infty}(\Sigma)$ is supported away from $\partial\Sigma$, we have
\begin{eqnarray*}
    -4\int_M\langle d\log(\phi),d\chi\rangle_{g_0}dV_{g_0}&=&-8\pi^2\int_{\Sigma}\langle \rho_0 d\phi^2,\phi^{-2}d\chi\rangle_{h_0}dA_{h_0}\\
    &=&-8\pi^2\int_{\Sigma}\langle d\rho,\phi^{-2}d\chi\rangle_{h_0}dA_{h_0}+8\pi^2\int_{\Sigma}\langle d\rho_0,d\chi\rangle_{h_0}dA_{h_0}.
\end{eqnarray*}
Further integration by parts, using that $\phi\in W^{1,2}_{loc}(\Sigma\setminus \partial\Sigma)$, $\inf \phi>0$, and \eqref{rho.wk}, then gives
\begin{eqnarray*}
    -4\int_M\langle d\log(\phi),d\chi\rangle_{g_0}dV_{g_0}&=&-8\pi^2\int_{\Sigma}\chi f_{\phi}\phi^2 \rho_0dA_{h_0}-16\pi^2\int_{\Sigma}\chi\langle d\rho,\phi^{-3}d\phi\rangle_{h_0}dA_{h_0}\\
    &&+8\pi^2\int_{\Sigma}(-\Delta_{h_0}\rho_0)\chi dA_{h_0}\\
    &=&-2\int_M \chi f_{\phi}\phi^2dV_{g_0}-8\int_M\chi|d\log(\phi)|_{g_0}^2dV_{g_0}\\
    &&-16\pi^2\int_{\Sigma}\chi\langle d\rho_0,d\log\phi\rangle_{h_0}dA_{h_0}+8\pi^2\int_{\Sigma}(-\Delta_{h_0}\rho_0)\chi dA_{h_0},
\end{eqnarray*}
and adding this back to the previous equation, we find
\begin{eqnarray*}
8\pi^2\int_{\Sigma}u_{\phi}(-\Delta_{h_0}(\chi\rho_0))&=&-8\pi^2\int_{\Sigma}\chi\langle d\log(\phi),d\rho_0\rangle_{h_0}dA_{h_0}+\int_M(-2|d\log(\phi)|_{g_0}^2\chi-R_{g_0}\chi)dV_{g_0}\\
    &&+8\pi^2\int_{\Sigma}(-\Delta_{h_0}\rho_0)\chi dA_{h_0}\\
    &&+8\pi^2\int_{\Sigma}\chi(\rho_0|dq_{\phi}|_{h_0}^2-\langle dq,\varphi_1d\varphi_2-\varphi_2d\varphi_1\rangle_{h_0})dA_{h_0}.
\end{eqnarray*}
Equivalently, expressing all integrals on $\Sigma$ and dividing through by $8\pi^2$ gives
\begin{eqnarray*}
    \int_{\Sigma}u_{\phi}(-\Delta_{h_0}(\chi\rho_0))&=&\int_{\Sigma}\chi[(-\Delta_{h_0}\rho_0)-\frac{R_{g_0}}{2}\rho_0-|d\log(\phi)|_{h_0}^2\rho_0-\langle d\log(\phi),d\rho_0\rangle_{h_0}]dA_{h_0}\\
    &&+\int_{\Sigma}\chi(\rho_0|dq_{\phi}|_{h_0}^2-\langle dq,\varphi_1d\varphi_2-\varphi_2d\varphi_1\rangle_{h_0})dA_{h_0}
\end{eqnarray*}
whenever $\chi\in C_c^{\infty}(\Sigma)$. In particular, we have the following.

\begin{lem}\label{lap.u}
    In the interior of $\Sigma$, the distributional Laplacian $\Delta_{h_0}u_{\phi}$ is given by
    \begin{eqnarray*}
        \Delta_{h_0}u_{\phi}&=&|d\log(\phi)|_{h_0}^2+\langle d\log(\phi),d\log(\rho_0)\rangle_{h_0}\\
        &&-|dq_{\phi}|^2-\langle dq_{\phi},d\log(\varphi_1/\varphi_2)\rangle_{h_0}+\frac{\Delta_{h_0}\rho_0}{\rho_0}+\frac{R_{g_0}}{2}.
    \end{eqnarray*}
\end{lem}
\begin{proof}
    First, recall that we already have $u_{\phi}\in L^1(M)^G$, or equivalently $u_{\phi}\in L^1(\Sigma,\rho_0dA_{h_0})$, so that integration against $u_{\phi}$ is indeed a well-defined distribution in the interior of $\Sigma$. 

    Then, for any $\eta\in C_c^{\infty}(\Sigma)$ supported in the interior of $\Sigma$, the quotient $\chi=\frac{\eta}{\rho_0}$ belongs to $C_c^{\infty}(\Sigma)$ as well, and the statement follows immediately from the preceding computations with this choice of $\chi$.
\end{proof}
\begin{remark}\label{lap.k.rel}
    Recall from Section \ref{sect:curv.formula} that
    \begin{eqnarray*}
        R_{g_0}&=&2\left(K_{h_0}-\frac{\Delta_{h_0}\varphi_1}{\varphi_1}-\frac{\Delta_{h_0}\varphi_2}{\varphi_2}-\frac{\langle d\varphi_1,d\varphi_2\rangle_{h_0}}{\rho_0}\right)\\
        &=&2K_{h_0}-2\rho_0^{-1}\Delta_{h_0}\rho_0+2\rho_0^{-1}\langle d\varphi_1,d\varphi_2\rangle_{h_0},
    \end{eqnarray*}
    so the final term in the preceding equation for $\Delta_{h_0}u_{\phi}$ can be rewritten as
    $$\frac{\Delta_{h_0}\rho_0}{\rho_0}+\frac{R_{g_0}}{2}=K_{h_0}+\frac{\langle d\varphi_1,d\varphi_2\rangle_{h_0}}{\rho_0}.$$
\end{remark}

Next, note that since $\log(\phi)\in W^{1,2}(M)$ and $q\in W^{1,2}(M)$, it follows immediately from Lemma \ref{lap.u} that
$$\Delta_{h_0}u\in L^1(\{\rho_0>\delta\})$$
for any $\delta>0$. As a consequence, we can apply the techniques of Brezis-Merle to obtain the following.

\begin{lem}\label{bm.lem}
    For each $\delta>0$ and $p\in [1,\infty)$, we have $e^{|u_{\phi}|}\in L^p(\{\rho_0>\delta\})$.
\end{lem}
\begin{proof}
First, cover $\{\rho_0>\delta\}$ by a finite collection of disks $D_{r_1}(x_1),\ldots, D_{r_k}(x_k)$ such that $D_{2r_j}(x_j)\subset \{\rho_0>\delta/2\}$. On each of the doubled disks $D_{2r_j}(x_j)$, we can then decompose
$$u_{\phi}=h_j+w_j\text{ on }D_{2r_j}(x_j),$$
where $h_j$ is harmonic on $D_{2r_j}(x_j)$, and $w_j$ is given by convolution of $\Delta_{h_0}u_{\phi}$ against the Green's function for the Laplacian on $D_{2r_j}(x_j)$ with Dirichlet boundary condition. 

Since $\Delta_{h_0}u_{\phi}\in L^1(\{\rho_0>\delta/2\})$, we can then appeal to Corollary 1 of \cite{BrezisMerle1991Uniform} to deduce that
$$e^{p|w_j|}\in L^1(D_{3r_j/2}(x_j))\text{ for every }p\in [1,\infty),$$
while standard interior estimates for harmonic functions yield
$$\|h_j\|_{C^0(D_{r_j}(x_j))}\leq Cr_j^{-2}\|h_j\|_{L^1(D_{3r_j/2}(x_j))}\leq C r_j^{-2}\int_{D_{3r_j/2}(x_j)}(|w_j|+|u_{\phi}|)<C_j.$$
Putting these together, we indeed have $e^{p|u_{\phi}|}\in L^1(D_{r_j}(x_j))$ for each of the disks $D_{r_1}(x_1),\ldots,D_{r_k}(x_k)$ in the cover, giving the desired estimate.
\end{proof}

With Lemma \ref{bm.lem} in place, we can now complete the proof of interior regularity with a relatively simple bootstrapping argument.

\begin{prop}
    For every $\delta>0$, we have $\phi,q_{\phi},u_{\phi}\in C^{\infty}(\{\rho_0>\delta\})$.
\end{prop}
\begin{proof}
    First, recall from Lemma \ref{rho.lap.lem} that $\rho\in W_0^{1,q}(\Sigma)$ for every $q\in [1,2)$, hence $\rho\in L^p(\Sigma)$ for every $p\in [1,\infty)$ by Sobolev embedding. Moreover, rewriting \eqref{rho.wk}, we have
    $$-\Delta_{h_0}\rho=C_{\phi}e^{2u_{\phi}}\rho.$$
    Combining this with Lemma \ref{bm.lem}, it follows that
    $$-\Delta_{h_0}\rho\in L^p(\{\rho_0>\delta\})$$
    for any $p\in [1,\infty)$ and $\delta>0$, which together with interior $L^p$ regularity for the Laplacian on $\{\rho_0>\delta\}$ implies 
    $$\rho\in W^{2,p}(\{\rho_0>2\delta\})$$
    for every $\delta>0$ and $p\in [1,\infty)$ as well. In particular, 
    $$\rho\in \bigcap_{p\in [1,\infty)}W^{2,p}_{loc}(\Sigma\setminus \partial\Sigma)\subset \bigcap_{\alpha\in (0,1)}C^{1,\alpha}_{loc}(\Sigma\setminus\partial\Sigma).$$
    And since $\rho_0^{-1}\in C^{\infty}(\{\rho_0>\delta\})$ for every $\delta>0$, the same qualitative regularity holds for $\phi^2=\rho/\rho_0$ on $\{\rho_0>\delta\}$.

    Next, by \eqref{q.wk}, we see that $q_{\phi}$ solves the equation
    $$div_{h_0}(2\rho \cdot dq_{\phi})=div_{h_0}(\phi^2(\varphi_1d\varphi_2-\varphi_2d\varphi_1))$$
    weakly on $\Sigma$. On any $\{\rho_0>\delta\}$, $\rho$ has a uniform lower bound $\rho\geq (\inf \phi^2)\delta>0$, and we've just seen that $\rho\in C^{1,\alpha}(\{\rho_0>\delta\})$ for every $\alpha\in (0,1)$. Thus, on each $\{\rho_0>\delta\}$, $q_{\phi}$ solves a uniformly elliptic equation with coefficients in $C^{0,\alpha}$, and classic interior Schauder estimates \cite[\S 6.1]{GilbargTrudinger2001Elliptic} give
    $$q_{\phi}\in C^{2,\alpha}_{loc}(\Sigma\setminus\partial\Sigma)\text{ for every }\alpha\in (0,1).$$

    Returning to Lemma \ref{lap.u}, it follows now that
    $$\Delta_{h_0}u_{\phi}\in C^{0,\alpha}(\{\rho_0>\delta\})$$
    for every $\delta>0$, and interior Schauder estimates again allow us to conclude that
    $$u_{\phi}\in C^{2,\alpha}_{loc}(\Sigma\setminus\partial\Sigma)$$
    for every $\alpha\in (0,1)$ as well. 

    Returning to \eqref{rho.wk}, we see now that 
    $$-\Delta_{h_0}\rho=C_{\phi}e^{2u_{\phi}}\rho\in C^{1,\alpha}(\{\rho_0>\delta\})$$
    for every $\delta>0$, so interior Schauder estimates let us upgrade regularity once again to conclude that $\rho\in C^{3,\alpha}_{loc}(\Sigma\setminus\partial\Sigma)$ for every $\alpha\in (0,1)$. From here, successively applying Schauder estimates with \eqref{q.wk}, Lemma \ref{lap.u}, and \eqref{rho.wk} shows that $u_{\phi},q_{\phi}$, and $\phi$ belong to $C^k(\{\rho_0>\delta\})$ for any $k\in \mathbb{N}$ and $\delta>0$, as claimed.   
\end{proof}

\subsection{The Einstein metric in the interior of $\Sigma$}\label{ein.sec}

Dropping the $\phi$ subscript from $q_{\phi}$, $v_{\phi}$, and $u_{\phi}$, we now define
\begin{equation}\label{fi.def}
f_1:=\phi e^q\varphi_1\text{, }f_2:=\phi e^{-q}\varphi_2,
\end{equation}
and consider the smooth metric defined on the interior region 
$$\mathring M:=M\setminus T^2(\partial\Sigma)$$
by 
$$g=g_{\phi}:=h+f_1^2d\theta_1^2+f_2^2d\theta_2^2,$$
where $h$ is the metric on $\mathring\Sigma:=\Sigma\setminus\partial\Sigma$ given by
$$h:=e^{2u}h_0=\phi^2e^{2v}h_0.$$

On $\mathring\Sigma$, it's then easy to see that \eqref{rho.wk} and \eqref{q.wk} can be expressed equivalently as
\begin{equation}\label{fib.ein}
    -div_h(f_1df_2)=-div_h(f_2df_1)=\lambda \rho
\end{equation}
for a common positive constant $\lambda>0$, ensuring $Ric_g(\partial_{\theta_i},\partial_{\theta_i})=\lambda f_i^2$.

Next, we compute the Gaussian curvature $K_h$ of $h$, using Lemma \ref{lap.u} and Remark \ref{lap.k.rel} to see that
\begin{eqnarray*}
    K_h&=&e^{-2u}(K_{h_0}-\Delta_{h_0}u)\\
    &=&-\frac{\langle d\varphi_1,d\varphi_2\rangle_h}{\rho_0}-|d\log(\phi)|_h^2-\langle d\log(\phi),d\log(\rho_0)\rangle_h\\
    &&+|dq|_h^2+\langle dq,d\log(\varphi_1/\varphi_2)\rangle_h,
\end{eqnarray*}
or equivalently,
\begin{equation}\label{base.k.eq}
    K_h=-\langle d\log(f_1),d\log(f_2)\rangle_h,
\end{equation}
Per the discussion in Section \ref{sect:curv.formula}, with \eqref{fib.ein} and \eqref{base.k.eq} established, we can check that $Ric_g=\lambda g$, once we know that the restriction $Ric_g|_{T\Sigma}$ of the Ricci tensor to $T\Sigma$ is pure trace. We turn now to the proof of this latter condition, which is ultimately a consequence of the gauge-invariance described in Proposition \ref{prop:Emin.indep.background}. 

\begin{thm}\label{ein.int}
    The metric $g$ solves $Ric_g=\lambda g$ on $\mathring M$.
\end{thm}
\begin{proof}

    Starting from the initial reference metric $g_0\in \cD$,  fix some $\delta>0$, and let $\chi\in C_c^{\infty}(\{\rho_0>\delta/2\})$ be a nonnegative cutoff function with $\chi\equiv 1$ on $\{\rho_0>\delta\}$. Then,  
    define a new reference metric $\tilde{g}_0\in \cD$ by setting
    $$\tilde{g}_0=\chi g_{\phi}+(1-\chi)g_0,$$
    so that
    $$\tilde{g}_0=g_{\phi}\text{ on }\{\rho_0>\delta\}.$$
    It follows from the proof of Proposition \ref{prop:Emin.indep.background} that the minimizer $\tilde{\phi}$ for $\cS_{\tilde{g}_0}$ then satisfies
    $$\tilde{\phi}\equiv 1\text{ on }\{\rho_0>\delta\},$$
    with
    $$v_{\tilde{\phi}}=q_{\tilde{\phi}}=0\text{ on }\{\rho_0>\delta\}$$
    as well.

    Next, let $\tilde{g}_t\in \cD$ be any smooth path such that
    $$\tilde{g}_t=\tilde{g}_0\text{ in }\{\rho_0<\delta\},$$
    $$\tilde{g}_t(\partial_{\theta_i},\partial_{\theta_j})=\tilde{g}_0(\partial_{\theta_i},\partial_{\theta_j}),$$
    and the associated area measures $dA_{\tilde{h}_t}$ on $\Sigma$ satisfy
    $$dA_{\tilde{h}_t}=dA_{\tilde{h}_0}$$
    for all $t$. Note that for any trace-free $\xi\in \Gamma(Sym(T*M))$, we can find such a perturbation with 
    $$\frac{d}{dt}|_{t=0}\tilde{h}_t=\xi.$$

    For any such path, it then follows from Proposition \ref{prop:Emin.indep.background} that
    \begin{equation}\label{ein.s.comp}
    \cS_{\tilde{g}_0}(\tilde{\phi})=\cI_{\tilde{g}_0}=\cI_{\tilde{g}_t}\leq \cS_{\tilde{g}_t}(\tilde{\phi}).
    \end{equation}
    Since $\tilde{\phi}$ is fixed and $\tilde{g}_t=\tilde{g}_0$ near $\partial\Sigma$, it is straightforward to check that the properties needed to ensure $\cS_{\tilde{g}_t}(\tilde{\phi})<\infty$ continue to hold for small values of $t$, so we can use Proposition \ref{alt.s.char} to recast \eqref{ein.s.comp} as
    \begin{eqnarray*}
        \frac{E_{\tilde{g}_0}(\tilde{\phi},0,0)+8\pi^2T_{\tilde{g}_0}(\tilde{\phi})+\int_M\log(f_{\tilde{\phi}}^0)f_{\tilde{\phi}}^0\tilde{\phi}^4 dV_{\tilde{g}_0}}{\int_M f_{\tilde{\phi}}^0\tilde{\phi}^4 dV_{\tilde{g}_0}}+\log\left(\int_Mf_{\tilde{\phi}}^0 \tilde{\phi}^4 dV_{\tilde{g}_0}\right)\\
        \leq \frac{E_{\tilde{g}_t}(\tilde{\phi},0,0)+8\pi^2T_{\tilde{g}_t}(\tilde{\phi})+\int_M\log(f_{\tilde{\phi}}^t)f_{\tilde{\phi}}^t\tilde{\phi}^4 dV_{\tilde{g}_t}}{\int_M f_{\tilde{\phi}}^t\tilde{\phi}^4 dV_{\tilde{g}_t}}+\log\left(\int_Mf_{\tilde{\phi}}^t \tilde{\phi}^4 dV_{\tilde{g}_t}\right),
    \end{eqnarray*}
    where $f^t_{\tilde{\phi}}:=\frac{-\Delta_{\tilde{h}_t}\tilde{\rho}}{\tilde{\phi}^4\rho_0}.$ Next, since $\tilde{g}_t=\tilde{g}_0$ on $\{\rho_0<\delta\}$ and $\tilde{\rho}=\tilde{\phi}^2\rho_0$ is fixed independent of $t$, we note that
    \begin{eqnarray*}   &&\int_Mf_{\tilde{\phi}}^t\tilde{\phi}^4dV_{\tilde{g}_t}=4\pi^2\int_{\Sigma}(-\Delta_{\tilde{h}_t}\tilde{\rho})dA_{\tilde{h}_t}\\
    &=&4\pi^2\int_{\{\rho_0<\delta\}}(-\Delta_{\tilde{h}_t}\tilde{\rho})dA_{\tilde{h}_t}-4\pi^2\int_{\partial\{\rho_0>\delta\}}\frac{\partial\tilde{\rho}}{\partial\nu_{\tilde{h}_t}}ds_{\tilde{h}_t}\\
    &=&4\pi^2\int_{\{\rho_0<\delta\}}(-\Delta_{\tilde{h}_0}\tilde{\rho})dA_{\tilde{h}_0}-4\pi^2\int_{\partial\{\rho_0>\delta\}}\frac{\partial\tilde{\rho}}{\partial\nu_{\tilde{h}_0}}ds_{\tilde{h}_0}
    \end{eqnarray*}
    is constant in $t$, so the comparison \eqref{ein.s.comp} simplifies further to
    \begin{eqnarray*}
        E_{\tilde{g}_0}(\tilde{\phi},0,0)+8\pi^2T_{\tilde{g}_0}(\tilde{\phi})+\int_M\log(f_{\tilde{\phi}}^0)f_{\tilde{\phi}}^0\tilde{\phi}^4 dV_{\tilde{g}_0}\\
        \leq E_{\tilde{g}_t}(\tilde{\phi},0,0)+8\pi^2T_{\tilde{g}_t}(\tilde{\phi})+\int_M\log(f_{\tilde{\phi}}^t)f_{\tilde{\phi}}^t\tilde{\phi}^4 dV_{\tilde{g}_t}.
    \end{eqnarray*}
    It remains to apply the first derivative test at $t=0$. 

    To treat the term $T_{\tilde{g}_t}(\tilde{\phi})$, we write
    $$T_{\tilde{g}_t}(\tilde{\phi})=\int_{\Sigma}\tilde{\rho}|dq^t|_{\tilde{h}_t}^2dA_{\tilde{h}_t},$$
    where $q^t\in Z_{\tilde{\phi}}$ is the unique weak solution of
    \begin{equation}\label{div.q.t}
        div_{\tilde{h}_t}[\tilde{\rho}(2dq^t+d\log(\tilde{\varphi}_1/\tilde{\varphi}_2))]=0
    \end{equation}
    and we recall that $q^0=0$ on $\{\rho_0>\delta\}$. It is straightforward to check that the path $t\mapsto q^t$ is differentiable in $Z_{\tilde{\phi}}$, and since $q^0=0$ on $\{\rho_0>\delta\}$, it follows that 
    $$\int_{\{\rho_0>\delta\}}\tilde\rho |dq^t|_{\tilde{h}_t}^2dA_{\tilde{h}_t}=O(t^2),$$
    as $t\to 0$. In particular, it follows that
    \begin{eqnarray*}
        T_{\tilde{g}_t}(\tilde{\phi})-T_{\tilde{g}_0}(\tilde{\phi})&=&\int_{\{\rho_0\leq \delta\}}\tilde{\rho}(|dq^t|_{\tilde{h}_t}^2-|dq^0|_{\tilde{h}_t}^2)dA_{\tilde{h}_t}+O(t^2)\\
        &=&\int_{\{\rho_0\leq \delta\}}\tilde{\rho}(|dq^t|_{\tilde{h}_0}^2-|dq^0|_{\tilde{h}_0}^2)dA_{\tilde{h}_0}+O(t^2),
    \end{eqnarray*}
    using that $\tilde{h}_t=\tilde{h}_0$ in $\{\rho_0\leq \delta\}$. Now, set $\dot{q}:=\frac{dq^t}{dt}|_{t=0}$, and note that differentiating \eqref{div.q.t} in $\{\rho_0\leq\delta\}$, where $\tilde{h}_t=\tilde{h}_0$, gives
    $$\int_{\Sigma}\tilde{\rho}\langle d\dot{q},d\eta\rangle_{\tilde{h}_0}dA_{\tilde{h}_0}=0$$
    for $\eta\in Z_{\tilde \phi}$ supported in $\{\rho_0\leq\delta\}$. Differentiating $T_{\tilde{g}_t}(\tilde{\phi})$ and recalling that $q^0$ is supported in $\{\rho_0\leq \delta\}$, we therefore have
    $$\frac{d}{dt}|_{t=0}T_{\tilde{g}_t}(\tilde{\phi})=2\int_{\{\rho_0\leq \delta\}}\tilde{\rho}\langle d\dot{q},dq^0\rangle_{\tilde{h}_0}dA_{\tilde{h}_0}=0.$$

    This reduces our first derivative test for \eqref{ein.s.comp} to
    \begin{equation}\label{eh.var}
    0=\frac{d}{dt}|_{t=0}\left(E_{\tilde{g}_t}(\tilde{\phi},0,0)+4\pi^2\int_{\Sigma}[\log(-\Delta_{\tilde{h}_t}\tilde{\rho})-\log(\tilde{\phi}^4\rho_0)](-\Delta_{\tilde{h}_t}\tilde{\rho})dA_{\tilde{h}_t}\right).
    \end{equation}
    Setting $\xi:=\frac{d}{dt}|_{t=0}\tilde{h}_t$ and recalling that $\xi$ is trace-free, we see that
    $$\frac{d}{dt}|_{t=0}(-\Delta_{\tilde{h}_t}\tilde{\rho})=div_{\tilde{h}_0}(\xi(d\tilde{\rho},\cdot)),$$
    which is clearly supported in $\{\rho_0>\delta\}$. On the other hand, in $\{\rho_0>\delta\}$, we know that
    $$-\Delta_{\tilde{h}_0}\tilde{\rho}=2\lambda \tilde{\rho}=2\lambda\tilde{\phi}^4\rho_0,$$
    and since $dA_{\tilde{h}_t}=dA_{\tilde{h}_0}$ by assumption, it follows that
    \begin{eqnarray*}
        &&\frac{d}{dt}|_{t=0}\int_{\Sigma}[\log(-\Delta_{\tilde{h}_t}\tilde{\rho})-\log(\tilde{\phi}^4\rho_0)](-\Delta_{\tilde{h}_t}\tilde{\rho})dA_{\tilde{h}_t}\\&=&\int_{\Sigma}div_{\tilde{h}_0}(\xi(d\tilde{\rho},\cdot))dA_{\tilde{h}_0}+\int_{\Sigma}\log(2\lambda)div_{\tilde{h}_0}(\xi(d\tilde{\rho},\cdot))dA_{\tilde{h}_0}=0,
    \end{eqnarray*}
    by applying the divergence theorem in $\{\rho_0>\delta\}$.

    Thus, \eqref{eh.var} reduces to the vanishing of $\frac{d}{dt}|_{t=0}E_{\tilde{g}_t}(\tilde{\phi},0,0)$. In particular, since $\tilde {g}_t$ is constant on $\{\rho_0<\delta\}$ and $\tilde{\phi}=1$ on $\{\rho_0\geq \delta\}$, we have
\begin{eqnarray*}
    \frac{d}{dt}|_{t=0}E_{\tilde{g}_t}(\tilde{\phi})&=&\frac{d}{dt}|_{t=0}\int_{\{\rho_0>\delta\}}R_{\tilde{g}_t}dV_{\tilde{g}_t}\\
    &=&\int_{\{\rho_0>\delta\}}\langle \frac{R_{\tilde{g}_0}}{2}\tilde{g}_0-Ric_{\tilde{g}_0},\xi\rangle_{\tilde{g}_0}dV_{\tilde{g}_0}\\
    &=&-\int_{\{\rho_0>\delta\}}\langle Ric_{\tilde{g}_0},\xi\rangle_{\tilde{g}_0}dV_{\tilde{g}_0},
\end{eqnarray*}
    so that on $\{\rho_0>\delta\}$, $Ric_{\tilde{g}_0}=Ric_g$ is orthogonal to every trace-free symmetric two-tensor field $\xi$ along $T\Sigma$. Hence, $Ric_g=\lambda g$ in $\{\rho_0>\delta\}$ for every $\delta>0$, as claimed.   
\end{proof}

Thus, we now know that $g$ defines a smooth, positive Einstein metric on $\mathring M$, leaving us to check that $g$ extends smoothly across the $G$-invariant spheres and fixed points making up $M\setminus \mathring M$. 

\section{Key identities and a priori estimates in $(\mathring M,g)$.}\label{sect:comp.sec}

As preparation for proving regularity of $g$ across $M\setminus \mathring M$, we first use this section to collect a number of useful identities and a priori estimates in $\mathring M$ that will see repeated use in the final regularity arguments.

\subsection{The Hessian identities and refined Bochner formulas}\label{sec.hess.id}

First, we introduce the notation
\begin{equation}
\alpha_i:=\log(f_i)\text{ for }i=1,2,
\end{equation}
and
\begin{equation}
    \psi:=\alpha_1+\alpha_2=\log(\rho);\text{ }\beta:=\alpha_1-\alpha_2=\log(f_1/f_2).
\end{equation}
Since the Laplacian $\Delta_g$ on $T^2$-invariant functions is given by
\begin{equation}\label{lap.cvt}
\Delta_g\chi=\rho^{-1}div_h(\rho d\chi)=\Delta_h\chi+\langle d\chi,d\psi\rangle_h,
\end{equation}
we see that
\begin{equation}\label{alpha.lap}
    -\Delta_g\alpha_1=-\Delta_g\alpha_2=\lambda,
\end{equation}
and consequently $-\Delta_g\psi=2\lambda$ and $-\Delta_g\beta=0$ on $\mathring M$.

Next, using the computations of Section \ref{sect:curv.formula}, we observe that the restriction of the Einstein condition $Ric_g=\lambda g$ to $T\Sigma$ is equivalent to the Hessian identity
\begin{equation}\label{hess.id.0}
    \frac{Hess_h f_1}{f_1}+\frac{Hess_h f_2}{f_2}=(K_h-\lambda) h=-(\langle d\alpha_1,d\alpha_2\rangle_h+\lambda)h.
\end{equation}
In particular, noting that
$$Hess_h\alpha_i=\frac{Hess_h f_i}{f_i}-d\alpha_i\otimes d\alpha_i,$$
we can rewrite \eqref{hess.id.0} as an identity
\begin{equation}\label{hess.eq}
    Hess_h(\psi)=-d\alpha_1\otimes d\alpha_1-d\alpha_2\otimes d\alpha_2-(\lambda+\langle d\alpha_1,d\alpha_2\rangle)h,
\end{equation}
for the Hessian of $\psi=\log(\rho)$. Equivalently, we have
$$Hess_h(\psi)=-\frac{1}{2}(d\psi\otimes d\psi+d\beta\otimes d\beta)-\frac{1}{4}(|d\psi|^2-|d\beta|^2+4\lambda)h,$$
where for the remainder of the section, all norms will be taken with respect to $h$ unless otherwise specified. As an initial application, we will use \eqref{hess.eq} to sharpen some estimates obtained via standard Bochner identities for $\psi$.

Next, note that by \eqref{alpha.lap} and the standard Bochner identities on $(\mathring M,g)$, we have
$$\frac{1}{2}\Delta_g|d\alpha_i|^2=|Hess_g(\alpha_i)|^2+\lambda |d\alpha_i|^2$$
and
$$\Delta_g\frac{1}{2}|d\psi|^2=|Hess_g(\psi)|^2+\lambda |d\psi|^2.$$

Moreover, for any $T^2$-invariant function $\chi:\mathring M\to \mathbb{R}$, we can expand the $4$-dimensional Hessian $Hess_g(\chi)$ in terms of the $2$-dimensional Hessian $Hess_h(\chi)$ on $\Sigma$ and the connection terms in the $T^2$ fibers, giving
$$Hess_g(\chi)=Hess_h(\chi)+f_1^2\langle d\alpha_1,d\chi\rangle d\theta_1^2+\langle f_2^2d\alpha_2,d\chi\rangle d\theta_2^2,$$
so that
$$|Hess_g(\chi)|^2=|Hess_h(\chi)|^2+\langle d\alpha_1,d\chi\rangle^2+\langle d\alpha_2,d\chi\rangle^2.$$
In particular, applying this in the preceding Bochner identities, we see that 

\begin{equation}\label{psi.boch}
    \Delta_g\frac{1}{2}|d\psi|^2=|Hess_h(\psi)|^2+\lambda |d\psi|^2+\langle d\alpha_1,d\psi\rangle^2+\langle d\alpha_2,d\psi\rangle^2
\end{equation}
and
\begin{equation}\label{alpha.boch}
    \frac{1}{2}\Delta_g|d\alpha_i|^2=|Hess_h(\alpha_i)|^2+\lambda |d\alpha_i|^2+\langle d\alpha_1,d\alpha_i\rangle^2+\langle d\alpha_2,d\alpha_i\rangle^2.
\end{equation}

For many of our applications, it will be useful to combine these identities with the following elementary 2-dimensional refined Kato-type identity. 

\begin{claim}\label{kato.lem}
    For any open subset $U\subset \Sigma$ and $\chi\in C^{\infty}(U)$, we have the pointwise identity
    $$|Hess_h(\chi)|^2-2|d|d\chi||^2=(\Delta_h\chi)[\Delta_h\chi-2|d\chi|^{-1}\langle d|d\chi|,d\chi\rangle]$$
    on $U\cap\{|d\chi|>0\}$.
\end{claim}
\begin{proof}
    The proof is a straightforward algebraic computation, but we include it for completeness. At an arbitrary point $p\in U$ where $d\chi(p)\neq 0$, let $A=Hess_h(\chi): T_p\Sigma\to T_p\Sigma$ be the symmetric linear map associated to $Hess_h(\chi)$, let $z=\frac{\nabla \chi}{|\nabla\chi|}$, and let $e_1,e_2\in T_p\Sigma$ be an orthonormal basis of eigenvectors for $A$ such that $Ae_i=\kappa_ie_i$. Then writing $z=z_1e_1+z_2e_2$, we check directly that
    \begin{eqnarray*}
        |Hess_h(\chi)|^2-2|d|d\chi||^2&=&|A|^2|z|^2-2|Az|^2\\
        &=&(\kappa_1^2+\kappa_2^2)|z|^2-2(\kappa_1^2z_1^2+\kappa_2^2z_2^2)\\
        &=&(\kappa_2^2-\kappa_1^2)z_1^2+(\kappa_1^2-\kappa_2^2)z_2^2\\
        &=&(\kappa_1+\kappa_2)((\kappa_2-\kappa_1)z_1^2+(\kappa_1-\kappa_2)z_2^2)\\
        &=&tr(A)(tr(A)|z|^2-2\langle Az,z\rangle),
    \end{eqnarray*}
    which gives the desired identity.
\end{proof}

As an application, we find a useful pointwise lower bound for the Laplacian of the logarithm $\log(|d\psi|_h^2+\lambda)$.

\begin{lem}\label{psi.boch.4d}
    On $\{\rho_0>0\}$, we have
    $$\frac{1}{2}\Delta_g\log(|d\psi|^2+\lambda)\geq \lambda.$$
\end{lem}
\begin{proof}
    To begin, we use \eqref{psi.boch} to compute
    \begin{eqnarray*}
        &&\frac{1}{2}\Delta_g\log(|d\psi|^2+\lambda)\\&=&\frac{(1/2)\Delta_g|d\psi|^2}{|d\psi|^2+\lambda}-\frac{2|d\psi|^2|d|d\psi||^2}{(|d\psi|^2+\lambda)^2}\\
        &=&\frac{(|d\psi|^2+\lambda)(|Hess_h(\psi)|^2+\lambda |d\psi|^2+\langle d\alpha_1,d\psi\rangle^2+\langle d\alpha_2,d\psi\rangle^2)}{(|d\psi|^2+\lambda)^2}
        -\frac{2|d\psi|^2|d|d\psi||^2}{(|d\psi|^2+\lambda)^2}.
    \end{eqnarray*}
    Next, we can use Claim \ref{kato.lem} together with the fact that
    $$\Delta_h\psi=-2\lambda-|d\psi|^2$$
    to expand 
$$|Hess_h(\psi)|^2=2|d|d\psi||^2+(2\lambda+|d\psi|^2)^2+2(2\lambda+|d\psi|^2)|d\psi|^{-2}Hess_h(\psi)(d\psi,d\psi),$$   
so that
\begin{eqnarray*}
    &&(|d\psi|^2+\lambda)|Hess_h(\psi)|^2-2|d\psi|^2|d|d\psi||^2\\
    &=&\lambda |Hess_h(\psi)|^2+(2\lambda+|d\psi|^2)^2|d\psi|^2
    +2(2\lambda+|d\psi|^2)Hess_h(\psi)(d\psi,d\psi)\\
    &=&\lambda|Hess_h(\psi)+d\psi\otimes d\psi|^2-\lambda|d\psi|^4+(2\lambda+|d\psi|^2)^2|d\psi|^2\\
    &&+2(\lambda+|d\psi|^2)Hess_h(\psi)(d\psi,d\psi),
\end{eqnarray*}
and applying this in the previous identity gives
\begin{eqnarray*}
    \frac{1}{2}\Delta_g\log(|d\psi|^2+\lambda)&=&\frac{(|d\psi|^2+\lambda)(\lambda |d\psi|^2+\langle d\alpha_1,d\psi\rangle^2+\langle d\alpha_2,d\psi\rangle^2)}{(|d\psi|^2+\lambda)^2}\\
    &&+\frac{(2\lambda+|d\psi|^2)^2|d\psi|^2+2(\lambda+|d\psi|^2)Hess_h(\psi)(d\psi,d\psi)}{(|d\psi|^2+\lambda)^2}\\
        &&+\frac{\lambda|Hess_h(\psi)+d\psi\otimes d\psi|^2-\lambda |d\psi|^4}{(|d\psi|^2+\lambda)^2}.
\end{eqnarray*}
Next, it follows from \eqref{hess.eq} that
$$Hess_h(\psi)(d\psi,d\psi)=-\langle d\alpha_1,d\psi\rangle^2-\langle d\alpha_2,d\psi\rangle^2-(\lambda+\langle d\alpha_1,d\alpha_2\rangle)|d\psi|^2,$$
and also
$$Hess_h(\psi)+d\psi\otimes d\psi=d\alpha_1\otimes d\alpha_2+d\alpha_2\otimes d\alpha_1-(\lambda+\langle d\alpha_1,d\alpha_2\rangle)h,$$
so that
\begin{equation}\label{nice.id}
|Hess_h(\psi)+d\psi\otimes d\psi|^2=2|d\alpha_1|^2|d\alpha_2|^2+2\lambda^2.
\end{equation}
Applying these in the Bochner computation above, we then see that
\begin{eqnarray*}
    \frac{1}{2}\Delta_g\log(|d\psi|^2+\lambda)&=&\frac{(|d\psi|^2+\lambda)(\lambda |d\psi|^2-\langle d\alpha_1,d\psi\rangle^2-\langle d\alpha_2,d\psi\rangle^2)}{(|d\psi|^2+\lambda)^2}\\
    &&+\frac{(2\lambda+|d\psi|^2)^2|d\psi|^2-2(\lambda+|d\psi|^2)|d\psi|^2(\lambda+\langle d\alpha_1,d\alpha_2\rangle)}{(|d\psi|^2+\lambda)^2}\\
        &&+\frac{\lambda(2|d\alpha_1|^2|d\alpha_2|^2+2\lambda^2)-\lambda |d\psi|^4}{(|d\psi|^2+\lambda)^2}\\
        &=&\frac{|d\psi|^2(|d\psi|^4-\langle d\alpha_1,d\psi\rangle^2-\langle d\alpha_2,d\psi\rangle^2-2|d\psi|^2\langle d\alpha_1,d\alpha_2\rangle)}{(|d\psi|^2+\lambda)^2}\\
        &&+\frac{\lambda(2|d\psi|^4-\langle d\alpha_1,d\psi\rangle^2-\langle d\alpha_2,d\psi\rangle^2-2|d\psi|^2\langle d\alpha_1,d\alpha_2\rangle)}{(|d\psi|^2+\lambda)^2}\\
        &&+\frac{2\lambda |d\alpha_1|^2|d\alpha_2|^2+3\lambda^2|d\psi|^2+2\lambda^3}{(|d\psi|^2+\lambda)^2}.
\end{eqnarray*}
Finally, recognizing that 
\begin{equation}\label{form.norm}
|d\psi|^4-\langle d\alpha_1,d\psi\rangle^2-\langle d\alpha_2,d\psi\rangle^2-2|d\psi|^2\langle d\alpha_1,d\alpha_2\rangle=|d\psi\wedge d\alpha_1|^2+|d\psi\wedge d\alpha_2|^2,
\end{equation}
we arrive at the identity
\begin{eqnarray*}
    \frac{1}{2}\Delta_g\log(|d\psi|^2+\lambda)&=&\frac{|d\psi\wedge d\alpha_1|^2+|d\psi\wedge d\alpha_2|^2}{(|d\psi|^2+\lambda)}\\
    &&+\frac{\lambda (|d\psi|^4+3\lambda|d\psi|^2+2\lambda^2+2|d\alpha_1|^2|d\alpha_2|^2)}{(|d\psi|^2+\lambda)^2}.
\end{eqnarray*}
Simplifying slightly, we see that
\begin{eqnarray*}
    \frac{1}{2}\Delta_g\log(|d\psi|^2+\lambda)&=&\frac{|d\psi\wedge d\alpha_1|^2+|d\psi\wedge d\alpha_2|^2}{(|d\psi|^2+\lambda)}\\
    &&+\frac{\lambda^2|d\psi|^2+\lambda^3+2\lambda |d\alpha_1|^2|d\alpha_2|^2}{(|d\psi|^2+\lambda)^2}+\lambda,
\end{eqnarray*}
or equivalently,
\begin{equation}\label{psi.boch.sharp}
    \frac{1}{2}\Delta_g\log(|d\psi|^2+\lambda)=\lambda+\frac{|d\psi\wedge d\alpha_1|^2+|d\psi\wedge d\alpha_2|^2+\lambda^2}{|d\psi|^2+\lambda}+\frac{2\lambda |d\alpha_1|^2|d\alpha_2|^2}{(|d\psi|^2+\lambda)^2}
\end{equation}
from which the desired bound follows immediately.
\end{proof}

\subsection{A global a priori bound for $|d\psi|_{h_0}^2$}\label{as.bd.sec}

By combining the estimates of the previous subsection with an Ahlfors-Schwarz-type comparison argument \cite{Ahlfors1938Schwarz} (see e.g. \cite[Theorem 9.5]{BeardonMinda2007HyperbolicMetricTheory}), we are able to obtain the following \emph{global} pointwise estimate for the $h_0$-norm $|d\psi|_{h_0}^2$, which will play a key role in the proof of regularity across $\partial \Sigma$.

\begin{prop}\label{dpsi.above}
    With respect to a fixed reference metric $h_0$ on $S$, we have 
    $$|d\psi|_{h_0}^2(x)+\lambda e^{2u(x)}\leq C dist_{h_0}(x,\partial\Sigma)^{-2}.$$
\end{prop}

The starting point for the proof is the following refinement of Lemma \ref{psi.boch.4d}. 

\begin{lem}\label{as.bd}
    Writing $h=e^{2u}h_0$ with respect to a flat reference metric $h_0$ on some domain $U\subset \{\rho_0>0\}$, we have the pointwise estimate
    $$\frac{1}{2}\Delta_h\log(|d\psi|_h^2+\lambda)+\Delta_h u\geq \frac{3}{4}(|d\psi|_h^2+\lambda).$$
    Equivalently,
    $$\frac{1}{2}\Delta_{h_0}\log(|d\psi|_{h_0}^2+\lambda e^{2u})\geq \frac{3}{4}(|d\psi|_{h_0}^2+\lambda e^{2u}).$$
\end{lem}
\begin{proof}
    Recall that for any $T^2$-invariant function $\chi$, we have 
    $$\Delta_g\chi=\rho^{-1}div_h(\rho d\chi)=\Delta_h\chi+\langle d\psi,d\chi\rangle.$$
    We apply this to $\chi=\frac{1}{2}\log(|d\psi|^2+\lambda)$ and use the identity preceding \eqref{psi.boch.sharp} to obtain
    \begin{eqnarray*}
        \frac{1}{2}\Delta_h\log(|d\psi|^2+\lambda)&=& \frac{1}{2}\Delta_g\log(|d\psi|^2+\lambda)-\frac{Hess_h(\psi)(d\psi,d\psi)}{|d\psi|^2+\lambda}\\
        &=&\lambda+\frac{\lambda^2}{|d\psi|^2+\lambda}+\frac{2\lambda|d\alpha_1|^2|d\alpha_2|^2}{(|d\psi|^2+\lambda)^2}\\
        &&+\frac{|d\psi\wedge d\alpha_1|^2+|d\psi\wedge d\alpha_2|^2-Hess_h(\psi)(d\psi,d\psi)}{|d\psi|^2+\lambda}.
    \end{eqnarray*}
    Now, in the proof of Lemma \ref{psi.boch.4d}, we have already computed
    $$Hess_h(\psi)(d\psi,d\psi)=-\langle d\alpha_1,d\psi\rangle^2-\langle d\alpha_2,d\psi\rangle^2-(\lambda+\langle d\alpha_1,d\alpha_2\rangle)|d\psi|^2,$$
    and using \eqref{form.norm}, it follows that
    $$|d\psi\wedge d\alpha_1|^2+|d\psi\wedge d\alpha_2|^2-Hess_h(\psi)(d\psi,d\psi)=|d\psi|^2(|d\psi|^2+\lambda-\langle d\alpha_1,d\alpha_2\rangle).$$
    Using this in our computation for $\frac{1}{2}\Delta_h\log(|d\psi|^2+\lambda)$, we arrive at the identity
    \begin{eqnarray*}
        \frac{1}{2}\Delta_h\log(|d\psi|^2+\lambda)
        &=&\lambda+\frac{\lambda^2}{|d\psi|^2+\lambda}+\frac{2\lambda|d\alpha_1|^2|d\alpha_2|^2}{(|d\psi|^2+\lambda)^2}\\
        &&+|d\psi|^2-\frac{|d\psi|^2\langle d\alpha_1,d\alpha_2\rangle}{|d\psi|^2+\lambda}.
    \end{eqnarray*}
    On the other hand, writing $h=e^{2u}h_0$ with respect to a flat local reference metric $h_0$, the conformal factor $u$ satisfies
    $$\Delta_hu=\langle d\alpha_1,d\alpha_2\rangle_h,$$
    so that adding this term to the previous equation gives
    \begin{eqnarray*}
        \frac{1}{2}\Delta_h\log(|d\psi|^2+\lambda)+\Delta_hu&=&|d\psi|^2+\lambda +\frac{\lambda^2}{|d\psi|^2+\lambda}+\frac{2\lambda|d\alpha_1|^2|d\alpha_2|^2}{(|d\psi|^2+\lambda)^2}\\
        &&+\frac{\lambda\langle d\alpha_1,d\alpha_2\rangle}{|d\psi|^2+\lambda}.
    \end{eqnarray*}
    Young's inequality gives the simple estimate
    $$\frac{|\langle d\alpha_1,d\alpha_2\rangle|}{|d\psi|^2+\lambda}\leq \frac{|d\alpha_1|^2|d\alpha_2|^2}{(|d\psi|^2+\lambda)^2}+\frac{1}{4}$$
    for the final term, giving us the simple estimate
    $$\frac{1}{2}\Delta_h\log(|d\psi|^2+\lambda)+\Delta_hu\geq |d\psi|^2+\frac{3\lambda}{4}+\frac{\lambda^2}{|d\psi|^2+\lambda}+\frac{\lambda|d\alpha_1|^2|d\alpha_2|^2}{(|d\psi|^2+\lambda)^2},$$
    and in particular
    $$\frac{1}{2}\Delta_h\log(|d\psi|^2+\lambda)+\Delta_hu\geq \frac{3}{4}(|d\psi|^2+\lambda),$$
    as claimed. Multiplying both sides by $e^{2u}$ gives the equivalent formulation in terms of the reference metric $h_0$.
\end{proof}

With Lemma \ref{as.bd} in hand, we can now finish the proof of Proposition \ref{dpsi.above}.

\begin{proof}[Proof of Proposition \ref{dpsi.above}]
The proof uses an Ahlfors-Schwarz-type comparison argument; see e.g. the proof of \cite[Theorem 9.5]{BeardonMinda2007HyperbolicMetricTheory}. We include some details for completeness.

First, note that there is some $r_0\in (0,1)$ depending on the reference surface $(S,h_0)$ such that for any $r<r_0$, the geodesic disk $D_{r_0}^{h_0}(p)$ can be conformally identified with the flat unit disk $D_1^2(0)$, such that the flat metric $\tilde{h}_0$ satisfies $C^{-1}h_0\leq \tilde{h}_0\leq Ch_0$, with $C$ independent of $p$. 

In particular, given $r\in [0,r_0)$ and $x\in \Sigma$ with $\dist(x,\partial\Sigma)>Cr$, it follows that the flat disk of radius $r$ about $x$ satisfies
$$D^{\tilde{h}_0}_r(x)\subset \Sigma\setminus \partial\Sigma$$
as well. 

With respect to the Euclidean coordinates centered at $p$, and for $a>0$, let
$$\zeta(x):=\log(a)-\log(r^2-|x|^2)\text{ for }x\in D_r(p),$$
and observe that
$$\Delta_{\tilde{h}_0}\zeta(x)=\frac{4r^2}{(r^2-|x|^2)^2}=4r^2e^{2(\zeta(x)-\log(a))}=(2r/a)^2 e^{2\zeta(x)},$$
equivalent to the familiar fact that the conformal metric $e^{2\zeta}\tilde{h}_0$ has constant curvature $-(2r/a)^2$. In particular, fixing $a=4r$, so that
$$\zeta(x)=\log(4r/(r^2-|x|^2)),$$
it follows that
$$\Delta_{\tilde{h}_0}\zeta=\frac{1}{4}e^{2\zeta}.$$

Next, letting 
$$\xi:=\frac{1}{2}\log(|d\psi|_{\tilde{h}_0}^2+\lambda e^{2\tilde{u}}),$$
where $h=e^{2\tilde{u}}\tilde{h}_0$, Lemma \ref{as.bd} gives
$$\Delta_{\tilde{h}_0}\xi\geq \frac{3}{4}e^{2\xi},$$
so that
$$\Delta_{\tilde{h}_0}(\xi-\zeta)\geq\frac{3}{4}e^{2\xi}-\frac{1}{4}e^{2\zeta}>\frac{3}{4}(e^{2\xi}-e^{2\zeta}).$$
Finally, since $\xi$ is bounded on $D_r(p)\subset \Sigma\setminus\partial\Sigma$ while $\zeta\to +\infty$ along $\partial D_r(p)$, we see that $\xi-\zeta$ must achieve its maximum in the interior of $D_r(p)$, where the preceding inequality forces
$$0>\frac{3}{4}(e^{2\xi}-e^{2\zeta}).$$
Thus, we must have $\xi\leq \zeta$ throughout $D_r(p)$; in particular,
$$\xi(p)\leq \log(4r)-\log(r^2)=-\log(r)+\log(4),$$
so exponentiating and letting $r\to dist(p,\partial\Sigma)/C$ gives the desired pointwise bound.
\end{proof}

\section{A bootstrap to regularity assuming $v\leq C$}\label{sect:boot.v.leq.C}

The starting point for the regularity arguments in the sections that follow is the observation that a pointwise upper bound for $v$ is sufficient for the full regularity of the triple $(\phi,v,q)$ used to define the Einstein metric $g$. Namely, regularity will be achieved by applying the following proposition, first to domains of the form $W=\{\varphi_i>\delta\}$ supported away from the fixed points, and finally to all of $M$.

\begin{prop}\label{v.bd.suff}
Let $\phi\in W^{1,2}(M)^G$, $v=v_{\phi}\in L^1(M)^G$, and $q=q_{\phi}\in Z_{\phi}\subset W^{1,2}(M)$ be as in Section \ref{sect:int.reg}, and let $U\subset M$ be a $G$-invariant domain. If $v$ has an upper bound of the form
$$v\leq C_U$$
almost everywhere on $U$, then $\phi$, $v$, and $q$ are smooth in the interior of $U$. 
\end{prop}

The proof is based on the regularity theory for the six-dimensional operator introduced in Section \ref{l0.sec}. Namely, since $\rho=\phi^2\rho_0$ satisfies
$$-\Delta_{h_0}\rho=2\lambda e^{2v}\phi^4\rho_0$$
with respect to the reference metric $h_0$, it follows that the lift of $\phi^2$ to $(P,G_0)$ satisfies
\begin{equation}\label{l0.phi}
-L_0(\phi^2)=2\lambda e^{2v}\phi^4.
\end{equation}
A priori, the arguments of Section \ref{l0.sec} only justify this in a strong sense away from $\partial\Sigma$, but it is not hard to see that \eqref{l0.phi} holds globally in the weak sense as well, by virtue of the following estimates.

\begin{lem}
    We have $\phi^2\in W^{1,p}(P,G_0)$ for every $p\in [1,2)$, and $\phi\in L^6(P,G_0)$.
\end{lem}
\begin{proof}
    For the $W^{1,p}$ bound, fix any $p\in [1,2)$, and note that
    \begin{eqnarray*}
        \int_P|d(\phi^2)|_{G_0}^pdV_{G_0}&=&2^p16\pi^2\int_{\Sigma}\phi^p|d\phi|^p\rho_0^2dA_{h_0}\\
        &=&C\int_{\Sigma}(|d\phi|^2\rho_0)^{p/2}\phi^p\rho_0^{2-p/2}dA_{h_0}.
    \end{eqnarray*}
    Now, since $p<2$, we have
    $$\phi^p\rho_0^{2-p/2}\leq C\phi^p\rho_0^{p/2}=C\rho^{p/2},$$
    and an application of H\"older's inequality above gives
    $$\int_P|d(\phi^2)|_{G_0}^pdV_{G_0}\leq C\left(\int_{\Sigma}|d\phi|^2\rho_0dA_{h_0}\right)^{p/2}\cdot\|\rho\|_{L^{2/(2-p)}}^{p/2},$$
    where the finiteness of the right-hand side follows from the facts that $\phi\in W^{1,2}(M)$ and $\rho\in L^r(S)$ for every $r\in [1,\infty)$.

    Together with the six-dimensional Sobolev embedding on $P$, the above    shows that $\phi\in L^p(P)$ for every $p<6$. To get the $L^6$ estimate, we first use the $4$-dimensional Sobolev embedding $\phi\in W^{1,2}(M)\subset L^4(M)$ and Cauchy-Schwarz to see that
    $$\int_M|d(\phi^3)|_{g_0}dV_{g_0}\leq C\|\phi\|_{L^4(M)}^2\|d\phi\|_{L^2(M)},$$
    so that $\phi^3\in W^{1,1}(M)$. In particular, we can then apply Lemma \ref{mtos.w11} with $\chi=\phi^3$ to see that $\phi^3\rho_0\in W^{1,1}(S)$, and therefore $\phi^3\rho_0\in L^2(S)$, by the two-dimensional Sobolev embedding. That is, $\int_{\Sigma}\phi^6\rho_0^2<\infty$, which is equivalent to $\phi\in L^6(P)$.
\end{proof}

The fact that $\phi^2\in W^{1,p}(P,G_0)$ for all $p\in [1,2)$, together with the validity of \eqref{l0.phi} away from the \emph{codimension-three} subset of $P$ lying over $\partial\Sigma$, then suffices to ensure that \eqref{l0.phi} holds globally on $P$ in the weak sense.  

Using this, we now show how an upper bound for $v$ leads to $C^{1,\alpha}$ estimates for $\phi$, to begin the regularity bootstrap.

\begin{lem}
    Under the hypotheses of Proposition \ref{v.bd.suff}, we have $\phi\in C^{1,\alpha}_{loc}(U)$ for every $\alpha\in (0,1)$.
\end{lem}
\begin{proof}
    Since 
    $$-L_0(\phi^2)=-\Delta_{G_0}(\phi^2)-\frac{\Delta_{h_0}\rho_0}{\rho_0}\phi^2=2\lambda e^{2v}\phi^4$$
    weakly on $P$, if $v\leq C$ on the domain $U\subset M$, then, identifying $U$ with its lift in $P$, we see immediately that
    $$|L_0(\phi^2)|\leq C\phi^4\in L^{3/2}(U),$$
    since $\phi\in L^6(P)$. By interior $L^p$ regularity for $L_0$, it follows that 
    $$\phi^2\in W^{2,3/2}(U_1)\subset W^{1,2}(U_1)$$
    on any smaller domain $U_1\subset\subset U$.

    Fixing some smaller domain $U_1\subset\subset U$, we now see that
    \begin{equation}\label{p.l0.weak.eq}
    \int_{U_1}\langle d(\phi^2),d\chi\rangle_{G_0} dV_{G_0}\leq C'\int_{U_1} \phi^4 \chi dV_{G_0}
    \end{equation}
    for every $0\leq \chi\in W_0^{1,2}(U_1)$. Now, fix a nonnegative cutoff function $\zeta\in C_c^{\infty}(U_1)$, and for $k<\infty$, set
    $$\chi_k:=\zeta^2 \phi^2\min\{\phi,k\}^2\in W_0^{1,2}(U_1).$$
    Testing against $\chi_k$ in \eqref{p.l0.weak.eq} then gives
    \begin{eqnarray*}
        &&\int_{\{\phi<k\}}\zeta^28\phi^4|d\phi|^2+\int_{\{\phi\geq k\}}\zeta^2 k^2 4\phi^2|d\phi|^2\\
        &\leq &C'\int_{U_1}\zeta^2\phi^6\min\{\phi,k\}^2-4\int_{U_1}\phi^3\min\{\phi,k\}^2\zeta\langle d\phi,d\zeta\rangle,
    \end{eqnarray*}
    which, with a little rearrangement, gives an estimate of the form
    $$\int_{U_1}\zeta^2\min\{\phi,k\}^2\phi^2|d\phi|^2\leq C\int_{U_1}\zeta^2\phi^6\min\{\phi,k\}^2+C\int_{U_1}|d\zeta|^2\phi^6.$$
    In particular, setting
    $$f_k:=\zeta \min\{\phi,k\}\phi^2,$$
    we see that this is equivalent to an estimate of the form
    \begin{equation}\label{dir.fk.est}
    \int |df_k|_{G_0}^2dV_{G_0}\leq C\int_{U_1}[\phi^2f_k^2+|d\zeta|^2\phi^6]dV_{G_0}.
    \end{equation}
    We now want to show boundedness of $\|f_k\|_{W^{1,2}(P)}$ as $k\to\infty$. To this end, first note that the results of Section \ref{sect:int.reg} already give $\phi\in L_{loc}^{\infty}(\Sigma\setminus \partial\Sigma)$, so that $\phi$ is uniformly bounded on the complement of a measure zero subset of $U_1$. Combining this with the fact that $\phi\in L^6(P)$, it follows that for any $\delta>0$, we can write
    $$U_1:=V_{\delta}\cup W_{\delta},$$
    where 
    $$\phi^8\leq C_{\delta}\text{ on }V_{\delta}\text{ and }\int_{W_\delta}\phi^6<\delta^6.$$
    Combining this with a simple application of H\"older's inequality, we then estimate the first term on the right-hand side of \eqref{dir.fk.est} by 
    \begin{eqnarray*}
        \int_{U_1}\phi^2f_k^2&\leq & C_{\delta}+C\int_{W_{\delta}}\phi^2f_k^2\\
        &\leq &C_{\delta}+C\|\phi\|_{L^6(W_{\delta})}^2\|f_k\|_{L^3(U_1)}^{2}\\
        &\leq &C_{\delta}+C\delta^2\|f_k\|_{L^3(U_1)}^{2}.
    \end{eqnarray*}
    In particular, using the Sobolev embedding $W_0^{1,2}(U_1)\to L^3(U_1)$ to estimate the final term and returning to \eqref{dir.fk.est}, we arrive at an estimate of the form
    $$\|f_k\|_{W^{1,2}(U_1)}^2\leq C_{\delta}+C\delta^2\|f_k\|_{W^{1,2}(U_1)}^2+C\|\zeta\|_{C^1}^2\|\phi\|_{L^6(P)}^6.$$
    Taking $\delta>0$ small enough such that $C\delta^2<\frac{1}{2}$, it follows that $\|f_k\|_{W^{1,2}(P)}$ is uniformly bounded as $k\to\infty$, and therefore
    $$\zeta \phi^3\in W_0^{1,2}(U_1).$$

    Thus, for any smaller domain $U_2\subset\subset U_1$, it follows that 
    $$\phi^3\in W^{1,2}(U_2)\subset L^3(U_2),$$
    so that
    $$|L_0(\phi^2)|\leq C\phi^4\in L^{9/4}(U_2),$$
    and interior $L^p$ regularity then gives
    $$\phi^2\in W^{2,9/4}(U_3)\subset L^9(U_3)$$
    on any smaller domain $U_3\subset \subset U_2$. With one more iteration of this argument, we see that $L_0(\phi^2)\in L_{loc}^{9/2}(U)$, and consequently
    $$\phi^2\in W^{2,9/2}_{loc}(U)\subset C^{0,2/3}_{loc}(U).$$
    From here, we then see that $L_0(\phi^2)\in L^{\infty}_{loc}(U)$, and therefore
    $$\phi^2\in \bigcap_{p\in [1,\infty)}W^{2,p}_{loc}(U)\subset \bigcap_{\alpha\in (0,1)}C^{1,\alpha}_{loc}(U),$$
    as desired.      \end{proof}

With the $C^{1,\alpha}$ regularity of $\phi$ in place, the $C^{1,\alpha}$ regularity of $q$ also follows without much additional effort.

\begin{lem}
    Under the hypotheses of Proposition \ref{v.bd.suff}, we also have $q\in C^{1,\alpha}_{loc}(U)$ for every $\alpha\in [0,1)$.
\end{lem}

\begin{proof}
We begin by improving the regularity of $q$. Revisiting Lemmas \ref{q.char} and \ref{zphi.density} and Remark \ref{ell.rk} in light of the $C^{1,\alpha}$ regularity of $\phi$, we now see that
$$2div_{g_0}(\phi^2 dq)=-div_{g_0}(\phi^2d\log(\varphi_1/\varphi_2))$$
holds in the weak sense on $U$. In particular, on any $U_1\subset\subset U$, the operator
$$\chi\mapsto div_{g_0}(\phi^2d\chi)$$
becomes a uniformly elliptic second-order operator with $C^{\alpha}$ coefficients for any $\alpha\in (0,1)$. Moreover, though $\log(\varphi_1/\varphi_2)$ is singular along $\partial\Sigma$, by the symmetries of $\varphi_1$ and $\varphi_2$, it's easy to see that the regular part of the divergence
$$-div_{g_0}(d\log(\varphi_1/\varphi_2))=\frac{div_{h_0}(\varphi_1d\varphi_2-\varphi_2d\varphi_1)}{\rho_0}$$
belongs to $C^{\infty}(M)$ and has the same symmetry as $q$ (recall \eqref{eq:L.sym.q}). Next, to estimate the term
$$\langle d\phi^2,d\log(\varphi_1/\varphi_2)\rangle=\frac{\langle d\phi^2,d\varphi_1\rangle_{h_0}}{\varphi_1}-\frac{\langle d\phi^2,d\varphi_2\rangle_{h_0}}{\varphi_2},$$ 
note that the $C^1$ regularity and $G$-invariance of $\phi^2\in C^1(M)^G$ forces 
$$\langle d\phi^2,d\varphi_i\rangle=0\text{ along the edge }\{\varphi_i=0\}\cap U_1,$$
and since $\langle d\phi^2,d\varphi_i\rangle \in C^{0,\alpha}(U_1)$, it follows that
$$|\langle d\phi^2,d\varphi_i\rangle| \leq C_{\alpha}\varphi_i^{\alpha}$$
in $U_1$ for any $\alpha\in (0,1)$. Putting this together, it follows that
$$|2div_{g_0}(\phi^2dq)|\leq |\phi^2\frac{div_{h_0}(\varphi_1d\varphi_2-\varphi_2d\varphi_1)}{\rho_0}|+C_{\alpha} \varphi_1^{\alpha-1}+C_{\alpha}\varphi_2^{\alpha-1}$$
in $U_1$ for any $\alpha\in (0,1)$, and therefore
$$2div_{g_0}(\phi^2dq)\in L^p(U_1)$$
for any $p\in [1,\infty)$. Since $\phi^2\in C^{0,\alpha}(U_1)$, the operator is sufficiently regular to apply standard $L^p$ regularity results, to conclude that
$$q\in\bigcap_{[1,\infty)} W^{2,p}(U_2)\subset \bigcap_{\alpha\in (0,1)}C^{1,\alpha}(U_2)$$
for any subdomain $U_2\subset\subset U_1\subset\subset U$. 
\end{proof}

\subsection{Continuity of $u$ via the Hessian identity}

We continue to assume the hypotheses of Proposition \ref{v.bd.suff} are in force. With $C^{1,\alpha}$ estimates for $\phi$ and $q$ established, we next want to improve the regularity for $v$.

As a starting point, we revisit the Hessian identity \eqref{hess.id.0}, noting in particular the vanishing of the trace-free part
$$(f_1Hess_h(f_2)+f_2Hess_h(f_1))^{TF}=0$$
in $\mathring \Sigma$. Writing $h=e^{2u}h_0$, note that the Hessian transforms conformally like
$$Hess_h(f_i)=Hess_{h_0}(f_i)-(df_i\otimes du+du\otimes df_i)+\langle df_i,du\rangle_{h_0}h_0,$$
so that the preceding identity gives
\begin{equation}\label{base.hess}
    (f_1Hess_{h_0}(f_2)+f_2Hess_{h_0}(f_1))^{TF}=(d\rho\otimes du+du\otimes d\rho)^{TF}
\end{equation}
with respect to the reference metric $h_0$. 

Now, by the six-dimensional $W^{2,p}$ estimates for $\phi$ and four-dimensional $W^{2,p}$ estimates for $q$ on $U$, for any $U_1\subset\subset U$, we know already that
$$\int_{U_1}|Hess_{h_0}(\phi)|^p\rho_0^2+|Hess_{h_0}(q)|^p\rho_0<\infty$$
for any $p\in [1,\infty)$. In particular, since $\rho_0^{-s}\in L^1(\Sigma)$ for $s\in (0,1)$, a simple application of H\"older's inequality implies that $Hess_{h_0}(\phi)$ and $Hess_{h_0}(q)$ belong to $L^p(U_1,dA_{h_0})$ as well, and consequently
$$\int_{U_1}|Hess_{h_0}(f_1)|^p+|Hess_{h_0}(f_2)|^p<\infty$$
for every $p\in [1,\infty)$. On the other hand, since
$$|(d\rho\otimes du+du\otimes d\rho)^{TF}|_{h_0}^2=2|d\rho|_{h_0}^2|du|_{h_0}^2,$$
it follows from \eqref{base.hess} that
$$|du|_{h_0}\leq \frac{f_1|Hess_{h_0}(f_2)|+f_2|Hess_{h_0}(f_1)|}{|d\rho|_{h_0}}.$$
Moreover, we know that 
$$\frac{|d\rho_0|_{h_0}^2}{\varphi_1^2+\varphi_2^2}(x)\to 1$$
uniformly as $x\to \partial\Sigma$, and using this with the previous bound and the $C^{1,\alpha}$ regularity of $\phi,q$ in $U_1$, it follows that 
$$|du|_{h_0}\leq C(|Hess_{h_0}(f_1)|+|Hess_{h_0}(f_2)|).$$
Together with the preceding $L^p$ bounds for $f_i$ on $U_1$, it follows that
\begin{equation}
    u\in \bigcap_{p\in [1,\infty)}W^{1,p}(U_1)\subset \bigcap_{\alpha\in [0,1)}C^{0,\alpha}(U_1).
\end{equation}

\subsection{Proof of Proposition \ref{v.bd.suff}}

So far, we've seen that, under the hypotheses of Proposition \ref{v.bd.suff}, we have $\phi,q\in C^{1,\alpha}_{loc}(U)$, and $u\in C^{0,\alpha}_{loc}(U)$ for every $\alpha\in [0,1)$. Now we're in a position to complete the proof.

\begin{proof}[Proof of Proposition \ref{v.bd.suff}]

Now that we have $\phi\in C^{1,\alpha}_{loc}(U)$ and $u\in C^{0,\alpha}_{loc}(U)$, Schauder theory for the equation
$$-L_0(\phi^2)=2\lambda e^{2u}\phi^2$$
gives $\phi\in C^{2,\alpha}_{loc}(U)$ for every $\alpha\in [0,1)$. Revisiting the defining equation
$$2div_{g_0}(\phi^2dq)=-div_{g_0}(\phi^2d\log(\varphi_1/\varphi_2))$$
in light of the $C^{2,\alpha}$ estimates for $\phi^2$, it likewise follows that
$$q\in C^{2,\alpha}_{loc}(U)$$
as well, where the respective symmetries of $\phi$ and $\log(\varphi_1/\varphi_2)$ are again necessary for controlling the cross-terms $\langle d\phi^2,d\log(\varphi_1/\varphi_2)\rangle_{h_0}$.

In particular, this gives $f_1,f_2\in C^{2,\alpha}_{loc}(U)$, and returning to \eqref{base.hess}, we then deduce that $u\in C^{1,\alpha}_{loc}(U)$ for every $\alpha\in [0,1)$, and in particular $\frac{\partial u}{\partial\nu}=0$ along $\partial\Sigma$. Combining this observation with Lemma \ref{lap.u}, it then follows that
\begin{eqnarray*}
        \Delta_{h_0}u&=&|d\log(\phi)|_{h_0}^2+\langle d\log(\phi),d\log(\rho_0)\rangle_{h_0}\\
        &&-|dq|^2-\langle dq,d\log(\varphi_1/\varphi_2)\rangle_{h_0}+\frac{\Delta_{h_0}\rho_0}{\rho_0}+\frac{R_{g_0}}{2}
    \end{eqnarray*}
globally on the domain $\hat{U}\subset S$ given by the $L$-orbit $\hat{U}=L\cdot U$. Now that we have this equation globally in $\hat{U}$, we see that $\Delta_{h_0}u\in C^{0,\alpha}_{loc}(\hat{U})$, and can appeal to Schauder theory to obtain $u\in C^{2,\alpha}_{loc}(U)$. 

From here, the remainder of the argument is a straightforward induction on Schauder estimates, applied to the equations for $L_0(\phi^2)$, $div_{g_0}(\phi^2dq)$, and $\Delta_{h_0}u$.
\end{proof}

\section{Regularity away from isolated fixed points}\label{sect.reg.rods}

For each $\delta>0$, introduce the set
$$U_{\delta}:=\{\varphi_2>\delta\}\subset M,$$
which we will also frequently identify with its projection $U_{\delta}\subset \Sigma$. Our object in this section is the following theorem, establishing regularity of $g$ on $U_{\delta}$, hence also on $\{\varphi_1>\delta\}$ by symmetry. In particular, letting $\delta>0$ vary, this gives smoothness of the Einstein metric $g$ away from the $2n+2$ isolated fixed points of $G$.

\begin{prop}\label{rod.smth}
    The metric $g$ extends smoothly across $U_{\delta}$.
\end{prop}

The content of Proposition \ref{rod.smth} essentially consists of two statements:
\begin{itemize}
\item smoothness of the triple $(\phi,v,q)$ on $U_{\delta}$, and
\item the admissibility criterion $v=q$ on $\Gamma_1$,
\end{itemize}
where the latter is necessary to rule out edge-cone type singularities.

In fact, our proof follows the opposite order: first, we establish a weak form of the admissibility criterion $v=q$ in the trace sense along $\Gamma_1$, then use this as a key ingredient in the proof of an upper bound for $v$ in $U_{\delta}$, which allows us to conclude by invoking Proposition \ref{v.bd.suff}.

\subsection{Initial improvements and $0$th order estimates in $U_{\delta}$}

One obvious improvement of the analysis in $U_{\delta}$ over the global analysis on $M$ is that the volume element $dV_{g_0}=4\pi^2\varphi_1\varphi_2dA_{h_0}$ is uniformly comparable on $U_{\delta}$ to the volume element $dV_{Y_2}=2\pi \varphi_1 dA_{h_0}$ of the three-dimensional quotient $Y_2$, so that
$$c\int_{U_{\delta}}\chi dV_{g_0}\leq\int_{U_{\delta}}\chi dV_{Y_2}\leq \frac{1}{2\pi\delta}\int_{U_{\delta}}\chi dV_{g_0}$$
for any $T^2$-invariant $0\leq \chi\in L^1(U_{\delta},dV_{g_0})$. Moreover, note that since $0$ is a regular value of $\varphi_2$, the image of $U_{\delta}$ in $Y_2$ is a smooth manifold with boundary for every $\delta\in (0,\delta_0)$ sufficiently small.

A particularly useful consequence is that the $T^2$-invariant subspaces of the $4$-dimensional and $3$-dimensional Sobolev spaces $W^{1,2}_{g_0}(U_{\delta})$ and $W^{1,2}_{g_0^2}(U_{\delta})$ coincide, and we can apply the $3$-dimensional Sobolev embedding theorem to obtain the following.

\begin{lem}\label{3d.edge.sob}
    For every $0<\delta<\delta_0$, there is a constant $C_{\delta}<\infty$ such that
    $$\|\chi\|_{L^6(U_{\delta},dV_{g_0})}\leq C_{\delta}\|\chi\|_{W^{1,2}_{g_0}(U_{\delta})}$$
    for any $T^2$-invariant $\chi\in W^{1,2}_{g_0}(U_{\delta})$. Moreover, this constant can be taken uniform over any compact subinterval of $(0,\delta_0)$.
\end{lem}

As a first application of Lemma \ref{3d.edge.sob}, we can upgrade the regularity for $\rho$ described in Lemma \ref{rho.lap.lem} to local $W^{2,1}$ control.

\begin{lem}\label{rho.w21}
    For any $\delta>0$, on the two-dimensional domain $(U_{\delta},h_0)$, we have $\rho\in W^{2,1}_{h_0}(U_{\delta})$. In particular, $\rho\in W_{h_0}^{1,2}(U_{\delta})\cap C^0(U_{\delta})$, with $\rho|_{\Gamma_1\cap U_{\delta}}=0$ pointwise.
\end{lem}
\begin{proof}
    To begin, recall from Lemma \ref{rho.lap.lem} that $\rho$ satisfies $0\leq -\Delta_{h_0}\rho\in L^1(\Sigma)$, with 
    $$|\int_{\Sigma}(-\Delta_{h_0}\rho)v|\leq C[1+\int_M\phi^4 e^{2v}dV_{g_0}]^{1/2}$$
    for every $0\leq v\in L^{\infty}(\Sigma)$. Using Lemma \ref{3d.edge.sob} to see that $\phi\in L^6(U_{\delta},dV_{g_0})$, we can further bound
    $$\int_{U_{\delta}}\phi^4 e^{2v}dV_{g_0}\leq C_{\phi,\delta}\left(\int_{U_{\delta}}e^{6v}dV_{g_0}\right)^{1/3}.$$
    In particular, this gives us a bound of the form
    $$|\int_{\Sigma}(-\Delta_{h_0}\rho)v|\leq C+\left(\int_{M\setminus U_{\delta}}\phi^4 e^{2v}dV_{g_0}\right)^{1/2}+C_{\delta}\|e^{v}\|_{L^6(\Sigma)},$$
    so that
    $$|\int_{\Sigma}(-\Delta_{h_0}\rho)v|\leq C[1+\|e^v\|_{L^6(\Sigma)}]$$
    for any bounded, nonnegative $v$ supported in $U_{\delta}$. Applying this to
    $$v_k:=\frac{1}{6}\log(|\Delta_{h_0}\rho|+e)\mathbb{1}_{U_{\delta}}$$
    and letting $k\to\infty$, we deduce that
    $$\int_{U_{\delta}}\log(e+|\Delta_{h_0}\rho|)|\Delta_{h_0}\rho|dA_{h_0}<\infty,$$
    so that the restriction of $(-\Delta_{h_0}\rho)$ to the reflected union $\hat{U}_{\delta}=U_{\delta}\cup R_1U_{\delta}\subset S$ belongs to the Orlicz space $Llog L(\hat U_{\delta})$. 

    By classic Calderon-Zygmund estimates (see, e.g., \cite[\S 1.2]{Mingione2010NonlinearCalderonZygmund}), the $LlogL$ estimate for $\Delta_{h_0}\rho$ on $\hat{U}_{\delta}$ yields the desired improved regularity $\rho\in W^{2,1}(U_{2\delta})$.

    Moreover, since, in this two-dimensional setting, $W^{2,1}(\hat{U}_{2\delta})$ embeds into $C^0(\hat{U}_{2\delta})\cap W^{1,2}(\hat{U}_{2\delta})$, and $\rho\circ R_1=-\rho$, the remaining conclusions follow immediately. 
\end{proof}

For later use, we record the following simple consequences of the $W^{2,1}$ estimate. In what follows, let $\Phi_t$ denote the flow of the vector field $\frac{\nabla \rho_0}{|\nabla \rho_0|^2}$ in a fixed neighborhood of $\Gamma_1\cap U_{\delta}$.

\begin{lem}\label{w21.consq}
    For any $\delta>0$, there is a function $\kappa\in L^1(\Gamma_1\cap U_{\delta})$ such that the trace $d\rho|_{\Gamma_1\cap U_{\delta}}$ satisfies
    $$d\rho|_{\Gamma_1\cap U_{\delta}}=\kappa d\rho_0$$
    and
    $$\lim_{t\to 0}\|\phi^2\circ \Phi_t-\kappa\|_{L^1(\Gamma_1\cap U_{\delta})}=0.$$
     Moreover, $\int_{U_{\delta}}|\langle d\phi^2,d\rho_0\rangle|_{h_0}dA_{h_0}<\infty$.
\end{lem}

\begin{proof}
The fact that the trace $d\rho|_{\Gamma_1\cap U_{\delta}}$ is purely normal, hence of the form $\kappa d\rho_0$ for some $\kappa\in L^1(\Gamma_1\cap U_{\delta})$, follows immediately from the facts that $\rho\in W^{2,1}(U_{\delta})\cap W^{1,2}(U_{\delta})$ and $\rho=0$ along $\Gamma_1\cap U_\delta$ in the trace sense; this can be shown, for instance, by smooth approximation and continuity of the relevant traces.

Next, let $\delta_0>0$ be small enough that the flow 
$$\Phi: \Gamma_1\cap U_{2\delta}\times [0,\delta_0]\to U_{\delta}$$ of the vector field $\frac{\nabla \rho_0}{|\nabla\rho_0|^2}$ is smooth and well-defined for $t\in [0,\delta_0]$. For any $p\in \Gamma_1\cap U_{2\delta}$, we then have
$$\phi^2(\Phi_t(p))=t^{-1}\rho(\Phi_t(p))=\frac{1}{t}\int_0^t\frac{\langle d\rho,d\rho_0\rangle_{h_0}}{|d\rho_0|^2}(\Phi_s(p))ds.$$
Subtracting $\kappa(p)$ and integrating over $p\in \Gamma_1\cap U_{2\delta}$, since $d\rho=\kappa d\rho_0$ in the trace sense along $\Gamma_1\cap U_{\delta}$, it follows that
\begin{eqnarray*}
    \lim_{t\to 0}\int_{\Gamma_1\cap U_{2\delta}}|\phi^2\circ \Phi_t-\kappa|&\leq &\lim_{t\to 0}\frac{1}{t}\int_0^t\int_{p\in\Gamma_1\cap U_{2\delta}}|\frac{\langle d\rho,d\rho_0\rangle_{h_0}}{|d\rho_0|^2}(\Phi_s(p))-\kappa(p)| ds\\
    &=&0,
\end{eqnarray*}
as claimed.

It remains to check that $\langle d\phi^2,d\rho_0\rangle_{h_0}\in L^1(U_{2\delta})$, or equivalently
$$\langle d\phi^2,\frac{d\rho_0}{|d\rho_0|^2}\rangle \in L^1(U_{3\delta})$$
for any $\delta>0$. By regularity in $\mathring\Sigma$, it suffices to obtain the estimate in the flow neighborhood $$W:=\Phi(\Gamma_1\cap U_{2\delta}\times [0,\delta_0]).$$
To this end, fix $p\in \Gamma_1\cap U_{2\delta}$, and write
\begin{eqnarray*}
    \langle d\phi^2,\frac{d\rho_0}{|d\rho_0|^2}\rangle(\Phi_t(p))&=&\frac{d}{dt}[\phi^2\circ \Phi_t](p)=\frac{d}{dt}[t^{-1}\rho\circ \Phi_t](p)\\
    &=&t^{-2}[t\frac{d}{dt}(\rho\circ \Phi_t)-\rho\circ\Phi_t](p).
\end{eqnarray*}
Applying the fundamental theorem of calculus to the last line then gives
\begin{eqnarray*}
    |\langle d\phi^2,\frac{d\rho_0}{|d\rho_0|^2}\rangle(\Phi_t(p))|&=&t^{-2}|\int_0^ts\frac{d^2}{ds^2}(\rho\circ \Phi_s)(p)ds|\\
    &\leq & Ct^{-2}\int_0^t s(|Hess_{h_0}(\rho)|+|d\rho|)(\Phi_s(p))ds
\end{eqnarray*}
Integrating over $p\in \Gamma_1\cap U_{2\delta}$, and repeatedly using the coarea formula with the uniform bounds $0<c_{\delta}\leq |d\rho_0|\leq C$, it follows that
$$\int_{W\cap\{\rho_0=t\}}|\langle d\phi^2,d\rho_0\rangle|\leq Ct^{-2}\int_0^t s\left(\int_{\{\rho_0=s\}\cap U_{\delta}}(|Hess_{h_0}(\rho)|+|d\rho|)\right)ds.$$
Finally, integrating over $t\in [0,\delta_0]$ and applying Fubini's theorem gives
\begin{eqnarray*}
\int_W|\langle d\phi^2,d\rho_0\rangle|&\leq &C\int_0^{\delta_0}t^{-2}\int_0^t s\left(\int_{\{\rho_0=s\}\cap U_{\delta}}(|Hess_{h_0}(\rho)|+|d\rho|)\right)dsdt\\
&=&C\int_0^{\delta_0}s\left(\int_{\{\rho_0=s\}\cap U_{\delta}}(|Hess_{h_0}(\rho)|+|d\rho|)\right)\int_s^{\delta_0}t^{-2}dtds\\
&\leq &C\|\rho\|_{W^{2,1}(U_{\delta})},
\end{eqnarray*}
giving the desired estimate.
\end{proof}

Next, we turn to the function
$$\alpha_2=\log(f_2)=\log(\phi)-q+\log(\varphi_2)$$
in $U_{\delta}$. Recall that $\alpha_2$ solves
$$-\Delta_g\alpha_2=\lambda\text{ on }\mathring M,$$
so that
\begin{equation}\label{wk.a2.eq}\int_{U_{\delta}}\phi^2\langle d\alpha_2,d\eta\rangle_{g_0}dV_{g_0}=\lambda\int_{U_{\delta}}\eta e^{2v}\phi^4 dV_{g_0}
\end{equation}
for any $\eta\in W_0^{1,2}(U_{\delta})$ such that $\eta=0$ on a neighborhood of $\Gamma_1$. On the other hand, the estimates of Section \ref{s.fin.sec} and the smoothness of $\log(\varphi_2)$ on $U_{\delta}$ give
\begin{equation}
    \int_{U_{\delta}} \phi^2|d\alpha_2|_{g_0}^2dV_{g_0}\leq C_{\delta},
\end{equation}
and in the proof of Lemma \ref{zphi.density}, we constructed a family of log cutoff functions $\chi_{\epsilon}\in C_c^{\infty}(\mathring M)$ with $0\leq \chi_{\epsilon}\leq 1$ such that
$$\chi_{\epsilon}=0\text{ on }\{\rho_0<\epsilon^2\},\text{ }\chi_{\epsilon}=1\text{ on }\{\rho_0>\epsilon\},$$
and
$$\lim_{\epsilon\to 0}\int_M \phi^2|d\chi_{\epsilon}|^2dV_{g_0}=0.$$
As a consequence, for any function $\eta\in [W_0^{1,2}\cap L^{\infty}](U_{\delta})$ such that $$\int_{U_{\delta}}\phi^2|d\eta|^2dV_{g_0}<\infty,$$
testing \eqref{wk.a2.eq} against $\eta_{\epsilon}=\chi_{\epsilon}\eta$ gives
\begin{eqnarray*}
    \lambda \int_{U_{\delta}}\eta e^{2v}\phi^4 dV_{g_0}&=&\lambda\lim_{\epsilon\to 0}\int_{U_{\delta}}\eta_{\epsilon}e^{2v}\phi^4dV_{g_0}\\
    &=&\lim_{\epsilon\to 0}\int_{U_{\delta}}\phi^2\langle d\alpha_2,\chi_{\epsilon}d\eta+\eta d\chi_{\epsilon}\rangle_{g_0} dV_{g_0}\\
    &=&\int_{U_{\delta}}\phi^2\langle d\alpha_2,d\eta\rangle_{g_0}dV_{g_0}.
\end{eqnarray*}
Thus, we see that \eqref{wk.a2.eq} holds for any test function $\eta\in W_0^{1,2}(U_{\delta})\cap L^{\infty}(U_{\delta})$ such that
$\int \phi^2 |d\eta|_{g_0}^2dV_{g_0}<\infty$. Combining this with a Moser-iteration type argument, we show next that $\alpha_2$ has a uniform lower bound on $U_{\delta}$. Note that some care is needed here, given the lack of a priori upper bounds on $\phi$.

\begin{lem}\label{f2.below}
    For each $\delta>0$, we have $\inf_{U_{\delta}}\alpha_2>-C_{\delta}$ for some $C_{\delta}<\infty$.
\end{lem}
\begin{proof}
Fixing an arbitrary $\epsilon>0$, we can choose a nonnegative cutoff function $\zeta=\zeta(\varphi_2)\in C_c^{\infty}(U_{\delta})$ such that $\zeta\equiv 1$ on $U_{\delta+\epsilon}$ and $|d\zeta|\leq \frac{C}{\epsilon}$ on $\{\delta\leq \varphi_2\leq \delta+\epsilon\}$. 

For any bounded, Lipschitz function $F: \mathbb{R}\to \mathbb{R}$, it follows that $\eta=\zeta^2 F(\alpha_2)$ meets the criteria for \eqref{wk.a2.eq} to hold, so that
$$\int_{U_{\delta}}\phi^2[ F'(\alpha_2)|d\alpha_2|_{g_0}^2\zeta^2+2\zeta F(\alpha_2)\langle d\alpha_2,d\zeta\rangle_{g_0}]dV_{g_0}=\lambda\int_{U_{\delta}}\zeta^2F(\alpha_2) e^{2v}\phi^4dV_{g_0}.$$

In particular, if $F$ satisfies $F\geq 0$ and $F'\leq 0$, it follows that
$$\int \phi^2 \zeta^2|F'(\alpha_2)||d\alpha_2|_{g_0}^2dV_{g_0}\leq \int 2\phi^2\zeta F(\alpha_2)\langle d\alpha_2,d\zeta\rangle_{g_0}dV_{g_0}.$$
Following the usual Moser iteration scheme, fix $p\in [1,\infty)$ and observe that for any $k\in (0,\infty)$, the bounded Lipschitz functions
$$F_{p,k,m}(t):=\min\{(-t)^+\min\{(-t)^+,k\}^{p-1},m\}$$
are nonnegative and decreasing, so we can test them in the preceding inequality and take $m\to\infty$ to obtain
\begin{eqnarray*}
    &&\int_{\{\alpha_2<0\}}\phi^2\zeta^2\min\{|\alpha_2|,k\}^{p-1}|d\alpha_2|_{g_0}^2dV_{g_0}\\
    &\leq& C\int 2\phi^2\zeta |\alpha_2|\min\{|\alpha_2|,k\}^{p-1}\langle d\alpha_2,d\zeta\rangle_{g_0}dV_{g_0}\\
    &\leq& \frac{1}{2}\int_{\{\alpha_2<0\}} \phi^2\zeta^2\min\{|\alpha_2|,k\}^{(p-1)}|d\alpha_2|_{g_0}^2dV_{g_0}\\
    &&+\frac{C'}{\epsilon^2}\int_{\{\alpha_2<0\}\cap U_{\delta}} \phi^2\min\{|\alpha_2|,k\}^{p-1}|\alpha_2|^2dV_{g_0}.
\end{eqnarray*}
Absorbing the term on the penultimate line into the left-hand side, we end up with 
\begin{eqnarray*}
    \int_{\{\alpha_2<0\}\cap U_{\delta+\epsilon}}\phi^2\min\{|\alpha_2|,k\}^{p-1}|d\alpha_2|_{g_0}^2dV_{g_0}\\
    \leq \frac{C''}{\epsilon^2}\int_{\{\alpha_2<0\}} \phi^2\min\{|\alpha_2|,k\}^{p-1}|\alpha_2|^2dV_{g_0},
\end{eqnarray*}
and letting $k\to\infty$ gives
$$\int_{\{\alpha_2<0\}\cap U_{\delta+\epsilon}}\phi^2|\alpha_2|^{p-1}|d\alpha_2|_{g_0}^2dV_{g_0}\leq \frac{C}{\epsilon^2}\int_{\{\alpha_2<0\}\cap U_{\delta}}\phi^2|\alpha_2|^{p+1}dV_{g_0}.$$
Together with the lower bound $\inf \phi>0$, this leaves us with an estimate of the form
$$\int_{U_{\delta+\epsilon}}|\alpha_2^-|^{p-1}|d\alpha_2|_{g_0}^2\leq \frac{C}{\epsilon^2}\int_{U_{\delta}}\phi^2|\alpha_2^-|^{p+1},$$
where of course $\alpha_2^-:=|\min\{\alpha_2,0\}|$. In particular, this is equivalent to an estimate of the form
$$\|(\alpha_2^-)^{\frac{p+1}{2}}\|_{W^{1,2}(U_{\delta+\epsilon},g_0)}^2\leq \frac{Cp^2}{\epsilon^2}\int_{U_{\delta}}\phi^2|\alpha_2^-|^{p+1}dV_{g_0},$$
and applying Lemma \ref{3d.edge.sob} with $\chi=(\alpha_2^-)^{\frac{p+1}{2}}$, it follows that
$$\left(\int |\alpha_2^-|^{3(p+1)}dV_{g_0}\right)^{1/3}\leq \frac{C_{\delta+\epsilon}p^2}{\epsilon^2}\int_{U_{\delta}}\phi^2|\alpha_2^-|^{p+1}dV_{g_0}.$$
Finally, applying Lemma \ref{3d.edge.sob} once more to $\phi$ to see that
$$\|\phi\|_{L^6(U_{\delta},dV_{g_0})}\leq C_{\delta},$$
an application of H\"older's inequality on the right-hand side gives
$$\int_{U_{\delta}}\phi^2|\alpha_2^-|^{p+1}dV_{g_0}\leq C_{\delta}\|\phi\|_{L^6(U_{\delta},dV_{g_0})}^2\|\alpha_2^-\|_{L^{3(p+1)/2}}^{p+1}.$$
Combining this with the previous estimates, we see finally that
$$\|\alpha_2^-\|^{p+1}_{L^{3(p+1)}(U_{\delta+\epsilon},dV_{g_0})}\leq \frac{C_{\delta}p^2}{\epsilon^2}\|\alpha_2^-\|_{L^{3(p+1)/2}(U_{\delta},dV_{g_0})}^{p+1}$$
for every $p\geq 1$, whenever $0<\delta<\delta+\epsilon<\delta_0$. In particular, starting from an arbitrary fixed $\delta_1<\frac{\delta_0}{10}$ and iterating with 
$$\epsilon_k=2^{-(k+1)\delta}\text{ and }\delta_{k+1}=\delta_k+\epsilon_k\in [\delta_1,2\delta_1]\text{ for }k\geq 1,$$
we note by Lemma \ref{3d.edge.sob} that the constant $C_{\delta_k}$ has a uniform upper bound $C_{\delta_k}=C(\delta_1)$ independent of $k$, so
the standard conclusion of the iteration argument (as in, e.g., \cite[\S 8.6]{GilbargTrudinger2001Elliptic}) then gives the desired $L^{\infty}$ control for $\alpha_2^-$ in $U_{2\delta_1}$. Since $\delta_1$ was arbitrary, this completes the proof.
\end{proof}

In particular, it follows that $f_2$ has a universal lower bound of the form
$$f_2\geq c_{\delta}>0\text{ in }U_{\delta}.$$
Together with the continuity of $\rho$ in $U_{\delta}$ and the vanishing $\rho|_{\Gamma_1}=0$, this also implies that 
$$f_1(x)=\frac{\rho}{f_2}(x)\to 0$$
as $x\to \Gamma_1$ in $U_{\delta}$.

\subsection{Lower bounds for $u+\alpha_2$ and weak admissibility statement along $\Gamma_1\cap U_{\delta}$.}

Next, consider the function
$$w:=u+\alpha_2=v+2\log \phi-q+\log(\varphi_2)$$
in $U_{\delta}$, and recall that the admissibility criterion corresponds to the equality $v=q$ on $\Gamma_1$, so we \emph{expect} to show that $w-2\log(\phi)-\log(\varphi_2)\to 0$ in some sense along $\Gamma_1$. In this subsection, we will confirm a version of this limiting behavior, and also show that $w$ has a uniform lower bound in $U_{\delta}$.

As a starting point, note that in $\mathring\Sigma$, we've already seen that
$$\Delta_{h_0}u=K_{h_0}-e^{2u}K_h=K_{h_0}+\langle d\alpha_1,d\alpha_2\rangle_{h_0},$$
while
$$\Delta_{h_0}\alpha_2=-\lambda e^{2u}-\langle d\alpha_2,d\psi\rangle_{h_0}=-\lambda e^{2u}-|d\alpha_2|_{h_0}^2-\langle d\alpha_1,d\alpha_2\rangle_{h_0}.$$
Thus, in the interior $\mathring\Sigma$, we have
$$-\Delta_{h_0}w=\lambda e^{2u}+|d\alpha_2|_{h_0}^2-K_{h_0}.$$
The following technical lemma gives a weak version of the statement that $w$ satisfies this equation on $\Sigma\cap U_{\delta}$ with Dirichlet data $w|_{\Gamma_1}=\log(\kappa)+\log(\varphi_2)=2\log(\phi)+\log(\varphi_2)$. 
\begin{remark}
    Note that the convergence statement $\phi^2\circ \Phi_t\to \kappa$ in $L^1(\Gamma_1\cap U_{\delta})$ from Lemma \ref{w21.consq} implies the corresponding convergence $2\log(\phi)\circ \Phi_t\to \log(\kappa)$ for the logarithms, by virtue of the positive lower bound $\inf_{\Sigma}\phi>0$ and the dominated convergence theorem.
\end{remark} 

\begin{lem}\label{wkadm.lem}
    For any $T^2$-invariant $\chi\in C_c^{1,\alpha}(U_{\delta})$ such that $-\frac{\Delta_{h_0}(\chi\rho_0)}{\rho_0}=\zeta\in L^{\infty}(M)$, we have
    $$\int_{\Sigma}w \zeta\rho_0 dA_{h_0}=\int_{\Sigma}(\lambda e^{2u}+|d\alpha_2|_{h_0}^2-K_{h_0})\chi \rho_0-\int_{\Gamma_1}(2\log\phi+\log(\varphi_2))\frac{\partial}{\partial\nu}(\chi\rho_0),$$
    where we identify $2\log\phi$ with $\log\kappa$ on $\Gamma_1\cap U_{\delta}$.
\end{lem}
\begin{proof}
We begin with the weak Euler-Lagrange equation: for $\chi$ and $\zeta$ as above, recall that
\begin{eqnarray*}
\int_{\Sigma} v\zeta \rho_0dA_{h_0}&=& \int_{\Sigma}2\lambda e^{2u}\chi \rho_0dA_{h_0}-\int_{\Sigma}(3\langle d\phi,d(\chi/\phi)\rangle_{h_0}+\frac{R_{g_0}}{2}\chi)\rho_0 dA_{h_0}\\
&&+\int_{\Sigma}(\chi\rho_0|dq|_{h_0}^2+q \cdot div_{h_0}(\chi(\varphi_1d\varphi_2-\varphi_2d\varphi_1)))dA_{h_0}.
\end{eqnarray*}
Using the curvature identity
$$\frac{R_{g_0}}{2}=K_{h_0}-\frac{\Delta_{h_0}\rho_0}{\rho_0}+\langle d\log(\varphi_1),d\log(\varphi_2)\rangle_{h_0},$$
this is equivalent to
\begin{eqnarray*}
\int_{\Sigma} v\zeta \rho_0dA_{h_0}&=& \int_{\Sigma}2\lambda e^{2u}\chi \rho_0 dA_{h_0}-\int_{\Sigma}3\langle d\phi,d(\chi/\phi)\rangle_{h_0}\rho_0dA_{h_0}\\
&&-\int_{\Sigma}(K_{h_0}-\frac{\Delta_{h_0}\rho_0}{\rho_0}+\langle d\log(\varphi_1),d\log(\varphi_2)\rangle_{h_0})\chi\rho_0 dA_{h_0}\\
&&+\int_{\Sigma}(\chi\rho_0|dq|_{h_0}^2+q \cdot div_{h_0}(\chi(\varphi_1d\varphi_2-\varphi_2d\varphi_1)))dA_{h_0}.
\end{eqnarray*}

Next, a straightforward integration by parts gives
\begin{eqnarray*}
    \int_{\Sigma}\log(\varphi_2)\zeta \rho_0 dA_{h_0}&=&-\int_{\Gamma_1}\log(\varphi_2)\frac{\partial(\chi\rho_0)}{\partial\nu}ds_{h_0}+\int_{\Sigma}\langle d\log(\varphi_2), d(\chi\rho_0)\rangle_{h_0}dA_{h_0}\\
    &=&-\int_{\Gamma_1}\log(\varphi_2)\frac{\partial(\chi\rho_0)}{\partial\nu}ds_{h_0}+\int_{\Sigma}\chi\rho_0|d\log(\varphi_2)|_{h_0}^2dA_{h_0}\\
    &&+\int_{\Sigma}\langle d(\varphi_1\chi),d\varphi_2\rangle_{h_0}dA_{h_0}.
\end{eqnarray*}
Since $\frac{\partial \varphi_2}{\partial\nu}=0$, one more integration by parts in the last term leaves
\begin{eqnarray*}
    \int_{\Sigma}\log(\varphi_2)\zeta \rho_0dA_{h_0}&=&-\int_{\Gamma_1}\log(\varphi_2)\frac{\partial(\chi\rho_0)}{\partial\nu}ds_{h_0}+\int_{\Sigma}\chi\rho_0|d\log(\varphi_2)|_{h_0}^2dA_{h_0}\\
    &&-\int_{\Sigma}(\varphi_1\chi)\Delta_{h_0}\varphi_2dA_{h_0}.
\end{eqnarray*}
For the $\log(\phi)$ term, integration by parts likewise gives
$$\int_{\Sigma}2\log\phi \zeta\rho_0 dA_{h_0}=-\int_{\Gamma_1}(2\log \phi)\frac{\partial(\chi\rho_0)}{\partial\nu}+2\int_{\Sigma}\langle d\log(\phi),d(\chi\rho_0)\rangle_{h_0}dA_{h_0};$$
to justify this carefully, we can first integrate by parts over $\{\rho_0>\epsilon\}$, and appeal to Lemma \ref{w21.consq} to pass to the limit $\epsilon\to 0$. Likewise, since $q\in W^{1,2}_{g_0}(M)$, we can justify the first-order integration by parts
$$-\int_{\Sigma}q\zeta \rho_0 dA_{h_0}=\int_{\Sigma}q\Delta_{h_0}(\chi\rho_0)dA_{h_0}=\int_{\Sigma}[-\rho_0\langle d\chi,d q\rangle+q div_{h_0}(\chi d\rho_0)]dA_{h_0}.$$
Summing the preceding computations gives the initial cumbersome identity
\begin{eqnarray*}
    \int_{\Sigma}w\zeta \rho_0 dA_{h_0}&=&-\int_{\Gamma_1}[2\log(\phi)+\log(\varphi_2)]\frac{\partial(\chi\rho_0)}{\partial\nu}ds_{h_0}-\int_{\Sigma}K_{h_0}\chi\rho_0 dA_{h_0}\\
    &&+\int_{\Sigma}2\lambda e^{2u}\chi \rho_0 dA_{h_0}-\int_{\Sigma}3\langle d\phi,d(\chi/\phi)\rangle_{h_0}\rho_0dA_{h_0}\\
&&+\int_{\Sigma}(\Delta_{h_0}\rho_0-\langle d\varphi_1,d\varphi_2\rangle_{h_0}-\varphi_1\Delta_{h_0}\varphi_2)\chi dA_{h_0}\\
&&+\int_{\Sigma}(\chi\rho_0|dq|_{h_0}^2+2q \cdot div_{h_0}(\chi[\varphi_1d\varphi_2]))dA_{h_0}\\
&&+\int_{\Sigma}\chi\rho_0|d\log(\varphi_2)|_{h_0}^2dA_{h_0}\\
    &&+2\int_{\Sigma}\langle d\log(\phi),d(\chi\rho_0)\rangle_{h_0}dA_{h_0}-\int_{\Sigma}\rho_0\langle d\chi,d q\rangle dA_{h_0}
\end{eqnarray*}
Note that the first line on the right-hand side already includes  the desired boundary terms and $-K_{h_0}$ contribution, so we're left to show that everything left over coincides with $\int_{\Sigma}(\lambda e^{2u}+|d\alpha_2|_{h_0}^2)\chi \rho_0 dA_{h_0}$.

Next, note that since $q\in W^{1,2}_{Y_2}(U_{\delta})$, we can write
$$\int_{\Sigma}2q div_{h_0}(\varphi_1\chi d\varphi_2)=-2\int_{\Sigma}\chi\langle dq,\varphi_1d\varphi_2\rangle_{h_0},$$
and consolidate the third, fifth and sixth lines on the right-hand side above, leaving, after a little rearrangement,
\begin{eqnarray*}
    \int_{\Sigma}w\zeta \rho_0 dA_{h_0}&=&-\int_{\Gamma_1}[2\log(\phi)+\log(\varphi_2)]\frac{\partial(\chi\rho_0)}{\partial\nu}ds_{h_0}-\int_{\Sigma}K_{h_0}\chi\rho_0 dA_{h_0}\\
&&+\int_{\Sigma}div_{h_0}(\varphi_2d\varphi_1)\chi dA_{h_0}+\int_{\Sigma}\chi\rho_0|dq-d\log(\varphi_2)|_{h_0}^2dA_{h_0}\\
    &&+\int_{\Sigma}2\lambda e^{2u}\chi \rho_0 dA_{h_0}+3\int_{\Sigma}|d\log(\phi)|_{h_0}^2\chi \rho_0 dA_{h_0}\\
    &&+\int_{\Sigma}\langle d\log(\phi),2\chi d\rho_0-\rho_0 d\chi\rangle_{h_0}dA_{h_0}-\int_{\Sigma}\rho_0\langle d\chi,d q\rangle dA_{h_0}.
\end{eqnarray*}
In particular, since
$$|d\alpha_2|^2=|dq-d\log(\varphi_2)|^2+|d\log\phi|^2+2\langle d\log\phi,d\log\varphi_2-dq\rangle,$$
we can rewrite this more suggestively as 
\begin{eqnarray*}
    \int_{\Sigma}w\zeta \rho_0 dA_{h_0}&=&-\int_{\Gamma_1}[2\log(\phi)+\log(\varphi_2)]\frac{\partial(\chi\rho_0)}{\partial\nu}ds_{h_0}\\
    &&+\int_{\Sigma}(\lambda e^{2u}+|d\alpha_2|_{h_0}^2-K_{h_0})\chi \rho_0 dA_{h_0}\\
&&+\int_{\Sigma}div_{h_0}(\varphi_2d\varphi_1)\chi dA_{h_0}\\
    &&+\int_{\Sigma}\lambda e^{2u}\chi \rho_0 dA_{h_0}+2\int_{\Sigma}|d\log(\phi)|_{h_0}^2\chi \rho_0 dA_{h_0}\\
    &&+\int_{\Sigma}\langle d\log(\phi),2\chi \varphi_2d\varphi_1-\rho_0 d\chi\rangle_{h_0}dA_{h_0}\\
    &&+\int_{\Sigma}\rho_0\langle dq,2\chi d\log(\phi)-d\chi\rangle_{h_0} dA_{h_0}.
\end{eqnarray*}
In this form, the lemma asserts that $\int_{\Sigma}w\zeta \rho_0 dA_{h_0}$ coincides exactly with the first two terms on the right-hand side; thus, \emph{our goal now is to show that}
\begin{eqnarray*}
    0&=&\int_{\Sigma}div_{h_0}(\varphi_2d\varphi_1)\chi dA_{h_0}\\
    &&+\int_{\Sigma}\lambda e^{2u}\chi \rho_0 dA_{h_0}+2\int_{\Sigma}|d\log(\phi)|_{h_0}^2\chi \rho_0 dA_{h_0}\\
    &&+\int_{\Sigma}\langle d\log(\phi),2\chi \varphi_2d\varphi_1-\rho_0 d\chi\rangle_{h_0}dA_{h_0}\\
    &&+\int_{\Sigma}\rho_0\langle dq,2\chi d\log(\phi)-d\chi\rangle_{h_0} dA_{h_0}.
\end{eqnarray*}
To check this, first note that
\begin{eqnarray*}
    -\frac{1}{2}div_{h_0}(\chi\phi^{-2}d\rho)&=&[\lambda e^{2u}+2|d\log\phi|^2]\chi\rho_0-\langle d\log(\phi),\rho_0d\chi\rangle\\
    &&+\langle d\rho_0,\chi d\log\phi-\frac{1}{2}d\chi\rangle,
\end{eqnarray*}
and using Lemma \ref{w21.consq} to justify the integration by parts
\begin{eqnarray*}
    &&\int_{\Sigma}[(\lambda e^{2u}+2|d\log\phi|^2)\chi\rho_0-\langle d\log(\phi),\rho_0d\chi\rangle]dA_{h_0}\\
    &=&-\frac{1}{2}\int_{\Gamma_1}\chi \frac{\partial\rho_0}{\partial\nu}+\int_{\Sigma}\langle \frac{1}{2}d\chi,d\rho_0\rangle_{h_0}dA_{h_0}
    -\int_{\Sigma}\chi \langle d\log\phi,d\rho_0\rangle_{h_0}dA_{h_0}\\
    &=&-\frac{1}{2}\int_{\Sigma}\chi\Delta_{h_0}\rho_0 dA_{h_0}-\int_{\Sigma}\chi\langle d\log\phi,d\rho_0\rangle_{h_0},
\end{eqnarray*}
we see that \emph{the desired identity is equivalent to}
\begin{eqnarray*}
    0&=&\int_{\Sigma}\frac{1}{2}div_{h_0}(\varphi_2d\varphi_1-\varphi_1d\varphi_2)\chi dA_{h_0}+\int_{\Sigma}\chi \langle d\log\phi,\varphi_2d\varphi_1-\varphi_1d\varphi_2\rangle_{h_0}dA_{h_0}\\
    &&+\int_{\Sigma}\rho_0\langle dq,2\chi d\log(\phi)-d\chi\rangle_{h_0}dA_{h_0}\\
    &=&\int_{\Sigma}\chi \phi^{-2}\frac{1}{2}div_{h_0}(\phi^2[\varphi_2d\varphi_1-\varphi_1d\varphi_2])dA_{h_0}-\int_{\Sigma}\rho\langle dq,d(\phi^{-2}\chi)\rangle_{h_0}dA_{h_0}
\end{eqnarray*}

Finally, recall that $q$ satisfies
$$div_{h_0}(\rho dq)+\frac{1}{2}div_{h_0}(\phi^2[\varphi_2d\varphi_1-\varphi_1d\varphi_2])=0$$
in $\mathring\Sigma$, and integrating against $\chi\phi^{-2}$ on $\{\rho_0>\epsilon\}$, we find that
\begin{eqnarray*}
    \int_{\partial\{\rho_0>\epsilon\}}\chi\langle \rho_0 dq,\nu\rangle&=&\int_{\{\rho_0>\epsilon\}}\chi \phi^{-2}div_{h_0}(\rho dq)+\rho\langle dq,d(\chi\phi^{-2})\rangle_{h_0}\\
    &=&-\int_{\{\rho_0>\epsilon\}}\chi \phi^{-2}\frac{1}{2}div_{h_0}(\phi^2[\varphi_2d\varphi_1-\varphi_1d\varphi_2])\\
    &&+\int_{\{\rho_0>\epsilon\}}\rho\langle dq,d(\phi^{-2}\chi)\rangle_{h_0}dA_{h_0}.
\end{eqnarray*}
As $\epsilon\to 0$, after flipping a sign, the right-hand side tends precisely to the integral whose vanishing we need to establish. Thus, it suffices to find a sequence $\epsilon_j\to 0$ along which
$$\lim_{\epsilon_j\to 0}\int_{\partial\{\rho_0>\epsilon_j\}}\chi\langle \rho_0 dq,\nu\rangle_{h_0}=0.$$
But this is straightforward: in each band $\{1/k<\rho_0<2/k\}$, we have
\begin{eqnarray*}
    \int_{\{1/k<\rho_0<2/k\}\cap U_{\delta}}\rho_0|dq|&\leq &\left(\int_{\{\rho_0<2/k\}}\rho_0|dq|^2\right)^{1/2}\left(\int_{\{\rho_0<2/k\}\cap U_{\delta}}\rho_0\right)^{1/2}\\
    &\leq &C_{\delta}(2/k)\left(\int_{\{\rho_0<2/k\}}\rho_0|dq|^2\right)^{1/2},
\end{eqnarray*}
so by the coarea formula for $\rho_0$, we can find $\epsilon_k\in [1/k,2/k]$ with
$$\int_{\{\rho_0=\epsilon_k\}\cap U_{\delta}}\rho_0|dq|\leq C\left(\int_{\{\rho_0<2/k\}}\rho_0|dq|^2\right)^{1/2}.$$
Since $q\in W^{1,2}_{g_0}(M)$, the right-hand side vanishes as $k\to\infty$, so the sequence $\epsilon_k\to 0$ has the desired properties, completing the proof.
\end{proof}

Note that once regularity is established for the triple $(\phi,v,q)$ in $U_{\delta}$, Lemma \ref{wkadm.lem} suffices to establish admissibility $v=q$ along $\Gamma_1$. At this stage, we don't yet have enough regularity to reach that conclusion, but we are ready to prove the following, weaker statement.

\begin{prop}\label{wbd.adm}
    On $U_{\delta}$, $w=u+\alpha_2$ has a lower bound $w\geq -C$, and we have the $L^1$ vanishing
    $$\lim_{t\to 0}\int_{\{\rho_0=t\}\cap U_{\delta}}|v-q|=0$$
    of $v-q$ along $\Gamma_1$.
\end{prop}

\begin{proof}

We've seen in Lemma \ref{wkadm.lem} that
\begin{eqnarray*}
    \int_{\Sigma}w(-\Delta_{h_0}\xi)dA_{h_0}&=&\int_{\Sigma}(\lambda e^{2u}+|d\alpha_2|_{h_0}^2-K_{h_0})\xi dA_{h_0}\\
    &&-\int_{\Gamma_1}(2\log \phi+\log(\varphi_2))\frac{\partial \xi}{\partial\nu}
\end{eqnarray*}
for $\xi$ of the form $\xi=\chi \rho_0\in C_c^{\infty}(U_{\delta})$. Fix some nonnegative, $R_1$-even cutoff $\eta\in C_c^{\infty}(U_{\delta})$, and let
$$G_p^0(x):=G_p(x)-G_p(R_1x)-G_p(R_2x)+G_p(R_1R_2x)$$
be the Dirichlet Green's function on $\Sigma$ centered at $p\in \mathring \Sigma$, with the convention $-\Delta_{h_0}G^0_p=\delta_p$. For any fixed $p\in \mathring\Sigma$, it's easy to see that $G^0_p$ is a smooth multiple of $\rho_0$ in $\Sigma$ away from the pole $p$, and since $w$ and the terms on the right-hand side are smooth on $\mathring\Sigma$, it follows that we can apply the preceding identity to a smooth approximation of $\eta G_p^0$ to obtain
\begin{eqnarray*}
\eta(p)w(p)&=&\int_{\Sigma}w(2\langle d\eta,dG^0_p\rangle_{h_0}+G^0_p\Delta_{h_0}\eta)dA_{h_0}\\
&&+\int_{\Sigma}(\lambda e^{2u}+|d\alpha_2|_{h_0}^2-K_{h_0})\eta G_p^0dA_{h_0}\\
&&-\int_{\Gamma_1}\eta (2\log \phi+\log(\varphi_2))\frac{\partial G_p^0}{\partial\nu}.
\end{eqnarray*}

Now, suppose that $\eta\equiv 1$ on $U_{2\delta}\subset U_{\delta}$; we will obtain the desired estimates on $U_{3\delta}$ to conclude. Indeed, since $\eta\equiv 1$ on $U_{2\delta}$ and $\eta\circ R_1=\eta\in C_c^{\infty}(U_{\delta})$, we see that the collection
$$\{\tau_p:=2\langle d\eta,dG_p^0\rangle_{h_0}+G_p^0\Delta_{h_0}\eta\mid p\in U_{3\delta}\}$$
defines a compact collection of smooth, $R_1$-odd functions on $\Sigma$, which therefore satisfy a uniform pointwise bound of the form
\begin{equation}\label{tau.p.pt}|\tau_p|(x)\leq C_{\delta}\rho_0(x)
\end{equation}
for all $p\in U_{3\delta}$ and $x\in \Sigma$, and satisfy
\begin{equation}\label{unif.van}
    \|\tau_p\|_{C^0(\Sigma)}\to 0\text{ as }p\to \Gamma_1\cap U_{\delta}.
\end{equation}

For any $p\in U_{3\delta}$, we then have
\begin{eqnarray*}
    w(p)&=&\int_{\Sigma}w\tau_p-\int_{\Gamma_1}\eta (2\log \phi+\log(\varphi_2))\frac{\partial G_p^0}{\partial\nu}\\
    &&+\int_{\Sigma}(\lambda e^{2u}+|d\alpha_2|_{h_0}^2-K_{h_0})\eta G_p^0dA_{h_0}
\end{eqnarray*}

Now, since $G_p^0\geq 0$ on $\Sigma$ and $\frac{\partial G_p^0}{\partial\nu}\leq 0$ along $\partial\Sigma$, using \eqref{tau.p.pt}, it follows that
\begin{eqnarray*}
    w(p)&\geq & -C\int_{\Sigma}|w|\rho_0 dA_{h_0}-\int_{\Sigma}K_{h_0}G_p^0dA_{h_0}\\
    &&-\inf_{U_{\delta}}(2\log\phi+\log \varphi_2)\int_{\Gamma_1}\eta\frac{\partial G_p^0}{\partial\nu}
\end{eqnarray*}
for any $p\in U_{3\delta}$. Since $\int_{U_{\delta}}|w|\rho_0<\infty$, $K_{h_0}$ is smooth, and $\inf \phi>0$, this gives the desired lower bound
$$w\geq -C_{\delta}\text{ on }U_{3\delta}.$$

Moreover, combining $w\in L^1(U_{\delta},\rho_0 dA_{h_0})$ with \eqref{tau.p.pt} and \eqref{unif.van}, we see that $\tau_pw\to 0$ pointwise almost everywhere while $|\tau_pw|\leq C_\delta|\rho_0w|$ as $p\to U_{3\delta}\cap\Gamma_1$, it follows from the dominated convergence theorem that
$$\lim_{p\to \Gamma_1\cap U_{3\delta}}\int_{\Sigma} w\tau_p=0.$$
This leaves
$$0=\lim_{p\to \Gamma_1\cap U_{3\delta}}[w(p)-h(p)-\omega(p)],$$
where
$$h(p):=-\int_{\Gamma_1}\eta (2\log \phi+\log(\varphi_2))\frac{\partial G_p^0}{\partial\nu}$$
and
$$\omega(p):=\int_{\Sigma}(\lambda e^{2u}+|d\alpha_2|_{h_0}^2-K_{h_0})\eta G_p^0 dA_{h_0}.$$
In suitable conformal coordinates near $\Gamma_1\cap U_{\delta}$, note that $h(p)$ can be identified with the usual harmonic extension of the $L^1$ function
$$\eta (2\log \phi+\log(\varphi_2))\in L^1(\Gamma_1)$$
via convolution with a Poisson kernel, so classical properties of the Poisson extension (cf. \cite[\S 1]{GarnettMarshall2005HarmonicMeasure}) give the $L^1$ convergence
$$\lim_{t\to 0}h\circ \Phi_t=\eta(2\log \phi+\log(\varphi_2))$$
in $L^1(\Gamma_1)$, where $\Phi_t$ is the flow of $\frac{\nabla \rho_0}{|\nabla\rho_0|^2}$ in $U_{3\delta}$. Thus, to obtain the desired $L^1$ vanishing of
$$v-q=w-2\log \phi-\log\varphi_2$$
along the level sets $\{\rho_0=\epsilon\}\cap U_{3\delta}$ as $\epsilon\to 0$, it suffices to show that
\begin{equation}\label{omega.van}
    \lim_{\epsilon\to 0}\int_{\{\rho_0=\epsilon\}\cap U_{3\delta}}|\omega|ds_{h_0}=0.
\end{equation}

To prove \eqref{omega.van}, we claim next that
\begin{equation}\label{green.clm}
\int_{\{\rho_0=\epsilon\}}|d\rho_0(p)|G^0(x,p)ds_{h_0}(p)\leq C\min\{\epsilon,\rho_0(x)\}
\end{equation}
Assuming this for a moment, Fubini's theorem and the definition of $\omega$ give
$$\int_{\{\rho_0=\epsilon\}\cap U_{3\delta}}|\omega|\leq C\int_{\Sigma}|\lambda e^{2u}+|d\alpha_2|_{h_0}^2-K_{h_0}|(x)\eta(x)\min\{\epsilon,\rho_0(x)\}dA_{h_0}(x),$$
and since $\lambda e^{2u}+|d\alpha_2|_{h_0}^2-K_{h_0}\in L^1(U_{\delta},\rho_0 dA_{h_0})$, the desired vanishing follows from an application of the dominated convergence theorem.

Finally, to check \eqref{green.clm}, consider the Lipschitz function
$$F_{\epsilon}(x):=\min\{\epsilon,\rho_0\},$$
and observe that for any $\chi\in C_c^{\infty}(\Sigma)$, the weak Laplacian of $F_{\epsilon}$ can be determined by 
$$\int\langle dF_{\epsilon},d\chi\rangle_{h_0}dA_{h_0}=-\int_{\{\rho_0<\epsilon\}}\chi\Delta_{h_0}\rho_0+\int_{\{\rho_0=\epsilon\}}\chi |d\rho_0|.$$
In particular, since $F_{\epsilon}|_{\partial\Sigma}=0$, it follows that
$$F_{\epsilon}(x)=\int_{\{\rho_0=\epsilon\}}|d\rho_0|G^0_x+\int_{\{\rho_0<\epsilon\}}(-\Delta_{h_0}\rho_0)G^0_x.$$
Either by using that $|\Delta_{h_0}\rho_0|\leq C\rho_0$, or simply assuming without loss of generality that $-\Delta_{h_0}\rho_0\geq 0$ in $\Sigma$, \eqref{green.clm} follows easily, completing the proof.

\end{proof}

\subsection{Proof of Proposition \ref{rod.smth} assuming a lower bound on $|d\rho|_h^2$}

Next, we complete the proof of Proposition \ref{rod.smth} under the assumption that there is a uniform lower gradient bound
\begin{equation}\label{drho.below}
    |d\rho|_h^2=e^{-2u}|d\rho|_{h_0}^2\geq c_{\delta}>0\text{ on }U_{\delta}\cap\{\rho_0<c_{\delta}\},
\end{equation}
for some $c_{\delta}>0$, before spending the remainder of the section proving \eqref{drho.below}.

\begin{proof}[Proof of Proposition \ref{rod.smth} assuming \eqref{drho.below}]

By the results of Proposition \ref{wbd.adm}, once we've proved that the triple $(\phi,v,q)$ is smooth, the admissibility criterion $v=q$ will indeed be satisfied along $\Gamma_1$, ensuring smoothness of the metric $g$.

Moreover, by Proposition \ref{v.bd.suff}, to prove smoothness of $(\phi,v,q)$ in $U_{3\delta}$, it suffices to show that $v$ is bounded above uniformly in $U_{2\delta}$. Assuming \eqref{drho.below}, we'll show that in fact 
$$u\leq C_{\delta}$$
on $U_{2\delta}$, which together with the lower bound $\log \phi\geq -C$ clearly suffices to bound $v=u-\log \phi$. And since we've already established regularity in the interior region $\{\rho_0\geq c_{\delta}\}$, we only need to prove this boundedness in the neighborhood
$$U_{2\delta}\cap \{\rho_0<c_{\delta}/2\}.$$

Clearly, \eqref{drho.below} is equivalent to an upper bound of the form
$$e^{2u}\leq C_{\delta}|d\rho|_{h_0}^2$$
for $e^{2u}$ on 
$$W_{\delta}:=U_{\delta}\cap\{\rho_0<c_{\delta}\},$$
which together with the equation
$$-\Delta_{h_0}\rho=2\lambda e^{2u}\rho$$
implies
$$0\leq -\Delta_{h_0}\rho\leq C |d\rho|_{h_0}^2\rho$$
on $W_{\delta}$. In particular, for any $F:\mathbb{R}\to\mathbb{R}$, it follows that
$$\Delta_{h_0}F(\rho)=F'(\rho)\Delta_{h_0}\rho+F''(\rho)|d\rho|^2\geq |d\rho|^2(F''(\rho)-CF'(\rho)\rho);$$
taking $F(t)=\int_0^te^{Cs^2/2}ds$, so that $F''(t)=Ce^{Ct^2/2}t=CF'(t)t$, this gives subharmonicity
$$\Delta_{h_0}F(\rho)\geq 0$$
on $W_{\delta}$. Recalling that $\rho\in (C^0\cap W^{1,2})(U_{\delta})$, by extending $F(\rho)$ in an $R_1$-odd way to $\hat{W}_{\delta}=W_{\delta}\cup R_1W_{\delta}$ and defining $\xi\in W^{1,2}(\hat{U}_{\delta})$ to be the harmonic extension of $F(\rho)|_{\partial \hat{W}_{\delta}}$, it then follows from the maximum principle that
$$F(\rho)\leq \xi$$
pointwise in $W_{\delta}$. On the other hand, since $\xi$ is smooth and $R_1$-odd in the interior of $W_{\delta}$, it must satisfy an estimate of the form
$$\xi\leq C_{\delta}\rho_0\text{ on }U_{2\delta}\cap \{\rho_0<c_{\delta}/2\},$$
and it follows that
$$\rho\leq CF(\rho)\leq C_{\delta}'\rho_0\text{ on }U_{2\delta}\cap\{\rho_0<c_{\delta}/2\}$$
as well. Consequently, $\phi^2=\frac{\rho}{\rho_0}$ is uniformly bounded on $U_{2\delta}$.

On the other hand, it follows from Proposition \ref{dpsi.above} that $\rho^{-2}|d\rho|_{h_0}^2=|d\psi|_{h_0}^2$ satisfies an upper bound of the form
$$\rho^{-2}|d\rho|_{h_0}^2\leq \frac{C}{\rho_0^2}\text{ on }U_{\delta},$$
which together with the upper bound for $\phi$ gives a uniform bound
$$|d\rho|_{h_0}^2\leq C\frac{\rho^2}{\rho_0^2}\leq C_{\delta}''\text{ on }U_{2\delta}.$$
Finally, combining this with \eqref{drho.below} gives the desired upper bound for $e^{2u}$, completing the proof.
\end{proof}

The remainder of the section is then devoted to confirming the lower bound \eqref{drho.below}.  

\subsection{An ODE for $f_1^2|d\psi|_h^2$}

Recalling the lower bound $f_2\geq c_{\delta}$ on $U_{\delta}$ supplied by Lemma \ref{f2.below}, we see that the quantity
$$J:=f_1^2|d\psi|_h^2=f_2^{-2}|d\rho|_h^2$$
is bounded above by
$$J\leq C_{\delta}|d\rho|_h^2\text{ on }U_{\delta}.$$
In particular, to prove \eqref{drho.below}, it suffices to show that
\begin{equation}\label{j.below}
    J\geq c_{\delta}\text{ on }U_{\delta}\cap\{\rho_0<c_{\delta}\}
\end{equation}
for some $c_{\delta}>0$, which will be our target in the remainder of the section.

The proof of Lemma \ref{j.below} is based on an ODE comparison argument for $J$ along the flow of 
$$X:=-\frac{\nabla^h\psi}{|d\psi|_h^2},$$
derived from the Hessian identities of Section \ref{sec.hess.id}. In what follows, all norms and inner products are \emph{taken with respect to the metric $h$}, rather than the reference metric $h_0$.

Namely, recalling from Section \ref{sec.hess.id} that
$$Hess_h(\psi)=-d\alpha_1\otimes d\alpha_1-d\alpha_2\otimes d\alpha_2-(\lambda+\langle d\alpha_1,d\alpha_2\rangle_h)h$$
in $\mathring\Sigma$, a straightforward computation gives
\begin{eqnarray*}
    d(f_1^2|d\psi|_h^2)&=&2f_1^2[Hess_h(\psi)(d\psi,\cdot)+|d\psi|_h^2d\alpha_1]\\  
    &=&2f_1^2[-(\lambda+\langle d\alpha_1,d\alpha_2\rangle)d\psi+\langle d\alpha_2,d\psi\rangle d\alpha_1-\langle d\alpha_2,d\psi\rangle d\alpha_2]\\
    &=&2f_1^2[(|d\alpha_2|^2-\lambda)d\psi-2\langle d\alpha_2,d\psi\rangle d\alpha_2],
\end{eqnarray*}
so that
\begin{equation}\label{dj.khaled}
dJ=2f_1^2(|d\alpha_2|^2-\lambda)d\psi-4f_1^2\langle d\alpha_2,d\psi\rangle d\alpha_2.
\end{equation}

In particular, on $\mathring\Sigma\cap\{|d\rho|>0\}$, pairing \eqref{dj.khaled} with the vector field
$$X:=-\frac{\nabla^h\psi}{|d\psi|_h^2}$$
generating the flow through level sets of $\psi$ gives
$$dJ(X)=-2f_1^2(|d\alpha_2|^2-\lambda)+4f_1^2\frac{\langle d\alpha_2,d\psi\rangle^2}{|d\psi|^2}.$$
Thus, the derivative of $J$ along the flow of $X$ has the two-sided bound
\begin{equation}\label{dj.odi}
|dJ(X)|\leq 2f_1^2|d\alpha_2|_h^2+2\lambda f_1^2.
\end{equation}

Over the next two subsections, we show that the term $f_1^2|d\alpha_2|_h^2$ is uniformly bounded in $U_{\delta}$, with some quantitative decay near $\Gamma_1$; this will then allow us to control $J$ near $\Gamma_1$ via \eqref{dj.odi} and a simple ODE comparison argument. As a first step, we obtain an a priori \emph{upper} bound for $J$, which we will use in the analysis of some key Bochner identities for $f_1^2|d\alpha_2|_h^2$.

\subsection{An a priori upper bound for $J$}
Writing 
$$J_{\lambda}:=J+\lambda f_1^2=f_1^2(|d\psi|_h^2+\lambda),$$
it follows from \eqref{dj.khaled} that
\begin{equation}\label{dj.jazzy.jeff}
dJ_{\lambda}=2f_1^2[|d\alpha_2|_h^2d\psi-(2\langle d\alpha_2,d\psi\rangle_h+\lambda)d\alpha_2],
\end{equation}
and as an immediate consequence, we have a pointwise bound of the form
$$|dJ_{\lambda}|\leq C\sqrt{J_{\lambda}}f_1(|d\alpha_2|_h^2+1).$$
Multiplying through by $e^u$ then gives
$$|d\sqrt{J_{\lambda}}|_{h_0}\leq Cf_1(e^{-u}|d\alpha_2|_{h_0}^2+e^u),$$
which together with the lower bounds $f_2\geq c_{\delta}$ and $e^uf_2\geq c_{\delta}$ supplied by Lemmas \ref{f2.below} and  Proposition \ref{wbd.adm} yield the control
$$|d\sqrt{J_{\lambda}}|_{h_0}\leq C_{\delta}\rho |d\alpha_2|_{h_0}^2+C\rho e^u$$
in $U_{\delta}$. In particular, since $\int_{\Sigma}\rho e^{2u}dA_{h_0}<\infty$ and $\alpha_2$ satisfies
$$\int_{U_{\delta}}\rho |d\alpha_2|_{h_0}^2dA_{h_0}\leq C\int_{U_{\delta}}\rho(|d\log\phi|_{h_0}^2+|dq|_{h_0}^2+|d\log\varphi_2|_{h_0}^2)dA_{h_0}<\infty,$$
it follows that
\begin{equation}\label{rt.j.w11}
\int_{U_{\delta}}|d\sqrt{J_{\lambda}}|_{h_0}dA_{h_0}<\infty,
\end{equation}
so that $\sqrt{J_{\lambda}}\in W_{h_0}^{1,1}(U_{\delta})$. 

We could argue similarly for $J$ itself; the utility of working with $J_{\lambda}$ comes from Lemma \ref{psi.boch.4d}, which, together with the equation $\Delta_g\alpha_1=-\lambda$, shows that $\log(J_{\lambda})$ solves
\begin{equation}\label{logj.sub0}
    \Delta_g\log(J_{\lambda})\geq 0\text{ on }\mathring M.
\end{equation}
In particular, combining this with the global control \eqref{rt.j.w11} and related estimates for $J_{\lambda}$, we obtain the following global upper bound in $U_{2\delta}$.

\begin{claim}\label{j.bd}
    We have $J_{\lambda}\in L^{\infty}(U_{2\delta})$ for any $\delta>0$.
\end{claim}
\begin{proof}

To begin, note that by combining \eqref{rt.j.w11} with the coarea formula for $\varphi_2$ on $U_{\delta}\setminus U_{2\delta}$, we can find a regular level set $\{\varphi_2=t\}$ for some $t\in [\delta,2\delta]$   along which
$$\int_{\{\varphi_2=t\}}|d\sqrt{J_{\lambda}}|\leq C_{\delta}.$$
In particular, since $J_{\lambda}$ is smooth in $\mathring\Sigma$, it follows by applying the fundamental theorem of calculus to $\sqrt{J_{\lambda}}$ along $\{\varphi_2=t\}$ that
$$\sup\{ J_{\lambda}(x)\mid \varphi_2(x)=t\}\leq C_{\delta}<\infty.$$
In other words, there is some $t\in [\delta,2\delta]$ such that
\begin{equation}\label{j.v.bd}
J_{\lambda}\leq T_{\delta}<\infty\text{ on }\partial U_t\setminus \Gamma_1.
\end{equation}

Next, by \eqref{logj.sub0}, we know that
\begin{equation}\label{logj.sub}
\int_{U_t}\rho\langle d\chi,d\log(J_{\lambda})\rangle_{h_0}dA_{h_0}\leq 0
\end{equation}
whenever $0\leq \chi\in C_c^{\infty}(U_t\setminus \Gamma_1)$, and more generally for $\chi\in W_0^{1,2}(K)$ for any $K\subset \subset U_t\setminus \Gamma_1$. We would like to test this against truncations of $\log(J_\lambda)$, but doing so requires the introduction of suitable cutoffs near $\Gamma_1$.

Since Lemma \ref{rho.w21} gives $\rho(x)\to 0$ as $x\to \Gamma_1$ in $U_{\delta}$, we can use cutoffs of the form $\zeta_{\epsilon}(\rho)$, where, for $\epsilon>0$, $0\leq \zeta_{\epsilon}\in C^{\infty}(\mathbb{R})$ is a smooth log cutoff function satisfying
$$\zeta_{\epsilon}(s)=0\text{ for }s<\epsilon^2; \text{ }\zeta_{\epsilon}(s)=1\text{ for }s\geq \epsilon,$$
and
$$|\zeta_{\epsilon}'(s)|\leq \frac{2}{s|\log\epsilon|}\text{ on }[\epsilon^2,\epsilon].$$
Note in particular that $\zeta_{\epsilon}(\rho)$ vanishes on a neighborhood of $\Gamma_1$ in $U_{\delta}$. Now, for $T_{\delta}$ as in \eqref{j.v.bd}, set $T=2T_{\delta}$, and for each $k\in (0,\infty)$, define
$$w_k:=\min\{k,\log(J_{\lambda}/T)^+\},$$
so that $w_k$ vanishes along $\partial U_t\setminus \Gamma_1$, and $dw_k=d\log(J_{\lambda})$ on $\{T<J_{\lambda}<Te^k\}$ and $0$ elsewhere. We can then apply \eqref{logj.sub} with
$$\chi=\zeta_{\epsilon}(\rho)w_k$$
in $U_t$, to find
\begin{equation}\label{logj.wk.sub}
\int_{U_t}(\zeta_{\epsilon}(\rho)\rho|dw_k|_{h_0}^2+\rho w_k\zeta_{\epsilon}'(\rho)\langle d\rho,d\log(J_{\lambda})\rangle_{h_0})dA_{h_0}\leq 0.
\end{equation}
Keeping $k$ fixed, we check now that
\begin{equation}\label{sub.err.van}
    \lim_{\epsilon\to 0}\int_{U_t}|\rho w_k\zeta_{\epsilon}'(\rho)\langle d\rho,d\log(J_{\lambda})\rangle_{h_0}|dA_{h_0}=0.
\end{equation}

To see this, note that $|\rho\zeta_{\epsilon}'(\rho)|\leq \frac{2}{|\log\epsilon|}$ by construction, and since $w_k\leq C_k$ and $w_k=0$ on $\{J_{\lambda}\leq T\}$, it follows that
$$\int_{U_t}|\rho w_k\zeta_{\epsilon}'(\rho)\langle d\rho,d\log(J_{\lambda})\rangle_{h_0}|dA_{h_0}\leq \frac{C_k}{|\log\epsilon|}\int_{U_t}|\langle d\rho,d\log(J_{\lambda})\rangle_{h_0}|dA_{h_0}.$$
On the other hand, pairing \eqref{dj.jazzy.jeff} with $d\psi$ with respect to the metric $h$, we see that
$$|\langle dJ_{\lambda},d\psi\rangle_h|\leq 2J_{\lambda}(|d\alpha_2|_h^2+\lambda),$$
and multiplying by $\frac{\rho}{J_{\lambda}}$ gives
$$|\langle d\log(J_{\lambda}),d\rho\rangle|_h\leq 2\rho(|d\alpha_2|_h^2+\lambda).$$
In particular, it follows that
\begin{eqnarray*}
    \int_{U_{\delta}}|\langle d\log(J_{\lambda}),d\rho\rangle|_{h_0}dA_{h_0}&=&\int_{U_{\delta}}|\langle d\log(J_{\lambda}),d\rho\rangle|_{h}dA_h\\
    &\leq &\int_{U_{\delta}}2\rho(|d\alpha_2|_h^2+\lambda)dA_h<\infty,
\end{eqnarray*}
and therefore
$$\lim_{\epsilon\to 0}\frac{C_k}{|\log\epsilon|}\int_{U_t}|\langle d\rho,d\log(J_{\lambda})\rangle_{h_0}|dA_{h_0}=0,$$
confirming \eqref{sub.err.van}.

Returning to \eqref{logj.wk.sub}, it then follows that
$$\int_{U_t}\rho |dw_k|^2_{h_0}dA_{h_0}\leq 0,$$
so that $w_k\equiv 0$ for every $k<\infty$. Thus, $J_{\lambda}\leq T$ on $U_t\supset U_{2\delta}$, as desired.
\end{proof}

\begin{remark}\label{j.trace}
As another consequence of the $W^{1,1}$ control $\sqrt{J},\sqrt{J_{\lambda}}\in W^{1,1}(U_{\delta})$, we also see that $\sqrt{J}$ and $\sqrt{J_{\lambda}}$ must have well-defined traces along $\Gamma_1\cap U_{\delta}$, which agree, since $f_1=\frac{\rho}{f_2}\to 0$ along $\Gamma_1\cap U_{\delta}$ by Lemmas \ref{rho.w21} and \ref{f2.below}. We can now describe that trace explicitly, which will be an important ingredient in the proof of lower bounds for $J$.

Namely, noting that
$$\log(J)=\log(f_2^{-2}e^{-2u}|d\rho|_{h_0}^2)=2(q-v)+\log(|d\rho|_{h_0}^2)-4\log(\phi)-2\log(\varphi_2),$$
it follows from Lemma \ref{w21.consq} and Proposition \ref{wbd.adm} that
$$\lim_{t\to 0}\log(J)\circ \Phi_t=\log(\kappa^2\varphi_2^2)-2\log(\kappa)-2\log(\varphi_2)=0,$$
pointwise almost everywhere in $\Gamma_1\cap U_{\delta}$ along the flow $\Phi_t$ of $\frac{\nabla \rho_0}{|\nabla\rho_0|^2}$.
In particular, it follows that $J\circ \Phi_t\to 1$ pointwise almost everywhere in $\Gamma_1\cap U_{\delta}$, and combining this with Claim \ref{j.bd}, it follows from the bounded convergence theorem that
$$J\circ \Phi_t\to 1\text{ and }\sqrt{J}\circ \Phi_t\to 1\text{ in }L^1(\Gamma_1\cap U_{\delta}).$$
In particular, $\sqrt{J}$ and $\sqrt{J_{\lambda}}$ have trace $1$ on $\Gamma_1\cap U_{\delta}$.

\end{remark}

\subsection{Bounds and decay for $f_1^2|d\alpha_2|^2$ near $\Gamma_1$}

To apply the inequality \eqref{dj.odi} effectively, we will need some boundedness and decay properties for the terms on the right-hand side, namely, for $f_1^2|d\alpha_2|_h^2$. In this section, we prove the needed boundedness and decay, via a pair of maximum principle arguments applied to suitable Bochner identities.

As a starting point, recall from \eqref{alpha.boch} that
$$\frac{1}{2}\Delta_g|d\alpha_2|^2=|Hess_h(\alpha_2)|^2+\lambda |d\alpha_2|^2+\langle d\alpha_1,d\alpha_2\rangle^2+|d\alpha_2|^4,$$
and since $\frac{1}{2}\Delta_gf_1^2=2f_1^2|d\alpha_1|^2-\lambda f_1^2$, it follows that
\begin{eqnarray*}
    \frac{1}{2}\Delta_g(f_1^2|d\alpha_2|^2)&=&f_1^2(|Hess_h(\alpha_2)|^2+\lambda |d\alpha_2|^2+\langle d\alpha_1,d\alpha_2\rangle^2+|d\alpha_2|^4)\\
    &&+(2f_1^2|d\alpha_1|^2-\lambda f_1^2)|d\alpha_2|^2+\langle d|d\alpha_2|^2,df_1^2\rangle\\
    &=&f_1^2(|d\alpha_2|^4+\langle d\alpha_1,d\alpha_2\rangle^2+2|d\alpha_1|^2|d\alpha_2|^2)\\
    &&+f_1^2(|Hess_h(\alpha_2)|^2+2\langle d|d\alpha_2|^2,d\alpha_1\rangle).
\end{eqnarray*}
Next, note that we can rewrite the terms in the last line as
\begin{eqnarray*}
    |Hess_h(\alpha_2)|^2+4Hess_h(\alpha_2)(d\alpha_2,d\alpha_1)&=&|Hess_h(\alpha_2)+d\alpha_1\otimes d\alpha_2+d\alpha_2\otimes d\alpha_1|^2\\
    &&-2|d\alpha_1|^2|d\alpha_2|^2-2\langle d\alpha_1,d\alpha_2\rangle^2,
\end{eqnarray*}
so the identity above becomes
\begin{eqnarray*}
    \frac{1}{2}\Delta_g(f_1^2|d\alpha_2|^2)&=&f_1^2(|d\alpha_2|^4-\langle d\alpha_1,d\alpha_2\rangle^2)\\
    &&+f_1^2|Hess_h(\alpha_2)+d\alpha_1\otimes d\alpha_2+d\alpha_2\otimes d\alpha_1|^2.
\end{eqnarray*}
Moreover, using the simple inequality
\begin{eqnarray*}
    |Hess_h(\alpha_2)+d\alpha_1\otimes d\alpha_2+d\alpha_2\otimes d\alpha_1|^2&\geq& \frac{1}{2}tr(Hess_h(\alpha_2)+d\alpha_1\otimes d\alpha_2+d\alpha_2\otimes d\alpha_1)^2\\
    &=&\frac{1}{2}(\Delta_h\alpha_2+2\langle d\alpha_1,d\alpha_2\rangle)^2
\end{eqnarray*}
and the fact that $\Delta_h\alpha_2=-\lambda-\langle d\psi,d\alpha_2\rangle,$
it follows after a little rearrangement that
\begin{eqnarray*}
    &&|Hess_h(\alpha_2)+d\alpha_1\otimes d\alpha_2+d\alpha_2\otimes d\alpha_1|^2\\
    &\geq&2\langle d\alpha_1,d\alpha_2\rangle^2-2(\lambda+\langle d\psi,d\alpha_2\rangle)\langle d\alpha_1,d\alpha_2\rangle+\frac{1}{2}(\lambda+\langle d\psi,d\alpha_2\rangle)^2\\
    &=&2\langle d\alpha_1,d\alpha_2\rangle^2-2\langle d\psi,d\alpha_2\rangle\langle d\alpha_1,d\alpha_2\rangle+\frac{1}{2}\langle d\psi,d\alpha_2\rangle^2-2\lambda\langle d\alpha_1,d\alpha_2\rangle\\
    && +\frac{1}{2}\lambda^2+\lambda\langle d\psi,d\alpha_2\rangle\\
    &\geq &\langle d\alpha_1,d\alpha_2\rangle^2-\frac{1}{2}\langle d\psi,d\alpha_2\rangle^2+2\lambda|d\alpha_2|^2-\lambda\langle d\psi,d\alpha_2\rangle+\frac{1}{2}\lambda^2\\
    &\geq &\langle d\alpha_1,d\alpha_2\rangle^2-\langle d\psi,d\alpha_2\rangle^2.    
\end{eqnarray*}
Finally, applying this in the preceding equation for $\frac{1}{2}\Delta_g(f_1^2|d\alpha_2|^2)$, we arrive at the estimate
\begin{equation}\label{b.sub.1}
    \frac{1}{2}\Delta_g(f_1^2|d\alpha_2|^2)\geq f_1^2|d\alpha_2|^2(|d\alpha_2|^2-|d\psi|^2).
\end{equation}
Using this, we can now show that $f_1^2|d\alpha_2|_h^2$ is bounded on $\Sigma\setminus \Gamma_2$.
\begin{lem}\label{b.bd}
    We have $f_1^2|d\alpha_2|_h^2\in L^{\infty}_{loc}(\Sigma\setminus \Gamma_2)$.
\end{lem}
\begin{proof}
We briefly introduce the auxiliary metric
$$\tilde{h}:=f_1^{-2}h=f_1^{-2}e^{2u}h_0$$
on $\Sigma\setminus \partial\Sigma$. With respect to this metric, note that
$$|d\alpha_2|_{\tilde{h}}^2=f_1^2|d\alpha_2|_h^2$$
and
$$|d\psi|_{\tilde{h}}^2=f_1^2|d\psi|_h^2.$$
As a consequence, setting $V_{\delta}:=U_{\delta}\setminus \Gamma_1$ for any fixed $\delta>0$, it follows from Claim \ref{j.bd} that 
\begin{equation}\label{psi.lip}
|d\psi|_{\tilde{h}}\leq C_{\delta}\text{ on }V_{\delta}.
\end{equation}
Next, we will show that there exists some $R_{\delta}>0$ such that
\begin{equation}\label{dist.below}
    dist_{\tilde{h}}(x,\partial V_{\delta})\geq R_{\delta}\text{ for all }x\in V_{2\delta};
\end{equation}
that is, $B_{R_{\delta}}^{\tilde{h}}(x)\subset V_{\delta}$ for every $x\in V_{2\delta}$.

To prove \eqref{dist.below}, let $x\in V_{2\delta}$, and let $\gamma:[0,a)\to V_{\delta}$ be a unit speed curve with respect to $\tilde{h}$ such that $\gamma(0)=x$ and 
$$\lim_{t\to a}\gamma(t)\in \partial V_{\delta}=\{\varphi_2=\delta\}\cup (\{\varphi_2>\delta\}\cap\Gamma_1).$$
It follows from \eqref{psi.lip} that
$$|\psi(\gamma(t))|\leq |\psi(x)|+C_{\delta}t$$
for every $t\in [0,a)$; in particular since $\psi(y)\to -\infty$ as $y\to \Gamma_1$ in $V_{\delta}$, we must have
$$\lim_{t\to a}\gamma(t)\in \partial V_{\delta}\setminus \Gamma_1=\{\varphi_2=\delta\}.$$
With respect to the reference metric $h_0$, the distance $c_{\delta}$ from $V_{2\delta}$ to $\{\varphi_2=\delta\}$ is positive, so we must have
$$\int_0^a|\gamma'(t)|_{h_0}dt\geq c_{\delta}>0.$$
On the other hand, since Lemma \ref{f2.below} and Proposition \ref{wbd.adm} give positive lower bounds for $f_2$ and $f_2e^u$ on $V_{\delta}$, we also have
$$1=|\gamma'(t)|_{\tilde{h}}^2=\rho^{-2}f_2^2e^{2u}|\gamma'(t)|_{h_0}^2\geq c_{\delta}'|\gamma'(t)|_{h_0},$$
so combining this with the lower bound for $\int_0^a|\gamma'|_{h_0}(t)dt$ yields
$$a=L_{\tilde{h}}(\gamma)\geq c_{\delta}c_{\delta}',$$
proving \eqref{dist.below} with $R_{\delta}=c_{\delta}c_{\delta}'$.

Next, we estimate the curvature 
$$K_{\tilde{h}}=e^{2\alpha_1}(K_h+\Delta_h\alpha_1).$$
Namely, using that $K_h=-\langle d\alpha_1,d\alpha_2\rangle_h$ and $\Delta_h\alpha_1=-\lambda-\langle d\alpha_1,d\psi\rangle_h,$ we see that
$$K_{\tilde{h}}=f_1^2(-\lambda-|d\psi|_h^2+|d\alpha_2|_h^2)\geq-C_{\delta}$$
on $V_{\delta}$, where the lower bound follows from Claim \ref{j.bd}. Thus, fixing an arbitrary $p\in V_{2\delta}$, it follows from the Laplacian comparison theorem  that the distance function
$$\chi_p(x):=dist_{\tilde{h}}(p,x)$$
satisfies a bound of the form
\begin{equation}\label{lap.comp}
\Delta_{\tilde{h}}\chi_p\leq C_{\delta}'/\chi_p,
\end{equation}
away from the cut locus of $p$ and globally in the barrier sense on $B^{\tilde{h}}_{R_{\delta}}(p)\subset V_{\delta}$. In particular, for each $x\in B^{\tilde{h}}_{R_{\delta}}(p)$, applying Calabi's construction from \cite[Lemma 7.1.9]{Petersen2016RiemannianGeometry} for sufficiently small $\epsilon>0$ gives a smooth function $\chi_{p,x}$ in a neighborhood $U_x$ of $x$ satisfying
$$0=\chi_{p,x}(x)-\chi_p(x)\leq \chi_{p,x}(y)-\chi_p(y)\text{ for }y\in U_x,$$
$$|d\chi_{p,x}|_{\tilde{h}}=1\text{ in }U_x,$$
and
$$\Delta_{\tilde{h}}\chi_{p,x}\leq \frac{C_{\delta}'}{\chi_p}+1.$$

Finally, defining the second-order operator $L$ by
$$L\zeta:=\Delta_{\tilde{h}}\zeta+\langle d\psi,d\zeta\rangle_{\tilde{h}}=f_1^2\Delta_g\zeta,$$
we see that \eqref{b.sub.1} is equivalent to 
$$\frac{1}{2}L|d\alpha_2|_{\tilde{h}}^2\geq |d\alpha_2|_{\tilde{h}}^2(|d\alpha_2|_{\tilde{h}}^2-|d\psi|_{\tilde{h}}^2)\geq |d\alpha_2|_{\tilde{h}}^2(|d\alpha_2|_{\tilde{h}}^2-C)$$
on $V_{\delta}$, where we've used \eqref{psi.lip} for the lower bound. Let $\eta=\eta_{\delta}: \mathbb{R}\to \mathbb{R}$ be a smooth, nonincreasing function such that $\eta(t)=0$ for $t\geq \frac{2R_{\delta}}{3}$, $\eta(t)=1$ for $t<\frac{R_{\delta}}{3}$ and
\begin{equation}\label{etader.bd}
(\eta')^2\leq C_{\delta}\eta,
\end{equation}
so that 
$$\xi:=\eta(\chi_p)|d\alpha_2|_{\tilde{h}}^2$$
must achieve its positive maximum at some interior point $x\in B_{2R_{\delta}/2}^{\tilde{h}}(p)$. If this $x$ is a smooth point for $\chi_p$, it follows directly that
$$d\xi(x)=[|d\alpha_2|_{\tilde{h}}^2\eta'(\chi_p)d\chi_p+\eta(\chi_p)d|d\alpha_2|_{\tilde{h}}^2](x)=0$$
and
\begin{eqnarray*}
    0&\geq &L\xi(x)\\
    &\geq &\eta(\chi_p)|d\alpha_2|_{\tilde{h}}^2(|d\alpha_2|_{\tilde{h}}^2-C)+2\langle d|d\alpha_2|_{\tilde{h}}^2,\eta'(\chi_p)d\chi_p\rangle_{\tilde{h}}\\
&&+|d\alpha_2|_{\tilde{h}}^2(\eta'(\chi_p)\Delta_{\tilde{h}}\chi_p+\eta''(\chi_p)+\eta'(\chi_p)\langle d\psi,d\chi_p\rangle_{\tilde{h}})\\
&=&\eta(\chi_p)|d\alpha_2|_{\tilde{h}}^2(|d\alpha_2|_{\tilde{h}}^2-C)-\frac{\eta'(\chi_p)^2}{\eta(\chi_p)}2|d\alpha_2|_{\tilde{h}}^2\\
&&+|d\alpha_2|_{\tilde{h}}^2(\eta'(\chi_p)\Delta_{\tilde{h}}\chi_p+\eta''(\chi_p)+\eta'(\chi_p)\langle d\psi,d\chi_p\rangle_{\tilde{h}}).
\end{eqnarray*}
Applying \eqref{lap.comp}, \eqref{psi.lip}, and \eqref{etader.bd}, we arrive at an estimate of the form
$$0\geq |d\alpha_2|_{\tilde{h}}^2(x)(\xi(x)-C_{\delta})$$
for a suitable constant $C_{\delta}$ at the maximum point $x$ for $\xi$, so that
$$\max\eta(\chi_p)|d\alpha_2|_{\tilde{h}}^2\leq C_{\delta}.$$
In the case where the maximum $x$ for $\xi$ occurs in the cut locus of $p$, we can use Calabi's barrier construction, as in the proof of \cite[Lemma 7.1.9]{Petersen2016RiemannianGeometry}, to find a smooth function $\zeta$ on a neighborhood $U_x$ of $x$ satisfying
$$0=\zeta(x)-\chi_p(x)\leq \zeta(y)-\chi_p(y),$$
$$|d\zeta|(y)=1,\text{ and }\Delta_{\tilde{h}}\zeta \leq C_{\delta}'/\chi_p+1$$
for all $y\in U_x$. Since $\zeta\geq \chi_p$ and $\eta$ is nonincreasing, it follows that
$$\eta(\zeta(y))|d\alpha_2|_{\tilde{h}}^2(y)\leq \eta(\chi_p(y))|d\alpha_2|_{\tilde{h}}^2(y)$$
for all $y\in U_x$, with equality when $y=x$. Thus, $x$ is also a maximum for the smooth function $\eta(\zeta)|d\alpha_2|_{\tilde{h}}^2$ on $U_x$, and we can argue exactly as before to arrive at the desired estimate
$$|d\alpha_2|_{\tilde{h}}^2(p)\leq \max \eta(\zeta)|d\alpha_2|_{\tilde{h}}^2\leq C_{\delta}'$$
in this case as well.
\end{proof}

Our goal now is to upgrade the boundedness of $f_1^2|d\alpha_2|_h^2$ to an estimate showing decay as $\rho\to 0$. As a starting point, we revisit the Bochner identity
$$\frac{1}{2}\Delta_g(|d\alpha_2|_h^2)=|Hess_h(\alpha_2)|^2+\lambda|d\alpha_2|^2+\langle d\alpha_1,d\alpha_2\rangle^2+|d\alpha_2|^4,$$
and use Lemma \ref{kato.lem} to rewrite the Hessian term, giving
\begin{eqnarray*}
    \frac{1}{2}\Delta_g|d\alpha_2|_h^2&=&\lambda |d\alpha_2|^2+\langle d\alpha_1,d\alpha_2\rangle^2+|d\alpha_2|^4\\
    &&+2|d|d\alpha_2||^2+(\Delta_h\alpha_2)[\Delta_h\alpha_2-2\langle \frac{d|d\alpha_2|}{|d\alpha_2|},d\alpha_2\rangle]\\
    &=&\lambda |d\alpha_2|^2+\langle d\alpha_1,d\alpha_2\rangle^2+|d\alpha_2|^4
    +2|d|d\alpha_2||^2\\
    &&+(\lambda+\langle d\psi,d\alpha_2\rangle)[\lambda+\langle d\psi,d\alpha_2\rangle+2\langle d|d\alpha_2|,\frac{d\alpha_2}{|d\alpha_2|}\rangle].
\end{eqnarray*}
Where $|d\alpha_2|>0$, it follows that
\begin{eqnarray*}
    \frac{1}{2}\Delta_g\log(|d\alpha_2|_h^2)&=&\frac{(1/2)\Delta_g|d\alpha_2|^2-2|d|d\alpha_2||^2}{|d\alpha_2|^2}\\
    &=&\lambda+|d\alpha_2|^2+\frac{\langle d\alpha_1,d\alpha_2\rangle^2}{|d\alpha_2|^2}\\
    &&+\frac{(\lambda+\langle d\psi,d\alpha_2\rangle)[\lambda+\langle d\psi,d\alpha_2\rangle+2\langle d\log(|d\alpha_2|),d\alpha_2\rangle]}{|d\alpha_2|^2}.
\end{eqnarray*}
In other words, defining an operator $\mathcal{L}$ for functions on $\{|d\alpha_2|>0\}\cap \{\rho_0>0\}$ by
$$\mathcal{L}(\chi):=\Delta_g\chi-2\frac{(\lambda+\langle d\psi,d\alpha_2\rangle)}{|d\alpha_2|^2}\langle d\alpha_2,d\chi\rangle,$$
we have
\begin{equation}\label{fancyl.boch}
    \mathcal{L}(\log |d\alpha_2|_h)=\lambda +|d\alpha_2|^2+\frac{(\lambda+\langle d\psi,d\alpha_2\rangle)^2+\langle d\alpha_1,d\alpha_2\rangle^2}{|d\alpha_2|^2}.
\end{equation}
On the other hand, since $\Delta_g\alpha_1=-\lambda$, we see that
$$\mathcal{L}\alpha_1=-\lambda-2\frac{(\lambda+\langle d\psi,d\alpha_2\rangle)\langle d\alpha_2,d\alpha_1\rangle}{|d\alpha_2|^2},$$
and adding this to \eqref{fancyl.boch} and simplifying yields the nice identity
\begin{equation}\label{lb.boch}
    \mathcal{L}(\log(f_1|d\alpha_2|_h))=2\lambda+2|d\alpha_2|^2+\frac{\lambda^2}{|d\alpha_2|^2}.
\end{equation}
We'll now use this to obtain improved decay for $f_1|d\alpha_2|_h$ near $\Gamma_1$.

\begin{prop}\label{b.decay}
    For each $\delta>0$, there exists $C_{\delta}<\infty$ such that
    $$f_1^2|d\alpha_2|^2\leq \frac{C_{\delta}}{\psi^2}\text{ on }U_{2\delta},$$
    where, again, $\psi=\log(\rho)$.
\end{prop}

\begin{proof}
First, since $\rho(x)\to 0$ as $x\to \Gamma_1$ in $U_{\delta}$, a bound of the desired form holds automatically on the interior region $\{\rho\geq 1/e\}\cap U_{\delta}$.

On the region $\Omega:=\{0<\rho<1/e\}\cap U_{\delta},$ fix a constant $A=A_{\delta}>0$, and, recalling that $\psi=\log\rho$, define 
$$\zeta:=\frac{A}{(\varphi_2-\delta)^2(1+\psi)^2}.$$
Since $f_1^2=f_2^{-2}\rho^2\leq C_{\delta}e^{2\psi}$ in $U_{\delta}$, by choosing $A_{\delta}>0$ sufficiently large, we can arrange that
$$\zeta f_1^{-2}\geq \frac{A_{\delta}C_{\delta}^{-1}e^{2|\psi|}}{\max(\varphi_2)^2|\psi|^2}\geq 1$$
on $\Omega$. 

Next, we observe that, in $\Omega$,
\begin{eqnarray*}
    &&\mathcal{L}[\log(|1+\psi|)]\\
    &=&div_g(\frac{d\psi}{1+\psi})-2\frac{(\lambda+\langle d\psi,d\alpha_2\rangle)\langle d\alpha_2,d\psi\rangle}{(1+\psi)|d\alpha_2|^2}\\
    &=&\frac{-2\lambda}{1+\psi}-\frac{|d\psi|^2}{(1+\psi)^2}-2\frac{\langle d\alpha_2,d\psi\rangle^2}{(1+\psi)|d\alpha_2|^2}-2\lambda \frac{\langle d\alpha_2,d\psi\rangle}{(1+\psi)|d\alpha_2|^2}\\
    &\geq & \frac{2\lambda}{|1+\psi|}-\frac{|d\psi|^2}{(1+\psi)^2}+2\frac{\langle d\alpha_2,d\psi\rangle^2}{|1+\psi||d\alpha_2|^2}-(\frac{|d\psi|^2}{(1+\psi)^2}+\frac{\lambda^2}{|d\alpha_2|^2})\\
    &\geq & -2\frac{|d\psi|^2}{(1+\psi)^2}-\frac{\lambda^2}{|d\alpha_2|^2},
\end{eqnarray*}
which together with Claim \ref{j.bd} gives
$$f_1^2\mathcal{L}[\log|1+\psi|]\geq -\frac{C}{(1+\psi)^2}-\frac{\lambda^2f_1^2}{|d\alpha_2|^2}\geq-C'\zeta-\frac{\lambda^2f_1^2}{|d\alpha_2|^2}.$$
Similarly, we compute
\begin{eqnarray*}
    &&\mathcal{L}[\log(\varphi_2-\delta)]\\
&=&\frac{\Delta_h\varphi_2+\langle d\varphi_2,d\psi\rangle_h}{\varphi_2-\delta}-\frac{|d\varphi_2|_h^2}{(\varphi_2-\delta)^2}-2\frac{(\lambda+\langle d\psi,d\alpha_2\rangle)\langle d\alpha_2,d\varphi_2\rangle}{(\varphi_2-\delta)|d\alpha_2|^2}\\
    &\geq &-\frac{|\Delta_h\varphi_2|}{\varphi_2-\delta}-\frac{|d\varphi_2|_h|d\psi|_h}{\varphi_2-\delta}-\frac{|d\varphi_2|_h^2}{(\varphi_2-\delta)^2}-2\lambda\frac{|d\varphi_2|_h}{(\varphi_2-\delta)|d\alpha_2|_h}-3\frac{|d\psi|_h|d\varphi_2|_h}{(\varphi_2-\delta)},
\end{eqnarray*}
and observe that 
$$|\Delta_h\varphi_2|+|d\varphi_2|_h^2=e^{-2u}(|\Delta_{h_0}\varphi_2|+|d\varphi_2|_{h_0}^2)\leq Ce^{-2u}=Ce^{-2u}f_2^{-2}\frac{\rho^2}{f_1^2}.$$
Since Proposition \ref{wbd.adm} gives $e^{-2u}f_2^{-2}\leq C_{\delta}$ on $\Omega$, it follows that
$$|\Delta_h\varphi_2|+|d\varphi_2|_h^2\leq C_{\delta}e^{-2|\psi|}f_1^{-2},$$
so that
\begin{eqnarray*}
    f_1^2\mathcal{L}[\log(\varphi_2-\delta)]&\geq &-\frac{C_{\delta}}{(\varphi_2-\delta)^2e^{2|\psi|}}\\
    &&-2\lambda\frac{f_1}{|d\alpha_2|_h(\varphi_2-\delta)e^{|\psi|}}-\frac{C_{\delta}f_1|d\psi|_h}{(\varphi_2-\delta)e^{|\psi|}},
\end{eqnarray*}
where the constant $C_{\delta}$ may differ from line to line. In particular, using Claim \ref{j.bd} to estimate the last term, and repeatedly using that $|1+\psi|\leq e^{|\psi|}$, we arrive at an estimate of the form
$$f_1^2\mathcal{L}[\log(\varphi_2-\delta)]\geq -C_{\delta}\zeta-\frac{C_{\delta}f_1\sqrt{\zeta}}{|d\alpha_2|_h}\geq -C_{\delta}\zeta(1+|d\alpha_2|_h^{-1}),$$
where we've used the estimate $f_1\leq \sqrt{\zeta}$ in the last line.

Combining this with the lower bound for $\mathcal{L}[\log(|1+\psi|)]$, we see that
$$f_1^2\mathcal{L}\log(\zeta)\leq C_{\delta}(\zeta+\frac{\zeta}{|d\alpha_2|_h})+\frac{2\lambda^2f_1^2}{|d\alpha_2|_h^2}.$$
Perturbing $\zeta$ slightly to $\zeta_{\epsilon}:=\zeta+\epsilon$, note moreover that
$$\mathcal{L}\log \zeta_{\epsilon}=\frac{\zeta}{\zeta+\epsilon}\mathcal{L}[\log \zeta]+\frac{\epsilon\zeta}{(\zeta+\epsilon)^2} |d\log(\zeta)|^2;$$
noting again that
$$f_1^2|d\log(\varphi_2-\delta)|_h^2\leq \frac{C\rho^2 f_2^{-2}e^{-2u}}{(\varphi_2-\delta)^2}\leq C\zeta$$
and
$$f_1^2|d\log(1+\psi)|_h^2=\frac{f_1^2|d\psi|^2}{(1+\psi)^2}\leq C\zeta,$$
it follows in particular that
$$f_1^2\mathcal{L}(\log(\zeta+\epsilon))\leq C_{\delta}'(\zeta+\frac{\zeta}{|d\alpha_2|_h})+\frac{2f_1^2\lambda^2}{|d\alpha_2|_h^2}).$$
for any $\epsilon\in (0,1)$. Finally, combining this with \eqref{lb.boch}, we see that
\begin{equation}\label{quo.sub}
f_1^2\mathcal{L}[\log(f_1^2|d\alpha_2|_h^2)-\log(\zeta+\epsilon)]\geq 2f_1^2|d\alpha_2|_h^2-C_{\delta}'(\zeta+\frac{\zeta}{|d\alpha_2|_h})
\end{equation}

Finally, for each $\epsilon\in (0,1)$, consider the function
$$w_{\epsilon}:=\log(f_1^2|d\alpha_2|^2)-\log(\zeta+\epsilon)+\epsilon\psi$$
on $\Omega$. By Lemma \ref{b.bd}, $\log(f_1^2|d\alpha_2|^2)$ is bounded from above on $\Omega$, while $-\log(\zeta+\epsilon)\leq -\log(\epsilon)$, so it follows that
$$\lim_{x\to \Gamma_1}w_{\epsilon}(x)\leq C_{\epsilon}+\epsilon\lim_{x\to \Gamma_1}\psi(x)=-\infty$$
in $\Omega$. Moreover, as 
$$x\to \partial\Omega\setminus \Gamma_1=\{\psi=-1\}\cup\{\varphi_2=\delta\},$$
we have $\zeta(x)\to +\infty$, so that $\log(\zeta+\epsilon)\to +\infty$, and it follows that
$$\lim_{x\to \partial\Omega}w_{\epsilon}(x)=-\infty$$
everywhere along $\partial\Omega$. In particular, $w_{\epsilon}$ must achieve its maximum at some interior point $x\in \Omega$, where
$$\mathcal{L}w_{\epsilon}(x)\leq 0.$$
To estimate $\mathcal{L}w_{\epsilon}$, first observe that
\begin{eqnarray*}
    f_1^2\mathcal{L}(\psi)&=&-2f_1^2\lambda-2\frac{f_1^2(\lambda+\langle d\psi,d\alpha_2\rangle)}{|d\alpha_2|^2}\langle d\alpha_2,d\psi\rangle\\
    &\geq &-C(1+\zeta/|d\alpha_2|_h),
\end{eqnarray*}
by repeated use of Claim \ref{j.bd}, and the boundedness of $f_1$ and $f_1^2/\zeta$ in $\Omega$. Combining this with \eqref{quo.sub}, at the maximum point $x=x_{\epsilon}$ for 
$$w_{\epsilon}=\log(\frac{f_1^2|d\alpha_2|^2}{\zeta+\epsilon})+\epsilon\psi,$$
we must have
\begin{eqnarray*}
    0&\geq &f_1^2\mathcal{L}[w_{\epsilon}](x)\\
    &\geq &2f_1^2|d\alpha_2|_h^2-C'(\zeta+\frac{\zeta}{|d\alpha_2|_h})-C\epsilon.    
\end{eqnarray*}
Setting 
$$\mu_{\epsilon}:=\max_{\Omega}w_{\epsilon}=w_{\epsilon}(x_{\epsilon}),$$
we then have $\log(\frac{f_1^2|d\alpha_2|^2}{\zeta+\epsilon})(x)> \mu_{\epsilon}$, so that
$$(f_1^2|d\alpha_2|^2)(x)\geq e^{\mu_{\epsilon}}(\zeta+\epsilon),$$
and applying this in the preceding estimate gives
$$0\geq 2 e^{\mu_{\epsilon}}(\zeta(x)+\epsilon)-C'(\zeta+e^{-\mu_{\epsilon}/2}f_1\sqrt{\zeta})(x)-C\epsilon.$$
Recalling that $f_1^2\leq \zeta$ pointwise in $\Omega$, we see finally that
$$0\geq (2e^{\mu_{\epsilon}}-C'-C'e^{-\mu_{\epsilon}/2})\zeta(x)+2(\epsilon e^{\mu_{\epsilon}}-C\epsilon).$$
If we had 
$$e^{\mu_{\epsilon}}>\max\{C,C',1\},$$
then the right-hand side above would be positive; thus, we must have an estimate of the form 
$$e^{\mu_{\epsilon}}\leq C_{\delta}.$$
In other words, at \emph{every} point $x\in \Omega$, we have
$$\log(\frac{f_1^2|d\alpha_2|^2}{\zeta+\epsilon})(x)+\epsilon\psi(x)\leq C_{\delta}$$
for every $\epsilon\in (0,1)$, and taking $\epsilon\to 0$ gives
$$f_1^2|d\alpha_2|_h^2\leq e^{C_{\delta}}\zeta$$
pointwise in $\Omega$. In particular, since
$$\zeta\leq \frac{4A}{\delta^2\psi^2}$$
in $\Omega\cap U_{2\delta}\cap\{\psi\leq -2\}$, this gives the desired estimate where $\{\psi\leq -2\}$; of course, in $U_{2\delta}\cap\{\psi>-2\}$, Lemma \ref{b.bd} already gives the desired bound.    
\end{proof}

\subsection{The lower bound for $|d\rho|_h^2$}\label{sec.drho.below}

We're now ready to prove the lower bound \eqref{drho.below}, thereby finishing the proof of Proposition \ref{rod.smth}. Since $f_2\geq c_{\delta}>0$ and $\rho(x)\to 0$ as $x\to \Gamma_1\cap U_{\delta}$, it suffices to prove the following variant.

\begin{lem}\label{j.below}
    There is a constant $\epsilon_{\delta}>0$ such that $J=f_1^2|d\psi|_h^2\geq \frac{1}{2}$ on $U_{\delta}\cap\{\rho<\epsilon_{\delta}\}$. 
\end{lem}

\begin{proof} 
First, recall from Remark \ref{j.trace} that $\sqrt{J}\in W^{1,1}(U_{\delta})$ has trace 
$$\sqrt{J}|_{\Gamma_1\cap U_{\delta}}\equiv 1.$$
Next, arguing as we did in the proof of Claim \ref{j.bd}, note that we can find $t\in (\delta,2\delta]$ such that $\sqrt{J}$ is continuous along the slice $\{\varphi_2=t\}$, and $\sqrt{J}=1$ at the endpoints $\Gamma_1\cap \{\varphi_2=t\}$, so that $\sqrt{J}(x)\to 1$ as $\{\varphi_2=t\}\ni x\to \Gamma_1$. Thus, by taking $\epsilon>0$ sufficiently small, we can arrange that
\begin{equation}\label{j.wall}
    J\geq \frac{9}{10}\text{ on }\partial U_t\cap\{\rho<\epsilon\}.
\end{equation}

Next, let $W_{\epsilon}$ be the connected component of $\{\rho<\epsilon\}\cap U_t$ containing one of the $n+1$ identical components of $\Gamma_1\cap U_t$, and suppose, for a contradiction, that
\begin{equation}
    \inf_{W_{\epsilon}}J<\frac{1}{2}.
\end{equation}
Since $\sqrt{J}\to 1$ along $\Gamma_1\cap U_{\delta}$ in the trace sense, it follows that $\sup_{W_{\epsilon}}J\geq 1$, so by connectedness, we can find some 
$$p\in W_{\epsilon}\subset \{\rho<\epsilon\}\text{ such that }J(p)=\frac{1}{2}.$$
After perturbing slightly, we may assume that $s_{\epsilon}:=\rho(p)\in (0,\epsilon)$ is a regular value for $\rho$ on $U_{\delta}$, and there exists an open interval $I\subset \{\rho=s_{\epsilon}\}$ such that
$$\frac{1}{3}\leq J\leq \frac{2}{3}\text{ on }I.$$

Now, for each $q\in I$, let $\gamma:[0,T)\to W_{\epsilon}$ be the maximal solution of the level set flow
$$\gamma(0)=q,\text{ }\gamma'(t)=-\frac{\nabla^h\psi}{|d\psi|_h^2}$$
for $\psi$, and observe that \eqref{dj.odi} gives
$$|(J\circ \gamma)'(t)|\leq (2 f_1^2|d\alpha_2|_h^2+2\lambda f_1^2)(\gamma(t)).$$
In particular, combining this with Proposition \ref{b.decay} and the bound $f_1^2\leq C_{\delta}\rho^2$, it follows that
$$|(J\circ \gamma)'(t)|\leq \frac{C}{\psi(\gamma(t))^2}+C e^{2\psi(\gamma(t))},$$
and since $\psi(\gamma(t))=\log(s_{\epsilon})-t\leq \log(\epsilon)-t$, it follows that
$$|J(\gamma(\tau))-J(q)|\leq C\int_0^{\tau}(\frac{1}{(t+|\log\epsilon|)^2}+\epsilon^2e^{-2t})dt\leq \frac{C}{|\log\epsilon|}+C\epsilon^2$$
for any $\tau\in [0,T)$. Since $\frac{1}{3}\leq J(q)\leq \frac{2}{3}$, for $\epsilon>0$ sufficiently small, it follows that
$$J(\gamma(t))\in [\frac{1}{4},\frac{3}{4}]$$
for all times $t\in [0,T)$. In particular, $|d\psi|_h(\gamma(t))\geq \frac{1}{2f_1}>c_{\delta}$ cannot vanish in finite time, and it follows from \eqref{j.wall} that $\gamma(t)$ can never reach $\partial U_t$, so we must have $T=\infty$. In particular, letting $F:I\times [0,\infty)\to W_{\epsilon}$ be the flow of $-\frac{\nabla^h\psi}{|d\psi|_h^2}$, and letting $Q$ be the box
$$Q:=F(I\times [0,\infty)),$$
it follows that $\frac{1}{4}\leq J\leq \frac{3}{4}$ throughout $Q$. 

Finally, let $\chi: \mathbb{R}\to [0,\infty)$ be a smooth function satisfying
$$\chi(t)=0\text{ for }t>\frac{4}{5}\text{ and }t<1/5$$
and 
$$\chi(t)=1\text{ for }t\in [1/4,3/4],$$
and consider the vector field 
$$Y:=\chi(J)\nabla^h\rho$$
on $U_{2\delta}$. By \eqref{j.wall}, we have that $Y$ vanishes near $\partial U_{\delta}\setminus \Gamma_1$; fixing $t>0$ and taking $s<<t$ so that $\{\rho_0<s\}\cap U_t\subset \{\rho<t\}$, integrating $div_h(Y)$ over the domains
$$\Omega_{t,s}:=\{\rho<t\}\cap\{\rho_0>s\}\cap U_t$$
gives
\begin{eqnarray*}
    \int_{\Omega_{t,s}}[\chi'(J)\langle dJ,d\rho\rangle_h+\chi(J)\Delta_h\rho]dA_h&=&\int_{\{\rho=t\}}\chi(J)|d\rho|_h ds_h\\
    &&-\int_{\{\rho_0=s\}}\chi(J)\langle d\rho,\frac{d\rho_0}{|d\rho_0|_h}\rangle_h ds_h,
\end{eqnarray*}
or equivalently,
\begin{eqnarray*}
    \int_{\{\rho=t\}\cap U_t}\chi(J)|d\rho|_h ds_h&=&\int_{\Omega_{t,s}}[\chi'(J)\langle dJ,d\rho\rangle_{h_0}+\chi(J)\Delta_{h_0}\rho]dA_{h_0}\\
    &&+\int_{\{\rho_0=s\}\cap U_t}\chi(J)\langle d\rho,\frac{d\rho_0}{|d\rho_0|_{h_0}}\rangle_{h_0}ds_{h_0}.
\end{eqnarray*}
Now, since $\chi(t)=0$ for $t>4/5$, it follows from Remark \ref{j.trace} and Lemma \ref{w21.consq} that
$$\lim_{s\to 0}\int_{\{\rho_0=s\}}\chi(J)\langle d\rho,d\rho_0\rangle_{h_0}ds_{h_0}=0.$$
Moreover, the argument from the proof of Claim \ref{j.bd} showing that $\langle d\log(J_{\lambda}),d\rho\rangle_{h_0}\in L^1(U_{\delta},dA_{h_0})$ also applies to $\langle d\log(J),d\rho\rangle_{h_0}$, and since $\chi(J)$ is supported on $\frac{1}{5}\leq J\leq \frac{4}{5}$, it follows that $\chi'(J)\langle dJ,d\rho\rangle_{h_0}\in L^1(U_{\delta},dA_{h_0})$ as well. And since $\chi(J)$ is bounded and $\Delta_{h_0}\rho\in L^1(\Sigma)$, it follows that we can take $s\to 0$ in the previous computation to get
$$\int_{\{\rho=t\}\cap U_t}\chi(J)|d\rho|_hds_h\leq\int_{\{\rho<t\}\cap U_t}|\chi'(J)\langle dJ,d\rho\rangle_{h_0}+\chi(J)\Delta_{h_0}\rho|dA_{h_0},$$
and in particular,
\begin{equation}\label{limt.flux.van}
\lim_{t\to 0}\int_{\{\rho=t\}\cap U_t}\chi(J)|d\rho|_hds_h=0.
\end{equation}
On the other hand, since $\chi(J)=1$ on $Q$, integrating the same vector field $Y$ over the box $Q\cap\{\rho>t\}$ for $t<s_\epsilon$ gives
\begin{eqnarray*}
    \int_I\chi(J)|d\rho|_hds_h-\int_{\{\rho=t\}\cap Q}\chi(J)|d\rho|_h ds_h&=&\int_{Q\cap\{\rho>t\}}\chi(J)\Delta_h\rho dA_h\\
    &\leq &0,
\end{eqnarray*}
so that
$$\int_I|d\rho|_hds_h\leq \int_{\{\rho=t\}\cap Q}\chi(J)|d\rho|_hds_h\leq \int_{\{\rho=t\}\cap U_t}\chi(J)|d\rho|_hds_h.$$
Letting $t\to 0$ and applying \eqref{limt.flux.van}, it then follows that
$\int_I|d\rho|_hds_h=0$; but $I$ is an open interval on which $|d\rho|_h>0$, giving a contradiction.

\end{proof}

\section{Removing singularities at the fixed points}\label{sect:full.reg}

By the results of Sections \ref{sect:int.reg} and \ref{sect.reg.rods}, we now see that the metric 
$$g=g_{\phi}:=h+f_1^2d\theta_1^2+f_2^2d\theta_2^2$$
obtained from our $\mathcal{S}_{g_0}$-minimizing $\phi$ defines a smooth Einstein metric on $M\setminus \cV$, where $\cV$ consists of the $2n+2$ isolated points fixed by the full $T^2$ action on $M$. It remains to check that $g$ is also smooth in a neighborhood of each of these points. In particular, Proposition \ref{v.bd.suff} reduces the desired global regularity statement to the following.
\begin{prop}\label{bdd.v}
    The function $v=v_{\phi}=u-\log\phi\in L^1(M,dV_{g_0})$ belongs to $L^{\infty}(M)$.
\end{prop}
Thus, our goal in this section is to prove Proposition \ref{bdd.v} by ruling out blow-up of $v$ near the fixed point set $\cV$. Fixing a point $p\in \cV$, we achieve this by working in the geodesic annuli 
$$A_{r,2r}(p):=B^{h_0}_{2r}(p)\setminus B^{h_0}_r(p)$$
with respect to the reference metric $h_0$ on $S$, obtaining an estimate for $\sup_{A_{r,2r}(p)}v$ which remains bounded as $r\to 0$. Namely, after identifying a constant $\sigma_r^2>0$ for which the lower bound
$$\sigma_r^2\leq C\inf_{A_{r,2r}(p)}\phi^2$$
can be proved relatively easily, the bulk of the work is devoted to establishing a complementary upper bound
$$\sup_{A_{r,2r}}e^{2u}\leq C \sigma_r^2,$$
which together give the desired uniform control
$$e^{2v}=e^{2u}\phi^{-2}\leq C \sigma_r^2\cdot \frac{C}{\sigma_r^2}=C^2\text{ on }A_{r,2r}(p).$$

\begin{remark}
In principle, much of the proof that follows could be repackaged as a blow-up argument for suitable rescaled metrics on a fixed reference annulus, but for clarity of presentation, we prefer to work directly in the shrinking annuli in the fixed reference space. 
\end{remark}

\begin{remark}\label{flat.simp}
Since all of our results thus far have been independent of the background reference metric $g_0$, we now have the freedom to fix $g_0$ in a neighborhood of $\cV$ for convenience. Namely, we may assume henceforth that there is a neighborhood $U$ of each $p\in \cV$ in which $(U,g_0)\cong B_1^4(0)$ coincides with the flat unit $4$-ball with the usual action of $S^1\times S^1\leq SO(4)$, identifying $p$ with the origin. In particular, for each $p\in \cV$, we can identify  $B_1^{h_0}(p)$ isometrically with $D_1^2(0)$ such that
\begin{equation}
\varphi_1(x,y)=x,\text{ }\varphi_2(x,y)=y,\text{ and }K_{h_0}=0.
\end{equation}
Moreover, we can arrange that these balls are disjoint.
\end{remark}

\subsection{The constants $\sigma_r$ and the lower bound for $\phi$}\label{scal.fix}

To begin, fix a vertex $p\in \Sigma$ of the fundamental polygon $\Sigma$, and fix $g_0$ as in Remark \ref{flat.simp}. Then, let
$$\gamma:[0,1]\to \Sigma$$
be a unit-speed parametrization of the segment in $(\Sigma,h_0)$ for which
$$\gamma(0)=p\text{ and }\varphi_1\circ \gamma=\varphi_2\circ \gamma,$$
so that $\gamma$ lies in the fixed point set of a reflection in $L$ exchanging $\Gamma_1$ and $\Gamma_2$, and observe that
$$dist_{h_0}(\gamma(r),\partial\Sigma)^2=\varphi_1(\gamma(r))^2=\varphi_2(\gamma(r))^2=\rho_0(\gamma(r))$$
for all $r\in [0,1]$. 

For $r\in (0,\delta_0)$, we then define
\begin{equation}\label{sig.def}
\sigma_r^2:=\phi^2(\gamma(r))=\frac{\rho(\gamma(r))}{\rho_0(\gamma(r))}.
\end{equation}
Since $\inf_{\Sigma}\phi>0$, clearly $\sigma_r$ has a positive lower bound independent of $r$. Moreover, along $\gamma(r)$, notice that Proposition \ref{dpsi.above} yields a bound of the form
$$|d\log(\rho)|_{h_0}(\gamma(r))\leq \frac{C}{r},$$
so that
$$|\log(\rho\circ\gamma)'(r)|\leq \frac{C}{r},$$
and integration gives an estimate of the form
$$|\log(\rho(\gamma(r))/\rho(\gamma(s)))|\leq C'|\log(r/s)|,$$
or equivalently
\begin{equation}\label{sig.cts}
    |\log(\sigma_r/\sigma_s)|\leq C''|\log(r/s)|
\end{equation}
for any $r,s\in (0,\delta_0)$. Moreover, since $dist_{h_0}(x,\partial\Sigma)\geq \frac{r}{6}$ for any $x\in B_{r/3}(\gamma(r))$, Proposition \ref{dpsi.above} also gives a local Lipschitz bound
$$Lip(\log(\rho); B_{r/3}(\gamma(r)))\leq \frac{C}{r},$$
so that 
$$c\leq \frac{\rho(x)}{\rho(\gamma(r))}\leq C\text{ for }x\in B_{r/3}(\gamma(r))$$
for some universal constants $0<c<C<\infty$. In particular, noting that $c\leq \frac{\rho_0(x)}{r^2}\leq C$ on $B_{r/3}(\gamma(r))$ as well, it follows that
\begin{equation}\label{mid.ball.cp}
    c'\sigma_r\leq \phi(x)\leq C'\sigma_r\text{ for all }x\in B_{r/3}(\gamma(r)).
\end{equation}
The following key lemma extends the lower bound in \eqref{mid.ball.cp} to the whole circle $\partial B_r(p)$.

\begin{lem}\label{ann.phi.below}
There is a constant $C<\infty$ independent of $r\in (0,1/2)$ such that if $dist_{h_0}(x,p)=r$, then
$$\sigma_r^2\leq C\phi^2(x).$$
Together with \eqref{sig.cts}, this implies moreover that
$$\sigma_r^2\leq C\phi^2(x)\text{ whenever }\frac{r}{2}<dist_{h_0}(x,p)<2r.$$
\end{lem}
\begin{proof}
    Let $\mathcal{A}_r\subset P$ be the lift of the quarter-annulus $[B_{2r}^{h_0}(p)\setminus B_r^{h_0}(p)]\cap\Sigma$ to the six-dimensional manifold $(P,G_0)$, and observe that the equation
    $$-L_0(\phi^2)=2\lambda e^{2v}\phi^4\geq 0$$
    now holds in the classical sense on $P\setminus \mathcal{V}$, where now
    $$L_0\chi:=\Delta_{G_0}\chi=\Delta_{\mathbb{R}^6}\chi\text{ for }\chi\in C^{\infty}(B_1^{G_0}(p)),$$
    by virtue of Remark \ref{flat.simp}. In particular, it follows from the mean value property for superharmonic functions in $B^6$ that
    \begin{equation}\label{phi.har}
    r^{-6}\int_{B_r^{G_0}(x)}\phi^2dV_{G_0}\leq C\phi^2(x)
    \end{equation}
    for every $x\in \mathcal{A}_r$ with $r\in (0,1/2)$. 
    
    Next, fix one of the totally geodesic copies
    $$\Sigma'\subset (P,G_0)$$
    of $(\Sigma,h_0)$ whose $SO(3)\times SO(3)$ orbit generates $P$, and consider the annulus 
    $$A_r:=\mathcal{A}_r\cap \Sigma'.$$
    It's easy to see that $A_r$ can be covered by $N$ geodesic balls $\{B^{G_0}_{r/6}(x_j)\}_{j=1}^N$  of radius $\frac{r}{6}$, with $x_j\in A_r$ and $N$ independent of $r$. At least one of these balls, say $B^{G_0}_{r/6}(x_1)$ must have the property that its projection to $\Sigma$ contains $\gamma(r)$, and therefore must lie inside  $B_{r/3}^{h_0}(\gamma(r))$. It then follows from \eqref{mid.ball.cp} that
    $$\inf_{B^{G_0}_{r/6}(x_1)}\phi^2\geq c_1 \sigma_r^2$$
    for some $c_1>0$ independent of $r$.

    Next, there must be another ball, say $B_{r/6}^{G_0}(x_2)$ which intersects $B_{r/6}^{G_0}(x_1)$ nontrivially, so that
    $$B_{r/6}^{G_0}(x_1)\subset B_{r}^{G_0}(x_2).$$
    Combining \eqref{phi.har} with the preceding estimate then gives
    \begin{equation*}
        \inf_{B_{r/4}^{G_0}(x_2)}\phi^2\geq  C^{-1}r^{-6}\int_{B_{r/6}^{G_0}(x_1)}\phi^2dV_{G_0}
        \geq  c_2\sigma_r^2,
    \end{equation*}
    with $c_2>0$ independent of $r$. Carrying on in this way, by induction, we eventually obtain a lower bound of the form
    $$\inf_{B_{r/6}^{G_0}(x_j)}\phi^2\geq c_N\sigma_r^2$$
    for all $1\leq j\leq N$. Since $N$ can be taken independent of $r$, and the union $\{B_{r/6}^{G_0}(x_j)\}_{j=1}^N$ covers $A_r$, this gives the desired lower bound. 
    \end{proof}

As an easy consequence of Lemma \ref{ann.phi.below}, note that the same pointwise bound $C\phi^2\geq \sigma_r^2$ must also hold in the four-dimensional annuli
$$B_{2r}^{g_0}(p)\setminus B_r^{g_0}(p)\subset M,$$
which, after squaring and integrating, gives an estimate of the form
$$r^4\sigma_r^4\leq C'\int_{B_{2r}^{g_0}(p)\setminus B_r^{g_0}(p)}\phi^4 dV_{g_0}.$$
In particular, since $\phi\in L^4(M,dV_{g_0})$, sending $r\to 0$ in this inequality implies
\begin{equation}\label{rsigr.van}
\lim_{r\to 0}r\sigma_r=0,
\end{equation}
which we will use repeatedly in the rest of the argument.

In the remainder of the section, we work toward a complementary upper bound of the form $e^{2u}\leq C\sigma_r^2$ in $B_{2r}^{h_0}(p)\setminus B_r^{h_0}(p)$. Before turning to the relevant localized estimates in these annuli, we first record the following relatively simple \emph{global} pointwise bound for the ratio $f_1/f_2$, which we will use throughout the argument.

\begin{lem}\label{fi.quo.lem}
Let $W\subset \Sigma$ be a fundamental triangle for the $\mathbb{Z}_{2n+2}$ action on $\Sigma$, such that $\partial W\cap \partial\Sigma$ is a connected component of $\Gamma_1$, and $\partial W\cap \mathring \Sigma$ consists of two edges, each fixed by a reflection exchanging $\Gamma_1$ and $\Gamma_2$. Then
$$f_1\leq f_2\text{ on }W.$$
\end{lem}
\begin{proof}
Recalling the notation
$$\beta:=\alpha_1-\alpha_2=\log(f_1/f_2)=2q+\log(\varphi_1/\varphi_2),$$
note that symmetries force
$$\beta=0\text{ on }\partial W\setminus \Gamma_1$$
and $\beta\to -\infty$ along $\Gamma_1\setminus \cV$, while
\begin{equation}\label{bet.harm}
\Delta_g\beta=0\text{ on }\mathring M.
\end{equation}

For any $k\in (0,\infty)$, setting
$$\beta^+_k:=\min\{\beta^+,k\},$$
it follows that
$$\chi\beta_k^+\in W_0^{1,2}(W)$$
for any smooth function $0\leq \chi\in C^{\infty}(W)$ vanishing in a neighborhood of $\cV$. In particular, integrating against \eqref{bet.harm} gives
\begin{eqnarray*}
    \int_W\chi|d\beta_k^+|_g^2dV_g&=& -\int_W\langle d\chi,\beta_k^+d\beta\rangle_gdV_g\\
    &=&-\int_W\phi^2\langle d\chi,\beta_k^+[2dq+d\log(\varphi_1/\varphi_2)]\rangle_{g_0}dV_{g_0}.
\end{eqnarray*}
Defining
\begin{equation}\label{rv.def}
    r_{\cV}(x):=\sqrt{\varphi_1^2(x)+\varphi_2^2(x)},
\end{equation} 
so that $r_{\cV}^2$ agrees with the squared distance $\dist_{\cV,g_0}^2$ in $B_1^{g_0}(\cV)$, we can let $\chi=\chi_{\delta}$ be a linear cutoff of the form $\chi_{\delta}=\zeta(r_{\cV})$, where 
$$\zeta(t)=0\text{ for }t<\delta,\text{ }\zeta(t)=1\text{ for }t>2\delta,\text{ and }0\leq \zeta'\leq  {C}/{\delta}\text{ on }[\delta,2\delta].$$
Testing the preceding identity against this $\chi_{\delta}$, and noting that
\begin{eqnarray*}
    |\langle \varphi_1d\varphi_1+\varphi_2d\varphi_2,d\log(\varphi_1/\varphi_2)\rangle_{g_0}|&=&||d\varphi_1|_{h_0}^2-|d\varphi_2|_{h_0}^2+\left(\frac{\varphi_2}{\varphi_1}-\frac{\varphi_1}{\varphi_2}\right)\langle d\varphi_1,d\varphi_2\rangle_{h_0}|\\
    &=&0\text{ where }r_{\cV}\leq 1,
\end{eqnarray*}
by Remark \ref{flat.simp}, we then see that
\begin{eqnarray*}
    \int_W\chi_{\delta}|d\beta_k^+|_g^2dV_g&\leq &\frac{Ck}{\delta}\int_{\{r_{\cV}<2\delta\}}\phi^2|dq|_{g_0}dV_{g_0}\\
    &\leq &\frac{C k}{\delta^2}\int_{\{r_{\cV}<2\delta\}}\phi^2dV_{g_0}+Ck\int_{\{r_{\cV}<2\delta\}}\phi^2|dq|_{g_0}^2dV_{g_0}.
\end{eqnarray*}
Since $\phi\in L^4(M,dV_{g_0})$, $Vol_{g_0}(r_{\cV}<2\delta)\leq C\delta^4$, and $\int \phi^2|dq|_{g_0}^2dV_{g_0}<\infty$, we see that the right-hand side $\to 0$ as $\delta\to 0$, leaving
$$\int_W|d\beta_k^+|^2_gdV_g=0$$
for any $k$; hence, $\beta^+\equiv 0$ in $W$, as claimed.    
\end{proof}

Of course, the opposite bound holds on a fundamental triangle $W$ with $\partial W\cap\partial\Sigma=\Gamma_2$, so the sets $\{f_1\leq f_2\}$ and $\{f_2\leq f_1\}$ are exactly the unions of the $n+1$ fundamental triangles with one edge in $\Gamma_1$ or $\Gamma_2$, respectively.

\subsection{Annular lower bounds for $e^uf_i$}

We next obtain localized lower bounds for $e^u\max\{f_1,f_2\}$ in the quarter-annuli 
\begin{equation}\label{ar.def}
    A_{r,2r}(p):=[B_{2r}^{h_0}(p)\setminus B_r^{h_0}(p)]\cap\Sigma
\end{equation}
about the vertex $p\in\cV$, which play a role analogous to that of Proposition \ref{wbd.adm} in the previous section.

As a starting point, note that under the assumptions of Remark \ref{flat.simp}, we now have
$$\Delta_{h_0}u=K_{h_0}+\langle d\alpha_1,d\alpha_2\rangle_{h_0}$$
and
$$\Delta_{h_0}f_i=-\lambda e^{2u}f_i-\langle d\alpha_1,d\alpha_2\rangle_{h_0} f_i,$$
where $\alpha_i=\log(f_i)$ on the whole punctured surface $S\setminus \cV$, and since $f_1+f_2>0$ on $S\setminus \cV$, it follows that
\begin{eqnarray*}
    &&\Delta_{h_0}(u+\log(f_1+f_2))\\&=&K_{h_0}+\langle d\alpha_1,d\alpha_2\rangle_{h_0}-\frac{(\lambda e^{2u}+\langle d\alpha_1,d\alpha_2\rangle_{h_0})(f_1+f_2)}{f_1+f_2}-|d\log(f_1+f_2)|_{h_0}^2,
\end{eqnarray*}
so that
\begin{equation}\label{ulogf.sub}
-\Delta_{h_0}(u+\log(f_1+f_2))=\lambda e^{2u}+|d\log(f_1+f_2)|_{h_0}^2-K_{h_0}\geq -K_{h_0}
\end{equation}
on $S\setminus \cV$ as well. Moreover, the admissibility conclusion from Proposition \ref{rod.smth} gives
$$u+\log(f_1+f_2)=v-q+2\log(\phi)+\log(\varphi_2)=\log(\varphi_2 \phi^2)\text{ on }\Gamma_1\setminus\cV,$$
and similarly $u+\log(f_1+f_2)=\log(\varphi_1\phi^2)$ on $\Gamma_2\setminus \cV$, so that
\begin{equation}\label{ulogf.bdry}
u+\log(f_1+f_2)=\log(r_{\cV}\phi^2)\geq \log(r_{\cV}\sigma_{r_{\cV}}^2)-C\text{ along }\partial\Sigma\setminus \cV,
\end{equation}
by Lemma \ref{ann.phi.below}, where $r_{\cV}$ is given by \eqref{rv.def}. We will then combine these with a comparison argument to prove the following.

\begin{lem}\label{eufi.below}
    There is a constant $C>0$ such that for any $t\in [0,1/2]$, we have 
    \begin{equation}
        t\sigma_t^2\leq Ce^u(f_1+f_2)\text{ on }A_{t,2t}(p).
    \end{equation}
\end{lem}
\begin{proof}
   First, using the Euclidean coordinates from Remark \ref{flat.simp} near $p$, consider the harmonic function
   $$h_t(z):=\log|\frac{z}{t}+\frac{it}{z}|$$
   on $B_1^{h_0}(p)\setminus \{p\}$, and along $B_1^{h_0}(p)\cap\partial\Sigma\setminus\{p\}$, where $z=Re(z)$ or $z=iIm(z)$, note that
   $$h_t(z)=\frac 12\log|(r_{\cV}/t)^2+(t/r_{\cV})^2|\geq |\log(r_{\cV}/t)|.$$

   On the other hand, it follows from \eqref{ulogf.bdry} and \eqref{sig.cts} that
   $$u+\log(f_1+f_2)\geq \log(r_{\cV}\sigma_{r_{\cV}}^2)-C\geq \log(t\sigma_t^2)-C|\log(r_{\cV}/t)|-C$$
   along $\partial\Sigma\setminus\{p\}$, so setting
   $$\xi_t:=u+\log(f_1+f_2)-\log(t\sigma_t^2)+C(1+h_t),$$
   we see that
   $$\xi_t\geq 0\text{ on }\partial\Sigma\cap (B_1^{h_0}(p)\setminus\{p\}),$$
   while it follows from \eqref{ulogf.sub} and Remark \ref{flat.simp} that
   $$-\Delta_{h_0}\xi_t\geq 0\text{ on }B_1^{h_0}(p)\setminus\{p\}.$$
   Moreover, since $u+\log(f_1+f_2)$ is bounded on the circle $\partial B_{h_0}^1(p)$ and $h_t\geq |\log(2)|$ there for $t\in (0,1/2)$, by fixing $C$ sufficiently large, we can arrange that
   $$\xi_t\geq 0\text{ on }\partial B_1^{h_0}(p)$$
   as well. That is, on the quarter disk $Q:=B_1^{h_0}(p)\cap\Sigma$, we have
   $$-\Delta_{h_0}\xi_t\geq 0\text{ on }Q\setminus\{p\}$$
   and
   $$\xi_t\geq 0\text{ on }\partial Q\setminus\{p\}.$$

   We would like to apply the maximum principle to conclude that $\xi_t\geq 0$ throughout $Q$, but some care is needed at the vertex. Instead, we will argue directly that
   \begin{equation}\label{xi.sgn}
   \int_{\Sigma}\rho_0\xi_t\zeta dA_{h_0} \geq 0
   \end{equation}
   for every nonnegative $0\leq \zeta\in C_c^{\infty}(Q)$, which clearly implies $\xi_t\geq 0$ in $Q$. To this end, fix some $0\leq \zeta\in C_c^{\infty}(Q)$, and let $\chi$ solve
   $$-\Delta_{G_0}\chi=-L_0\chi=\zeta\text{ on }B_1^{G_0}(p)\subset P,$$
   with Dirichlet boundary condition $\chi=0$ on $\partial B_1^{G_0}(p)$. Recalling Section \ref{l0.sec}, it follows that $\chi$ is smooth, nonnegative, and $SO(3)\times SO(3)$-invariant in $B_1^{G_0}(p)$, and therefore descends to a smooth, nonnegative $\chi\in C^{\infty}(B_1^{h_0}(p))$ solving
   $$-\Delta_{h_0}(\chi\rho_0)=\zeta\rho_0\text{ on }B_1^{h_0}(p),$$
   with $\chi=0$ on $\partial B_1^{h_0}(p)$, extending $\chi$ and $\zeta$ from $Q$ to $B_1^{h_0}(p)$ in an $R_i$-even way. Next, let $0\leq \eta\in C_c^{\infty}(S\setminus \cV)$ be a smooth, nonnegative cutoff such that $\eta\equiv 1$ on the support of $\zeta$. We can then integrate $\xi_t$ against $-\Delta_{h_0}(\eta\chi\rho_0)$, and safely integrate by parts to obtain
   \begin{eqnarray*}
       &&\int_Q\xi_t(-\Delta_{h_0}[\eta\chi\rho_0])dA_{h_0}\\
       &=&-\int_{\partial Q}\xi_t\langle d(\eta\chi\rho_0),\nu\rangle_{h_0}ds_{h_0}+\int_Q\langle d\xi_t,d([\eta\chi\rho_0])\rangle_{h_0}dA_{h_0}\\
       &=&\int_{\partial \Sigma\cap Q}\xi_t\eta \chi|\frac{\partial\rho_0}{\partial\nu}|_{h_0}ds_{h_0}-\int_{\partial B_1^{h_0}(p)\cap Q}\xi_t\eta \rho_0\partial_\nu\chi ds_{h_0}+\int_Q\eta \chi\rho_0(-\Delta_{h_0}\xi_t) dA_{h_0}\\
       &\geq & 0,
   \end{eqnarray*}
   using that $-\Delta_{h_0}\xi_t\geq 0$, $\xi_t\geq 0$ along $\partial Q\setminus\{p\}$, and appealing to the nonnegativity and boundary conditions for $\chi$, and the vanishing of $\rho_0$ along $\partial\Sigma$. Expanding 
   $$-\Delta_{h_0}(\eta\chi\rho_0)=\zeta \rho_0-2\langle d\eta,d\chi\rangle_{h_0}\rho_0-\chi\Delta_{h_0}(\eta\rho_0),$$
   it follows that
   $$\int_Q\xi_t\zeta\rho_0\geq\int_Q\xi_t[2\langle d\eta,d\chi\rangle_{h_0}\rho_0+\chi\Delta_{h_0}(\eta\rho_0)]dA_{h_0},$$
   so to prove \eqref{xi.sgn}, it suffices to find $\eta=\eta_{\delta}\in C_c^{\infty}(S\setminus \cV)$ as above for which
   \begin{equation}\label{eta.suff}
   \lim_{\delta\to 0}\int_Q\xi_t[2\langle d\eta_{\delta},d\chi\rangle_{h_0}\rho_0+\chi\Delta_{h_0}(\eta_{\delta}\rho_0)]dA_{h_0}=0.
   \end{equation}
Taking $\eta_{\delta}=\eta_{\delta}(r_{\cV})$ to be a simple radial cutoff such that $\eta_{\delta}\equiv 1$ on $S\setminus B^{h_0}_{2\delta}(\cV)$, $\eta_{\delta}\equiv 0$ on $B_{\delta}^{h_0}(\cV)$, and $|d\eta_{\delta}|_{h_0}^2+|Hess_{h_0}(\eta)|_{h_0}\leq \frac{C}{\delta^2}$ on $B_{2\delta}^{h_0}(\cV)\setminus B_{\delta}^{h_0}(\cV)$, we claim that \eqref{eta.suff} holds. 

   To see this, first note that, since $v\in L^1(M,dV_{g_0})$ and $\log(\phi),q\in W^{1,2}_{g_0}(M,g_0)$, it's easy to see that
   \begin{equation}\label{xit.int}
   \xi_t\in L^1(Q,\rho_0dA_{h_0}).
   \end{equation}
   Moreover, it follows from the $R_i$-evenness of $\chi$ that $|\langle d\chi,d\varphi_i\rangle_{h_0}|\leq C\varphi_i$, and consequently 
   $$|\langle d\chi,dr_{\cV}\rangle_{h_0}|\leq C r_{\cV}.$$
   In particular, for the radial cutoff $\eta_{\delta}$ chosen above, this gives
   $$\int_Q|\xi_t\langle d\eta_{\delta},d\chi\rangle_{h_0}|\rho_0dA_{h_0}\leq C\int_{B_{2\delta}(p)\cap Q}|\xi_t|\rho_0dA_{h_0},$$
   which together with the integrability \eqref{xit.int} forces
   $$\lim_{\delta\to 0}\int_Q|\xi_t\langle d\eta_{\delta},d\chi\rangle_{h_0}|\rho_0dA_{h_0}=0.$$

   Thus, to prove \eqref{eta.suff}, and consequently \eqref{xi.sgn}, it remains to check that
   \begin{equation}\label{xi.lap.van}
       \lim_{\delta\to 0}\int_Q\xi_t\chi\Delta_{h_0}(\eta_{\delta}\rho_0)dA_{h_0}=0.
   \end{equation}
   As another consequence of the bound $|\langle d\chi,dr_{\cV}\rangle|\leq Cr_{\cV}$, note that we have
   $$|\chi-\chi(0)|\leq Cr_{\cV}^2,$$
   while it's easy to see that $|\Delta_{h_0}(\eta_{\delta}\rho_0)|\leq \frac{C}{\delta^2}\rho_0$ on $B_{2\delta}^{h_0}(p)\setminus B_{\delta}^{h_0}(p)$, which together with \eqref{xit.int} implies
   $$\lim_{\delta\to 0}\int_Q|(\chi-\chi(0))\xi_t\Delta_{h_0}(\eta_{\delta}\rho_0)|dA_{h_0}=0.$$
   Thus, replacing $\chi$ with $\chi(0)$ in \eqref{xi.lap.van}, all that remains to prove \eqref{xi.sgn} is to show that the cutoffs $\eta_{\delta}$ satisfy
   \begin{equation}\label{cut.eta.van}
       \lim_{\delta\to 0}\int_Q\xi_t\Delta_{h_0}(\eta_{\delta}\rho_0)dA_{h_0}=0.
   \end{equation}
   By inspection, clearly $h_t\in L^p(Q,dA_{h_0})$ for any $p\in [1,\infty)$, and writing 
   $$u+\log(f_1+f_2)=v+2\log\phi+\log(\varphi_1 e^q+\varphi_2e^{-q}),$$
   since $q,\log\phi \in W^{1,2}_{g_0}(M)$, it's easy to see that
   $$|\xi_t-v|\in L^4(Q,\rho_0 dA_{h_0}),$$
   so that
   \begin{eqnarray*}
       \int_Q|\xi_t-v||\Delta_{h_0}(\eta_{\delta}\rho_0)|dA_{h_0}&\leq &\left(\int_{B_{2\delta}^{h_0}(p)}|\xi_t-v|^4\rho_0\right)^{1/4}\left(\int_{B_{2\delta}^{h_0}(p)}(C/\delta^2)^{4/3}\rho_0\right)^{3/4}\\
       &\leq &\left(\int_{B_{2\delta}^{h_0}(p)}|\xi_t-v|^4\rho_0\right)^{1/4}\cdot C'\delta\\
       &\to&0
   \end{eqnarray*}
   as $\delta\to 0$. 
   
   Finally, to estimate the integral of $v\Delta_{h_0}(\eta_{\delta}\rho_0)$ over $Q$, note that $\Delta_{h_0}(\eta_{\delta}\rho_0)=-\Delta_{h_0}([1-\eta_{\delta}]\rho_0)$ in $B_1^{h_0}(p)$, and since $1-\eta_{\delta}=0$ on $\Sigma\setminus Q$, it follows that 
   \begin{eqnarray*}
    \int_Q v \Delta_{h_0}(\eta_{\delta}\rho_0)&=&-\int_Qv\Delta_{h_0}([1-\eta_{\delta}]\rho_0)\\
    &=&-\int_{\Sigma}v\Delta_{h_0}([1-\eta_{\delta}]\rho_0).
   \end{eqnarray*}

   Now, we apply Proposition \ref{weak.el} with $\chi=1-\eta_{\delta}$ to find
   \begin{eqnarray*}
        -\int_{\Sigma}v\Delta_{h_0}([1-\eta_{\delta}]\rho_0)dA_{h_0}&= & \int_{\Sigma}[1-\eta_{\delta}]\phi^{-2}(-\Delta_{h_0}\rho) dA_{h_0}\\
        &&-\int_{\Sigma}(3\langle d\phi,d([1-\eta_{\delta}]/\phi)\rangle_{h_0}+\frac{R_{g_0}}{2}[1-\eta_{\delta}])\rho_0 dA_{h_0}\\
        &&+\int_{\Sigma}([1-\eta_{\delta}]\rho_0|dq_{\phi}|_{h_0}^2+q \cdot div_{h_0}([1-\eta_{\delta}](\varphi_1d\varphi_2-\varphi_2d\varphi_1)))dA_{h_0}.
\end{eqnarray*}
Note that all of the terms on the right which are $0$th order in $1-\eta_{\delta}$ have the form $\int_{\Sigma}[1-\eta_{\delta}]fdA_{h_0}$ for $f\in L^1(Q)$, and therefore vanish as $\delta\to 0$. Moreover, since $\eta_{\delta}$ is radial and  $d\eta_{\delta}$ is supported in the region where the assumptions of Remark \ref{flat.simp} are in place, we see that $\langle d\eta_{\delta},\varphi_1d\varphi_2-\varphi_2d\varphi_1\rangle_{h_0}=0$ identically, leaving
$$\lim_{\delta\to 0}-\int_{\Sigma}v\Delta_{h_0}([1-\eta_{\delta}]\rho_0)dA_{h_0}=\lim_{\delta\to 0}\int_{\Sigma}3\langle d\log(\phi),d\eta_{\delta}\rangle \rho_0 dA_{h_0}.$$
Using that $\log(\phi)\in W^{1,2}_{g_0}(M)$, it's likewise easy to check that the limit on the right-hand side vanishes. Thus, \eqref{xi.lap.van} holds, and we can indeed conclude that
$$\xi_t\geq 0\text{ on }Q.$$

In particular, on the annulus $A_{t,2t}(p)$, where
$$|h_t(z)|=|\log|\frac{z}{t}+\frac{it}{z}||\leq C,$$
it follows that
$$u+\log(f_1+f_2)-\log(t\sigma_t^2)\geq -C(1+h_t)\geq -C'$$
with $C'$ independent of $t$, giving the desired bound.
\end{proof}

\subsection{Annular counterparts of the edge and interior estimates}

By Lemmas \ref{ann.phi.below} and \ref{eufi.below} and \eqref{sig.cts}, we now have pointwise lower bounds of the form
\begin{equation}\label{big.ann.est}
    \sigma_r^2\leq C\phi^2\text{ and }r\sigma_r^2\leq C e^u(f_1+f_2)\text{ on }A_{r/2,4r}(p),
\end{equation}
and on the wedge $W=\{\varphi_1\leq \varphi_2\}\cap B_1^{h_0}(p)$, Lemma \ref{fi.quo.lem} gives $f_2\geq f_1$, so that the lower bound for $e^u(f_1+f_2)$ is equivalent to
\begin{equation}\label{euf2.low}
    e^{u}f_2\geq c r\sigma_r^2\text{ on }W\cap A_{r/2,4r}(p).
\end{equation}
In this form, \eqref{euf2.low} has a nice geometric interpretation in terms of the auxiliary metric
$$\tilde{h}_r:=(\sigma_r^{-2}r^{-2}f_2)^2h\text{ on }W.$$
Namely, setting $A_{r/2,4r}^+:=W\cap A_{r/2,4r}(p)$ and 
$$\hat{A}_{r/2,4r}^+:=A_{r/2,4r}^+\cup R_1A_{r/2,4r}^+\subset S,$$
for any path $\gamma:[a,b]\to \hat{A}_{r/2,4r}^+$, it follows from \eqref{euf2.low} that the length of $\gamma$ with respect to $\tilde h_r$ satisfies
$$L_{\tilde{h}_r}(\gamma)=\sigma_r^{-2}r^{-2}\int_a^b[(e^uf_2)\circ \gamma]|\gamma'|_{h_0}\geq c r^{-1}L_{h_0}(\gamma),$$
and it follows that
\begin{equation}\label{dist.cp}
\dist_{\tilde{h}_r}(x,\partial \hat{A}_{r/2,4r}^+)\geq \frac{c}{r}\dist_{h_0}(x,\partial \hat{A}_{r/2,4r}^+)
\end{equation}
for any $x\in A_{r/2,4r}^+$. 

The distance comparison \eqref{dist.cp} becomes particularly useful when coupled with the observation (see e.g. \cite[Lemma 3.1]{cortierMinerbe2016stationaryVacuum}) that the 3-dimensional metric
$$\tilde{g}_r^2:=\tilde{h}_r+\sigma_r^{-4}r^{-4}\rho^2d\theta_1^2=(r^{-2}\sigma_r^{-2}f_2)^2(h+f_1^2d\theta_1^2)$$
on $Y_2$ has Ricci curvature
\begin{equation}
    Ric_{\tilde{g}^2}=2d\alpha_2\otimes d\alpha_2+2\lambda r^4\sigma_r^4 f_2^{-2}\tilde{g}^2>0,
\end{equation}
while $\alpha_2$ solves
$$-\Delta_{\tilde{g}^2}\alpha_2=-\rho^{-1}div_{\tilde{h}_r}(f_1df_2)=\lambda r^4\sigma_r^4f_2^{-2},$$
so that the Bochner formula for $|d\alpha_2|_{\tilde{g}^2}^2$ with respect to $\tilde{g}^2$ yields an estimate of the form
\begin{equation}\label{tilg.a.boch}
\Delta_{\tilde{g}^2}(|d\alpha_2|_{\tilde{g}^2}^2+\lambda r^4\sigma_r^4f_2^{-2})\geq 2(|d\alpha_2|_{\tilde{g^2}}^2+\lambda r^4\sigma_r^4f_2^{-2})^2.
\end{equation}
By combining this with \eqref{dist.cp} and a standard maximum principle argument, we obtain the following gradient estimate, which can be compared with Lemma \ref{b.bd} of the previous section.

\begin{lem}\label{da2.edge.ann}
    If $x\in A_{r,2r}(p)$ with $\varphi_2(x)>2\varphi_1(x)$, then there is a constant $C$ independent of $r$ such that
    \begin{equation}\label{da2.edge}
    f_2^{-2}|d\alpha_2|_h^2+\lambda f_2^{-2}\leq \frac{C}{\sigma_r^4r^4}.
    \end{equation}
\end{lem}

\begin{proof} 

For the metric $\tilde{g}^2=\tilde{g}^2_r$ as above, if $t<1$ and $x$ is any point for which the geodesic distance $\dist_{\tilde{g}^2}(x,\partial Y_2)>t$, then letting $\chi(y):=(\dist_{\tilde{g}^2}(y,x)^2-t^2)^2,$  another application of the Laplacian comparison theorem as in \cite[Lemma 7.1.9]{Petersen2016RiemannianGeometry}  gives
$$\Delta_{\tilde{g}^2}\chi(y) \geq -12\sqrt{\chi(y)}+8\dist_{\tilde{g}^2}(y,x)^2$$
on $B_t^{\tilde{g}^2}(x)$. Setting $$\xi:=|d\alpha_2|_{\tilde{g}^2}^2+\lambda r^4\sigma_r^4f_2^{-2}$$
and appealing to \eqref{tilg.a.boch}, it follows that, at the point $y\in B_t^{\tilde{g}^2}(x)$ where the product $\chi\xi$ achieves its maximum, we must have 
\begin{eqnarray*}
    0&\geq &\xi(y)[-12\sqrt{\chi(y)}+8\dist_{\tilde{g}^2}(y,x)^2]+2\langle d\xi,d\chi\rangle_{\tilde{g}^2}+2\chi(y)\xi^2(y)\\
    &=&\xi(y)[-12\sqrt{\chi(y)}+8\dist_{\tilde{g}^2}(y,x)^2-2\frac{|d\chi|^2(y)}{\chi(y)}+2\chi(y)\xi(y)]\\
    &\geq &\xi(y)[2\chi(y)\xi(y)-Ct^2],
\end{eqnarray*}
so that $\chi(y)\xi(y)\leq C'\sqrt{\chi(y)}\leq C't^2$ at the maximum. In particular, evaluating at $x$ gives an estimate of the form
$$|d\alpha_2|_{\tilde{g}_r^2}^2(x)+\lambda r^4\sigma_r^4f_2^{-2}\leq Ct^{-2}.$$

In particular, if $x\in A_{r,2r}(p)\cap\{\varphi_2>2\varphi_1\}\subset A_{r/2,4r}^+$, then we see that
$$dist_{h_0}(x,\partial \hat{A}_{r/2,4r}^+)>cr$$
with $c>0$ independent of $r$, and it follows from \eqref{dist.cp} that the lift of $x$ to $Y_2$ satisfies
$$\dist_{\tilde{g}_r^2}(x,\partial Y_2)\geq \dist_{\tilde{h}}(x,\partial \hat{A}_{r/2,4r}^+)\geq c .$$
Thus, we can apply the preceding estimate with $t=c$, to obtain a bound of the form
$$|d\alpha_2|_{\tilde{g}_r^2}^2+\lambda r^4\sigma_r^4f_2^{-2}\leq C,$$
which is clearly equivalent to \eqref{da2.edge}.
\end{proof}

\begin{remark}\label{u.big.me.small}
    More finely, if, for some $K<\infty$ we knew that
    \begin{equation}
        e^uf_2\geq K\sigma_r^2r\text{ on }A^+_{r/2,4r},
    \end{equation}
    the estimate \eqref{dist.cp} would take the form
    $$\dist_{\tilde{h}_r}(x,\partial \hat{A}^{+}_{r/2,4r})\geq {K}{r}\dist_{h_0}(x,\partial \hat{A}^+_{r/2,4r}),$$
    and it is easy to see that the same argument, changing only the distance comparison in the final paragraph, gives a bound of the form
    \begin{equation}
        f_2^{-2}|d\alpha_2|_h^2+\lambda f_2^{-2}\leq \frac{C K^{-2}}{\sigma_r^4r^4}
    \end{equation}
    on $A_{r,2r}(p)\cap \{\varphi_1<\varphi_2\}$.
\end{remark}

Next, we record a useful localized version of Lemma \ref{j.below} near the $\Gamma_1$ wall of the annulus, where regularity across $\Gamma_1\setminus \cV$ now allows for a simpler proof.

\begin{lem}\label{ann.drho.below}
    There is some $\delta_0>0$ such that for each $r\in (0,1/2)$, we have
    $$f_1^2|d\psi|_h^2\geq \frac{1}{2}\text{ on }A_{r,2r}(p)\cap\{3\varphi_1<\varphi_2\}\cap\{\rho<\delta_0\sigma_r^2r^2\}.$$
\end{lem}

\begin{proof}

As in the proof of Lemma \ref{j.below}, we first argue that
\begin{equation}\label{no.half}
    f_1^2|d\psi|_h^2(x)\neq \frac{1}{2}\text{ for any }x\in A_{r,2r}(p)\cap\{3\varphi_1<\varphi_2\}\cap\{\rho<\delta_0\sigma_r^2r^2\}
\end{equation}
for some $\delta_0>0$ independent of $r$. As in the previous section, we let
$$J:=f_1^2|d\psi|_h^2=f_2^{-2}|d\rho|_h^2,$$
and note that we still have the differential inequality \eqref{dj.odi}, which we rewrite as
$$|dJ(Y)|\leq 2\frac{f_1}{f_2}|d\alpha_2|_h^2+2\lambda \frac{f_1}{f_2},$$
where now $Y:=-\frac{\nabla \rho}{|\nabla\rho|^2}=-\frac{\nabla \psi}{\rho|\nabla\psi|_h^2}$. Multiplying \eqref{da2.edge} by $\rho$ and using \eqref{sig.cts}, it follows that
\begin{equation}\label{dj.mustard}
| dJ(Y)|\leq \frac{C\rho}{\sigma_r^4r^4},
\end{equation}
on $A_{r/2,4r}(p)\cap\{2\varphi_1<\varphi_2\}$.

Fixing a small $\delta_0>0$ to be chosen shortly, suppose for a contradiction that \eqref{no.half} fails, and let 
$$x\in A_{r,2r}(p)\cap\{3\varphi_1<\varphi_2\}$$
be a point where $\rho(x)<\delta_0\sigma_r^2 r^2$ and $J(x)=\frac{1}{2}$, and let 
$$\gamma:[0,T)\to A_{r/2,4r}(p)\cap\{2\varphi_1<\varphi_2\}$$
solve $\gamma'(t)=Y(\gamma(t))$ starting at $\gamma(0)=x$, where $T$ is the first time such that one of the following occurs:
$$\lim_{t\to T}\gamma(t)\in \Gamma_1,$$
$$\lim_{t\to T}\gamma(t)\in \partial (A_{r/2,4r}(p)\cap\{2\varphi_1<\varphi_2\})\setminus \Gamma_1,$$
or
$$\lim_{t\to T}|\gamma'(t)|=\infty.$$
Since $\rho$ satisfies
$$(\rho\circ \gamma)(t)=\rho(x)-t<\delta_0\sigma_r^2 r^2-t,$$
it follows that $\gamma(t)$ can run for time at most
\begin{equation}\label{t.up}
    T<\delta_0\sigma_r^2 r^2
\end{equation}
before reaching $\Gamma_1$.

Now, it follows from \eqref{dj.mustard} that
$$|(J\circ\gamma)'(t)|\leq \frac{C\rho}{\sigma_r^4r^4}\leq \frac{C\delta_0}{\sigma_r^2 r^2}$$
for all $t\in [0,T)$, so for any 
$$t\leq T,$$
we must have
$$|J(\gamma(t))-\frac{1}{2}|\leq C\delta_0^2.$$
In particular, taking $\delta_0^2<\frac{1}{4C}$, this forces
$$J(\gamma(t))\in [1/4,3/4]\text{ for all }t\in [0,T'].$$
On the other hand, with respect to the reference metric $h_0$, the velocity $\gamma'(t)=Y(\gamma(t))$ has length
$$|Y|_{h_0}=e^{-u}|d\rho|_h^{-1}=\frac{e^{-u}f_2^{-1}}{\sqrt{J}},$$
so the bound $J\circ \gamma\geq \frac{1}{4}$ together with \eqref{euf2.low}
gives
$$|\gamma'(t)|_{h_0}\leq \frac{C}{r\sigma_r^2}$$
for all $t\in [0,T)$. Thus, the maximal time $T$ must be characterized by $\gamma(T)=\lim_{t\to T}\gamma(t)\in \partial (A_{r/2,4r}(p)\cap\{2\varphi_1<\varphi_2\})$.

On the other hand, since $T<\delta_0\sigma_r^2 r^2$, we see that
$$\dist_{h_0}(x,\gamma(T))\leq \frac{CT'}{r\sigma_r^2}\leq C\delta_0r,$$
and since $\dist_{h_0}(x,\partial (A_{r/2,4r}(p)\cap\{2\varphi_1<\varphi_2\})\setminus \Gamma_1)>cr$, by taking $\delta_0<\frac{c}{C}$, we ensure that $\gamma(T)\in \Gamma_1$. In particular, we now have that
$$J(\gamma(t))\leq \frac{3}{4}$$
for all $t\in [0,T),$ but $\gamma(T)\in \Gamma_1$, so that
$$J(\gamma(T))=f_2^{-2}|d\rho|_h^2=|df_1|_h^2=1$$
by admissibility of $g$ along $\Gamma_1$, giving a contradiction. 

We've now confirmed \eqref{no.half}, and it is easy to deduce that the same statement holds with $A_{r,2r}$ replaced by the larger annulus $A_{r/2,4r}(p)$. As a consequence, for any connected component $U$ of $A_{r/2,4r}(p)\cap\{3\varphi_1<\varphi_2\}\cap\{\rho<\delta_0\sigma_r^2r^2\}$, we must have either 
$$f_1^2|d\psi|_h^2<\frac{1}{2}\text{ on }U,\text{ or }f_1^2|d\psi|_h^2>\frac{1}{2}\text{ on }U,$$
and it remains to exclude the first possibility for $\delta_0$ sufficiently small.

If $U$ is a component of $A_{r/2,4r}(p)\cap\{3\varphi_1<\varphi_2\}\cap\{\rho<\delta_0\sigma_r^2r^2\}$ with $f_1^2|d\psi|_h^2<\frac{1}{2}$ on $U$, then clearly $U\cap \Gamma_1=\varnothing$, so that
$$\partial U\subset \{\rho=\delta_0\sigma_r^2r^2\}\cup[\partial(A_{r/2,4r}(p)\cap\{3\varphi_1<\varphi_2\})]\setminus \Gamma_1.$$
Moreover, since we already have a lower bound of the form $\rho\geq c_0\sigma_r^2r^2$ on $\{3\varphi_1=\varphi_2\}$, $U$ must avoid this segment if $\delta_0<c_0$, so that 
$$\partial U\subset \{\rho=\delta_0\sigma_r^2r^2\}\cup\{r_{cV}=r/2\}\cup\{r_{\cV}=4r\}.$$
Now, letting $\xi: A_{r/2,4r}(p)\cap\{0\leq\varphi_1\leq\varphi_2\}$ solve
$$\Delta_{h_0}\xi=0\text{ on }\{r_{\cV}=r/2\text{ or }4r\}$$
and
$$\xi=1\text{ on }A_{r/2,4r}\cap (\Gamma_1\cup \{\varphi_1=\varphi_2\}),$$
we clearly have $0\leq \xi\leq 1$ pointwise, and 
$$\xi\geq c_0r^2\text{ on }A_{r,2r}(p)\cap\{0<\varphi_1<\varphi_2\}$$
for some $c_0>0$ independent of $r$. Along $\partial U$, we then have either
$$\rho=\delta_0\sigma_r^2r^2\geq \delta_0\sigma_r^2r^2\xi$$
or
$$\delta_0\sigma_r^2r^2\xi=0\leq \rho,$$
so that $\rho\geq \delta_0\sigma_r^2r^2\xi$ throughout $\partial U$, and since $\Delta_{h_0}(\rho-\delta_0\sigma_r^2r^2\xi)\leq 0$, the maximum principle forces
$$\rho\geq \delta_0\sigma_r^2r^2\xi\text{ throughout }U.$$

In particular, on $A_{r,2r}(p)\cap U$, it follows that $\rho\geq c_0\delta_0\sigma_r^2r^2$ for any such component $U$. Thus, the set $A_{r,2r}(p)\cap \{\rho<c_0\delta_0\sigma_r^2r^2\}$ must be contained in the union of components where $f_1^2|d\psi_h|^2>\frac{1}{2}$, giving the desired conclusion with $\delta_0'=c_0\delta_0$.

\end{proof}

\begin{remark}\label{u.big.rho.big}
    Per Remark \ref{u.big.me.small}, if we knew that
    \begin{equation}
        e^uf_2\geq K\sigma_r^2r\text{ on }A^{+}_{r/2,4r}:=A_{r/2,4r}(p)\cap \{\varphi_1<\varphi_2\},
    \end{equation}
    then the estimate \eqref{dj.mustard} would improve to
    $$|dJ(Y)|\leq \frac{C K^{-2}\rho}{\sigma_r^4r^4}$$
    on $A_{2r/3,3r}(p)\cap \{2\varphi_1<\varphi_2\}$. Keeping track of this estimate in the preceding proof, it follows in this scenario that 
    $$f_1^2|d\psi|_h^2\geq \frac{1}{2}\text{ on }A_{r,2r}(p)\cap\{3\varphi_1<\varphi_2\}\cap \{\rho<K^2\delta_0\sigma_r^2r^2\};$$
    that is, the smallness condition on $\rho$ can be relaxed proportionally to the lower bound on $e^{2u}f_2^2$.
\end{remark}

Since the roles of $\varphi_1$ and $\varphi_2$ can of course be exchanged by symmetry, Lemmas \ref{da2.edge.ann} and \ref{ann.drho.below} supply most of the estimates we will need in the regions
$$A_{r,2r}(p)\cap (\{\varphi_2>3\varphi_1\}\cup \{\varphi_1>3\varphi_2\}).$$
Next, we turn to the inner wedges
$$A_{r,2r}(p)\cap \{\delta<\frac{\varphi_1}{\varphi_2}<\delta^{-1}\}.$$

For any fixed $\delta>0$, throughout the wedge
$$V_{\delta}:=\{\delta<\frac{\varphi_1}{\varphi_2}<\delta^{-1}\}\cap B_1^{h_0}(p),$$
note that Proposition \ref{dpsi.above} gives the pointwise bounds
\begin{equation}\label{wedge.psi.grad}
|d\psi(x)|\leq \frac{C_{\delta}}{r_{\cV}(x)}\text{ and }|\Delta_{h_0}\psi(x)|=|2\lambda e^{2u}+|d\psi|_{h_0}^2|\leq \frac{C_{\delta}}{r_{\cV}^2},
\end{equation}
for all $x\in V_\delta$, and as discussed in Section \ref{scal.fix}, this also leads to the pointwise control
\begin{equation}\label{wedge.rho}
    \rho\leq C_{\delta}\sigma_r^2r^2 \text{ on }V_{\delta}.
\end{equation}
Combining this with rescaled versions of the standard interior $L^p$ estimates for the Laplacian on the disk, we also obtain estimates of the form
\begin{equation}\label{psi.w22}
\int_{V_{\delta}\cap A_{r/4,2r}(p)}|Hess_{h_0}(\psi)|^2dA_{h_0}\leq C_{\delta}r^{-2}
\end{equation}
for any  $r,\delta>0$. Building on these estimates, we next prove the following pointwise bounds.

\begin{lem}\label{eufi.up}
For any $r\in [0,1/2]$ and $\delta\in (0,1)$, on $A_{r,2r}(p)\cap V_{\delta}$, we have an upper bound of the form
$$e^u(f_1+f_2)\leq C_{\delta}\sigma_r^2r,$$
complementing the lower bound of Lemma \ref{eufi.below}. Moreover, we have the Dirichlet bound
\begin{equation}\label{log.dir.bd}
    \int_{V_{\delta}\cap A_{r,2r}(p)}|d\log(f_1+f_2)|_{h_0}^2\leq C_{\delta}.
\end{equation}
\end{lem}

\begin{proof}

 As a starting point, consider the function
 $$w:=u+\log(f_1+f_2)-\log(r\sigma_r^2)+C,$$
 where, appealing to Lemma \ref{eufi.below}, we can choose $C$ large enough, independent of $r$, such that
 $$w\geq 1\text{ on }A_{r/2,4r}(p).$$
 Then, recall from the proof of Lemma \ref{eufi.below} that
 $$-\Delta_{h_0}w=|d\log(f_1+f_2)|_{h_0}^2+\lambda e^{2u}$$
 on $\mathring \Sigma$; in particular, $w$ is superharmonic, so that
 \begin{equation}\label{mean.below}
    w(x)\geq \frac{1}{\pi s^2}\int_{B_s(x)}w
 \end{equation}
 whenever $B_s(x)\subset \mathring\Sigma$. Next, we will show that a mean value bound opposite to \eqref{mean.below} holds away from $\Gamma_1\cup \Gamma_2$, to obtain a pointwise Harnack estimate for $w$ in the wedge $V_{\delta}\cap A_{r,2r}(p)$.

 To this end, recall the expression \eqref{base.hess} for the trace-free Hessian identity with respect to $h_0$, which we rewrite as
 \begin{eqnarray*}
     Hess_{h_0}(\psi)^{TF}&=&(d\psi\otimes du+du\otimes d\psi-d\alpha_1\otimes d\alpha_2-d\alpha_2\otimes d\alpha_1)^{TF}\\
     &=&(d\psi\otimes du+du\otimes d\psi-\frac{1}{2}d\psi\otimes d\psi-\frac{1}{2}d\beta\otimes d\beta)^{TF},
 \end{eqnarray*}
 where again $\beta=\alpha_1-\alpha_2$. Expanding the $du$ terms as
 $$du=dw-d\log(f_1+f_2)=dw-\frac{1}{2}d\psi-\frac{1}{2}\frac{f_1-f_2}{f_1+f_2}d\beta,$$
 it then follows that, for some universal constant $C>0$, 
$$|(d\beta\otimes d\beta)^{TF}|\leq C(|Hess_{h_0}(\psi)|+|d\psi||dw|+|d\psi|^2+|d\psi||d\beta|),$$
and consequently
\begin{equation}\label{dbeta.band}
|d\beta|_{h_0}^2\leq C(|Hess_{h_0}(\psi)|+2|d\psi|^2+|d\psi||dw|)
\end{equation}
On the other hand, we also have
$$|\Delta_{h_0}w|\leq |d\log(f_1+f_2)|_{h_0}^2+\lambda e^{2u}\leq \frac{1}{2}(|d\psi|^2+|d\beta|^2)+\lambda e^{2u},$$
and combining these estimates with \eqref{wedge.psi.grad} gives
\begin{equation}\label{lap.w.bd}
-\Delta_{h_0}w=|\Delta_{h_0}w|\leq \frac{C_{\delta}}{r^2}+C|Hess_{h_0}\psi|+\frac{C_{\delta}}{r}|dw|
\end{equation}
on the wedge $V_{\delta}\cap A_{r/2,4r}(p)$. Using \eqref{psi.w22} to control $Hess_{h_0}(\psi)$, applying a simple Moser iteration to \eqref{lap.w.bd} leads to a pointwise upper bound of the form
\begin{equation}\label{mean.above}
    w(x)\leq \frac{C_{\delta}}{s^2}\int_{B_s(x)}w
\end{equation}
whenever $B_s(x)\subset V_{2\delta}\cap A_{r/2,4r}(p).$ 

Combining \eqref{mean.above} and \eqref{mean.below} leads to a Harnack estimate of the form $\frac{w(x)}{w(y)}\leq C_{\delta}$ whenever $x,y\in V_{\delta}\cap A_{r,2r}(p)$ with $|x-y|<\delta r/4$, and since any two points in $V_{\delta}\cap A_{r,2r}(p)$ can be connected by $\leq k_{\delta}$ segments in $V_{\delta}\cap A_{r,2r}(p)$ of length $\delta r/4$, it follows that
\begin{equation}\label{w.harn}
\frac{w(x)}{w(y)}\leq C'_{\delta}\text{ for any }x,y\in V_{\delta}\cap A_{r,2r}(p).
\end{equation}
Moreover, an identical bound clearly holds with $A_{r,2r}(p)$ replaced with the larger annulus $A_{r/2,4r}(p)$.

Now, let 
$$Q_{r,\delta}:=\frac{1}{\sigma_r^2r}\inf\{e^uf_2(x)\mid x\in A_{r/2,4r}(p)\cap V_{\delta}\cap \{\varphi_1<\varphi_2\}\}$$
with $\delta<\frac{1}{10}$, and recall that the function $u+\alpha_2\in C^{\infty}(S\setminus \Gamma_2)$ solves
$$-\Delta_{h_0}(u+\alpha_2)=|d\alpha_2|_{h_0}^2+\lambda e^{2u}.$$
In particular, $u+\alpha_2$ is superharmonic on $S\setminus \Gamma_2$, so that
$$(u+\alpha_2-\log(\sigma_r^2r))(x)\geq \frac{1}{\pi s^2}\int_{B_s^{h_0}(x)}(u+\alpha_2-\log(\sigma_r^2r))$$
whenever $B_s^{h_0}(x)\subset S\setminus \Gamma_2$. If $x\in A_{r/2,4r}(p)\cap\{\varphi_1<\varphi_2\}$, then we can apply this with $s=\frac{r}{4}$, so that $|B_{r/4}(x)\cap V_{\delta}\cap A_{r/2,4r}(p)|\geq c_0r^2$, and obtain a pointwise estimate of the form
$$u+\alpha_2-\log(\sigma_r^2r)\geq c_0\log(Q)\text{ on }A_{r/2,4r}(p).$$

In particular, setting $K=K_{\delta,r}:=Q^{c_0}$, it follows that $e^uf_2\geq K\sigma_r^2r$ on $A_{r/2,r}(p)$, so by Remark \ref{u.big.rho.big}, we have the bound
$$f_1^2|d\psi|_h^2\geq \frac{1}{2}\text{ on }A_{r,2r}(p)\cap \{3\varphi_1<\varphi_2\}\cap\{\rho<K^2\delta_0\sigma_r^2r^2\}.$$
On the other hand, fixing any point $x_r\in A_{r,2r}(p)\cap \{4\varphi_1=\varphi_2\}$, we can apply  \eqref{wedge.rho} to see that
$$\rho(x_r)\leq C_0\sigma_r^2r^2,$$
so \emph{if} $K_{\delta,r}^2>\delta_0^{-1}C_0$, it follows that
$$\frac{1}{2}\leq f_1^2(x_r)|d\psi|_h^2(x_r)=\rho(x_r)^2f_2(x_r)^{-2}e^{-2u(x_r)}|d\psi|_{h_0}^2(x_r).$$
Then, applying \eqref{wedge.rho}, \eqref{wedge.psi.grad}, and $e^{-u}f_2^{-1}\leq \frac{1}{K\sigma_r^2r}$, we arrive at a bound of the form
$$\frac{1}{2}\leq C_1(\sigma_r^2r^2)^2\cdot (K\sigma_r^2r)^{-2}\cdot \frac{C_1}{r^2}=\frac{C_1^2}{K^2}.$$
That is, either $K_{\delta,r}^2\leq \delta_0^{-1}C_0$, or $K_{\delta,r}^2\leq 2C_1^2$, and in either case we see that
$$K_{\delta,r}\leq C_2$$
with $C_2$ independent of $r$. 

In particular, it follows that $Q_{\delta,r}=c_0^{-1}\log(K)\leq C'$ as well, so there exists some $y_0$ in $A_{r/2,4r}(p)\cap \{\varphi_1<\varphi_2\}$ with 
$$(u+\alpha_2)(y_0)-\log(\sigma_r^2r)\leq C,$$
with $C$ independent of $r$. Since $f_1(y_0)\leq f_2(y_0)$ holds automatically in this domain, this then gives
$$w(y_0)\leq u(y_0)+\log(2f_2(y_0))-\log(\sigma_r^2r)\leq C+\log(2),$$
and finally it follows from the Harnack estimate \eqref{w.harn}, with $y=y_0$, that
$$w(x)\leq C_{\delta}''\text{ for any }x\in V_{\delta}\cap A_{r,2r}(p),$$
In particular, this gives the desired upper bound for $e^u(f_1+f_2)$.

To prove \eqref{log.dir.bd}, first note that by \eqref{sig.cts}, and since $\delta>0$ was arbitrary, we have
$$w\leq C_{\delta}''\text{ on }V_{\delta}\cap A_{r/2,4r}(p)$$
as well. Now, let $\chi\in C_c^{\infty}(V_{\delta}\cap A_{r/2,4r}(p))$ be a cutoff such that 
$$\chi\equiv 1\text{ on }V_{2\delta}\cap A_{r,2r}(p)$$
and $|\Delta_{h_0}\chi|\leq \frac{C_{\delta}}{r^2}$. Integrating
$$-\Delta_{h_0}w=|d\log(f_1+f_2)|_{h_0}^2+\lambda e^{2u}$$
against $\chi$ then gives a bound of the form
$$\int_{V_{2\delta}\cap A_{r,2r}(p)}|d\log(f_1+f_2)|_{h_0}^2\leq \int |\Delta_{h_0}\chi||w|\leq C_{\delta}.$$

\end{proof}

Next, note that we can conveniently write 
$$\log(f_1+f_2)=\frac{\psi}{2}+\log(2\cosh(\beta/2)),$$
and since $|d\psi|_{h_0}^2\leq \frac{C_{\delta}}{r^2}$ pointwise, it follows from \eqref{log.dir.bd} that
$$\int_{V_{\delta}\cap A_{r,2r}(p)}|d\log(2\cosh(\beta/2))|_{h_0}^2\leq C_{\delta}$$
as well. In particular, since $\beta=0$ along the central geodesic $\{\varphi_1=\varphi_2\}=\{f_1=f_2\}$, it follows from a simple rescaling of a Poincar\'e-type inequality for functions on $V_{\delta}\cap A_{1/2,1}(p)$ vanishing along $\gamma$ that 
$$\int_{V_{\delta}\cap A_{r,2r}(p)}|\log(\cosh(\beta/2))|^2dA_{h_0}\leq C_{\delta}r^2.$$
On the other hand, a straightforward calculation gives the pointwise bound 
$$|\beta|\leq 2\log(\cosh(\beta/2))+2\log 2,$$
which together with the previous bound gives the rescaled $L^2$ bound
\begin{equation}\label{beta.l2}\int_{V_{\delta}\cap A_{r,2r}(p)}|\beta|^2dA_{h_0}\leq C_{\delta}'r^2
\end{equation}
for $\beta=\alpha_1-\alpha_2$. Starting from \eqref{beta.l2}, we can now amass a large collection of sharp pointwise bounds in $V_{\delta}.$

\begin{lem}\label{wedge.bonanza}
    Let $V_{\delta}:=\{\delta<\frac{\varphi_2}{\varphi_1}<\delta^{-1}\}\cap B_1^{h_0}(p)$ as above. Then for any $x\in V_{\delta}$ with $r_{\cV}(x)=r\in (0,1/2),$ we have bounds of the form
    \begin{equation}\label{wedge.al.bd}
    |\alpha_1(x)-\log(r\sigma_r)|+|\alpha_2(x)-\log(r\sigma_r)|\leq C_{\delta}
    \end{equation}
    \begin{equation}\label{wedge.u.bd}
    |u(x)-\log(\sigma_r)|\leq C_{\delta}
    \end{equation}
    and
    \begin{equation}\label{wedge.dal.bd}
        |d\alpha_i|_{h_0}(x)\leq \frac{C_{\delta}}{r}.
    \end{equation}
\end{lem}
\begin{proof}
Let $x\in V_{2\delta}$ with $r_{\cV}(x)=r$, so that
$$B_{\delta r/10}^{h_0}(x)\subset V_{\delta}\cap A_{2r/3,4r/3}.$$
Recall that $\beta=\alpha_1-\alpha_2$ satisfies
$$\Delta_g\beta=0\text{ on }\mathring M,$$
or equivalently
$$\Delta_{h_0}\beta+\langle d\psi,d\beta\rangle_{h_0}=0$$
on $\mathring\Sigma$. In particular, since 
$$|d\psi|\leq \frac{C_{\delta}}{r}\text{ on }V_{\delta}\cap A_{2r/3,4r/3},$$
it follows that $|\Delta_{h_0}\beta|_{h_0}\leq \frac{C_{\delta}}{r}|d\beta|_{h_0}$ on the ball $B_{\delta r/10}^{h_0}(x)$, and a simple Moser iteration argument gives
$$|\beta(x)|^2\leq \frac{C_{\delta}}{r^2}\int_{B_{\delta r/10}(x)}|\beta|^2.$$
On the other hand, \eqref{beta.l2} bounds the right-hand side by another constant $C_{\delta}'$, giving the pointwise bound
$$|\beta|\leq C_{\delta}'\text{ in }V_{2\delta}.$$
Moreover, we've already seen in Section \ref{scal.fix} that 
$$|\psi-\log(\sigma_r^2r^2)|\leq C_{\delta}\text{ on }V_{2\delta},$$
so putting these together gives \eqref{wedge.al.bd}.

Next, combining \eqref{wedge.al.bd} with Lemmas \ref{eufi.below} and \ref{eufi.up} controlling $|u+\log(f_1+f_2)-\log(r\sigma_r^2)|$ from above and below immediately gives \eqref{wedge.u.bd}. 

With \eqref{wedge.u.bd} in place, we now have the uniform equivalence of metrics
\begin{equation}\label{met.equiv}
c_{\delta} h\leq \sigma_r^2 h_0\leq C_{\delta}h\text{ on }V_{\delta}\cap A_{r/2,4r}(p),
\end{equation}
and in particular we can translate estimates for the Einstein metric $g$ to estimates in the reference metric $h_0$. In particular, recalling from Section \ref{sec.hess.id} the refined Bochner inequalities
$$\Delta_g\frac{1}{2}|d\alpha_i|_g^2\geq |d\alpha_i|_g^4,$$
the same standard maximum principle and Laplacian comparison argument used in the proof of Lemma \ref{da2.edge.ann} gives the comparison
$$|d\alpha_i|_h(x)^2\leq \frac{C}{\dist_h(x,\partial\Sigma)^2}.$$
If $x\in A_{r,2r}(p)$, it then follows from \eqref{met.equiv} that
$$\dist_h(x,\partial\Sigma)\geq c_{\delta}'\sigma_r r$$
and
$$|d\alpha_i|_h(x)^2\geq c_{\delta}\sigma_r^{-2}|d\alpha_i|_{h_0}(x)^2,$$
so that the gradient estimate in $g$ is equivalent to an estimate of the form
$$|d\alpha_i|_{h_0}(x)^2\leq \frac{C_{\delta}'}{r^2},$$
giving \eqref{wedge.dal.bd}.
    
\end{proof}

Note that if the estimate \eqref{wedge.u.bd} held throughout $A_{r,2r}(p)$, then we could immediately conclude, but a little more work is needed to extend the desired bound $e^{2u}\leq C\sigma_r^2$ to a neighborhood of $\Gamma_i\cap A_{r,2r}(p)$.
With this in mind, we record a few more simple consequences of the estimates above that play a role in the final argument of the next subsection. 

First, fix a nonnegative cutoff function $\chi$ such that $\chi\in C_c^{\infty}(A_{r/2,4r}(p)\cap \{\varphi_2>\varphi_1\})$ and
$$\chi\equiv 1\text{ on }A_{r,2r}(p)\cap \{\varphi_2>2\varphi_1\}$$
and $|d\chi|_{h_0}^2+|\Delta_{h_0}\chi|\leq Cr^{-2}$. Integrating the equation 
$$-\Delta_{h_0}(u+\alpha_2-\log(r\sigma_r^2))=|d\alpha_2|_{h_0}^2+\lambda e^{2u}$$
against $\chi$ gives
\begin{equation}\label{edge.dir.1}
\int_{A_{r,2r}(p)\cap\{\varphi_2>2\varphi_1\}}(|d\alpha_2|_{h_0}^2+\lambda e^{2u})\leq \frac{C}{r^2}\int_{A_{r/2,4r}(p)\cap \{\varphi_2>\varphi_1\}}|u+\alpha_2-\log(r\sigma_r^2)|.
\end{equation}
To control the right-hand side, recall that \eqref{eufi.below} gives a lower bound of the form
\begin{equation}\label{u.a.log.low}
u+\alpha_2-\log(r\sigma_r^2)\geq -C_1
\end{equation}
on $A_{r/4,8r}(p)\cap\{\varphi_2>\varphi_1\}$ for some constant $C_1$, while Lemma \ref{wedge.bonanza} supplies a corresponding uniform upper bound in the overlapping inner wedge $A_{r/4,8r}(p)\cap\{\varphi_1<\varphi_2<\delta^{-1}\varphi_1\}$ for any fixed $\delta>0$. In particular, taking $\delta_0>0$ sufficiently small independent of $r$, we can find a collection of points $x_1,\ldots,x_k\in A_{r/2,4r}(p)\cap\{\varphi_1=\delta_0\varphi_2\}$, with $k$ independent of $r$, such that
$$A_{r/2,4r}\cap\{\varphi_1<\delta_0\varphi_2\}\subset \bigcup_{j=1}^kB_{r/4}^{h_0}(x_j)\subset \Sigma\cup R_1\Sigma.$$
Since $u+\alpha_2$ is superharmonic throughout $S\setminus \Gamma_2$, it follows that
$$\frac{c}{r^2}\int_{B_{r/4}^{h_0}(x_j)}(u+\alpha_2-\log(r\sigma_r^2))\leq (u+\alpha_2-\log(r\sigma_r^2))(x_j)\leq C_{\delta_0},$$
which together with the lower bound \eqref{u.a.log.low} gives
$$\int_{B_{r/4}^{h_0}(x_j)}|u+\alpha_2-\log(r\sigma_r^2)|\leq C_{\delta_0}r^2.$$
Summing over $j=1,\ldots,k$, it follows that
$$\int_{A_{r/2,4r}(p)\cap\{\varphi_1<\delta_0\varphi_2\}}|u+\alpha_2-\log(r\sigma_r^2)|\leq C_{\delta_0}r^2,$$
which together with the pointwise bounds of Lemma \ref{wedge.bonanza} gives
\begin{equation}\label{u.a2.l1}
\int_{A_{r/2,4r}(p)\cap\{\varphi_2>\varphi_1\}}|u+\alpha_2-\log(r\sigma_r^2)|\leq C_{\delta_0}r^2.
\end{equation}

Applying this to the right-hand side of \eqref{edge.dir.1}, it follows that
\begin{equation}\label{rhs.ann.ctrl}
    \int_{A_{r,2r}(p)\cap \{\varphi_2>2\varphi_1\}}(|d\alpha_2|_{h_0}^2+\lambda e^{2u})\leq C,
\end{equation}
with $C$ independent of $r\in [0,1/2]$. Moreover, since the pointwise bound \eqref{wedge.al.bd} for $|\alpha_2-\log(r\sigma_r)|$ holds on the intersection of $A_{r,2r}(p)\cap\{\varphi_2>2\varphi_1\}$ with the wedge $\{\varphi_2<10\varphi_1\}$, a simple Poincar\'e-type inequality then upgrades the preceding Dirichlet estimate to a bound of the form
\begin{equation}\label{a2.w12}
\int_{A_{r,2r}(p)\cap\{\varphi_1<\varphi_2\}}(r^{-2}|\alpha_2-\log(r\sigma_r)|_{h_0}^2+|d\alpha_2|_{h_0}^2)dA_{h_0}\leq C.
\end{equation}
In particular, by rescaling the standard Moser-Trudinger inequality \cite{Trudinger1967OrliczSpaceImbeddings} on the region $\{\varphi_1<\varphi_2\}\cap A_{1/2,1}(p)$ to bound
$$e^{\alpha_2-\log(r\sigma_r)}=f_2/(r\sigma_r),$$ it follows moreover that
\begin{equation}\label{f2.lps}
    \int_{A_{r,2r}(p)\cap\{\varphi_1<\varphi_2\}}[(r\sigma_r)^{-k}f_2^k+(r\sigma_r)^kf_2^{-k}]\leq C_kr^2
\end{equation}
for any $k\in [0,\infty).$ 

Recalling that $\rho\leq f_2^2$ on $\{\varphi_2\geq\varphi_1\}$ and arguing symmetrically on $\{\varphi_1\geq \varphi_2\}$, it follows from \eqref{f2.lps} that
\begin{equation}\label{rho.ann.lps}
    \int_{A_{r,2r}(p)}(r^{-2}\sigma_r^{-2}\rho)^k\leq C_k'r^2
\end{equation}
for every $k\in [0,\infty)$ as well, and by \eqref{sig.cts}, each of the preceding estimates formulated on $A_{r,2r}(p)$ hold, with a different constant, on the larger annulus $A_{r/8,8r}(p)$. 

Finally, noting that the estimates \eqref{rho.ann.lps} and \eqref{a2.w12} extend to the annulus $A_{r/8,8r}(p)\cap \{\varphi_1<\varphi_2\}$, since the equation
$$\lambda e^{2u}=-div_{h_0}(\rho d\alpha_2)$$
holds classically on $S\setminus \Gamma_2$, multiplying by
by a nonnegative cutoff supported in $A_{r/8,8r}(p)\cap \{\varphi_1<\varphi_2\}$ and $\equiv 1$ on $A_{r/4,4r}(p)\cap\{2\varphi_1<\varphi_2\}$ leads to a bound of the form
\begin{eqnarray*}
    \int_{A_{r/4,4r}(p)\cap\{2\varphi_1<\varphi_2\}}e^{2u}\rho&\leq &\frac{C}{r}\int_{A_{r/4,4r}(p)\cap\{\varphi_1<\varphi_2\}}\rho|d\alpha_2|_{h_0}dA_{h_0}\\
    &\leq & \frac{C}{r}\|\rho\|_{L^2(A_{r/4,4r}(p))}\|d\alpha_2\|_{L^2(A_{r/4,4r}(p)\cap \{\varphi_1<\varphi_2\}}\\
    &\leq &C'\sigma_r^{2}r^2
\end{eqnarray*}
At the same time, since $-\Delta_{h_0}\rho=2\lambda e^{2u}\rho$ globally on $S$, for any $a\in [1,2)$, a rescaling of standard Green's function estimates on the disk gives
$$\left(s^{a-2}\int_{B_s(x_0)}|d\rho|_{h_0}^a\right)^{1/a}\leq C_a\int_{B_{2s}(x_0)}2\lambda e^{2u}\rho+C_as^{-2}\int_{B_{2s}(x_0)}|\rho|$$
on any disk $B_s^{h_0}(x_0)\subset S$. In particular, applying this to $x_0\in A_{r,2r}(p)\cap\{3\varphi_1<\varphi_2\}$ with $s=r/10$, we can use the preceding estimate and \eqref{rho.ann.lps} to control the right-hand side, finding
$$\left(r^{a-2}\int_{B_{r/10}(x_0)}|d\rho|_{h_0}^a\right)^{1/a}\leq C_a\sigma_r^2r^2$$
for $s=r/10$ and $a\in [1,2)$. A simple Vitali covering argument then gives a bound of the same form for the $\|d\rho\|_{L^a(A_{r,2r}(p)\cap\{3\varphi_1<\varphi_2\}}$, and using Lemma \ref{wedge.bonanza} to supply comparable estimates in the wedge $A_{r,2r}(p)\cap V_{1/3}$, we see finally that
\begin{equation}\label{drho.la}
    \|d\rho\|_{L^a(A_{r,2r}(p)}\leq C_a\sigma_r^2r^2r^{(2-a)/a}
\end{equation}
for any $a\in [1,2)$.

\subsection{An upper curvature bound and the proof of Proposition \ref{bdd.v}}

To finish proving the desired bound $e^{2u}\leq C\sigma_r^2$ on $A_{r,2r}(p)$, we first prove one more useful estimate: a uniform upper bound on the curvature $K_h$ on the whole surface $S$.

\begin{lem}
    We have $K_h\leq \lambda/3$ on $S$.
\end{lem}
\begin{proof}
To begin, recall that, in $\mathring\Sigma$,
$$\Delta_h\psi=-2\lambda-|d\psi|_h^2$$
and
$$\Delta_h\beta=-\langle d\psi,d\beta\rangle_h,$$
so applying the Bochner formula to each gives
$$\frac{1}{2}\Delta_h|d\psi|_h^2=|Hess_h(\psi)|_h^2+K_h|d\psi|_h^2-\langle d\psi,d|d\psi|_h^2\rangle_h$$
and
$$\frac{1}{2}\Delta_h|d\beta|_h^2=|Hess_h(\beta)|_h^2+K_h|d\beta|_h^2-\langle d\beta,d\langle d\psi,d\beta\rangle_h\rangle_h.$$
Writing $\langle d\alpha_1,d\alpha_2\rangle=\frac{1}{4}(|d\psi|^2-|d\beta|^2)$, subtracting the second Bochner formula from the first then gives
\begin{eqnarray*}
\Delta_h\langle d\alpha_1,d\alpha_2\rangle_h&=&\frac{1}{2}(|Hess_h(\psi)|_h^2-|Hess_h(\beta)|_h^2)+2K_h\langle d\alpha_1,d\alpha_2\rangle_h\\
&&-\frac{1}{2}\langle d\psi,d|d\psi|^2\rangle+\frac{1}{2}\langle d\beta,d\langle d\psi,d\beta\rangle\rangle\\
&=&2\langle Hess_h(\alpha_1),Hess_h(\alpha_2)\rangle_h+2K_h\langle d\alpha_1,d\alpha_2\rangle_h\\
&&-\langle d\psi,d\langle d\alpha_1,d\alpha_2\rangle\rangle
-2Hess_h(\psi)(\nabla_h\alpha_1,\nabla_h\alpha_2)
\end{eqnarray*}
and with a little more rearrangement, recalling that $K_h=-\langle d\alpha_1,d\alpha_2\rangle_h,$ we arrive at the convenient identity
$$\Delta_hK_h+3\langle d\psi,dK_h\rangle_h=4K_h^2-2\rho^{-1}\langle Hess_h(f_1),Hess_h(f_2)\rangle_h.$$
Now, recalling \eqref{hess.id.0}, we see that
\begin{eqnarray*}
    \frac{2}{\rho}\langle Hess_h(f_1),Hess_h(f_2)\rangle_h&\leq &\frac{1}{2}|f_1^{-1}Hess_h(f_1)+f_2^{-1}Hess_h(f_2)|_h^2\\
    &=&(K_h-\lambda)^2,
\end{eqnarray*}
and using this in the identity for $\Delta_hK_h$ yields
\begin{equation}\label{lap.k.miracle}
\Delta_hK_h+3\langle d\psi,dK_h\rangle_h\geq 4K_h^2-(K_h-\lambda)^2.
\end{equation}

If we knew that $K_h$ achieved its maximum in the interior of $\Sigma$, then we could arrive at the desired bound by applying the maximum principle directly to \eqref{lap.k.miracle}. Of course, we cannot argue directly in this way, but using \eqref{lap.k.miracle}, we see that applying $\Delta_h$ to $\rho^{3/2}(3K_h-\lambda)$ gives
\begin{eqnarray*}
    \Delta_h(\rho^{3/2}(3K_h-\lambda))&\geq &(3K_h-\lambda)[-3\lambda \rho^{3/2}+\frac{3}{4}\rho^{3/2}|d\psi|^2]\\
    &&+3\rho^{3/2}(3K_h^2+2K_h\lambda-\lambda^2)\\
    &\geq &\rho^{3/2}(3K_h-\lambda)(3K_h+\frac{3}{4}|d\psi|^2)
\end{eqnarray*}
on $\mathring\Sigma$. In particular, setting 
$$w:=\rho^{3/2}(3K_h-\lambda)^+,$$
it follows that $w$ is weakly subharmonic in $\mathring\Sigma$, with $w=\frac{\partial w}{\partial\nu}=0$ along $\partial\Sigma\setminus \cV$, so it can be extended $R_i$-evenly to a weakly subharmonic function 
$$w: S\setminus \cV\to [0,\infty).$$
Thus, once we have established continuity
\begin{equation}\label{corner.van}
\lim_{x\to p}w(x)=0
\end{equation}
at the vertex $p\in \cV$, we can conclude by the maximum principle, as desired.

For \eqref{corner.van}, first note that an application of the mean value inequality gives an estimate of the form
$$w(x)\leq Cr_{\cV}(x)^{-2}\int_{B_{r_{\cV}(x)/2}(x)}w,$$
so the pointwise vanishing \eqref{corner.van} is equivalent to the integral vanishing
\begin{equation}\label{w.ann.van}
\lim_{r\to 0}r^{-2}\int_{A_{r,2r}(p)}w=0.
\end{equation}
For any fixed $\delta>0$, in the central wedge $V_{\delta}\cap A_{r,2r}(p)$, \eqref{wedge.dal.bd} and \eqref{wedge.u.bd} give
$$|K_h|=|e^{-2u}\langle d\alpha_1,d\alpha_2\rangle_{h_0}|\leq C_{\delta}\sigma_r^{-2}r^{-2},$$
while
$$\rho\leq C_{\delta}\sigma_r^2r^2,$$
so that
\begin{equation*}
    \frac{1}{r^2}\int_{V_{\delta}\cap A_{r,2r}(p)}wdA_{h_0}\leq  C_{\delta}'(\sigma_r^2r^2)^{3/2}\sigma_r^{-2}r^{-2}
    =C_{\delta}'\sigma_rr\to 0
\end{equation*}
as $r\to 0$, by \eqref{rsigr.van}. 

It remains to check vanishing on the complementary regions $\{\varphi_1<\delta\varphi_2\}$ and $\{\varphi_2<\delta\varphi_1\}$ for some fixed $\delta\in (0,1)$. To this end, recall that Lemma \ref{da2.edge.ann} gives 
$$|d\alpha_2|_h^2\leq C\frac{f_2^2}{\sigma_r^4r^4}\text{ on }\{\varphi_2>2\varphi_1\}\cap A_{r,2r}(p),$$
and writing $K_h=|d\alpha_2|_h^2-\langle d\psi,d\alpha_2\rangle_h,$ it follows that 
$$|K_h|\leq |d\psi|_h|d\alpha_2|_h+|d\alpha_2|_h^2\leq \frac{Cf_2}{\sigma_r^2 r^2}(|d\psi|_h+f_2\sigma_r^{-2}r^{-2})$$
on $\{\varphi_2>2\varphi_1\}\cap A_{r,2r}(p)$. Since $\rho\leq f_2^2$ on $\{\varphi_2>\varphi_1\}$, we therefore see that
\begin{eqnarray*}
    \int_{A_{r,2r}(p)\cap \{\varphi_2>2\varphi_1\}}\rho^{3/2}|K_h|&\leq &C\int_{A_{r,2r}(p)\cap\{\varphi_2>2\varphi_1\}}\sigma_r^{-2}r^{-2}f_2^2|d\rho|_h\\
    &&+C\int_{A_{r,2r}(p)\cap\{\varphi_2>2\varphi_1\}}\sigma_r^{-4}r^{-4}f_2^5,
\end{eqnarray*}
and since Lemma \ref{eufi.below} gives the lower bound $Ce^uf_2\geq \sigma_r^2r$ on this domain, we can use
$$|d\rho|_h=e^{-u}f_2^{-1}f_2|d\rho|_{h_0}\leq C\sigma_r^{-2}r^{-1}f_2|d\rho|_{h_0}$$
on the right-hand side to obtain
\begin{eqnarray*}
    \int_{A_{r,2r}(p)\cap \{\varphi_2>2\varphi_1\}}\rho^{3/2}|K_h|&\leq &C'\int_{A_{r,2r}(p)\cap\{\varphi_2>2\varphi_1\}}\sigma_r^{-4}r^{-3}f_2^3|d\rho|_{h_0}\\
    &&+C\int_{A_{r,2r}(p)\cap\{\varphi_2>2\varphi_1\}}\sigma_r^{-4}r^{-4}f_2^5.
\end{eqnarray*}
Finally, applying Young's inequality to the first term on the right-hand side, and appealing to \eqref{drho.la} with $a=3/2$ and \eqref{f2.lps} with $k=9$ yields
$$\int_{A_{r,2r}(p)\cap\{\varphi_2>\varphi_1\}}(\sigma_r^{-4}r^{-3}f_2^3|d\rho|_{h_0})\leq C\sigma_r r^3,$$
and similarly, by applying \eqref{f2.lps} to the second term with $k=6$, we arrive at the estimate
$$\int_{A_{r,2r}(p)\cap \{\varphi_2>2\varphi_1\}}\rho^{3/2}|K_h|\leq C \sigma_r r^3.$$
Applying \eqref{rsigr.van} once again, it follows that
$$\lim_{r\to 0}r^{-2}\int_{A_{r,2r}(p)\cap \{\varphi_2>2\varphi_1\}}\rho^{3/2}|K_h|=0.$$
Together with the previous estimate on the wedge $V_{\delta}\cap A_{r,2r}(p)$, this proves \eqref{w.ann.van} and therefore \eqref{corner.van}, so we can apply the maximum principle for $w$ in $\mathring\Sigma$ to conclude.
\end{proof}

With the global curvature bound $K_h\leq \lambda/3$ in place, we can now prove that $e^{2u}\leq C\sigma_r^2$ on $A_{r,2r}(p)$, to complete the proof of regularity.

\begin{prop}
    There is a constant $C<\infty$ such that
    $$e^{2u}\leq C\sigma_r^2$$
    on each annulus $A_{r,2r}(p)$ with $r\in [0,1/2]$. In particular, we have a uniform bound $v\leq C'$ on $\Sigma$.
\end{prop}
\begin{proof}
By Lemma \ref{wedge.bonanza}, we already have a bound of the form
$$e^{2u}\leq C_{\delta}\sigma_r^2\text{ on }A_{r,2r}(p)\cap V_{\delta}$$
for each $\delta>0$, so it remains to prove the bound along the edge, say on the domain $\Omega_r:=A_{r,2r}(p)\cap \{\varphi_2>2\varphi_1\}.$

Now, since $\Delta_{h_0}u=-e^{2u}K_h$ and 
$$\Delta_{h_0}\alpha_2=-\lambda e^{2u}-\langle d\psi,d\alpha_2\rangle_{h_0}=-\lambda e^{2u}-|d\alpha_2|_{h_0}^2+e^{2u}K_h,$$
we see that
$$\Delta_{h_0}(2u-\alpha_2-\log(\sigma_r/r))=|d\alpha_2|_{h_0}^2+(\lambda-3K_h)e^{2u}\geq 0$$
on $S\setminus \Gamma_2$, by virtue of the curvature bound $K_h\leq \lambda/3$. For any $x\in \Omega_r$, the mean value inequality then gives
$$(2u-\alpha_2-\log(\sigma_r/r))(x)\leq \frac{1}{\pi s^2}\int_{B_s^{h_0}(x)}(2u-\alpha_2-\log(\sigma_r/r)),$$
for $s\in (0,r)$. Taking $s=r/4$ and applying \eqref{u.a2.l1} and \eqref{a2.w12} to estimate the right-hand side, we then see that
$$2u-\alpha_2-\log(\sigma_r/r)\leq C\text{ on }\Omega_r.$$

Recall next that Lemma \ref{da2.edge.ann} gives the pointwise gradient bound
$$|d\alpha_2|_{h_0}^2=e^{2u}|d\alpha_2|_h^2\leq \frac{Cf_2^2e^{2u}}{\sigma_r^4r^4}\text{ on }\Omega_r,$$
and combining this with the upper bound for $2u-\alpha_2-\log(\sigma_r/r)$, it follows that
$$|d\alpha_2|_{h_0}^2\leq \frac{C'}{r^2}(\sigma_r^{-1}r^{-1}f_2)^3\text{ on }\Omega_r.$$
Using \eqref{f2.lps} to control the right-hand side, integrating this estimate to the power $k/2$ gives an estimate of the form
$$r^{1-2/k}\|d\alpha_2\|_{L^k(\Omega_r,dA_{h_0})}\leq C_k',$$
and taking $k=4$, a rescaling of the usual Morrey estimates on $A_{1/2,1}(p)\cap \{\varphi_2>2\varphi_1\}$ gives the bound
$$r^{1/2}[\alpha]_{C^{0,1/2}(\Omega_r)}\leq C''$$
for the $C^{1/2}$-seminorm. Finally, using this H\"older control together with \eqref{a2.w12} gives a pointwise bound
$$\alpha_2-\log(r\sigma_r)\leq C$$
with $C$ independent of $r$. 

Combining these estimates, we see finally that
$$2u\leq C+\alpha_2+\log(\sigma_r/r)\leq C_2+\log(\sigma_r^2),$$
which is clearly equivalent to the desired bound $e^{2u(x)}\leq C_3\sigma_r^2.$ Now that $e^{2u}\leq C\sigma_r^2$ throughout $A_{r,2r}(p)$, and Lemma \ref{ann.phi.below} gives $\phi^2\geq c\sigma_r^2$ in $A_{r,2r}(p)$, we obtain the desired uniform bound
$$e^{2v}=e^{2u}\phi^{-2}\leq C''\text{ on }A_{r,2r}(p)$$
with $C''$ independent of $r$, giving the desired $L^{\infty}$ bound for $v$ on $\Sigma$.
\end{proof}

\appendix

\printbibliography

@article{WangZiller1986,
  author  = {Wang, McKenzie Y. and Ziller, Wolfgang},
  title   = {Existence and non-existence of homogeneous {Einstein} metrics},
  journal = {Invent. Math.},
  volume  = {84},
  number  = {1},
  year    = {1986},
  pages   = {177--194},
  doi     = {10.1007/BF01388738}
}

@article{BohmWangZiller2004,
  author  = {B{\"o}hm, Christoph and Wang, McKenzie and Ziller, Wolfgang},
  title   = {A variational approach for compact homogeneous {Einstein} manifolds},
  journal = {Geom. Funct. Anal.},
  volume  = {14},
  number  = {4},
  year    = {2004},
  pages   = {681--733},
  doi     = {10.1007/s00039-004-0471-x}
}

@article{Bohm2004,
  author  = {B{\"o}hm, Christoph},
  title   = {Homogeneous {Einstein} metrics and simplicial complexes},
  journal = {J. Differential Geom.},
  volume  = {67},
  number  = {1},
  year    = {2004},
  pages   = {79--165},
  doi     = {10.4310/jdg/1099587730}
}

@article{Jensen1971,
  author  = {Jensen, Gary R.},
  title   = {The scalar curvature of left-invariant {Riemannian} metrics},
  journal = {Indiana Univ. Math. J.},
  volume  = {20},
  number  = {12},
  year    = {1971},
  pages   = {1125--1144},
  doi     = {10.1512/iumj.1971.20.20104}
}

@article{ChenLeBrunWeber2008,
  author  = {Chen, Xiuxiong and LeBrun, Claude and Weber, Brian},
  title   = {On conformally {K{\"a}hler}, {Einstein} manifolds},
  journal = {J. Amer. Math. Soc.},
  volume  = {21},
  number  = {4},
  year    = {2008},
  pages   = {1137--1168},
  doi     = {10.1090/S0894-0347-08-00594-8}
}

@article{BermanBoucksomGuedjZeriahi2013,
  author  = {Berman, Robert J. and Boucksom, S{\'e}bastien
             and Guedj, Vincent and Zeriahi, Ahmed},
  title   = {A variational approach to complex {Monge-Amp{\`e}re} equations},
  journal = {Publ. Math. Inst. Hautes {\'E}tudes Sci.},
  volume  = {117},
  year    = {2013},
  pages   = {179--245},
  doi     = {10.1007/s10240-012-0046-6}
}

@article{BermanBerndtsson2013,
  author  = {Berman, Robert J. and Berndtsson, Bo},
  title   = {Real {Monge-Amp{\`e}re} equations and {K{\"a}hler-Ricci}
             solitons on toric log {Fano} varieties},
  journal = {Ann. Fac. Sci. Toulouse Math. (6)},
  volume  = {22},
  number  = {4},
  year    = {2013},
  pages   = {649--711},
  doi     = {10.5802/afst.1386}
}

@article{BermanBoucksomJonsson2021,
  author  = {Berman, Robert J. and Boucksom, S{\'e}bastien
             and Jonsson, Mattias},
  title   = {A variational approach to the {Yau--Tian--Donaldson} conjecture},
  journal = {J. Amer. Math. Soc.},
  volume  = {34},
  number  = {3},
  year    = {2021},
  pages   = {605--652},
  doi     = {10.1090/jams/964}
}

@article{Lauret2002,
  author  = {Lauret, Jorge},
  title   = {Finding {Einstein} solvmanifolds by a variational method},
  journal = {Math. Z.},
  volume  = {241},
  year    = {2002},
  pages   = {83--99},
  doi     = {10.1007/s002090100407}
}

@article{liu2025toricEinstein,
  title={Toric Einstein 4-manifolds with non-negative sectional curvature},
  author={Liu, Tianyue},
  journal={Annals of Global Analysis and Geometry},
  volume={67},
  number={3},
  pages={14},
  year={2025},
  publisher={Springer}
}

@article{DarvasRubinstein2017,
  author  = {Darvas, Tam{\'a}s and Rubinstein, Yanir A.},
  title   = {{Tian}'s properness conjectures and {Finsler} geometry
             of the space of {K{\"a}hler} metrics},
  journal = {J. Amer. Math. Soc.},
  volume  = {30},
  number  = {2},
  year    = {2017},
  pages   = {347--387},
  doi     = {10.1090/jams/873}
}

@article{BermanBoucksomEyssidieuxGuedjZeriahi2019,
  author  = {Berman, Robert J. and Boucksom, S{\'e}bastien
             and Eyssidieux, Philippe and Guedj, Vincent
             and Zeriahi, Ahmed},
  title   = {{K{\"a}hler--Einstein} metrics and the
             {K{\"a}hler--Ricci} flow on log {Fano} varieties},
  journal = {J. Reine Angew. Math.},
  volume  = {751},
  year    = {2019},
  pages   = {27--89},
  doi     = {10.1515/crelle-2016-0033}
}

@article{BermanGuenancia2014,
  author  = {Berman, Robert J. and Guenancia, Henri},
  title   = {{K{\"a}hler--Einstein} metrics on stable varieties
             and log canonical pairs},
  journal = {Geom. Funct. Anal.},
  volume  = {24},
  number  = {6},
  year    = {2014},
  pages   = {1683--1730},
  doi     = {10.1007/s00039-014-0301-8}
}

@article{BergerCherrier1988,
  author  = {Berger, M. S. and Cherrier, P.},
  title   = {A new variational method for finding {Einstein}
             metrics on compact {K{\"a}hler} manifolds, {I}},
  journal = {J. Funct. Anal.},
  volume  = {79},
  number  = {1},
  year    = {1988},
  pages   = {103--135},
  doi     = {10.1016/0022-1236(88)90032-8}
}

@article{Yamabe1960Deformation,
  author  = {Yamabe, Hidehiko},
  title   = {On a deformation of {Riemannian} structures on compact manifolds},
  journal = {Osaka Math. J.},
  volume  = {12},
  number  = {1},
  pages   = {21--37},
  year    = {1960},
  doi     = {10.18910/8577}
}

@article{Page1978,
  author  = {Page, Don N.},
  title   = {A compact rotating gravitational instanton},
  journal = {Phys. Lett. B},
  volume  = {79},
  number  = {3},
  pages   = {235--238},
  year    = {1978},
  doi     = {10.1016/0370-2693(78)90231-9}
}

@article{TianYau1987,
  author  = {Tian, Gang and Yau, Shing-Tung},
  title   = {{K{\"a}hler--Einstein} metrics on complex surfaces with {$C_1>0$}},
  journal = {Comm. Math. Phys.},
  volume  = {112},
  number  = {1},
  pages   = {175--203},
  year    = {1987},
  doi     = {10.1007/BF01217685}
}

@article{Tian1990,
  author  = {Tian, Gang},
  title   = {On {Calabi}'s conjecture for complex surfaces with positive first {Chern} class},
  journal = {Invent. Math.},
  volume  = {101},
  number  = {1},
  pages   = {101--172},
  year    = {1990},
  doi     = {10.1007/BF01231499}
}

@book{GilbargTrudinger2001Elliptic,
  author    = {Gilbarg, David and Trudinger, Neil S.},
  title     = {Elliptic Partial Differential Equations of Second Order},
  edition   = {2},
  series    = {Classics in Mathematics},
  publisher = {Springer-Verlag},
  address   = {Berlin},
  year      = {2001},
  note      = {Reprint of the 1998 edition},
  doi       = {10.1007/978-3-642-61798-0}
}

@article{BrezisMerle1991Uniform,
  author  = {Brezis, Ha{\"i}m and Merle, Frank},
  title   = {Uniform estimates and blow-up behavior for solutions of
             {$-\Delta u=V(x)e^u$} in two dimensions},
  journal = {Communications in Partial Differential Equations},
  volume  = {16},
  number  = {8-9},
  pages   = {1223--1253},
  year    = {1991},
  doi     = {10.1080/03605309108820797}
}

@article{Ahlfors1938Schwarz,
  author  = {Ahlfors, Lars V.},
  title   = {An extension of {Schwarz}'s lemma},
  journal = {Transactions of the American Mathematical Society},
  volume  = {43},
  number  = {3},
  pages   = {359--364},
  year    = {1938},
  doi     = {10.1090/S0002-9947-1938-1501949-6}
}

@article{Aub78,
  author  = {Aubin, Thierry},
  title   = {{\'E}quations du type {Monge--Amp{\`e}re} sur les vari{\'e}t{\'e}s k{\"a}hl{\'e}riennes compactes},
  journal = {Bull. Sci. Math. (2)},
  volume  = {102},
  pages   = {63--95},
  year    = {1978}
}

@article{cortierMinerbe2016stationaryVacuum,
  title={On complete stationary vacuum initial data},
  author={Cortier, Julien and Minerbe, Vincent},
  journal={Journal of Geometry and Physics},
  volume={99},
  pages={20--27},
  year={2016},
  publisher={Elsevier}
}

@incollection{BeardonMinda2007HyperbolicMetricTheory,
  author    = {Beardon, Alan F. and Minda, David},
  title     = {The Hyperbolic Metric and Geometric Function Theory},
  booktitle = {Quasiconformal Mappings and Their Applications},
  editor    = {Ponnusamy, S. and Sugawa, T. and Vuorinen, M.},
  publisher = {Narosa Publishing House},
  address   = {New Delhi},
  pages     = {9--56},
  year      = {2007}
}

@book{GarnettMarshall2005HarmonicMeasure,
  author    = {Garnett, John B. and Marshall, Donald E.},
  title     = {Harmonic Measure},
  series    = {New Mathematical Monographs},
  volume    = {2},
  publisher = {Cambridge University Press},
  address   = {Cambridge},
  year      = {2005},
  doi       = {10.1017/CBO9780511546617}
}

@article{Mingione2010NonlinearCalderonZygmund,
  author  = {Mingione, Giuseppe},
  title   = {Nonlinear Aspects of {Calder{\'o}n--Zygmund} Theory},
  journal = {Jahresbericht der Deutschen Mathematiker-Vereinigung},
  volume  = {112},
  number  = {3},
  pages   = {159--191},
  year    = {2010},
  doi     = {10.1365/s13291-010-0004-5}
}

@book{Petersen2016RiemannianGeometry,
  author    = {Petersen, Peter},
  title     = {Riemannian Geometry},
  edition   = {3},
  series    = {Graduate Texts in Mathematics},
  volume    = {171},
  publisher = {Springer},
  address   = {Cham},
  year      = {2016},
  doi       = {10.1007/978-3-319-26654-1}
}

@article{Trudinger1967OrliczSpaceImbeddings,
  author  = {Trudinger, Neil S.},
  title   = {On Imbeddings into {Orlicz} Spaces and Some Applications},
  journal = {Journal of Mathematics and Mechanics},
  volume  = {17},
  number  = {5},
  pages   = {473--483},
  year    = {1967},
  doi     = {10.1512/iumj.1968.17.17028}
}

@misc{AlaeeKhuriKunduri2026MassToricManifolds,
  author        = {Alaee, Aghil and Khuri, Marcus and Kunduri, Hari},
  title         = {A Comparison Theorem for the Mass of {ALE} and {ALF}
                   Toric $4$-Manifolds},
  year          = {2026},
  eprint        = {2605.12352},
  archivePrefix = {arXiv},
  primaryClass  = {math.DG},
  doi           = {10.48550/arXiv.2605.12352}
}

@article{schoen1984conformal,
  title={Conformal deformation of a Riemannian metric to constant scalar curvature},
  author={Schoen, Richard},
  journal={Journal of Differential Geometry},
  volume={20},
  number={2},
  pages={479--495},
  year={1984},
  publisher={Lehigh University}
}

@article{HebeyVaugon1993Equivariant,
  author  = {Hebey, Emmanuel and Vaugon, Michel},
  title   = {Le probl{\`e}me de {Yamabe} {\'e}quivariant},
  journal = {Bulletin des Sciences Math{\'e}matiques},
  volume  = {117},
  number  = {2},
  pages   = {241--286},
  year    = {1993}
}

@article{Lions1985LimitI,
  author  = {Lions, Pierre-Louis},
  title   = {The concentration-compactness principle in the calculus
             of variations. The limit case, {Part 1}},
  journal = {Revista Matem{\'a}tica Iberoamericana},
  volume  = {1},
  number  = {1},
  pages   = {145--201},
  year    = {1985},
  doi     = {10.4171/RMI/6}
}

@article{Aubin1976Sobolev,
  author  = {Aubin, Thierry},
  title   = {Probl{\`e}mes isop{\'e}rim{\'e}triques et espaces de {Sobolev}},
  journal = {Journal of Differential Geometry},
  volume  = {11},
  number  = {4},
  pages   = {573--598},
  year    = {1976},
  doi     = {10.4310/jdg/1214433725}
}

@article{Talenti1976Sobolev,
  author  = {Talenti, Giorgio},
  title   = {Best constant in {Sobolev} inequality},
  journal = {Annali di Matematica Pura ed Applicata},
  volume  = {110},
  pages   = {353--372},
  year    = {1976},
  doi     = {10.1007/BF02418013}
}

@article{Hawking1977GravitationalInstantons,
  author  = {Hawking, S. W.},
  title   = {Gravitational instantons},
  journal = {Physics Letters A},
  volume  = {60},
  number  = {2},
  pages   = {81--83},
  year    = {1977},
  doi     = {10.1016/0375-9601(77)90386-3}
}

@article{aubin1976yamabe,
  title={{\'E}quations diff{\'e}rentielles non lin{\'e}aires et probl{\`e}me de Yamabe concernant la courbure scalaire},
  author={Aubin, Thierry},
  journal={J. Math. Pures Appl.(9)},
  volume={55},
  pages={269--296},
  year={1976}
}

@article{liSun2026geometryInstantons,
  title={Geometry of gravitational instantons},
  author={Li, Mingyang and Sun, Song},
  journal={arXiv preprint arXiv:2606.10443},
  year={2026}
}

@article{Kobayashi1987,
  author  = {Kobayashi, Osamu},
  title   = {Scalar curvature of a metric with unit volume},
  journal = {Mathematische Annalen},
  volume  = {279},
  number  = {2},
  year    = {1987},
  pages   = {253--265},
  doi     = {10.1007/BF01461722}
}

@misc{Ozuch2026PositiveEinstein,
  author        = {Ozuch, Tristan},
  title         = {Positive {Einstein} metrics on infinitely many four-manifolds},
  year          = {2026},
  eprint        = {2610.08052},
  archivePrefix = {arXiv},
  primaryClass  = {math.DG},
  url           = {https://arxiv.org/abs/2610.08052}
}

@misc{broderLafuente2026t2Einstein,
      title={$T^2$-invariant Einstein metrics on closed 4-manifolds}, 
      author={Kyle Broder and Ramiro A. Lafuente},
      year={2026},
      eprint={2610.18219},
      archivePrefix={arXiv},
      primaryClass={math.DG},
      url={https://arxiv.org/abs/2610.18219}, 
}

@article{Yau1978,
  author  = {Yau, Shing-Tung},
  title   = {On the {Ricci} curvature of a compact {K{\"a}hler} manifold and the complex {Monge--Amp{\`e}re} equation. {I}},
  journal = {Comm. Pure Appl. Math.},
  volume  = {31},
  number  = {3},
  pages   = {339--411},
  year    = {1978},
  doi     = {10.1002/cpa.3160310304}
}

@article{Anderson2006,
  author        = {Anderson, Michael T.},
  title         = {{Dehn} filling and {Einstein} metrics in higher dimensions},
  journal       = {J. Differential Geom.},
  volume        = {73},
  number        = {2},
  pages         = {219--261},
  year          = {2006},
  doi           = {10.4310/jdg/1146169911},
  eprint        = {math/0303260},
  archivePrefix = {arXiv},
  primaryClass  = {math.DG}
}

@article{Bamler2012,
  author        = {Bamler, Richard H.},
  title         = {Construction of {Einstein} metrics by generalized {Dehn} filling},
  journal       = {J. Eur. Math. Soc. (JEMS)},
  volume        = {14},
  number        = {3},
  pages         = {887--909},
  year          = {2012},
  doi           = {10.4171/JEMS/321},
  eprint        = {0911.4730},
  archivePrefix = {arXiv},
  primaryClass  = {math.DG}
}

@article{FinePremoselli2020,
  author        = {Fine, Joel and Premoselli, Bruno},
  title         = {Examples of compact {Einstein} four-manifolds with negative curvature},
  journal       = {J. Amer. Math. Soc.},
  volume        = {33},
  number        = {4},
  pages         = {991--1038},
  year          = {2020},
  doi           = {10.1090/jams/944},
  eprint        = {1802.00608},
  archivePrefix = {arXiv},
  primaryClass  = {math.DG}
}

@misc{HamenstadtJackel2024,
  author        = {Hamenst{\"a}dt, Ursula and J{\"a}ckel, Frieder},
  title         = {Negatively curved {Einstein} metrics on {Gromov--Thurston} manifolds},
  year          = {2024},
  eprint        = {2411.12956},
  archivePrefix = {arXiv},
  primaryClass  = {math.DG},
  note          = {To appear in J. Eur. Math. Soc. (JEMS)},
  url           = {https://arxiv.org/abs/2411.12956}
}

@misc{GuenanciaHamenstadt2025,
  author        = {Guenancia, Henri and Hamenst{\"a}dt, Ursula},
  title         = {{K{\"a}hler--Einstein} metrics of negative curvature},
  year          = {2025},
  eprint        = {2503.02838},
  archivePrefix = {arXiv},
  primaryClass  = {math.DG},
  url           = {https://arxiv.org/abs/2503.02838}
}

@article{Trudinger1968Conformal,
  author  = {Trudinger, Neil S.},
  title   = {Remarks concerning the conformal deformation of {Riemannian} structures on compact manifolds},
  journal = {Ann. Scuola Norm. Sup. Pisa Cl. Sci. (3)},
  volume  = {22},
  number  = {2},
  pages   = {265--274},
  year    = {1968},
  url     = {https://numdam.org/item/ASNSP_1968_3_22_2_265_0/}
}

@misc{liSun2026ToricEinstein,
      title={Compact toric Einstein four-manifolds}, 
      author={Mingyang Li and Song Sun},
      year={2026},
      eprint={2610.10425},
      archivePrefix={arXiv},
      primaryClass={math.DG},
      url={https://arxiv.org/abs/2610.10425}, 
}

@inproceedings{viaclovsky2010monopole,
  title={Monopole metrics and the orbifold Yamabe problem},
  author={Viaclovsky, Jeff A},
  booktitle={Annales de l'institut Fourier},
  volume={60},
  number={7},
  pages={2503--2543},
  year={2010}
}

@book{evans1990weakConv,
  title={Weak convergence methods for nonlinear partial differential equations},
  author={Evans, Lawrence C},
  volume={74},
  year={1990},
  publisher={American Mathematical Soc.}
}

@incollection{schoen1987variational,
  title={Variational theory for the total scalar curvature functional for Riemannian metrics and related topics},
  author={Schoen, Richard M},
  booktitle={Topics in Calculus of Variations: Lectures given at the 2nd 1987 Session of the Centro Internazionale Matematico Estivo (CIME) held at Montecatini Terme, Italy, July 20--28, 1987},
  pages={120--154},
  year={1987},
  publisher={Springer}
}

@article{dammerman2006metricsWithSymmetry,
  title={Metrics of special curvature with symmetry},
  author={Dammerman, Brandon},
  journal={arXiv preprint math/0610738},
  year={2006}
}

@article{anderson2008survey,
  title={A survey of Einstein metrics on 4-manifolds},
  author={Anderson, Michael T},
  journal={arXiv preprint arXiv:0810.4830},
  year={2008}
}

@article{chuStern2025minimalDoublings,
  title={Minimal surface doublings and electrostatics for Schr\"odinger operators},
  author={Chu, Adrian Chun-Pong and Stern, Daniel},
  journal={arXiv preprint arXiv:2509.18630},
  year={2025}
}

@book{besse2007einstein,
  title={Einstein manifolds},
  author={Besse, Arthur L},
  year={2007},
  publisher={Springer}
}

@article{obata1971conformal,
  title={The conjectures on conformal transformations of Riemannian manifolds},
  author={Obata, Morio},
  journal={Bulletin of the American Mathematical Society},
  volume={77},
  number={2},
  pages={265--270},
  year={1971}
}

@article{leeParker1987yamabe,
  title={The Yamabe problem},
  author={Lee, John M and Parker, Thomas H},
  journal={Bulletin of the American Mathematical Society},
  volume={17},
  number={1},
  pages={37--91},
  year={1987}
}
\end{document}